\RequirePackage[2021-11-14]{latexrelease}
\documentclass{amsart}

\makeatletter
\DeclareRobustCommand*\cal{\@fontswitch\relax\mathcal}
\makeatother

\usepackage[margin=2cm]{geometry}
\usepackage{enumitem}
\setlist{noitemsep,topsep=2pt plus 2pt,parsep=0pt,partopsep=0pt,itemsep=1pt plus 1pt}

\usepackage{silence}
\usepackage[T1]{fontenc}
\usepackage[utf8]{inputenc}
\usepackage{newunicodechar}

\usepackage{amsmath}
\usepackage{amssymb}
\usepackage{amsthm}
\usepackage[mathscr]{euscript}

\usepackage[all,cmtip]{xy}
\SilentMatrices

\usepackage{tikz}
\usetikzlibrary{cd}
\usetikzlibrary{decorations.markings}
\usetikzlibrary{decorations.pathmorphing}
\usetikzlibrary{decorations.pathreplacing,calc,tikzmark}
\usetikzlibrary{positioning,arrows.meta,automata,shapes}

\tikzstyle{none}=[]
\tikzstyle{new style 0}=[fill=white, draw=black, shape=circle]
\tikzstyle{new style 1}=[fill=black, draw=black, shape=circle]

\tikzstyle{new edge style 0}=[->]
\tikzstyle{new edge style 1}=[<-]
\tikzstyle{new edge style 2}=[-]
\tikzstyle{new edge style 3}=[-, dashed]
\tikzstyle{new edge style 4}=[-, dashed, orange]
\tikzstyle{new edge style 5}=[-, orange]

\usepackage{pgfplots}
\pgfplotsset{width=7cm,compat=1.18}
\pgfdeclarelayer{edgelayer}
\pgfdeclarelayer{nodelayer}
\pgfsetlayers{edgelayer,nodelayer,main}

\usepackage{xr-hyper}
\def\gcref#1{Grady–Pav\-lov \cite[\cref*{GL-#1}]{GradyPavlov.GCH}}
\def\scref#1{Grady–Pav\-lov \cite[\cref*{SL-#1}]{GradyPavlov.Str}}

\newif\ifpdf \pdffalse \ifx\pdfoutput\undefined\else\ifx\pdfoutput\relax\else\ifnum\pdfoutput<1 \else\pdftrue\fi\fi\fi
\ifpdf
\usepackage[pdfencoding=auto,psdextra,backref=section,linktoc=all]{hyperref}

\else
\usepackage[hypertex,backref=section,linktoc=all]{hyperref}
\fi

\usepackage[capitalize,nosort]{cleveref}

\crefformat{equation}{(#2#1#3)}
\crefname{figure}{Figure}{Figures}

\catcode`@=11
\let\over\@@over
\let\atop\@@atop
\let\above\@@above
\let\overwithdelims\@@overwithdelims
\let\atopwithdelims\@@atopwithdelims
\let\abovewithdelims\@@abovewithdelims
\def\eqalign#1{\null\,\vcenter{\openup\jot\m@th\ialign{\strut\hfil$\displaystyle{##}$&$\displaystyle{{}##}$\hfil\crcr#1\crcr}}\,}
\newskip\xcentering
\def\eqalignno#1{\displ@y \tabskip\xcentering
	\halign to\displaywidth{\hfil$\@lign\displaystyle{##}$\tabskip\z@skip
		&$\@lign\displaystyle{{}##}$\hfil\tabskip\xcentering
		&\llap{$\@lign##$}\tabskip\z@skip\crcr
		#1\crcr}}
\def\eqlabel#1{\refstepcounter{equation}\label{#1}\ifmmode\ifinner\else\eqno\fi\fi\hbox{\@eqnnum}} 

\catcode`@=12

\makeatletter
\@addtoreset{equation}{section}
\@addtoreset{equation}{subsection}
\let\proof@qed\displaymath@qed

\let\c@subsubsection\c@equation

\makeatother

\theoremstyle{definition}
\newtheorem{theorem}[equation]{Theorem}
\newtheorem{definition}[equation]{Definition}
\newtheorem{proposition}[equation]{Proposition}
\newtheorem{conjecture}[equation]{Conjecture}
\newtheorem{construction}[equation]{Construction}
\newtheorem{remark}[equation]{Remark}

\newtheorem{notation}[equation]{Notation}

\newtheorem{example}[equation]{Example}

\numberwithin{equation}{subsection}
\numberwithin{subsubsection}{subsection}

\newunicodechar{α}{\alpha} \newunicodechar{β}{\beta} \newunicodechar{γ}{\gamma} \newunicodechar{δ}{\delta} \newunicodechar{ϵ}{\epsilon} \newunicodechar{ζ}{\zeta} \newunicodechar{η}{\eta} \newunicodechar{θ}{\theta} \newunicodechar{ι}{\iota} \newunicodechar{κ}{\kappa} \newunicodechar{λ}{\lambda} \newunicodechar{μ}{\mu} \newunicodechar{ν}{\nu} \newunicodechar{ξ}{\xi} \newunicodechar{ο}{o} \newunicodechar{π}{\pi} \newunicodechar{ρ}{\rho} \newunicodechar{σ}{\sigma} \newunicodechar{τ}{\tau} \newunicodechar{υ}{\upsilon} \newunicodechar{ϕ}{\phi} \newunicodechar{χ}{\chi} \newunicodechar{ψ}{\psi} \newunicodechar{ω}{\omega} \newunicodechar{ε}{\varepsilon} \newunicodechar{ϑ}{\vartheta} \newunicodechar{ϖ}{\varpi} \newunicodechar{ϱ}{\varrho} \newunicodechar{ς}{\varsigma} \newunicodechar{φ}{\varphi} \newunicodechar{Γ}{\Gamma} \newunicodechar{Δ}{\Delta} \newunicodechar{Θ}{\Theta} \newunicodechar{Λ}{\Lambda} \newunicodechar{Ξ}{\Xi} \newunicodechar{Π}{\Pi} \newunicodechar{Σ}{\Sigma} \newunicodechar{Υ}{\Upsilon} \newunicodechar{Φ}{\Phi} \newunicodechar{Ψ}{\Psi} \newunicodechar{Ω}{\Omega} \newunicodechar{ℵ}{\aleph} \newunicodechar{ℏ}{\hbar} \newunicodechar{ℓ}{\ell} \newunicodechar{℘}{\wp} \newunicodechar{ℜ}{\Re} \newunicodechar{ℑ}{\Im} \newunicodechar{∂}{\partial} \newunicodechar{∞}{\infty} \newunicodechar{∅}{\emptyset} \newunicodechar{∇}{\nabla} \newunicodechar{√}{\surd} \newunicodechar{⊤}{\top} \newunicodechar{⊥}{\bot} \newunicodechar{∠}{\angle} \newunicodechar{△}{\triangle} \newunicodechar{∀}{\forall} \newunicodechar{∃}{\exists} \newunicodechar{♭}{\flat} \newunicodechar{♮}{\natural} \newunicodechar{♯}{\sharp} \newunicodechar{♣}{\clubsuit} \newunicodechar{♢}{\diamondsuit} \newunicodechar{♡}{\heartsuit} \newunicodechar{♠}{\spadesuit} \newunicodechar{∐}{\coprod} \newunicodechar{⋁}{\bigvee} \newunicodechar{⋀}{\bigwedge} \newunicodechar{⨄}{\biguplus} \newunicodechar{⋂}{\bigcap} \newunicodechar{⋃}{\bigcup} \newunicodechar{∫}{\int} \newunicodechar{∏}{\prod} \newunicodechar{∑}{\sum} \newunicodechar{⨂}{\bigotimes} \newunicodechar{⨁}{\bigoplus} \newunicodechar{⨀}{\bigodot} \newunicodechar{⨆}{\bigsqcup} \newunicodechar{◁}{\triangleleft} \newunicodechar{▷}{\triangleright} \newunicodechar{▽}{\bigtriangledown} \newunicodechar{∧}{\wedge} \newunicodechar{∨}{\vee} \newunicodechar{∩}{\cap} \newunicodechar{∪}{\cup} \newunicodechar{⊓}{\sqcap} \newunicodechar{⊔}{\sqcup} \newunicodechar{⊎}{\uplus} \newunicodechar{⨿}{\amalg} \newunicodechar{⋄}{\diamond} \newunicodechar{∙}{\bullet} \newunicodechar{≀}{\wr} \newunicodechar{⊙}{\odot} \newunicodechar{⊘}{\oslash} \newunicodechar{⊗}{\otimes} \newunicodechar{⊖}{\ominus} \newunicodechar{⊕}{\oplus} \newunicodechar{∓}{\mp} \newunicodechar{∘}{\circ} \newunicodechar{○}{\Orb} \newunicodechar{∖}{\setminus} \newunicodechar{⋅}{\cdot} \newunicodechar{∗}{\ast} \newunicodechar{⨯}{\times} \newunicodechar{⋆}{\star} \newunicodechar{∝}{\propto} \newunicodechar{⊑}{\sqsubseteq} \newunicodechar{⊒}{\sqsupseteq} \newunicodechar{∥}{\parallel} \newunicodechar{∣}{\divides} \newunicodechar{⊣}{\dashv} \newunicodechar{⊢}{\vdash} \newunicodechar{↗}{\nearrow} \newunicodechar{↘}{\searrow} \newunicodechar{↖}{\nwarrow} \newunicodechar{↙}{\swarrow} \newunicodechar{⇔}{\Leftrightarrow} \newunicodechar{⇐}{\Leftarrow} \newunicodechar{⇒}{\Rightarrow} \newunicodechar{≠}{\neq} \newunicodechar{≤}{\leq} \newunicodechar{≥}{\geq} \newunicodechar{≻}{\succ} \newunicodechar{≺}{\prec} \newunicodechar{≈}{\approx} \newunicodechar{≽}{\succeq} \newunicodechar{≼}{\preceq} \newunicodechar{⊃}{\supset} \newunicodechar{⊂}{\subset} \newunicodechar{⊇}{\supseteq} \newunicodechar{⊆}{\subseteq} \newunicodechar{∈}{\in} \newunicodechar{∋}{\ni} \newunicodechar{≫}{\gg} \newunicodechar{≪}{\ll} \newunicodechar{↔}{\leftrightarrow} \newunicodechar{↦}{\mapsto} \newunicodechar{∼}{\sim} \newunicodechar{≃}{\simeq} \newunicodechar{≡}{\equiv} \newunicodechar{≍}{\asymp} \newunicodechar{⌣}{\smile} \newunicodechar{⌢}{\frown} \newunicodechar{↼}{\leftharpoonup} \newunicodechar{↽}{\leftharpoondown} \newunicodechar{⇀}{\rightharpoonup} \newunicodechar{⇁}{\rightharpoondown} \newunicodechar{↪}{\hookrightarrow} \newunicodechar{↩}{\hookleftarrow} \newunicodechar{⋈}{\bowtie} \newunicodechar{⊨}{\models} \newunicodechar{⟹}{\Longrightarrow} \newunicodechar{⟶}{\longrightarrow} \newunicodechar{⟵}{\longleftarrow} \newunicodechar{⟸}{\Longleftarrow} \newunicodechar{⟼}{\longmapsto} \newunicodechar{⟷}{\longleftrightarrow} \newunicodechar{⟺}{\Longleftrightarrow} \newunicodechar{⋯}{\cdots} \newunicodechar{⋮}{\vdots} \newunicodechar{⋱}{\ddots} \newunicodechar{↕}{\updownarrow} \newunicodechar{⇑}{\Uparrow} \newunicodechar{⇓}{\Downarrow} \newunicodechar{⇕}{\Updownarrow} \newunicodechar{⌉}{\rceil} \newunicodechar{⌈}{\lceil} \newunicodechar{⌋}{\rfloor} \newunicodechar{⌊}{\lfloor} \newunicodechar{≅}{\cong} \newunicodechar{∉}{\notin} \newunicodechar{⇌}{\rightleftharpoons} \newunicodechar{≐}{\doteq} \newunicodechar{∄}{\not\exists} \newunicodechar{∌}{\not\ni} \newunicodechar{∔}{\dot+} \newunicodechar{∕}{/} \newunicodechar{∤}{\not|} \newunicodechar{∦}{\not\|} \newunicodechar{∬}{\int\!\!\!\int} \newunicodechar{∭}{\int\!\!\!\int\!\!\!\int} \newunicodechar{∮}{\oint} \newunicodechar{∸}{\dot-} \newunicodechar{≁}{\not\sim} \newunicodechar{≄}{\not\simeq} \newunicodechar{≆}{\not\cong} \newunicodechar{≇}{\not\cong} \newunicodechar{≉}{\not\approx} \newunicodechar{≔}{:=} \newunicodechar{≕}{=:} \newunicodechar{≢}{\not\equiv} \newunicodechar{≭}{\not\asump} \newunicodechar{≮}{\not<} \newunicodechar{≯}{\not<} \newunicodechar{≰}{\not\le} \newunicodechar{≱}{\not\ge} \newunicodechar{⊀}{\not\prec} \newunicodechar{⊁}{\not\succ} \newunicodechar{⊄}{\not\subset} \newunicodechar{⊅}{\not\supset} \newunicodechar{⊈}{\not\subseteq} \newunicodechar{⊉}{\not\supseteq} \newunicodechar{⊦}{\vdash} \newunicodechar{⊧}{\models} \newunicodechar{⊬}{\not\vdash} \newunicodechar{⊭}{\not\models} \newunicodechar{⊲}{\triangleleft} \newunicodechar{⊳}{\triangleright} \newunicodechar{⋠}{\not\preceq} \newunicodechar{⋡}{\not\succeq} \newunicodechar{⋤}{\not\sqsubseteq} \newunicodechar{⋥}{\not\sqsupseteq} \newunicodechar{⋪}{\not\triangleleft} \newunicodechar{⋫}{\not\triangleright} \newunicodechar{◻}{\square} 

\newunicodechar{↑}{\uparrow}
\newunicodechar{↓}{\downarrow}
\newunicodechar{→}{\to}
\newunicodechar{←}{\gets}
\newunicodechar{⟨}{\langle}
\newunicodechar{⟩}{\rangle}
\newunicodechar{…}{\ldots}
\newunicodechar{∞}{\ifmmode\infty\else$\infty$\fi}
\newunicodechar{≔}{\colon=}
\let\oldO\O
\newunicodechar{Ø}{\oldO}
\newunicodechar{⊠}{\boxtimes}
\def\ppet{\mathbin{\bar\boxtimes}} 

\def\hom{\mathop{\rm Hom}\nolimits}
\def\map{\mathop{\rm Map}\nolimits}
\def\rmap{\mathop{\rm map}\nolimits}
\def\Fun{\mathop{\rm Fun}\nolimits}
\def\Funmon{\mathop{\rm Hom}\nolimits}
\def\Funmonglob{\Fun^\otimes}
\def\Funmonuple{\Fun^\otimes_{\uple}}
\def\Funmontheta{\Funmon_\Theta}
\def\hocolim{\mathop{\rm hocolim}}
\def\holim{\mathop{\rm holim}}
\def\colim{\mathop{\rm colim}}
\def\sing{\mathop{\rm sing}}
\def\core{\mathop{\rm core}\nolimits}
\def\Sd{\mathop{\rm Sd}\nolimits}
\def\Ex{\mathop{\rm Ex}\nolimits}
\def\Exi{\Ex^\infty}

\def\FEmb{{\sf FEmb}}
\def\frakFEmb{\mathfrak{FEmb}}
\def\Struct{{\sf Struct}}
\def\frakStruct{\mathfrak{Struct}}

\let\into\hookrightarrow
\def\cC{{\cal C}}

\def\GL{{\rm GL}}
\def\O{{\rm O}}

\def\stu{\Re}

\def\RR{{\bf R}}
\def\NN{{\bf N}}
\def\ZZ{{\bf Z}}

\def\proj{{\sf proj}}
\def\inj{{\sf inj}}
\def\glob{{\sf glob}}
\def\local{{\sf local}}
\def\uple{{\sf uple}}
\def\tdeloop{{\rm B}}
\def\rdf{{\bf R}}

\def\gs{{\cal S}} 

\def\set{\mathscr{S}\mathrm{et}}
\def\sm{{\rm C}^\infty}
\def\sset{\mathrm{s}\set}
\def\poset{\mathscr{P}{\rm o}\set}
\def\smset{\sm\set}
\def\smposet{\sm\poset}
\def\smsset{\sm\sset}
\def\PSh{\mathscr{PS}\mathrm{h}}
\def\sPSh{\PSh_\Delta}
\def\smPSh{{\sm\PSh}}
\def\smsPSh{{\sm\sPSh}}
\def\cat{{\rm Cat}^\otimes}
\def\catuple{{\rm Cat}^{\otimes,\uple}}
\def\smcat{\sm\cat}
\def\smcatuple{\sm\catuple}
\def\fraksmcat{\sm\mathfrak{Cat}^{\otimes}}
\def\fraksmcatuple{\sm\mathfrak{Cat}^{\otimes,\uple}}

\def\Man{{\sf Man}}
\def\SMan{{\sf SMan}}

\def\stcart{{\sf StCart}}
\def\cart{{\sf Cart}}
\def\stman{{\sf StMan}}
\def\man{{\sf Man}}

\def\Ob{{\sf Ob}}
\def\Mor{{\sf Mor}}
\def\ev{{\rm ev}}

\def\Cut{{\sf Cut}}
\def\tCut{{\sf Cut}_\pitchfork^\uple}
\def\tCutglob{{\sf Cut}_\pitchfork}

\def\Yo#1{{\cal Y}_{#1}}

\def\Emb{{\rm Emb}}
\def\id{{\rm id}}
\def\FRiem{{\rm FRiem}}
\def\Vect{{\rm Vect}}
\def\Bord{{\sf Bord}}
\def\BBord{\mathfrak{Bord}}
\def\EBord{{\sf E}}
\def\frakEBord{\mathfrak{E}}
\def\gBord{{\sf gBord}}
\def\frakB{\mathfrak{B}}
\def\frakD{\mathfrak{D}}
\def\FFT{{\sf FFT}}
\def\op{{\rm op}}
\def\csp{{\sf B_\smallint}}
\def\deloop{{\bf B}}
\def\hq{{/\!/}}
\def\simcon{\sim_{\rm con}}
\def\Ei{{\rm E}_{\infty}}
\def\Nec{{\sf Nec}}
\def\Cnec{\mathfrak{C}^{\rm nec}}
\def\Cech{{\sf \check Cech}}
\def\frakA{\mathfrak{A}}
\def\fraka{\mathfrak{a}}
\def\frakb{\mathfrak{b}} 

\def\gpop{{\bf P}} 
\def\tpop{{\rm P}} 

\def\gsim{{\bf\Delta}} 
\def\gbou{{\bf\partial\Delta}} 

\def\frake{\mathfrak{e}} 

\def\cC{{\sf C}} 
\def\cD{{\sf D}} 

\def\ax#1{(A#1)}
\def\frakax#1{($\mathfrak{A}$#1)}

\def\ltoarr#1{\mathop{\count0=#1 \loop\ifnum\count0>0 \smash-\mkern-7mu \advance\count0 -1 \repeat \mathord\rightarrow}\limits} 
\def\lto#1#2{\mathrel{\ltoarr{#1}^{#2}}} 
\def\lgetsarr#1{\mathop{\mathord\leftarrow \count0=#1 \loop\ifnum\count0>0 \mkern-7mu\smash-\advance\count0 -1 \repeat}\limits} 
\def\lgets#1#2{\mathrel{\lgetsarr{#1}\limits^{#2}}} 
\def\frakFRiem{\mathfrak{FRiem}}
\def\tglob{\otimes_{\glob}}
\mathcode`:="603A 
\mathchardef\colon="303A 

\DeclareRobustCommand\and{\end{tabular}\hskip 1em plus.17fil \begin{tabular}[t]{c}}
\def\maketitle{\null\vskip 2em \begin{center}{\LARGE\atitle\par}\vskip 1.5em {\large\lineskip .5em \begin{tabular}[t]{c}\authors\end{tabular}\par}\end{center}\par\vskip 1.5em }
\def\typesetabstract{\begingroup
	\skip0=20pt \advance\skip0 -\lastskip \advance\skip0 -\baselineskip \vskip\skip0
	\box\abstractbox
	\prevdepth0pt
	\normalsize
	\dimen0=34pt \advance\dimen0 -\baselineskip \vskip\dimen0 \relax
	\endgroup}
\def\authors{Daniel Grady\\\href{http://www.gradydaniel.com/}{gradydaniel.com}
\\Wichita State University
\and
Dmitri Pavlov\\\href{https://dmitripavlov.org/}{dmitripavlov.org}
\\Texas Tech University}
\def\atitle{Extended field theories are local and have classifying spaces}

\begin{document}

\maketitle

\begin{abstract}
We show that all extended functorial field theories, both topological and nontopological, are local.
We define the smooth $(\infty,d)$-category of bordisms with geometric data, such as Riemannian metrics or geometric string structures,
and prove that it satisfies codescent with respect to the target $\gs$, which implies the locality property.
We apply this result to construct a classifying space for concordance classes of functorial field theories with geometric data,
solving a conjecture of Stolz and Teichner about the existence of such a space.
We use our classifying space construction to develop a geometric theory of power operations, following the recent work of Barthel, Berwick-Evans, and Stapleton.
\end{abstract}

\typesetabstract

\tableofcontents

\section{Introduction}

A key feature of functorial field theory (FFT) that is usually expected, or demanded, is that it is \emph{local} in the sense that global phenomena are entirely determined by local phenomena. 
All influences should propagate through space (or spacetime) at some finite speed and there is no `action at a distance'. This is a familiar property for classical field theories, but there is still nontrivial locality even for quantum field theories (where one has nontrivial moduli for bordisms). 
Although the axiomatization of field theories via functorial field theories (see the original paper by Segal \cite{Segal.CFT}
or Atiyah's paper \cite{Atiyah} for the topological case) provides a beautiful framework that captures many necessary ingredients, 
it has the unfortunate side effect of allowing field theories that are not local. 
This is almost evident from the axioms, for suppose we are given a bordism~$\Sigma$ (i.e., a $d$-dimensional smooth manifold with boundary)
between two $(d-1)$-manifolds $\Sigma_0$ and~$\Sigma_1$ and suppose we have a smooth map $f:\Sigma\to X$ to some target manifold~$X$.
To such data, a functorial field theory associates a linear map between Hilbert spaces
$$Z_{(\Sigma,f)}:H_{(\Sigma_0,f\vert_{\Sigma_0})}\to H_{(\Sigma_1,f\vert_{{\Sigma_1}})}.$$
Now in this setting `local' means local with respect to the target manifold~$X$.
In applications, $\Sigma$ should represent the `worldline' of a particle, or `worldvolume' of a string
and $X$~is the ambient space in which it propagates.
If this field theory were local, we should be able to restrict it to some small region $U\subset X$. 
Moreover, we should be able to reconstruct all of the field theory from such restrictions, once $U$ runs over elements of an open cover of $X$. 
Here we have a problem. 
If the bounding manifolds $\Sigma_0$ and $\Sigma_1$ have images in~$X$ (under~$f$) that cover a large region in~$X$,
we have no hope for being able to assemble the field theory from its restriction to a cover $\{U_\alpha\}$ of~$X$. 
Indeed, there is no reason to expect that the maps $f_i:\Sigma_i\to X$, $i=0,1$, factor through some~$U_\alpha$
and we are unable to cut~$\Sigma_i$ into smaller pieces so that it factors. 

This defect of functorial field theories can be resolved by allowing ourselves to further cut~$\Sigma_i$ into smaller pieces by codimension~1 submanifolds. 
This naturally leads to the fully-extended setting of Freed \cite{Freed.Extended}, Lawrence \cite{Lawrence}, Baez–Dolan \cite{BaezDolan}, and Lurie \cite{Lurie.TFT},
where we have a $d$-category of bordisms, starting with points as objects, 1-manifolds as morphisms, 2-manifolds as morphisms between morphisms etc. 
It has long been expected that fully-extended field theories ought to be local, but no actual proof has emerged in the literature. 
One part of the issue for this is that a rigorous definition (in the topological case) of the fully-extended bordism category has not appeared until Lurie \cite{Lurie.TFT} and Calaque–Scheimbauer \cite{CalaqueScheimbauer}.
Moreover, to talk about the nontopological case, one really needs to parameterize the bordism category over smooth manifolds,
and this adds another layer of complexity that has not yet been explored in the literature in the fully extended case (see Stolz and Teichner \cite{StolzTeichner.Elliptic,StolzTeichner.SUSY} for the nonextended case). 

The main goal of the present work is to establish several results, which culminate in a proof of the locality property for the bordism category.
Specifically, we provide the following.
\begin{itemize}
\item An axiomatization of the fully extended $(\infty,d)$-category of bordisms that allows for geometric structures on bordisms and, optionally, incorporates isotopies as higher morphisms (i.e., $n$-morphisms for $n>d$). 
\item A proof of existence and uniqueness (up to contractible space of choices) of bordism categories that satisfy the aforementioned axioms (in both cases, with and without isotopies).
\item A proof that any bordism category satisfying the axioms (in both cases, with and without isotopies) satisfies locality in the sense of \cref{local}.
\end{itemize}

To incorporate a variety of geometric structures on smooth manifolds,
such as Riemannian metrics, orientations, spin and string structures,
as well as more complicated higher geometric structures such as geometric string structures,
we generalize the usual definition of tangential structures, i.e., the space of lifts
$$\xymatrix{
& Y\ar[d]\\
M\ar[r]_-{\tau_M}\ar[ru]^-{\ell} & \deloop\GL(d),
}$$ 
where $\tau_M$ denotes the tangent bundle and $Y$ is a space encoding the tangential structure.
We are interested in a variety of structures that do not fit into the above formalism. 
For example, we can consider Riemannian metrics with constraints on the sectional curvature (e.g., positive, negative, nonpositive, nonnegative).

We adopt an alternative perspective.
Tangential structures can be pulled back along open embeddings.
They are also local in that structures can be glued from local data.  
Thus, a tangential structure naturally gives rise to an $\infty$-sheaf on the site of $d$-dimensional manifolds, with open embeddings between them, and we might as well consider any $\infty$-sheaf on this site as a legitimate geometric structure.
This perspective is explained in  \cref{geostr} and in more detail in our separate companion paper \cite{GradyPavlov.Str}.
The following list provides examples of geometric structures that can be encoded in our setting.
\begin{itemize}
\item Smooth maps to some fixed target manifold.
\item Super-Euclidean structures of Stolz--Teichner \cite{StolzTeichner.SUSY}.
\item Principal $G$-bundles with connection.
\item Conformal, complex, symplectic, contact, K\"ahler structures.
\item Bundle $n$-gerbes with connection.
\item Embeddings into a target manifold $W$.
\item Foliations, possibly equipped with additional structures such as transversal metrics.
\item Geometric spin-c, string, fivebrane, and ninebrane structures.
\item Topological tangential structures, such as orientation, spin, string, framing, etc. 
\item Riemannian and pseudo-Riemannian metrics, possibly with restrictions on Ricci or sectional curvature.
\end{itemize}

For a geometric structure $\gs$, we get a smooth symmetric monoidal $(\infty,d)$-category of bordisms $\Bord_d^\gs$.  
For $T$ is a smooth symmetric monoidal $(\infty,d)$-category, a $T$-valued fully-extended field theory $\mathcal{F}$ with $\gs$-structure
is a symmetric monoidal functor $\mathcal{F}:\Bord_d^\gs\to T $. To say that $\mathcal{F}$ is local roughly means that
for any open cover of the target~$\gs$
we can reconstruct $\mathcal{F}$ from its restriction to bordisms $M→\gs$ that factor through some element of the open cover. Thus, a natural formulation of locality involves descent on a general smooth stack $\gs\in \sPSh(\FEmb_d)$,
as in \cref{geometric.structure}. 

\begin{theorem}
\label{local}
The smooth symmetric monoidal $(\infty,d)$-category of symmetric monoidal functors $\Funmonglob(\Bord_d^\gs,T)$
satisfies descent with respect to the target $\gs$, for all smooth symmetric monoidal $(\infty,d)$-categories $T$.
That is to say, the functor
$$\FFT_{d,T}:\Struct_d^{\op}→\smcat_{\infty,d},\qquad \gs↦\FFT_{d,T}(\gs)≔\Funmonglob(\Bord_d^\gs,T)$$
is an ∞-sheaf. 
\end{theorem}

\begin{proof}
Combine \cref{a1a2a3} with \cref{derivedhom}.
\end{proof}

\begin{remark}
We use the term $\infty$-sheaf as shorthand for simplicial presheaves that satisfy the homotopy descent condition.
We emphasize that we do not use quasicategories in the present paper.
Instead, we use model categories for all constructions.
\end{remark}

The strategy of the proof roughly goes as follows.
First, we show that it suffices to prove \cref{local} when the functor $\FFT_{d,T}$ is restricted to the full subcategory of representable presheaves in $\sPSh(\FEmb_d)$.
This follows by showing that the bordism category depends homotopy cocontinuously on the geometric structure at the level of \emph{presheaves}.
This part of the proof follows easily from our definition of the bordism category. 

In the case of representable presheaves, it turns out that our bordism category simplifies to a very concrete and easy to understand category of \emph{embedded} bordisms.
To prove \cref{local} for representables, we must show that for every bordism embedded in a fixed ambient manifold~$W$, and for every open cover of~$W$,
we can decompose the bordism into pieces in such a way that each piece in the decomposition is contained in some element of the open cover.
Decomposing bordisms into such pieces turns out to be quite delicate.
Furthermore, the entire space of such decompositions must be contractible (see \cref{contractible.cuts}), which is the most difficult part of the proof. 

This reduction to representable presheaves allows us to treat any bordism category that satisfies the following axioms. 
\begin{definition}
\label{axioms}
Fix $d\geq 0$.
Let $\gs$ be a geometric structure in the sense of \cref{geometric.structure}.
Let $\Bord_d^\gs$ be a corresponding object in $\smcat_{\infty,d}$ that satisfies the following axioms:
\begin{enumerate}
\item[\ax1] $\Bord_d^\gs$ is functorial in the geometric structure $\gs$. 
\item[\ax2] $\Bord_d^\gs$ depends homotopy cocontinuously on $\gs$ in the injective model structure on simplicial presheaves. 
\item[\ax3] Restricting to the subcategory of representable presheaves $W\to U$ (\cref{fembdef}),
the corresponding bordism category $\Bord_d^{W\to U}$ is naturally weakly equivalent to the category of embedded bordisms (\cref{embeddedbordcat.unenriched}).
\end{enumerate}
We call the functor $\gs\mapsto \Bord_d^\gs$ a \emph{smooth symmetric monoidal $(\infty,d)$-category of bordisms} (\cref{bordstr}).
\end{definition}

The main theorem thus naturally splits up into two parts.
\begin{theorem}[\cref{proofof104}]
\label{existencebord}
There exists a smooth symmetric monoidal $(\infty,d)$-category of bordisms.
\end{theorem}

\begin{theorem}[\cref{axioms.contractible.codescent}]
\label{axtheorem}
\cref{local} is true for any smooth symmetric monoidal $(\infty,d)$-category of bordisms. 
\end{theorem}

Breaking the main theorem up this way has the following advantage.
If, for whatever reason, the reader does not prefer our smooth symmetric monoidal $(\infty,d)$-category of bordisms,
the reader is free to replace this category by any category satisfying the axioms of \cref{axioms} and the main theorem will still hold.
In particular, this provides some flexibility in the proof.
The reader may be comforted to know that the structure of the proof is not overly rigid; it does not depend crucially on complicated definitions.
The axioms \ax1, \ax2, and \ax3 uniquely determine the smooth symmetric monoidal $(\infty,d)$-category of bordisms up to weak equivalence.
Hence any bordism category that satisfies the axioms is in fact equivalent to our bordism category.

In the sequel \cite{GradyPavlov.GCH} to this paper we make use of a related category of bordisms,
which we denote by $\BBord_d^\gs$, where $\gs$ is a geometric structure in the sense of \cref{geometric.structure.isotopy}.
This bordism category satisfies an analogous set of axioms as \cref{axioms}.

\begin{definition}
\label{multiple.isotopy.axiom}
Fix $d\geq 0$.
Let $\gs$ be a geometric structure in the sense of \cref{geometric.structure.isotopy}.
Let $\BBord_d^\gs$ be a corresponding object in $\fraksmcat_{\infty,d}$ that satisfies the following axioms.
\begin{enumerate}
\item[\frakax1] $\BBord_d^\gs$ is $\smsset$-enriched functorial in the geometric structure $\gs$. 
\item[\frakax2] $\BBord_d^\gs$ depends homotopy cocontinuously on $\gs$ in the injective model structure on simplicial presheaves. 
\item[\frakax3] Restricting to the subcategory of representable presheaves $W\to U$ (\cref{fembdelta}),
the corresponding bordism category $\BBord_d^{W\to U}$ is naturally weakly equivalent to the category of embedded bordisms with isotopies (\cref{embeddedbordcat.enriched}).
\end{enumerate}
We call the functor $\gs\mapsto \BBord_d^\gs$ a \emph{smooth symmetric monoidal $(\infty,d)$-category of bordisms with isotopies} (\cref{enrichedbordstr}). 
\end{definition}

In parallel with our proof of \cref{axtheorem}, we also prove the following.

\begin{theorem}[\cref{proofof106}]
\label{existencebord2}
There exists a smooth symmetric monoidal $(\infty,d)$-category of bordisms with isotopies.
\end{theorem}

\begin{theorem}[\cref{axioms.contractible.codescent}]
\label{axtheorem2}
\cref{local} is true for any smooth symmetric monoidal $(\infty,d)$-category of bordisms \emph{with isotopies}.
\end{theorem}

\cref{local} is closely related to the cobordism hypothesis of Baez–Dolan \cite{BaezDolan} and Lurie \cite{Lurie.TFT}, as witnessed by the following result.

\begin{theorem}[Geometric cobordism hypothesis, Part I]
\label{cob.part1}
Fix $d\geq 0$ and $\gs\in \Struct_d$.
Fix a fibrant $T\in \smcat_{\infty,d}$ in the model structure of \cref{globular.model.structure}.
Define the presheaf $T_d$ on $\FEmb_d$ with values in $\smcat_{\infty,d}$ by 
$$T_d(M\to U)≔\Funmonglob(\Bord_d^{M\to U},T).$$
Then $T_d$ is an $\infty$-sheaf and we have a weak equivalence of $\infty$-sheaves
\begin{equation}\label{cob1}\Funmonglob(\Bord_d^\gs,T)\simeq  \Funmonglob(\gs,T_d),\end{equation}
where we regard the presheaf of simplicial sets $\gs$, i.e., presheaf of $\infty$-groupoids, as a presheaf of $(\infty,d)$-categories.
The analogous statement for $\BBord_d$, with $\frakFEmb_d$ replacing $\FEmb_d$, also holds.
\end{theorem}

\begin{proof}
Write $\gs\in \sPSh(\FEmb_d)$ as a homotopy colimit of representables, use  \cref{local} and the definition of~$T_d$.
\end{proof}

To compare theorem \cref{cob.part1} (for $\BBord_d$) with the result of Lurie \cite{Lurie.TFT} in the special case of topological structures,
observe that every presheaf on $\frakFEmb_d$ canonically gives rise to an $\O(d)$-equivariant space
by evaluating a presheaf on the object given by the canonical map $(\RR^d\to \RR^0)\in \frakFEmb_d$, whose simplicial monoid of endomorphisms
is homotopy equivalent to (the singular complex of) $\O(d)$.
The $\O(d)$-action on this space is supplied by the $\O(d)$-action on $\RR^d$.
Natural transformations between presheaves on $\frakFEmb_d$ induce equivariant maps between the corresponding $\O(d)$-equivariant spaces, so the right side of \cref{cob1} can be thought of as a geometric refinement of the equivariant maps of spaces that appears in Lurie \cite{Lurie.TFT}.
See \scref{fiberwiseshape} for details.

Equation \cref{cob1} does not quite yield a full proof of the cobordism hypothesis,
because the target category~$T$ has changed to the somewhat mysterious sheaf of $(\infty,d)$-categories $T_d$.
In order to complete the proof of (the geometric refinement of) the cobordism hypothesis,
it remains to relate $T_d$ to the original target category~$T$.
More precisely, one must show that when $T$ has all duals,
there is an equivalence of $(\infty,d)$-categories $T_d(\RR^d\to \RR^0)\simeq T^{\times}$, where $T^{\times}$ denotes the core of~$T$. 
This also proves that the right side of \cref{cob1} is an $\infty$-groupoid (i.e., we can replace $\Funmonglob$ by the mapping simplicial set functor).
This is the focus of the sequel paper \cite{GradyPavlov.GCH} and the statement is only true in the case of the bordism category with isotopies $\BBord_d$.
We expand on this in the sequel.

In \cref{applications}, we offer some applications.
The locality property (\cref{local}) puts us in a position to apply the results of Berwick-Evans, Boavida de Brito, and the second author \cite{BEBdBP} (see also \cite{Pavlov.Diffeo})
to the $\infty$-sheaf $\FFT_{d,T}$. 
More precisely, we can construct (rather explicitly) a space $\csp\FFT_{d,T,\gs}$ such that the following holds. 

\begin{theorem}
\label{classcnc}
For all smooth manifolds $X$ and geometric structures $\gs$, we have a bijective correspondence 
$$\left\{\begin{array}{c}
\text{Concordance classes of }\\
\text{$d$-dimensional, $T$-valued }\\
\text{field theories with structure $\gs\times X$}
\end{array}
\right\}
\cong [X,\csp\FFT_{d,T,\gs}],$$
where on the right we take homotopy classes of maps.
\end{theorem}

This resolves a long-standing conjecture of Stolz and Teichner \cite{StolzTeichner.SUSY}
that postulates the existence of such a classifying space for extended functorial field theories with any geometry.
Our setup is sufficiently general to cover supersymmetric field theories required by the Stolz–Teichner program (see \cref{def.cartsp}).
For a more general version, which allows for twisted structures on $X$, see \cref{concordance.twisted}.
We also use this classifying space construction to provide a geometric model for power operations,
building on the work of Barthel, Berwick-Evans, and Stapleton \cite{BarthelBerwickEvansStapleton}.

In \cref{symminfn}, we introduce the notion of a smooth, symmetric monoidal $(\infty,d)$-category.
In \cref{bordcts}, we define geometric structures on bordisms and introduce two smooth variants of the bordism category, analogous to the distinction between bicategories and double categories.
In \cref{homotopy.theory}, we develop necessary tools from homotopy theory that will be used in the proof of \cref{local}.
In \cref{axioms.section}, we prove that the bordism categories $\Bord_{d}$ and $\BBord_{d}$ satisfy the axioms introduced above.
In \cref{codescent}, we prove that the various bordism categories we have introduced satisfy codescent, which is the main theorem of the paper.
In \cref{applications} we establish several applications. After reviewing background in \cref{shapebordism}, we construct an explicit classifying space of field theories (in the sense of \cref{classcnc}) in \cref{classifying.space}. In \cref{power.operations}, we construct a geometric model for power operations.

\subsection*{Prerequisites}

We assume familiarity with the language of simplicial homotopy theory and model categories, and, more specifically, the following topics.
\begin{itemize}
\item
Simplicial homotopy theory, including simplicial sets, simplicial maps, simplicial weak equivalences, Quillen's Theorem~A,
Kan's $\Exi$ functor, and the simplicial Whitehead theorem.
See Goerss–Jardine \cite{GoerssJardine} and Dugger–Isaksen \cite{DuggerIsaksen}.
\item
Simplicial categories, Dwyer–Kan weak equivalences, the homotopy coherent nerve functor and its left adjoint.
See Bergner \cite[Sections 4 and 7]{Bergner}.
\item
Model categories, including model structures, left Quillen functors, injective and projective model structures on presheaves with values in a model category,
Reedy model structures, combinatorial model categories, left Bousfield localizations.
See Hovey \cite{Hovey}, Hirschhorn \cite{Hirschhorn}, Barwick \cite{Barwick.Model},
as well as Lurie \cite[Appendix~A]{Lurie.HTT} and the survey of Balchin \cite{Balchin}.
\item
Homotopy limits and colimits.
See Bousfield–Kan \cite{BousfieldKan}, Hirschhorn \cite{Hirschhorn}, Shulman \cite{Shulman}, Riehl \cite{Riehl}.
\item
Segal's $Γ$-objects.
See Segal \cite{Segal.Gamma}, Bousfield–Friedlander \cite{BousfieldFriedlander}, Schwede \cite{Schwede}.
\item
Rezk's complete Segal spaces \cite{Rezk.CSS} and their generalization to $n$-fold Segal spaces by Barwick \cite{Barwick.CSS}.
\item
Simplicial presheaves and descent.
See Dugger–Hollander–Isaksen \cite{DHI} and Jardine \cite{Jardine}.
\item
Elementary Morse theory.
See Milnor \cite{Milnor}.
\end{itemize}

\subsection*{Acknowledgments}

We thank Peter Teichner and Stephan Stolz for many discussions about bordism categories and field theories.
We thank Urs Schreiber, Hisham Sati, and Vincent Braunack-Mayer for discussions on the physics applications and for providing feedback on the first version of this article. 
We thank André Henriques, Thomas Nikolaus, Christopher Schommer-Pries, and David Reutter for discussions about this paper.
We thank Daniel Brügmann, Nino Scalbi, and Alexander Zahrer for a careful reading of this paper.
We thank Nils Carqueville, Domenico Fiorenza, and Konrad Waldorf for organizing a workshop on this paper and \cite{GradyPavlov.GCH}, which led to many improvements in the exposition.

\clearpage
\section*{Table of notation}
\begin{center}
{\small
\begin{tabular}{c|c|l}
Notation & Definition & Description\\
\hline
$\Delta^{\times d}$ & \cref{threecats} & The $d$-fold product of the simplex category\\
$\Gamma$ & \cref{threecats} & Segal's Gamma category\\
$\cart$ & \cref{def.eucl} & Smooth cartesian spaces (e.g., $\RR^n$, for some $n\in \NN$)\\
$\man$ & \cref{def.cartsp} & Smooth manifolds\\
$\stcart$ & \cref{def.cartsp} & Structured cartesian spaces (e.g., $\RR^n$, $\RR^{n\mid m}$, for some $n,m\in \NN$)\\
$\stman$ & \cref{def.cartsp} & Structured manifolds\\
$\stu$ & \cref{def.cartsp} & The reduction functor $\stman→\man$\\
$\PSh(\mathscr{C},{\sf V})$ & \cref{presheaf.notation} & Presheaves on a small category $\mathscr{C}$ with values in ${\sf V}$\\
$\PSh(\mathscr{C})$ & \cref{presheaf.notation}  & Presheaves of sets on a small category $\mathscr{C}$\\
$\sPSh(\mathscr{C})$ & \cref{presheaf.notation}  & Simplicial presheaves on a small category $\mathscr{C}$\\
$⊗$, $\hom$ & \cref{presheaf.notation} & Day convolution and internal hom on presheaves\\
$-⊠-$& \cref{external.product} & External product of presheaves\\
$\Yo{}$& \cref{external.product} & The Yoneda embedding\\
$\smcatuple_{\infty,d}$ & \cref{multiple.model.structure} &  $\sPSh(\stcart\times \Gamma\times \Delta^{\times d})$, with the uple model structure\\
$\smcat_{\infty,d}$ & \cref{globular.model.structure} &$\sPSh(\stcart\times \Gamma\times \Delta^{\times d})$, with the globular model structure\\
$\catuple_{\infty,d}$ & \cref{globular.multiple.model.structure} &  $\sPSh(\Gamma\times \Delta^{\times d})$, with the uple model structure\\
$\cat_{\infty,d}$ & \cref{globular.multiple.model.structure} &$\sPSh(\Gamma\times \Delta^{\times d})$, with the globular model structure\\
$⊗$, $\Funmonuple$& \cref{def.funmonuple} & (noncartesian) tensor product and internal hom in $\smcatuple$\\
$\Theta_d$ & \cref{functor.categories} & Joyal's cell category\\
$\tglob$, $\Funmonglob$& \cref{globular.hom} & Globular product and functor object in $\smcat_{∞,d}$\\
$\smset$ & \cref{smsets} & Presheaves of sets on $\cart$\\
$\smsset$ & \cref{smsets} & Simplicial presheaves on $\cart$, with the cartesian $\RR$-local model structure\\
$\smsPSh(\mathscr{C})$ & \cref{smsets} & $\smsset$-valued presheaves on $\mathscr{C}$ with the injective model structure\\
$\fraksmcatuple_{\infty,d}$  & \cref{multiple.model.structure.smooth} &$\smPSh(\stcart\times \Gamma\times \Delta^{\times d})$, with the uple model structure\\
$\fraksmcat_{\infty,d}$  & \cref{globular.model.structure.smooth} &$\smPSh(\stcart\times \Gamma\times \Delta^{\times d})$, with the globular model structure\\
$B^{\mathfrak{p}}_ε$, $\mathfrak{p}$ & \cref{stalks} & A ball of radius~$ε$ inside an object of $\stcart$ indexing a point~$\mathfrak{p}$\\
$\Cnec$ & \cref{Cnec.def} & The Dugger–Spivak simplicial category of necklaces\\
$\FEmb_d$ & \cref{fembdef} & Submersions with $d$-dimensional fibers and fiberwise embeddings\\
$\Struct_d$ & \cref{geometric.structure} & $\sPSh(\FEmb_d)$, with the Čech-local injective model structure\\
$\frakFEmb_d$ & \cref{fembdelta} & $\FEmb_d$ enriched in $\smset$\\
$\frakStruct_d$ & \cref{geometric.structure.isotopy} & $\smsPSh(\frakFEmb_d,\smsset)$, with the Čech-local injective model structure\\
$\Bord_{d,\uple}$ & \cref{bord} & The $d$-uple category of bordisms\\ 
$\Bord_d$ & \cref{bord} & The globular $d$-category of bordisms\\ 
$\Bord_{d,\uple}^\gs$ & \cref{bordstr} & The $d$-uple category of bordisms with geometric structure~$\gs$\\ 
$\Bord_d^\gs$ & \cref{bordstr} & The globular $d$-category of bordisms with geometric structure~$\gs$\\ 
$\BBord_{d,\uple}^\gs$ & \cref{enrichedbordstr} & The $d$-uple category of bordisms with geometric structure~$\gs$ and isotopies\\ 
$\BBord_d^\gs$ & \cref{enrichedbordstr} & The globular $d$-category of bordisms with geometric structure~$\gs$ and isotopies\\ 
$\core(M,C,P)$ & \cref{bord} & The core of a bordism $(M,C,P)$\\
$\EBord_d^p$ & \cref{embeddedbordcat.unenriched} & Bordisms embedded fiberwise in $(p:W\to U)\in \FEmb_d$\\
$\frakEBord_d^p$ & \cref{embeddedbordcat.enriched}& Bordisms with isotopies embedded fiberwise in $(p:W\to U)\in \frakFEmb_d$\\
$\core[f,N,C,P]$ & \cref{bracketcore}& The embedded core of an embedded bordism $[f,N,C,P]$\\
$\frake_d^p$ & \cref{emapconstr}& The comparison map for bordisms and bordisms embedded into $\mathfrak{p}$\\
$d$ & \cref{filtration.notation} & The dimension of bordisms\\
$\mathfrak{p}$ & \cref{filtration.notation} & A point of $\stcart$, indexing a stalk\\
$\mathfrak{l}$ & \cref{filtration.notation} & A point of $\cart$, indexing a stalk\\
$p:W\to U$ & \cref{filtration.notation} & An object of $\frakFEmb_d$\\
${\cal W}=\{p_a\}_{a∈A}$ & \cref{filtration.notation} & A covering family of~$p$ in $\frakFEmb_d$\\
$\frakB$ & \cref{filtration.notation} & The stalk $\mathfrak{l}^*\mathfrak{p}^*\frakEBord_d^p$\\
$\frakB_{-1}$ & \cref{filtration.notation} & The colimit over the Čech nerve of ${\cal W}$: $\colim_{[n]\in \Delta^\op}\coprod_{\alpha:[n]\to A}\mathfrak{l}^*\mathfrak{p}^*\frakEBord_d^{p_\alpha}$\\
$\frakB_i$ & \cref{subcov} & The filtration on $\frakB$\\
$⟨\ell⟩$ & \cref{simplex.notation} & An object of $Γ$\\
$i$ & \cref{simplex.notation} & An integer indexing a direction in $Δ^{⨯d}$\\
${\bf m}$ & \cref{simplex.notation} & A multisimplex in $\Delta^{\{1,\hdots, \hat \imath,\hdots d\}}=\Delta^{\times d-1}$\\
$x$, $y$ & \cref{simplex.notation} & Elements of $\frakB_i(\langle \ell\rangle,{\bf m})([0])=\frakB_{i-1}(\langle \ell\rangle,{\bf m})([0])$\\
$\frakB^i_{i-1}$, $\frakB^i_i$ & \cref{simplex.notation} & A presheaf of sets on~$Δ$: $\frakB_j(\langle \ell\rangle,{\bf m})$\\
$\frakB_{i-1}^{i,x,y}$, $\frakB_i^{i,x,y}$ & \cref{etalemapsnow} & The category of necklaces from~$x$ to~$y$ in $\frakB^i_{i-1}$ or $\frakB^i_i$\\
\end{tabular}
}
\end{center}
\clearpage

\section{Smooth symmetric monoidal $(\infty,d)$-categories}
\label{symminfn}

In this section we develop a model for smooth symmetric monoidal $(\infty,d)$-categories.
The encoding of symmetric monoidal structures follows Segal \cite{Segal.Gamma}.
The encoding of categorical composition follows the same idea of Segal, as developed by Rezk \cite{Rezk.CSS} (in the $(\infty,1)$-categorical case) and Barwick \cite{Barwick.CSS} (in the $(\infty,d)$-categorical case).
The encoding of smoothness follows Stolz–Teichner \cite{StolzTeichner.Elliptic, StolzTeichner.SUSY}.
The use of the cartesian site in this context was proposed by Urs Schreiber (see, e.g., Fiorenza–Schreiber–Stasheff \cite{FiorenzaSchreiberStasheff}). 

\subsection{Homotopy theory and higher categories}
\label{threecats}

Recall the simplex category~$\Delta$, whose objects are nonempty, totally ordered sets $[d]=\{0,\ldots,d\}$, and whose morphisms are order-preserving maps. 
Recall also Segal's category~$Γ$ \cite{Segal.Gamma} (see also Bousfield–Friedlander \cite{BousfieldFriedlander} and Schwede \cite{Schwede}),
whose opposite category has objects which are finite pointed sets $\langle \ell\rangle=\{\ast,1,\ldots,\ell\}$ and basepoint-preserving functions as morphisms.
Finally, recall the category $\cart$ (see, e.g., Fiorenza–Schreiber–Stasheff \cite[Definition~3.1.1]{FiorenzaSchreiberStasheff})
whose objects are cartesian spaces~$\RR^n$ and morphisms are smooth functions between them.
More generally, we can use any category $\stcart$ that satisfies certain properties making it similar to $\cart$ (\cref{def.cartsp}),
e.g., we could take $\stcart=\cart$, other choices are discussed in \cref{def.supercart,def.infinitesimal}.
The basic building blocks for smooth symmetric monoidal $(\infty,d)$-categories lie in the threefold product category
$$\stcart\times \Gamma\times \Delta^{\times d}.$$
One should think of each category as capturing an independent nature of the objects we are constructing. 
We provide some indications of how to think of each piece conceptually. 

\subsubsection{\bf The category $\Delta^{\times d}$}
Here an object is a multisimplex ${\bf m}=([m_1],[m_2],\ldots, [m_d])\in \Delta^{\times d}$,
which we can think of as indexing a composable grid of morphisms, with $m_i$~composed morphisms in the $i$th direction.
For example, for $d=2$, we have two composition directions.
If $m_1=m_2=4$, this gives a grid
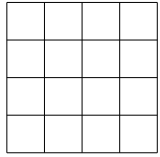
\begin{center}
\def\xMin{0}%
\def\xMax{4}%
\def\yMin{0}%
\def\yMax{4}%
\begin{tikzpicture}[scale=.5]
\foreach \i in {\xMin,...,\xMax}{\draw [very thin] (\i,\yMin) -- (\i,\yMax);}
\foreach \i in {\yMin,...,\yMax}{\draw [very thin] (\xMin,\i) -- (\xMax,\i);}
\end{tikzpicture}
\end{center}
A single square represents a 2-morphism.
The edges represent 1-morphisms and vertices represent objects.
These 2-morphisms can be composed in two directions, giving rise to a grid.
One can construct a model structure on simplicial presheaves on $\Delta^{\times d}$, whose fibrant objects are $d$-fold Segal spaces,
see \cref{multiple.model.structure} and \cref{globular.model.structure} for a construction including symmetric monoidal and smooth structures.
For an expository account of $d$-fold Segal spaces, see
Barwick \cite{Barwick.CSS}, Lurie \cite{Lurie.TFT}, Barwick–Schommer-Pries \cite{BarwickSchommerPries}, Calaque–Scheimbauer \cite{CalaqueScheimbauer}.

\subsubsection{\bf The category $\Gamma$} 
This category captures the symmetric monoidal structure. 
An object $\langle\ell\rangle\in\Gamma$ is simply a finite set $\{\ast,1,2,\ldots,\ell\}$ and a morphism is just a function in the opposite direction satisfying $f(\ast)=\ast$. 
This category encodes the symmetric monoidal structure as follows.
Consider a $Γ$-object $X:Γ^\op→\sset$.
The function $\phi_{⟨\ell⟩}:\langle\ell\rangle \to \langle1\rangle$ sending $i$ to~$1$ (for all $i$ such that $1\leq i\leq \ell$),
gives, by functoriality, a map
$$X(\phi_{⟨\ell⟩}):X\langle\ell\rangle→X\langle1\rangle.$$
We also have maps $\delta_i:\langle \ell \rangle\to \langle 1\rangle$ that send $i$ to $1$ and $j$ to $*$ for $j\neq i$. 
They induce a map $$δ_{\langle \ell\rangle}:X\langle \ell\rangle \to X\langle 1\rangle^{\times \ell}.$$
The multiplicative structure is given by the following zigzag, where the left leg is a weak equivalence, so can be (formally) inverted:
$$X\langle 1\rangle^{\times \ell}\lgets5{\delta_{⟨\ell⟩}} X\langle \ell\rangle \lto5{X(\phi_{⟨\ell⟩})} X\langle 1\rangle.$$
One can think of $X\langle\ell\rangle$ as the space of $\ell$-tuples that can be multiplied.
The map $X(\phi_{⟨\ell⟩})$ then performs the multiplication.
The map $δ_{⟨\ell⟩}$ extracts the individual components of an $\ell$-tuple
and is a weak equivalence because any $\ell$-tuple can be deformed to an $\ell$-tuple that can be multiplied.
The elegance of Segal's method is that it circumvents the need for explicitly keeping track of coherent homotopies in the symmetric monoidal structure,
although these are of course still hidden in the above equivalence.
For more information, see Segal \cite{Segal.Gamma}.

\subsubsection{\bf The category $\stcart$}
The categories $\cart$ and, more generally, $\stcart$ capture the smooth structure.

\begin{definition}
\label{def.eucl}
The site $\cart$ is defined as follows.
Objects are open subsets of $\RR^n$ diffeomorphic to $\RR^n$, where $n≥0$ is arbitrary.
Morphisms are smooth maps.
The covering families are given by {\it good open covers}, defined as open covers in which every finite intersection is empty or diffeomorphic to some $\RR^n$.
\end{definition}

The functor category $\Fun(\cart^\op,\mathscr{D})$ can be thought of as the category of smoothly parametrized objects of $\mathscr{D}$ over cartesian spaces.
As a basic example, one can consider a smooth manifold $M$ as a functor $M:\cart^\op\to \set$, 
by sending a cartesian space~$U$ to the set of smooth maps $\sm(U,M)$. 
The category $\cart$ allows us to form
smooth families of symmetric monoidal $(\infty,d)$-categories
as well as isotopies of bordisms (\cref{isotopybordsec}).
For more information about sheaves on this site, see Fiorenza–Schreiber–Stasheff \cite{FiorenzaSchreiberStasheff} or Sati–Schreiber \cite{SatiSchreiber.POC}.

For applications (including the Stolz–Teichner program), we find it convenient to generalize the site $\cart$ slightly and allow for “enhancements”
such as those provided by various flavors of supermanifolds and smooth infinitesimal analysis. We will denote this generalized site by $\stcart$.
Of course, the site $\cart$ of cartesian spaces remains our main example throughout the text.
The following definition follows Pavlov \cite[Definition~1.5]{Pavlov.Numerable}.
It is given for arbitrary manifolds, since we need manifolds to encode bordisms,

\begin{definition}
\label{def.cartsp}
(See also Pavlov \cite[Definition~1.5]{Pavlov.Numerable}.)
Throughout the paper, $\stman$ denotes a fixed small subcanonical site with finite products
equipped with a \emph{reduction functor} $\stu:\stman→\man$ that preserves and reflects covering families,
preserves finite products,
and admits cartesian lifts for all open embeddings in $\man$,
and an \emph{embedding functor} $ι:\man→\stman$ (typically omitted from the notation) that preserves finite products
and satisfies $\stu ι=\id_\man$.
We say that a morphism in $\stman$ is an open embedding if it is a cartesian morphism and its image under the functor~$\stu$ is an open embedding.
We say that a morphism in $\stman$ is a local diffeomorphism if it is a cartesian morphism and its image under the functor~$\stu$ is a local diffeomorphism.
We say that a morphism $f:X→Y$ in $\stman$ is a submersion with $d$-dimensional fibers
if locally on~$X$ and~$Y$ the morphism~$f$ is isomorphic to a projection map of the form $ι(\RR^d)⨯U→U$ for some $U∈\stcart$.

We define a full subcategory $\stcart⊂\stman$ by taking objects $X∈\stman$ such that $\stu(X)∈\man$ is isomorphic to an object in $\cart$.
\end{definition}

A canonical example of $\stman$ is given by $\stman=\man$, with $\stu$ and $ι$ being the identity functors.
Since covering families are created by the functor~$\stu$, we can omit their description from the examples below.

\begin{example}
\label{def.supercart}
The site of supermanifolds with the forgetful functor~$\stu$ being the reduction functor and $ι$ being the inclusion functor satisfies the conditions of \cref{def.cartsp}.
\end{example}

\begin{example}
\label{def.infinitesimal}
As a variation on the previous example, we could also take various sites from smooth infinitesimal analysis, e.g., the Dubuc site.
That is, supermanifolds can be replaced by the opposite category of the full subcategory of commutative real algebras
isomorphic to the tensor product of a Weil algebra and the algebra of functions on some smooth manifold.
See Moerdijk–Reyes \cite{MoerdijkReyes}.
\end{example}

\subsubsection{\bf Presheaves}
Throughout this paper, we will make use of internal homs and simplicial enrichments of presheaves on monoidal categories
(like $\stcart$ or $Δ$ with its cartesian structure or $Γ$ with its smash product structure).
We introduce the following systematic notation.

\begin{notation}
\label{presheaf.notation}
Throughout the remainder of the paper,
presheaves with values in another category ${\sf V}$ will be denoted by
$$\PSh(\mathscr{C},{\sf V})≔\Fun(\mathscr{C}^\op,{\sf V}).$$
If ${\sf V}$ is a closed symmetric monoidal category,
then $\PSh(\mathscr{C},{\sf V})$ is tensored, powered, and enriched over ${\sf V}$.
(In many cases, the symmetric monoidal structure on ${\sf V}$ is the cartesian monoidal structure.)
Furthermore, if $\mathscr{C}$ is equipped with a symmetric monoidal structure,
then $\PSh(\mathscr{C},{\sf V})$ is equipped with a closed symmetric monoidal structure given by the Day convolution
(Day \cite{Day}, see also the $n$Lab \cite[Section~3]{nLab.Day}).
We denote the Day convolution product by $\otimes$ and the internal hom by $\hom(-,-)$. 
We also write
$$\PSh(\mathscr{C})=\PSh(\mathscr{C},\set),\qquad \sPSh(\mathscr{C})=\PSh(\mathscr{C},\sset)$$
and refer to the latter as {\it simplicial presheaves\/} on~$\mathscr{C}$.
We equip $\sPSh(\mathscr{C})$ with
the injective model structure, turning it into a simplicial model category $\sPSh(\mathscr{C})_\inj$ (Lurie \cite[Remark A.3.3.4]{Lurie.HTT}). 
\end{notation}

\begin{definition}
\label{external.product}
Given small categories $\mathscr{C}$ and $\mathscr{D}$ and a cocomplete closed symmetric monoidal category ${\sf V}$ (such as $\mathscr{D}=\sset$), we denote by
$$⊠:\PSh(\mathscr{C},{\sf V})⨯\PSh(\mathscr{D},{\sf V})→\PSh(\mathscr{C}\times \mathscr{D},{\sf V})$$
the unique ${\sf V}$-enriched separately cocontinuous functor that sends $(\Yo{a},\Yo{b})$ to $\Yo{(a,b)}$,
where $\Yo{}$ is the Yoneda embedding.
\end{definition}

\subsection{Simplicial presheaves and their left Bousfield localizations}
\label{smooth.multiple.categories}

We will present the $(\infty,1)$-category of sheaves of symmetric monoidal $(\infty,d)$-categories
by a particular model structure,
following Barwick \cite{Barwick.CSS} in the case of $\Delta^{\times d}$
and Toën \cite{Toen} for symmetric monoidal $(∞,d)$-categories.
We begin by taking the injective model structure on $\sPSh(\stcart\times \Gamma\times \Delta^{\times d})$.
To encode the Segal conditions for $Δ$ and $Γ$, completeness conditions for~$Δ$, and descent conditions for $\stcart$ on these objects,
we perform a left Bousfield localization of this model category.

\begin{definition}
\label{left.Bousfield.localization}
(Hirschhorn \cite[Definition~3.1.1.(1)]{Hirschhorn}, Barwick \cite[Definition~4.2]{Barwick.Model}, Lurie \cite[§A.3.7]{Lurie.HTT}.)
Suppose $M$ is a model category and $S$ is a set of morphisms in~$M$.
The {\it left Bousfield localization\/} of~$M$ at~$S$ is a model category $L_S M$ together with a left Quillen functor $F:M→L_S M$
that satisfies the following universal property:
\begin{itemize}
\item Composing with~$F$ maps left Quillen functors $L_S M→N$
bijectively to left Quillen functors $G:M→N$ such that the left derived functor of~$G$ sends elements of $S$ to weak equivalences in~$N$.
\end{itemize}
An object $X∈M$ is {\it $S$-local\/} if for any morphism $f:Y→Z$ in~$S$ the induced map $\rmap(Z,X)→\rmap(Y,X)$ is a weak equivalence.
A morphism $f:Y→Z$ is an {\it $S$-local equivalence\/} if for any $S$-local object~$X$ the induced map $\rmap(Z,X)→\rmap(Y,X)$ is a weak equivalence.
Here $\rmap(-,-)$ denotes the homotopy function complex in a model category, see Hirschhorn \cite[Definition~17.4.1]{Hirschhorn}.
\end{definition}

\begin{proposition}
\label{left.Bousfield.localization.exists}
(Barwick \cite[Theorem~4.7]{Barwick.Model}.)
Suppose $M$ is a left proper combinatorial model category and $S$ is a set of morphisms in~$M$.
Then the left Bousfield localization (\cref{left.Bousfield.localization}) of~$M$ at~$S$ exists and is a left proper combinatorial model category. Moreover:
\begin{itemize}
\item its underlying category is $M$;
\item its weak equivalences are precisely the $S$-local equivalences;
\item its cofibrations and acyclic fibrations coincide with those of~$M$;
\item its fibrant objects are precisely the $S$-local fibrant objects of~$M$;
\item if $M$ is simplicial, then so is $L_S M$;
\item if $M$ is $V$-enriched, then $L_S M$ is $V$-enriched if and only if the derived monoidal product of an element of~$S$
with an object of~$V$ (which can be taken from any class of homotopy generators of~$V$) is a weak equivalence in $L_S M$.
The case $V=M$ yields a criterion for monoidal model structures.
(Barwick \cite[Proposition~4.47]{Barwick.Model}, Gorchinskiy–Guletskiĭ \cite[Lemma~31 (journal), Lemma~28 (arXiv v4)]{GorchinskiyGuletskii}).
\end{itemize}
\end{proposition}

\begin{definition}
\label{bousloc}
Fix $d\geq 0$.
We define the following morphisms in $\sPSh(\Delta)$, $\sPSh(\Gamma)$, and $\sPSh(S)$ (where $S$ is a site such as $\stcart$, $\cart$, $\FEmb_d$),
which will be used later to construct a left Bousfield localization. 
\begin{enumerate}
\item[(i)] (Segal's special $Δ$-condition.)
For $[a],[b]\in \Delta$, the maps 
\begin{equation}\label{segal1}\phi^{a,b}:\Yo{[a]}\sqcup_{\Yo{[0]}}\Yo{[b]}\to \Yo{[a+b]},\end{equation}
where the pushout uses the morphism $[0]→[a]$ that picks out the terminal vertex and $[0]→[b]$ that picks out the initial vertex.
Local objects are precisely Segal's special $Δ$-objects~$X$, defined by the condition
that the map $X_{a+b}→X_a⨯_{X_0}X_b$ is a weak equivalence.
\item[(ii)] (Completeness condition.)
The map
\begin{equation}\label{segal2}x:E\to \Yo{[0]}\end{equation}
where $E$ is obtained by evaluating the functor
$$\sset=\PSh(Δ)→\sPSh(Δ)$$
on the nerve of the groupoid with two objects $p$ and~$q$ and two nonidentity morphisms $p\to q$ and $q\to p$,
where the functor converts presheaves of sets into presheaves of discrete simplicial sets.
Local objects with respect to (i) and (ii) are Rezk's complete Segal objects,
defined by Segal's special $Δ$-condition
and the completeness condition, which says that the map $X_0→X_1$ that sends objects to their identity morphisms
is a weak equivalence onto the subobject of invertible 1-morphisms.
\item[(iii)] (Segal's special $Γ$-condition.)
For $\langle \kappa\rangle,\langle \ell\rangle\in\Gamma$, the maps
\begin{equation}\label{monoidal1}t^{\kappa,\ell}:\Yo{\langle \kappa\rangle}\sqcup_{\Yo{\langle 0\rangle}}\Yo{\langle\ell\rangle}\to \Yo{\langle \kappa+\ell\rangle}.\end{equation}
Also, the map $$τ:∅→\Yo{\langle0\rangle},\eqlabel{monoidal2}$$
which forces the 0th space to be contractible.
Local objects are precisely Segal's special $Γ$-objects,
for which the maps $X_{\kappa+\ell}→X_\kappa⨯X_\ell$ and $X_0→*$ are weak equivalences.
\item[(iv)] (Sheaf condition.)
Suppose $S$ is a site, e.g., $S=\stcart$ (\cref{def.cartsp}) with its Grothendieck topology of good open covers,
other examples of interest include $\cart$ (\cref{def.eucl}), $\FEmb_d$ (\cref{fembdef}), $\frakFEmb_d$ (\cref{fembdelta}).
For ${\cal U}=\{U_\alpha\}_{α∈I}$ a covering family of $V$ in $S$,
the Čech nerve $c_{\cal U}∈\sPSh(S)$
has as its presheaf of $m$-simplices the presheaf of sets
$$\coprod_{α_0,\ldots,α_m: U_{α_0}∩⋯∩U_{α_m}≠∅} \Yo{U_{α_0}∩⋯∩U_{α_m}}.$$
Here we abbreviate
$$U_{α_0}∩⋯∩U_{α_m} = U_{α_0} ⨯_V ⋯ ⨯_V U_{α_m},$$
which boils down to the usual intersection for sites of interest to us.
The $k$th face map deletes $U_{α_k}$, the $k$th degeneracy map duplicates $U_{α_k}$.
We have a canonical map $$i^{{\cal U},V}:c_{\cal U}→\Yo{V}\eqlabel{cech}$$
that in simplicial degree~$m$ is induced by the inclusions $U_{α_0}∩⋯∩U_{α_m}→V$.
The local objects for the maps $i^{{\cal U},V}$ are ∞-sheaves (alias ∞-stacks), for which
the restriction map 
\begin{equation}\label{homotopy.descent}X(V)→\holim_α X(U_{α_0}∩⋯∩U_{α_m})\end{equation}
is a weak equivalence.
For more information on the Čech descent, see Dugger–Hollander–Isaksen \cite[Appendix~A]{DHI}.
\end{enumerate}
\end{definition}

\begin{remark}
Although it may seem as if the nature of the two maps in (iii) are different, both of them arise from inverting maps of the form
$$\Yo{\langle s_1\rangle}⊔\Yo{\langle s_2\rangle}⊔⋯⊔\Yo{\langle s_m\rangle}→\Yo{\langle s_1+⋯+s_m\rangle}.$$
Taking $m=0$ produces the map~$τ$.
Taking $m=2$ produces a map weakly equivalent to the map $t^{s_1,s_2}$ because the object $\Yo{\langle0\rangle}$ in the pushout was already forced to be contractible by the map~$τ$.
In fact, it is clear that we can alternatively localize at maps of the above form, which gives a somewhat more uniform description.
\end{remark}

\begin{definition}
\label{gammasp}
The model structure on $\Gamma$-spaces is given by taking the injective model structure on $\sPSh(\Gamma)$ and performing left Bousfield localization at the morphisms \cref{monoidal1} and \cref{monoidal2}. We denote this model structure by $\sPSh(\Gamma)_{\local}$.
\end{definition}

This model structure on $\Gamma$-spaces has more cofibrations than the Bousfield–Friedlander model structure \cite{BousfieldFriedlander}, which, in turn, has more cofibrations than Schwede's Q-model structure \cite{Schwede}.
All three model structures have the same weak equivalences.

\begin{definition}
\label{rezkdelta}
The Rezk model structure \cite{Rezk.CSS} on $Δ$-spaces is
given by taking the Reedy model structure on $\sPSh(\Delta)$, which coincides with the injective model structure, and performing left Bousfield localization at the morphisms \cref{segal1} and \cref{segal2}. We denote this model category by $\sPSh(\Delta)_{\local}$.
\end{definition}

\subsection{Smooth multiple and globular $(\infty,d)$-categories}

To combine symmetric monoidal structure, higher categorical structure, and smooth structure, we need to promote the morphisms in \cref{bousloc} to morphisms in the category of simplicial presheaves on the threefold product $\stcart\times \Gamma\times \Delta^{\times d}$. We do this by using the external product~$⊠$ of \cref{external.product}. The resulting model category has fibrant objects that are local with respect to each one of the morphisms in \cref{bousloc}
after evaluating on an arbitrary fixed choice of all objects in factors other than the one under consideration. 
For instance, fixing an arbitrary object in $Γ⨯\Delta^{\times d}$ should yield an ∞-sheaf of simplicial sets on $\stcart$,
meaning the homotopy descent condition \cref{homotopy.descent} is satisfied.

\begin{notation}
\label{localmapsnot}
Fix $d\geq 0$.
We consider the following sets of maps in $\sPSh(\stcart\times\Gamma\times\Delta^{\times d})$. 
\begin{itemize}
\item Let $S_\Delta=\{\Yo{c}⊠\phi^{a,b}, \Yo{c}⊠x\mid c\in\stcart\times\Gamma\times Δ^{\{1,\ldots,k-1,k+1,\ldots,d\}}, 1≤k≤d, a,b∈Δ\}$,
where $\phi^{a,b}$ are the morphisms \cref{segal2} in \cref{bousloc} and $x$ is the morphism \cref{segal1}. 
\item Let $S_\Gamma=\{\Yo{c}⊠t^{\kappa,\ell}, \Yo{c}⊠\tau\mid c\in\stcart\times\Delta^{\times d}, \kappa,\ell∈Γ\}$,
where $t^{\kappa,\ell}$ are the morphisms \cref{monoidal1} and $\tau$ is the morphism \cref{monoidal2}.
\item Let $S_\stcart=\{\Yo{c}⊠i^{{\cal U},V}\mid c\in\Gamma\times\Delta^{\times d}, V∈\stcart\}$,
where $i^{{\cal U},V}$ are the morphisms \cref{cech}, with ${\cal U}$ ranging over good covers of each~$V$.
\end{itemize}
\end{notation}

\begin{definition}[The multiple model structure]
\label{multiple.model.structure}
We define a model structure 
$$\smcatuple_{\infty,d}≔\sPSh(\stcart\times \Gamma\times \Delta^{\times d})_{\uple}$$
by performing the left Bousfield localization (\cref{left.Bousfield.localization}) of the model category $\sPSh(\stcart\times \Gamma\times \Delta^{\times d})_{\inj}$
at the set of morphisms $S_{\stcart}\cup S_{\Gamma}\cup S_{\Delta}$ (\cref{localmapsnot}).
The localization exists by \cref{left.Bousfield.localization.exists}.
\end{definition}

\begin{proposition}
\label{multiple.monoidal}
The model structure of \cref{multiple.model.structure} is a monoidal model structure.
\end{proposition}

\begin{proof}
According to the criterion for monoidality stated in \cref{left.Bousfield.localization.exists},
it suffices to show that the monoidal product of
a simplicial presheaf~$F$ on
$\stcart\times \Gamma\times \Delta^{\times d}$
and one of the morphisms in $S_{\stcart}\cup S_{\Gamma}\cup S_{\Delta}$
is a weak equivalence in $\smcatuple_{\infty,d}$.
This monoidal product is derived because all objects in the model category under consideration are cofibrant.
The simplicial presheaf~$F$ can be assumed to be representable because the monoidal product is separately homotopy cocontinuous in each argument
and $S$-local weak equivalences are closed under homotopy colimits in the original model category.

Observe that replacing the representable presheaf~$\Yo{c}$ in \cref{localmapsnot} with an arbitrary (nonrepresentable) simplicial presheaf~$F$
on the same category produces a weak equivalence in $\smcatuple_{\infty,d}$,
since simplicial presheaves can be presented as homotopy colimits of representable presheaves
and $S$-local weak equivalences in~$C$ are closed under homotopy colimits in~$C$.
More formally, we can say that the external product functor $F⊠-$
is a left Quillen functor
originating in the category of simplicial presheaves on $\stcart$ respectively $Γ$ respectively~$Δ$
(left Bousfield localized at the corresponding morphisms in \cref{bousloc})
and landing in the model category $\smcatuple_{\infty,d}$.

Now recall that morphisms in $S_{\stcart}\cup S_{\Gamma}\cup S_{\Delta}$ have the form $\Yo{c}⊠f$,
where $f$ is one of the maps of \cref{bousloc} (i.e., a morphism of simplicial presheaves on $\stcart$, $Γ$, or $Δ$ respectively),
whereas $c$ is an object in the product of remaining factors,
as explained in \cref{localmapsnot}.
Taking the monoidal product $(\Yo{c}⊠f)⊗\Yo{G}$ of the morphism $\Yo{c}⊠f$ and the representable presheaf of an arbitrary object $G∈\stcart⨯Γ⨯Δ^{⨯d}$
produces a morphism of the form $(\Yo{c}⊗\Yo{G_1})⊠(f⊗\Yo{G_2})$, where $G_1$ and $G_2$ denote the corresponding projections of~$G$.
Since $(\Yo{c}⊗\Yo{G_1})⊠-$ is a left Quillen functor, it remains to show that $f⊗\Yo{G_2}$ is a weak equivalence.

Thus, the problem reduces to showing
that the category of simplicial presheaves on $\stcart$ respectively $Γ$ respectively~$Δ$
left Bousfield localized at the morphisms of \cref{bousloc} is a monoidal model category. 
For $\stcart$, observe that the (cartesian) product of the representable presheaf of $W∈\stcart$ with a map \cref{cech}
$$i^{{\cal U},V}:c_{\cal U}→\Yo{V}$$
is again a map of the same form, for the covering family $\{U_α⨯W\}_{α∈I}$ of $V⨯W$ in $\stcart$. 
For $Γ$, observe that the (smash) product of $⟨r⟩∈Γ$ with a map \cref{monoidal1}
$$t^{\kappa,\ell}:\Yo{\langle \kappa\rangle}\sqcup_{\Yo{\langle 0\rangle}}\Yo{\langle\ell\rangle}\to \Yo{\langle \kappa+\ell\rangle}$$
is again a map of the same form, namely, $t^{⟨κ⟩∧⟨r⟩,⟨\ell⟩∧⟨r⟩}$.
Likewise, the smash product of $τ:∅→\Yo{\langle0\rangle}$ \cref{monoidal2}
with any $⟨r⟩∈Γ$ is $τ$ again.
For $Δ$, the corresponding statement was proved by Rezk \cite[Theorem~7.2]{Rezk.CSS}.
\end{proof}

\cref{smooth.multiple.categories} gives the correct smooth variant of $d$-fold complete Segal spaces.
However, these are not quite a model for $(\infty,d)$-categories until we impose the globular condition of Barwick–Schommer-Pries \cite[Notation~12.1]{BarwickSchommerPries}.
For $d=2$ this amounts to passing from double categories to bicategories.
Recall that double categories have two distinct notions of 1-morphisms: horizontal and vertical.
Both can be composed, and 2-cells are squares involving two vertical and two horizontal morphisms.
In other words, our model so far describes the $d$-fold analog of double categories.
In order to get rid of the extra 1-morphisms, we further localize the functor category $\sPSh(\stcart\times \Gamma\times \Delta^{\times d})_{\uple}$
at the morphisms given by taking
external tensorings of the map \cref{glob} with representable presheaves on $\stcart\times\Gamma$. 

\begin{definition}
\label{globular.maps}
Fix $d\geq 0$.
We define the following maps, which we include in our left Bousfield localization for the globular model structure.
\begin{enumerate}
\item[(v)] (Globular maps.)
For an object ${\bf m}=([m_1],\ldots,[m_d])\in \Delta^{\times d}$, let $\hat{\bf m}$ be the object with $j$th component given by 
$$[\hat{m}_j]=\begin{cases}
\left[0\right],& \text{if there is $i<j$ with $m_i=0$,}\cr
\left[m_j\right],& \text{otherwise.}\cr
\end{cases}$$
There is a canonical map from ${\bf m}$ to $\hat{\bf m}$ given by identities or unique maps to $[0]$ for each index~$j$.
The globular maps are defined as the following morphisms in $\sPSh(Δ^{⨯d})$: 
\begin{equation}\label{glob}\psi^{{\bf m}}:\Yo{{\bf m}}\to \Yo{\bf \hat{m}}.\end{equation}
\end{enumerate}
\end{definition}

The local objects are multisimplicial spaces~$X$ such that $X_0$,
which we interpret as an object in $(d-1)$-fold simplicial spaces,
is homotopy constant and $X_k$ is a local object in $(d-1)$-fold simplicial spaces for any $k≥0$.
For $d=2$, the locality condition boils down to forcing the degeneration maps $X_{0,0}\to X_{0,b}$ to be equivalences.
This makes all vertical morphisms homotopic to identities.

\begin{notation}
Fix $d\geq 0$. We consider the following set of morphisms in $\sPSh(\stcart\times \Gamma\times \Delta^{\times d})$.
\begin{itemize}
\item Let $S_\glob=\{c⊠\psi^{\bf m}\mid c\in\stcart\times\Gamma, {\bf m}∈Δ^{⨯d}\}$,
where $\psi^{\bf m}$ are the morphisms \cref{glob}. 
\end{itemize}
\end{notation}

\begin{definition}[The globular model structure]
\label{globular.model.structure}
Fix $d\geq 0$.
We define the model structure 
$$\smcat_{∞,d}≔\sPSh(\stcart\times \Gamma\times \Delta^{\times d})_{\glob}$$
as the left Bousfield localization (\cref{left.Bousfield.localization})
of the model category $\sPSh(\stcart\times \Gamma\times \Delta^{\times d})_{\inj}$
at the set of morphisms $S_{\stcart}\cup S_{\Gamma}\cup S_{\Delta}\cup S_\glob$,
which exists by \cref{left.Bousfield.localization.exists}.
Alternatively, $\smcat_{∞,d}$ is the left Bousfield localization of $\smcatuple_{\infty,d}$ at the set of morphisms $S_{\glob}$.
\end{definition}

Again, the existence of the localization is established in Barwick \cite[Theorem 4.7]{Barwick.Model}. 
The following definition builds upon the work of Toën \cite{Toen}, although his model for $(∞,d)$-categories uses Segal categories, not complete Segal spaces.

\begin{definition}
\label{globular.multiple.model.structure}
Fix $d\geq 0$.
We define the model structures 
$$\catuple_{∞,d}≔\sPSh(\Gamma\times \Delta^{\times d})_\uple,
\qquad\cat_{∞,d}≔\sPSh(\Gamma\times \Delta^{\times d})_\glob$$
as the left Bousfield localizations (\cref{left.Bousfield.localization})
of the model category $\sPSh(\Gamma\times \Delta^{\times d})_\inj$
at the sets of morphisms defined like $S_\Gamma\cup S_\Delta$ and $S_\Gamma\cup S_\Delta\cup S_\glob$, respectively, but dropping $\stcart$ from all constructions.
The localizations exist by \cref{left.Bousfield.localization.exists}.
\end{definition}

\subsection{Functor categories}
\label{functor.categories}

Given $X,Y∈\smcatuple_{∞,d}$ (respectively $\smcat_{∞,d}$),
all functors $X→Y$ should themselves naturally organize into an object in $\smcatuple_{∞,d}$ respectively $\smcat_{∞,d}$.
By \cref{multiple.monoidal}, the multiple injective model structure on $\smcatuple_{∞,d}$ is a symmetric monoidal model category,
which allows for a convenient formalization of such a construction in the uple case.

\begin{definition}
\label{def.funmonuple}
Fix $d\geq 0$.
Let $X,Y\in \smcatuple_{∞,d}$.
We define the \emph{uple functor object} in $\smcatuple_{\infty,d}$ as the derived internal hom in simplicial presheaves on $\stcart\times \Gamma\times \Delta^{\times d}$.
Explicitly, we define
$$\Funmonuple(X,Y)≔\Funmon(X,R(Y))\in \smcatuple,$$
where $R$ denotes the fibrant replacement functor.  
\end{definition}

One can easily see that including the globular condition yields a model structure that does not satisfy the pushout product axiom,
so functor objects can no longer be computed by deriving the internal hom.
For example, working in complete globular 2-fold Segal spaces (i.e., simplicial presheaves on $Δ⨯Δ$),
the cartesian product of the cofibrant object $\Yo{[1],[0]}$ and the acyclic cofibration with cofibrant source $\Yo{[0],[0]}→\Yo{[0],[1]}$
is the map $\Yo{[1],[0]}→\Yo{[1],[1]}$ that is not a weak equivalence in the globular model structure:
its domain represents the $\infty$-groupoid of 1-morphisms with invertible 2-morphisms and higher homotopies,
whereas its codomain represents the $\infty$-groupoid of arbitrary 2-morphisms and higher homotopies.

To define functor objects for globular $d$-fold Segal spaces,
we transfer the derived internal hom of Rezk's $\Theta_d$-spaces (which form a cartesian model category)
along the Quillen equivalence between globular complete $d$-fold Segal spaces and Rezk's $\Theta_d$-spaces.
We emphasize that the resulting functor object is not computed as the left or right derived functor of some Quillen bifunctor.

Recall from Rezk \cite{Rezk.Theta} the notion of a $\Theta_d$-space and from Bergner–Rezk \cite{BergnerRezk} the relationship between $\Theta_d$-spaces and $d$-fold Segal spaces.
There is a functor $g:\Delta^{\times d}\to \Theta_d$ that arises as a composition
$$g:\Delta^{\times d}\lto3{g_1} \Delta^{\times d-1}\times \Theta_1 \lto3{g_2} \Delta^{\times d-2}\times \Theta_2 \lto3{g_3} \cdots \lto3{g_{d-1}} \Delta\times \Theta_{d-1}\lto3{g_d} \Theta_d,$$ 
where each $g_i$ is defined by 
$$g_i(([m_1],\ldots,[m_{d-i+1}],c)=([m_1],\ldots,[m_{d-i}],[m_{d-i+1}](c,\ldots, c )),$$
where $c\in \Theta_{i-1}$.
If $d=1$, then $g=g_1$ is the identity functor. If $d=2$, the functor $g=g_2$ is given by
$g([m_1],[m_2])=[m_1]([m_2],[m_2],\ldots,[m_2]),$
where $[m_2]$ is repeated $m_1$~times. 
The functor~$g$ is described explicitly on morphisms in Bergner–Rezk \cite{BergnerRezk}. 

By left and right Kan extension, we have an adjoint triple $g^{\#}\dashv g^*\dashv g_*$:
$$\xymatrix{
\sPSh(\Delta^{\times d})\ar@<.25cm>[r]^-{g^{\#}}\ar@<-.25cm>[r]_-{g_*} & \ar[l]|-{g^*} \sPSh(\Theta_d),
}$$
where the bottom adjunction $g^*\dashv g_*$ is a Quillen equivalence with respect to the globular model structure on $d$-fold Segal spaces and the Rezk model structure on $\Theta_d$-spaces (Bergner–Rezk \cite[Corollary 7.3]{BergnerRezk}).
By Rezk \cite[Theorem 8.1]{Rezk.Theta}, the category $\sPSh(\Theta_d)$ is a cartesian model category, when equipped with the Rezk model structure.
The following proposition enhances this result to the setting of smooth symmetric monoidal $(∞,d)$-categories.

\begin{proposition}
Fix $d\geq 0$.
We define the model structure 
$$\sPSh(\stcart\times\Gamma\times\Theta_d)_\local$$
as the left Bousfield localization (\cref{left.Bousfield.localization})
of the model category $\sPSh(\stcart\times \Gamma\times \Theta_d)_\inj$
at the set of morphisms $S_{\stcart}\cup S_{\Gamma}\cup S_\Theta$,
where $S_\stcart$ and $S_Γ$ are described in \cref{localmapsnot} (with $\Theta_d$ replacing $Δ^{⨯d}$)
and $S_\Theta$ is obtained by applying the functors $c⊠-$ ($c∈\stcart⨯Γ$) to the set $\mathscr{S}_\Theta$ defined inductively by Rezk \cite[§8]{Rezk.Theta}.
The resulting model structure is monoidal.
Furthermore,
the functor
$$\tilde g=\id_{\stcart}\times \id_{\Gamma}\times g: \stcart⨯Γ⨯\Delta^{\times d} → \stcart⨯Γ⨯\Theta_d$$
induces a Quillen equivalence
$$\xymatrix{
\sPSh(\stcart⨯Γ⨯\Delta^{\times d})\ar@<-.125cm>[r]_-{\tilde g_*} & \ar@<-.125cm>[l]_{\tilde g^*} \sPSh(\stcart⨯Γ⨯\Theta_d).
}$$
\end{proposition}

\begin{proof}
The left Bousfield localization exists by \cref{left.Bousfield.localization.exists}.
It is cartesian by the argument of \cref{multiple.monoidal}, using Rezk \cite[Theorem~8.1]{Rezk.Theta} for monoidality in the case of $\Theta_d$.
Promoting the argument of Bergner–Rezk \cite[Corollary 7.3]{BergnerRezk}
shows that the adjoint pair $\tilde g_*⊣\tilde g^*$ is a Quillen equivalence.
\end{proof}

Using this cartesian presentation, we can transfer the derived internal hom from $\Theta_d$-spaces
to $d$-fold Segal spaces through the Quillen equivalence $\tilde g^*\dashv \tilde g_*$ as follows. 

\begin{definition}
\label{globular.hom}
Fix $d\geq 0$ and let $\Funmontheta(-,-)$ denote the internal hom in $\sPSh(\stcart\times\Gamma\times\Theta_d)$.
Denote by~$R$ the fibrant replacement functor for $\sPSh(\stcart\times\Gamma\times\Delta^{⨯d})_\glob$.
We define the {\it globular product\/}
$$(Y,Z)\mapsto Y \tglob Z ≔ \tilde g^* (\tilde g_* R Y ⊗ \tilde g_* R Z)$$
and the {\it globular functor object\/}
$$(Y,Z)\mapsto \Funmonglob(Y,Z) ≔ \tilde g^* \Funmontheta(\tilde g_* R Y,\tilde g_* R Z).$$
\end{definition}

Denote by $\map(-,-)$ the simplicial enrichment functor in a category of simplicial presheaves and by $\rdf\map$ its right derived functor.

\begin{proposition}
\label{derivedhom}
For all $X,Y,Z∈\smcat_{∞,d}$
we have a natural weak equivalence of derived mapping spaces 
$$\rdf\map(X \tglob Y,Z)\simeq \rdf\map(X,\Funmonglob(Y,Z)).$$
In particular, $\tglob$ is separately homotopy cocontinuous in each argument and $\Funmonglob$ is separately homotopy continuous in each argument.
\end{proposition}

\begin{proof}
Using $R$ for the fibrant replacement functor in
$\sPSh(\stcart\times \Gamma\times \Delta^{\times d})_\glob$,
we have a chain of natural isomorphisms and weak equivalences
\begin{align*}
\rdf\map(X \tglob Y,Z)
&=\map(X \tglob Y,R Z)
=\map(\tilde g^* (\tilde g_* R X ⊗ \tilde g_* R Y),R Z)\cr
&≅\map(\tilde g_* R X ⊗ \tilde g_* R Y,\tilde g_* R Z)
≅\map(\tilde g_* R X, \Funmontheta(\tilde g_* R Y,\tilde g_* R Z))\cr
&≃\map(\tilde g^* \tilde g_* R X, \tilde g^* \Funmontheta(\tilde g_* R Y,\tilde g_* R Z))
≃\map(R X, \tilde g^* \Funmontheta(\tilde g_* R Y,\tilde g_* R Z))\cr
&=\map(R X, \Funmonglob(Y,Z))
≃\rdf\map(X, \Funmonglob(Y,Z)).\cr
\end{align*}
The first map labeled $≃$ is a weak equivalence because $\tilde g^*⊣\tilde g_*$ is a Quillen equivalence,
the second map labeled $≃$ is a weak equivalence because the derived counit $\tilde g^* \tilde g_* R X → R X$ is a weak equivalence,
the third map labeled $≃$ is a weak equivalence because derived mapping simplicial sets preserve weak equivalences.
\end{proof}

\subsection{Smooth sets}

There are three ways in which simplicial sets enter the construction of the bordism category:
\begin{itemize}
\item Higher gauge transformations of a geometric structure form an $\infty$-groupoid, which is encoded as a Kan complex, i.e., a simplicial set (\cref{geometric.structure}).
\item Gluing bordisms together is performed using open neighborhoods of the core (\cref{bord}).
Such open neighborhoods and their open embeddings form a category.
Taking the nerve of this category amounts to modding out by the equivalence relation that identifies two bordisms if they coincide on a smaller open neighborhood,
i.e., we pass to the germ of the core.
\item The isotopy space of cuts encodes the homotopy type of diffeomorphism groups of bordisms, responsible for the $(∞,d)$-category structure (\cref{frakbord}).
\end{itemize}

The first type of simplicial sets naturally persists throughout the entire paper.
The second type of simplicial sets is contractible and is quickly disposed of in \cref{emaplociso,bordisms.with.isotopies.representable}. 
The third type of simplicial sets is obtained as the smooth singular complex of a presheaf of sets on the site $\cart$.
The latter presheaf can be naturally seen as encoding the {\it smooth space\/} of isotopies of cuts.

\begin{definition}
Given a category~$C$, we refer to objects of the category
$$\sm C=\PSh(\cart,C)$$
as \emph{smooth $C$-objects}.
In particular, we have \emph{smooth sets} ($\smset$; $C=\set$),
\emph{smooth posets} ($\smposet$; $C=\poset$),
and \emph{smooth simplicial sets} ($\smsset$; $C=\sset$).
\end{definition}

\begin{definition}
Let $l\in \NN$.
We define the \emph{extended $l$-simplex} $\gsim^l$ as the smooth manifold given by equipping the subspace
$$\gsim^l=\biggl\{t\in \RR^{l+1} \Bigm| \sum_{i}t_i=1\biggr\}$$
with the canonical smooth structure.
\end{definition}

Working with presheaves on~$\cart$ has the following advantage, which will be exploited in the paper:
if a map of smooth simplicial sets is a stalkwise weak equivalence, it is also a weak equivalence in the model structure of \cref{smsets}.
Rather than constantly invoke the smooth singular complex functor, we work directly with presheaves and simplicial presheaves on the site $\cart$.
The following result can be seen as the formal justification of such an approach.

\begin{proposition}
\label{smsets}
(Pavlov \cite[Theorem~12.7]{Pavlov.Diffeo}.)
The category $$\smsset=\sPSh(\cart)$$
of smooth simplicial sets
admits a cartesian simplicial model structure
whose weak equivalences are transferred along the smooth singular complex functor
$$\smsset→\sset, \qquad F↦(l↦F(\gsim^l)_l)$$
and generating cofibrations
are given by the maps $(∂Δ^m→Δ^m)\ppet(\gbou^n→\gsim^n)$,
where $\gbou^n=|∂Δ^n|$ with $|{-}|$ denoting the left adjoint of the singular functor $\smset→\sset$
and $\ppet$ denotes the pushout product associated to the external tensor product~$⊠$.
The resulting Quillen adjunction 
$$\xymatrix{
\sset\ar@<.1cm>[r] &\ar@<.1cm>[l]\smsset.}$$
is a Quillen equivalence.
Given a small category~$\mathscr{C}$, we set
$$\smsPSh(\mathscr{C})=\PSh(\mathscr{C},\smsset)$$
and equip it with the injective model structure.
If $\mathscr{C}$ is monoidal, then $\smsPSh(\mathscr{C})$ is equipped with the Day convolution monoidal structure.
\end{proposition}

\begin{remark}
By Pavlov \cite[Proposition~12.5]{Pavlov.Diffeo},
a stalkwise weak equivalence of simplicial presheaves is a weak equivalence in the model structure on $\smsset$ in \cref{smsets}.
The converse is false, since the map $\RR^n→\RR^0$ is not a stalkwise weak equivalence, but its singular complex is a weak equivalence of simplicial sets.
Thus, the two model structures on $\smsset$ are not Quillen equivalent.
\end{remark}

\begin{remark}
\label{enrichedbous}
By Hirschhorn \cite[Theorem~4.1.1(4)]{Hirschhorn},
any left Bousfield localization of an $\sset$-enriched model category is automatically $\sset$-enriched.
Likewise, any left Bousfield localization of a $\smsset$-enriched model category is automatically $\smsset$-enriched.
Indeed, the terminal object generates $\smsset$ under homotopy colimits and therefore is a homotopy generator for $\smsset$.
Now apply the criterion in \cref{left.Bousfield.localization.exists}.
\end{remark}

\begin{proposition}
\label{multiple.model.structure.smooth}
There is a monoidal model structure 
$$\fraksmcatuple_{\infty,d}≔\smsPSh(\stcart\times\Gamma\times\Delta^{\times d})_{\uple}$$
given by performing the left Bousfield localization (\cref{left.Bousfield.localization}) of the model category $$\smsPSh(\stcart\times\Gamma\times\Delta^{\times d})_{\inj}$$
at the set of morphisms defined like $S_{\stcart}\cup S_{\Gamma}\cup S_{\Delta}$, but using the Yoneda embedding for presheaves valued in $\smset$ instead of $\sset$. 
The Quillen equivalence $\sset \rightleftarrows \smsset$ induces a Quillen equivalence
$$\xymatrix{
\smcatuple_{∞,d}\ar@<.1cm>[r] &\ar@<.1cm>[l]\fraksmcatuple_{∞,d}.}$$
\end{proposition}

\begin{proof}
The localization exists and is enriched over $\smsset$ by \cref{left.Bousfield.localization.exists} and \cref{enrichedbous}.
It is monoidal by the same argument as in \cref{multiple.monoidal}, substituting $\smsset$ for $\sset$. 
The Quillen equivalence $\sset \rightleftarrows \smsset$ induces a Quillen equivalence of injective model structures.
By Hirschhorn \cite[Theorem~3.3.20]{Hirschhorn}, this Quillen equivalence descends to a Quillen equivalence of localized model structures.
\end{proof}

\begin{proposition}
\label{globular.model.structure.smooth}
There is a model structure 
$$\fraksmcat_{\infty,d}≔\smsPSh(\stcart\times\Gamma\times\Delta^{\times d})_{\glob}$$
given by performing the left Bousfield localization (\cref{left.Bousfield.localization}) of the model category $$\smsPSh(\stcart\times\Gamma\times\Delta^{\times d})_{\inj}$$
at the set of morphisms defined like $S_{\stcart}\cup S_{\Gamma}\cup S_{\Delta}\cup S_\glob$, but using the Yoneda embedding for presheaves valued in $\smset$ instead of $\sset$.
The Quillen equivalence $\sset \rightleftarrows \smsset$ induces a Quillen equivalence
$$\xymatrix{
\smcat_{∞,d}\ar@<.1cm>[r] &\ar@<.1cm>[l]\fraksmcat_{∞,d}.}$$
\end{proposition}

\begin{proof}
Same as the proof of \cref{multiple.model.structure.smooth}, adding $S_\glob$ to the set of localizing morphisms.
\end{proof}

\begin{remark}
\label{companion}
The 1-bordism $τ$ is not quite what we need: its source and target were induced from~$ρ$,
meaning the ambient manifold is $M_2=\RR⊔\RR$, and from each copy of~$\RR$ we cut out a single point.
On the other hand, the target of~$η$ is given by the ambient manifold $M_0=\RR$, from which we cut out two points $0$ and~$2$.
These two vertices in   can be connected by a zigzag of 1-simplices
that restricts from $M_0=\RR$ to $M_1=(-1,1)∪(1,3)$ and then embeds into $M_2=\RR⊔\RR$ via the inclusion maps
into the first respectively second summand.

Such a pair of 1-simplices is also not a 1-bordism, but can be converted to a zigzag of 1-bordisms
using a lifting construction that is entirely analogous to the one for~$ρ$,
replacing the category $\cart$ used for isotopies with the category~$Δ$ used in $\smsset$,
and replacing $\RR^1$ with $Δ^1$.
Thus, instead of $(\RR^0,⟨1⟩,[0],\RR^1)$ we use $(\RR^0,⟨1⟩,[0],\RR^0,[0])$, etc.,
where the last component now corresponds to the category $Δ$ used in $\smsset$.

To convert~$ρ$ into a 1-bordism, consider the acyclic following cofibration~$ι$ in $\smcat_{∞,d}$:
$$(\RR^0,⟨1⟩,[0],\RR^1) ⊔_{(\RR^0,⟨1⟩,[0],\RR^0)} (\RR^0,⟨1⟩,[1],\RR^1)⊔_{(\RR^0,⟨1⟩,[0],\RR^0)}(\RR^0,⟨1⟩,[0],\RR^1) → (\RR^0,⟨1⟩,[1],\RR^1).$$
The domain of~$ι$ maps into via the map with three components, where the third component is~$ρ$
and the first respectively second component is uniquely determined by the condition that it is simplicially degenerate in the direction of~$Δ$
respectively given by a pullback along a map $\RR^1→\RR^0$.
Extending the resulting map along~$ι$ and then restricting along the inclusion
$$(\RR^0,⟨1⟩,[1],\RR^0≅\{-2\}) → (\RR^0,⟨1⟩,[1],\RR^1).$$
produces an element~$τ$ in hat can be seen as the homotopy coherent analogue of the companion of~$ρ$.
\end{remark}

\section{Criteria for weak equivalences of higher categorical structures}
\label{homotopy.theory}

\subsection{Cartesian spaces and stalks}

First, we show how to reduce the problem of showing that a morphism in $$\smcat_{∞,d}=\sPSh(\stcart⨯Γ⨯Δ^{⨯d})_\glob$$ (\cref{globular.model.structure})
is a weak equivalence to the same problem in the model category $$E=\cat_{∞,d}=\sPSh(Γ⨯Δ^{⨯d})_\glob$$ (\cref{globular.multiple.model.structure}).

\begin{proposition}
\label{stalks}
Suppose $E$ is a combinatorial model category,
$\stcart$ is the site from \cref{def.cartsp},
and $\PSh(\stcart,E)$ denotes the Čech localization of the injective model structure,
which exists by \cref{left.Bousfield.localization.exists}.
Then the following hold.
\begin{enumerate}
\item A morphism in $\PSh(\stcart,E)$ is a weak equivalence
if and only if all of its stalks at all points $\mathfrak{p}$ of $\stcart$ are weak equivalences in~$E$.
\item A point $\mathfrak{p}$ of $\stcart$ is given by a pair $\mathfrak{p}=(T,ρ)$, where $T∈\stcart$ and $ρ$ is an isomorphism $\stu T\to \RR^n$ for some $n≥0$.
Denote by $B^{\mathfrak{p}}_ε$ the open subobject of~$T$
given by the unique cartesian lift of the open ball $B_ε^n⊂\RR^n≅\stu T$ of radius~$ε$ centered at~0.
The associated stalk functor
$$\mathfrak{p}^*:\PSh(\stcart,E)\to E$$
sends $X∈\PSh(\stcart,E)$ to the stalk
$$\mathfrak{p}^*X≔\colim_{ε\to 0}X(B^{\mathfrak{p}}_ε)∈E,\eqlabel{stalk}$$
where the (homotopy) colimit is taken over the inclusions of open balls.
Different choices of $\rho$ and isomorphic choices of $T$ lead to isomorphic stalks. 
\end{enumerate}
\end{proposition}

\begin{proof}
The proof is standard, see, for example, Amabel–Debray–Haine \cite[Proposition A.5.4]{AmabelDebrayHaine} for the case $\stcart=\cart$ and $\stu=\id$.
Recall that a family of relative functors $\{F_i:C→D_i\}_{i∈I}$ is {\it jointly conservative\/} if for every morphism~$f$ in~$C$
the following property holds: if for all~$i$ the morphism $F_i(f)$ is a weak equivalence in~$D_i$, then $f$ is a weak equivalence in~$C$.
In particular, the family of functors 
$$\{\rdf ι_T^*:\PSh(\stcart,E)→\PSh(T,E)\}_{T∈\stcart}$$
given by the right derived functors of functors that restrict a presheaf from the site $\stcart$ to the small site of~$T$
is a jointly conservative family for the model category $\PSh(\stcart,E)$.

Since the small site of $T∈\stcart$ coincides with the small site of the ordinary manifold $\stu(T)∈\cart$,
the family of ordinary stalk functors 
$$\{p^*:\PSh(T,E)→E\}_{p∈\stu(T)}$$
is a jointly conservative family for $\PSh(T,E)$.
Thus, the family of functors 
$$\{p^*\rdf ι_T^*:\PSh(\stcart,E)→E\}_{T∈\stcart,p∈\stu(T)}$$
is a jointly conservative family for $\PSh(\stcart,E)$.

Any point $p∈\stu(T)$ can be translated to $0∈\stu(T)$ using an isomorphism in $\stcart$.
Identifying $\stu(T)$ with some~$\RR^n$ and using the fact that stalks at $0∈\RR^n$
can be computed by taking the sequential colimit over $ε$-balls around~$0$,
we arrive at the stated formula for stalks.
\end{proof}

The following proposition will be used in the description of the double stalk for the isotopy embedded bordism category. 

\begin{proposition}
\label{doublestalk}
Fix a point $\mathfrak{l}$ of $\cart$ and a point $\mathfrak{p}$ of $\stcart$. Let $X\in \PSh(\stcart\times \cart)$ be a presheaf of sets. Then we have a natural isomorphism 
$$\mathfrak{l}^*\mathfrak{p}^*X\cong \colim_{\epsilon\to 0}X(B_{\epsilon}^{\mathfrak{l}},B_{\epsilon}^{\mathfrak{p}}),$$
where the colimit runs over the totally ordered set of positive real numbers and $B_{\epsilon}^{\mathfrak{l}}$ and $B_{\epsilon}^{\mathfrak{p}}$ are open balls of radius~$\epsilon$, as defined in \cref{stalks}. 
\end{proposition}
\begin{proof}
Let $\mathcal{P}_{\mathfrak{l}}$ and $\mathcal{P}_{\mathfrak{p}}$ be the filtered posets of open balls and their inclusions, which are both isomorphic to the totally ordered set of positive real numbers $\RR_{>0}$. The functor
$$\RR_{>0}\to \mathcal{P}_{\mathfrak{l}}\times \mathcal{P}_{\mathfrak{p}}$$
that sends $\epsilon$ to the pair $(B_{\epsilon}^{\mathfrak{l}},B_{\epsilon}^{\mathfrak{p}})$ and a morphism $\epsilon<\epsilon'$ to the corresponding pair of inclusions of open balls is an initial functor. The claim follows. 
\end{proof}

\subsection{Reduction of multiple to single}

Next, we explain how the problem of showing that a morphism in $\cat_{\infty,d}$ or $\catuple_{\infty,d}$ (\cref{globular.multiple.model.structure}) is a weak equivalence
can be reduced to the same problem in the model categories $\sPSh(Δ)_\local$ (\cref{gammasp}) or $\sPSh(Γ)_\local$ (\cref{rezkdelta}).

\begin{proposition}
\label{multiple.single}
If $f$ is a morphism in $\cat_{\infty,d}$ such that for each ${\bf m}∈Δ^{⨯d}$
the functor
$$\cat_{\infty,d}→\sPSh(Γ)_\local$$
induced by
$$Γ→Γ⨯Δ^{⨯d},\qquad ⟨\ell⟩↦(⟨\ell⟩,{\bf m})$$
sends $f$ to a weak equivalence,
then $f$ is a weak equivalence.

Given $i∈\{1,…,d\}$, if $f$ is a morphism in $\cat_{\infty,d}$ such that for each $⟨\ell⟩∈Γ$, and ${\bf m}∈Δ^{d-1}$
the functor
$$\cat_{\infty,d}→\sPSh(Δ)_\local$$
induced by
$$Δ→Γ⨯Δ^{⨯d},\qquad [\omega]↦(m_1,…,m_{i-1},\omega,m_i,…,m_{d-1},⟨\ell⟩)$$
sends $f$ to a weak equivalence,
then $f$ is a weak equivalence. The analogous statement for the $\catuple_{\infty,d}$ also holds. 
\end{proposition}

\begin{proof}
We will prove the second claim in the globular case, the first claim and the uple case are analogous. 
If $f:F→G$ is a morphism with indicated properties,
then the induced natural transformation in $\PSh(\Gamma\times Δ^{\times d-1},\sPSh(\Delta)_\local)$
with components $f_{\langle \ell\rangle,{\bf m}}:F_{\langle \ell\rangle,{\bf m}}→G_{\langle \ell\rangle,{\bf m}}$
is an objectwise weak equivalence. 
Such objectwise weak equivalences are precisely weak equivalences
in the left Bousfield localization of $\sPSh(\Gamma\times\Delta^{\times d})_\inj$
with respect to the set of maps
$$S_{\Delta,i}=\{\Yo{c}⊠\phi^{a,b}, \Yo{c}⊠x\mid c\in\Gamma\times Δ^{\{1,\ldots,i-1,i+1,\ldots,d\}}, a,b∈Δ\},$$
where $\phi^{a,b}$ are the morphisms \cref{segal2} in \cref{bousloc} and $x$ is the morphism \cref{segal1}.
Since the set $S_{Δ,i}$ is a subset of the set $S_Γ∪S_Δ∪S_\glob$ used to construct $\cat_{∞,d}$ in \cref{globular.multiple.model.structure},
objectwise weak equivalences in
$\PSh(\Gamma\times Δ^{\times d-1},\sPSh(\Delta)_\local)$
are also weak equivalences in $\cat_{∞,d}$.
\end{proof}

\subsection{The simplicial Whitehead theorem}

We will need the following special case of the simplicial Whitehead theorem.

\begin{proposition}
\label{new.whitehead}
Suppose $C$ is a small category such that for any full subcategory $F⊂C$
with finitely many objects
there is a zigzag of natural transformations
that connects the inclusion functor $F→C$
to the composition $F→1→C$, where $1$ is the terminal category and $1→C$ picks out some object of~$C$.
Then the nerve of~$C$ is weakly contractible.
\end{proposition}

\begin{proof}
Set $X$ to the nerve of~$C$.
It suffices to show that any map
$$g:∂Δ^n → \Exi X$$
can be extended to a map
$$Δ^n → \Exi X.$$
The map~$g$ factors as $g:\partial \Delta^n\to \Ex^k X\to \Exi X$, for some $k\geq 0$. The first map in the factorization is adjoint to a map for the form $f:\Sd^k ∂Δ^n→X$.

Take $F$ to be the full subcategory of~$C$ on objects in the image of $f_0$.
By assumption,
the inclusion $F→C$ is connected by a zigzag of natural transformations to some constant functor $F→1→C$.
Taking nerves of all natural transformations,
we get a simplicial homotopy $h_0:H⨯\Sd^k ∂Δ^n→X$
from $f$ to a constant simplicial map $\Sd^k ∂Δ^n→Δ^0→X$,
where $H$ is an arbitrary zigzag of 1-simplices.

The adjoint map to~$h_0$ has the form $H→\hom(\Sd^k ∂Δ^n,X)=\hom(∂Δ^n,\Ex^k X)$,
and we postcompose it with the map $\hom(∂Δ^n,\Ex^k X)→\hom(∂Δ^n,\Exi X)$
induced by the inclusion $\Ex^k X→\Exi X$. 
The latter composition is adjoint to the map $h_1:B=H⨯∂Δ^n→\Exi X$,
which is a simplicial homotopy from~$g$ to a constant simplicial map~$p:∂Δ^n→Δ^0→\Exi X$.
We extend $h_1$ to a map $h_1':A→X$,
where $A=B⊔_{∂Δ^n}Δ^n$ (the map $\partial \Delta^n\to B$ is supplied by the last vertex in $H$)
and $h'_1$ is defined on $Δ^n$ by extending~$p$ to~$Δ^n$ as a constant map.

The inclusion $A→H⨯Δ^n$ is an acyclic cofibration,
thus the map $A→\Exi X$ can be extended to a map $H⨯Δ^n→\Exi X$.
Precomposing the latter map with the map $\Delta^n=\Delta^0\times Δ^n\to H\times \Delta^n$ that picks out the initial vertex of $H$ provides the desired extension
of the map $g:∂Δ^n→\Exi X$ to~$Δ^n$.
\end{proof}

\subsection{Dugger and Spivak's rigidification of quasicategories}

Another tool needed for the proof of the main theorem is a rigidification of quasicategories.
The reason we will need this rigidification is the following.
In \cref{codescent}, we will have a monomorphism of simplicial sets $f:X\to Y$ which we claim is a weak equivalence in the Joyal model structure.
In our case, the objects $X$ and $Y$ are not fibrant in the Joyal model structure and it is unclear how to prove that this map is a weak equivalence directly.
On the other hand, if we work with simplicial categories, then the induced map on hom spaces turns out to be much easier to understand.  

A morphism of simplicial sets $X\to Y$ is a weak equivalence in the Joyal model structure if it induces a Dwyer–Kan equivalence on corresponding simplicial categories $\mathfrak{C}(X)\to \mathfrak{C}(Y)$.
Here $\mathfrak{C}$ is left adjoint to the homotopy coherent nerve functor.
Dugger and Spivak have a concrete model \cite{Dugger.Spivak} for a functor $\Cnec$ weakly equivalent to $\mathfrak{C}$,
where the mapping simplicial sets are described using necklaces.
We now review this description. 

\begin{definition}
The \emph{monoidal category of bipointed simplicial sets} $\sset_{*,*}$ is defined as follows.
Objects are triples $(X,x,y)$, where $X$ is a simplicial set and $x$ and $y$ are vertices in~$X$.
Morphisms are simplicial maps $X→X'$ that preserve $x$ and~$y$.
The monoidal product is defined by
$$(X,x,y)∨(X',x',y')=(X⊔_{y,Δ^0,x'}X',x,y').$$
A simplex~$Δ^n$ is turned into an object of~$\sset_{*,*}$ by declaring $x$ to be the initial vertex and $y$ to be the final vertex.
In particular, $Δ^0$ is the monoidal unit.
The \emph{monoidal category of necklaces} $\Nec$ is defined as the full monoidal subcategory of $\sset_{*,*}$
generated by the bipointed simplices~$Δ^n$ ($n≥0$).
\end{definition}

Thus, an object of $\Nec$ is a \emph{necklace} $\Delta^{n_1}\vee \cdots \vee \Delta^{n_k}$ ($k≥0$).
Each $\Delta^{n_i}$ is called a \emph{bead} and the initial/final vertex of each $\Delta^{n_i}$ is called a \emph{joint}.
Morphisms of necklaces are compositions of morphisms induced by the canonical inclusions $\Delta^{n_i}\vee \Delta^{n_{i+1}}\into \Delta^{n_i+n_{i+1}}$,
and the face and degeneracy maps of each bead.

The following alternative description of the category $\Nec$ will be used later.

\begin{proposition}
\label{alternative.nec}
The category $\Nec$ is equivalent to the following category~$\Nec'$.
Objects are pairs $([m],Υ)$, where $[m]∈Δ$ and $Υ⊂[m]$ is a subset that contains $\{0,m\}$, whose elements are called \emph{joints}.
Morphisms $([m],Υ)→([m'],Υ')$ are morphisms of simplices $σ:[m]→[m']$ such that $σ(Υ)⊃Υ'$.
(Hence, $σ(0)=0$ and $σ(m)=m'$.)
Composition and identities are induced from~$Δ$ via the forgetful maps $([m],Υ)↦[m]$, $σ↦σ$.
\end{proposition}

\begin{proof}
We define a functor $\Nec→\Nec'$ as follows.
Send an object~$X$ of $\Nec$ to the pair $(X_0,Υ)$,
where the set~$X_0$ of vertices of~$X$ is equipped with the unique total order such that $x≤x'$ whenever there is a 1-simplex from~$x$ to~$x'$,
and $Υ$ denotes the set of all joints in~$X_0$.
Send a morphism $f:X→Y$ in $\Nec$ to the induced map~$f_0$ on vertices.

We define a functor $\Nec'→\Nec$ as follows.
Send an object $([n],Υ)$ of $\Nec'$ to the necklace $⋁_{\{υ,υ'\}⊂Υ}Δ^{υ,…,υ'}$, where $υ<υ'$ run over consecutive elements of~$Υ$.
Send a morphism $σ:([n],Υ)→([n'],Υ')$ to the unique simplicial map whose vertex map is~$σ$.
By inspection, the two functors are mutually inverse to each other, which completes the proof.
\end{proof}

\begin{definition}
\label{Cnec.def}
Suppose $X$ is a simplicial set.
The simplicial category $\Cnec(X)$ is defined as follows.
Objects are vertices of~$X$.
Given two objects $x$ and $y$,
the simplicial set of morphisms $\Cnec(X)(x,y)$ is the nerve of the \emph{category of necklaces from~$x$ to~$y$ in~$X$}, defined as the comma category $\Nec/(X,x,y)$.
Composition in $\Cnec(X)$ is induced by the monoidal structure on necklaces.
\end{definition}

\begin{proposition}
\label{necklaces.prop}
The following three conditions on a simplicial map $f:X\to Y$ are equivalent:
\begin{enumerate}
\item $f$ is a weak equivalence in the Joyal model structure;
\item promoting $f$ to a simplicial object of (discrete) simplicial sets yields a weak equivalence in the Rezk model structure (\cref{rezkdelta});
\item the simplicial functor
$$\Cnec(f):\Cnec(X)\to \Cnec(Y)$$
is a Dwyer--Kan weak equivalence of simplicial categories. 
\end{enumerate}
\end{proposition}

\begin{proof}
We first prove (1) $\Leftrightarrow$ (2).
By Joyal–Tierney \cite[Theorem~4.11]{JoyalTierney},
the functor $p_1^*:\PSh(Δ)→\sPSh(Δ)$ that takes the constant simplicial set objectwise is a left Quillen equivalence that preserves and reflects weak equivalences.
Thus, $f$ is a Joyal weak equivalence if and only if $p_1^*f$ is a Rezk weak equivalence.

Next, we prove (3) $\Leftrightarrow$ (1). 
By Dugger–Spivak \cite[Theorem 5.2 (arXiv); 5.3 (journal)]{Dugger.Spivak}, the functor~$\Cnec$ is connected to the functor~$\mathfrak{C}$ by a zigzag of weak equivalences.
By Joyal's theorem (Lurie \cite[Theorem~2.2.5.1]{Lurie.HTT}), the functor $\mathfrak{C}$ is a left Quillen equivalence that preserves and reflects weak equivalences.
Therefore, the functor $\Cnec$ preserves and reflects weak equivalences.
Thus, $f$ is a Joyal weak equivalence if and only if $\Cnec(f)$ is a Dwyer–Kan weak equivalence. 
\end{proof}

\subsection{$\Gamma$-spaces}
In this section, we give a sufficient condition to detect acyclic cofibrations in the local injective model structure on $Γ$-spaces.
Recall that the category~$Γ$ is defined as the opposite category of finite pointed sets and basepoint-preserving maps of sets.
We denote the basepoint by~$*$.

\begin{definition}
(Lurie \cite[Definitions 2.1.1.8, 2.1.2.1]{Lurie.HA}.)
A morphism~$γ$ in the category~$Γ$ corresponding to a map of finite pointed sets $f:⟨m⟩→⟨n⟩$ is
\begin{itemize}
\item \emph{inert} if $f|_A:A→B$ is a bijection, where $B=⟨n⟩∖\{*\}$ and $A=f^{-1}B$;
\item \emph{active} if $f^{-1}\{*\}=\{*\}$.
\end{itemize}
\end{definition}

Given $⟨m⟩∈Γ$, the space $X(⟨m⟩)$ can be interpreted as the space of multipliable $m$-tuples in a $Γ$-space~$X$.
The structure map $X(γ)$ for a morphism $γ$ in~$Γ$ removes elements from a tuple if $γ$ is an inert morphism and multiplies elements if $γ$ an active morphism.
Any morphism in~$Γ$ factors as an active morphism followed by an inert morphism.

\begin{definition}
\label{boundarygamma}
Let $\Yo{}:\Gamma\into \PSh(\Gamma)$ denote the Yoneda embedding.
We define the boundary of a representable $\Yo{\langle \ell\rangle}\in \PSh(\Gamma)$ as the subpresheaf $\partial\langle \ell\rangle\into \Yo{\langle \ell\rangle}$,
whose value at $\langle m\rangle$
consists of all maps of finite pointed sets $\langle \ell \rangle\to \langle m\rangle$ that are not active,
i.e., at least one element $k\in \langle \ell\rangle$ with $k\neq \ast$ is sent to $\ast\in  \langle m\rangle$.
\end{definition}

\begin{proposition}
\label{gammacyc}
For $\ell>1$, the inclusion $\partial \langle \ell\rangle\into \Yo{\langle \ell\rangle}$ (\cref{boundarygamma}), viewed as a morphism of discrete simplicial presheaves,
is an acyclic cofibration in the local injective  model structure on $\Gamma$-spaces.
\end{proposition}

\begin{proof}
We prove by induction on $\ell$ that $∂⟨\ell⟩→\Yo{⟨\ell⟩}$ is a weak equivalence.
The map in $\PSh(\Gamma)$,
$$S_\ell≔\Yo{⟨1⟩}⊔_{\Yo{⟨0⟩}}\Yo{⟨1⟩}⊔_{\Yo{⟨0⟩}}⋯⊔_{\Yo{⟨0⟩}}\Yo{⟨1⟩}→\Yo{⟨\ell⟩},$$
is an acyclic cofibration for any $⟨\ell⟩$, which follows by repeatedly composing cobase changes of maps of the form \cref{monoidal1}.
The resulting subpresheaf picks out those maps of pointed sets $⟨\ell⟩→⟨m⟩$ that send all elements to~$\ast$ except for some $k∈⟨\ell⟩$, $k≠\ast$.

For $\ell=2$, the map $\Yo{⟨1⟩}⊔_{\Yo{⟨0⟩}}\Yo{⟨1⟩}→\Yo{⟨2⟩}$ coincides with $∂⟨2⟩→\Yo{⟨2⟩}$, which establishes the base of induction.

Suppose $∂⟨m⟩→\Yo{⟨m⟩}$ is a weak equivalence for all $m<\ell$.
To show that $∂⟨\ell⟩→\Yo{⟨\ell⟩}$ is a weak equivalence,
we present the inclusion $S_\ell→∂⟨\ell⟩$ as a composition
of cobase changes of coproducts of maps of the form $∂⟨m⟩→\Yo{⟨m⟩}$,
and then use the 2-out-of-3 property for the maps $S_\ell→∂⟨\ell⟩$ and $S_\ell→\Yo{⟨\ell⟩}$
to conclude that $∂⟨\ell⟩→\Yo{⟨\ell⟩}$ is a weak equivalence.

We introduce the following filtration on $S_\ell→∂⟨\ell⟩$:
$$S_\ell=W_1→W_2→⋯→W_{\ell-1}=∂⟨\ell⟩.$$
Here $W_k$ is a subpresheaf of $\Yo{⟨\ell⟩}$
comprising precisely those maps of pointed finite sets $⟨\ell⟩→⟨m⟩$
such that the preimage of $⟨m⟩∖\{\ast\}$ has $k$ or fewer elements.
By construction, $S_\ell=W_1$ and $W_{\ell-1}=∂⟨\ell⟩$.

It remains to present the inclusion $W_{k-1}→W_k$ for any $1<k<\ell$
as a cobase change of a coproduct of maps of the form  $∂⟨k⟩→\Yo{⟨k⟩}$.
The indexing set of the coproduct is the set of all inert maps $f$ of pointed finite sets $⟨\ell⟩→⟨k⟩$
such that $f$ is (strictly) increasing when restricted to $⟨\ell⟩∖f^{-1}\{*\}$.
(Here we use the natural order on $⟨\ell⟩∖\{*\}=\{1,…,\ell\}$.)
The map~$f:\langle \ell \rangle \to \langle k\rangle$ corresponds to a map $\Yo{⟨k⟩}→W_k$ under the (contravariant) Yoneda embedding,
and the attaching map $∂⟨k⟩→W_{k-1}$ is the composition $∂⟨k⟩→\Yo{⟨k⟩}→W_k$, which factors through $W_{k-1}$ by construction.

To show that the resulting commutative square is a pushout square,
pick an arbitrary element of $W_k$ that is not in $W_{k-1}$,
i.e., a map of pointed finite sets $g:⟨\ell⟩→⟨m⟩$ such that the preimage~$A$ of $⟨m⟩∖\{\ast\}$ has exactly $k$ elements.
We have to show that this element comes from a unique element in the coproduct,
i.e., there is a unique pair $(f,σ)$, where $f:⟨\ell⟩→⟨k⟩$ is an element in the indexing set of the coproduct
and $σ:⟨k⟩→⟨m⟩$ is a map of pointed finite sets such that $σf=g$.
Indeed, $f$ must be the unique map $f:⟨\ell⟩→⟨k⟩$ of pointed finite sets such that the preimage of $⟨k⟩∖\{\ast\}$ equals~$A$ and $f$ is strictly increasing.
For such an $f$ there is a unique $σ:⟨k⟩→⟨m⟩$ for which $σf=g$ because $f$ is injective away from the preimage of~$\{\ast\}$.
\end{proof}

\begin{definition}
\label{indecomposable.element}
Suppose $X∈\PSh(Γ)$, $⟨k⟩∈Γ$, and $x∈X(⟨k⟩)$.
We say that $x$ is \emph{indecomposable} if for any active morphism~$γ:⟨k⟩→⟨m⟩$ and $y∈X(⟨m⟩)$,
the relation $x=X(γ)(y)$ implies that there is an active morphism $δ:⟨m⟩→⟨k⟩$
such that $y=X(δ)(x)$ and $δγ=\id_{⟨k⟩}$.
\end{definition}

\begin{remark}
\label{indecomposable.bordism}
In the context of bordism categories, elements $x∈X(⟨k⟩)$ involve maps $M→⟨k⟩$ that encode a $k$-tuple of bordisms via the fibers over various points of $⟨k⟩=\{*,1,…,k\}$,
where the fiber over~$*$ is the trash bin.
The $Γ$-structure maps act by composition of maps of sets: $M→⟨k⟩→⟨\ell⟩$.
Inert maps can move some components to the trash bin and do nothing else.
Active maps can take disjoint unions of some fibers and do not move anything to the trash bin.
An element~$x$ is indecomposable if every fiber has a single connected component.
To see this, observe that $x=X(γ)(y)$ means that the fibers of~$x$ are decomposed into disjoint unions according to the active morphism~$γ$.
The relation $δγ=\id_{⟨k⟩}$ forces the underlying map of sets of~$γ$ to be surjective 
and the underlying map of sets of~$δ$ to be injective.
Then $y=X(δ)(x)$ means that the fibers of~$y$ are precisely the fibers of~$x$ together with some empty fibers. Since this property holds for all active $\gamma$, this implies that all fibers of~$x$ are connected.
\end{remark}

\begin{definition}
\label{EZ.property}
Suppose $X∈\PSh(Γ)$.
We say that $X$ has the \emph{modified Eilenberg–Zilber property}
if for any $⟨k⟩∈Γ$, $x∈X(⟨k⟩)$ there is an active morphism $σ:⟨k⟩→⟨\ell⟩$ and an indecomposable $g∈X(⟨\ell⟩)$ such that $x=X(σ)(g)$
and for any other such pair $(g',σ':⟨k⟩→⟨\ell'⟩)$ the following two conditions are satisfied:
(1) there is a unique isomorphism $δ:⟨\ell⟩→⟨\ell'⟩$ such that $g=X(δ)(g')$;
(2) for this unique~$δ$ we also have $δσ=σ'$.
\end{definition}

\begin{remark}
\label{indecomposable.free}
Condition~(2) of \cref{EZ.property} is equivalent to saying that the action of automorphisms of $⟨k⟩∈Γ$ (i.e., the symmetric group~$Σ_k$)
on indecomposable elements in $X(⟨k⟩)$ is free.
\end{remark}

\begin{remark}
Continuing \cref{indecomposable.bordism}, in the context of bordism categories the modified Eilenberg–Zilber property
says that any $k$-tuple of bordisms $M→⟨k⟩$ decomposes as a composition of maps $M→⟨\ell⟩→⟨k⟩$,
where $M→⟨\ell⟩$ is indecomposable, i.e., its fibers are connected.
Furthermore, this decomposition is unique up to a unique permutation,
and the action of the symmetric group $Σ_k$ on $k$-tuples $M→⟨k⟩$ is free.
\end{remark}

\begin{proposition}
\label{EZ.criterion}
Suppose $ι:X→Y$ is an inclusion in $\PSh(Γ)$ such that
$X$ and $Y$ have the modified Eilenberg–Zilber property,
inert structure maps of~$Y$ preserve indecomposable elements of~$Y$,
and indecomposable elements of~$Y$ in degree~$⟨1⟩$ belong to~$X$.
Then $ι$ is a weak equivalence in the model category $\sPSh(\Gamma)_{\local}$ (\cref{gammasp}).
\end{proposition}

\begin{proof}
{\bf Step 1.} (A filtration on $Y$.)
Set $Y_0=X$.
For $k>0$, denote by $Y_k⊂Y$ the subobject given by the union of $Y_{k-1}$
and the subobject of~$Y$ generated by indecomposable elements of $Y(⟨k⟩)$.
Hence, we have a filtration 
$$X=Y_0⊂Y_1⊂Y_2⊂⋯⊂Y.$$
By assumption, $Y_0=Y_1$.

Any morphism $γ$ in $Γ$ factors as an active morphism~$a$ followed by an inert morphism~$i$.
Acting by~$γ$ on an indecomposable element~$x$ in~$Y$ yields $y=Y(γ)(x)=Y(a)(Y(i)(x))$.
Since inert structure maps preserve indecomposables by assumption,
the element $Y(i)(x)$ is indecomposable and together with the morphism~$a$ provides the required data for the modified Eilenberg–Zilber property of~$y$ in~$Y$.

By induction on~$k$, if $Y_{k-1}$ has the modified Eilenberg–Zilber property, then so does~$Y_k$.
Indeed, uniqueness follows from the uniqueness for~$Y$.
To establish existence, observe that any newly added element of~$Y_k$ has the form
$Y(γ)(x)$ for some morphism~$γ$ in~$Γ$ and indecomposable~$x∈Y(⟨k⟩)$,
and the required decomposition for such an element was constructed above.

We now show that each inclusion $Y_{k-1}\into Y_k$ is a weak equivalence.
We claim that it is a cobase change of a coproduct of boundary inclusions $∂⟨k⟩\into \Yo{⟨k⟩}$ (\cref{boundarygamma}), which are weak equivalences by \cref{gammacyc}. 

\medskip
\paragraph{\bf Step 2} (The proposed pushout.)
Denote by $Z_k$ the set of indecomposable elements in $Y_k(⟨k⟩)∖Y_{k-1}(⟨k⟩)$.
The symmetric group $\Sigma_k$ acts on $Y_k(⟨k⟩)$ and $Y_{k-1}(⟨k⟩)$ via automorphisms of $⟨k⟩∈Γ$.
Therefore, it also acts on $Y_k⟨(k⟩)∖Y_{k-1}(⟨k⟩)$ and hence on $Z_k$.
The action of $Σ_k$ on $Z_k$ is free by \cref{indecomposable.free}.
For $g∈Z_k$, the composition $∂⟨k⟩→\Yo{⟨k⟩}\lto3{g}Y_k$ factors through $Y_{k-1}$.
Indeed, elements of $∂⟨k⟩(⟨\ell⟩)$ are maps of finite pointed sets $⟨k⟩→⟨\ell⟩$ that have at most $k-1$ elements of $⟨k⟩$ map to non-basepoint elements of~$⟨\ell⟩$,
hence have the form $γδ$, where $δ$ is an inert map (of finite pointed sets) that is not an isomorphism.
Acting by such a map on~$g$ yields $Y(γ)(Y(δ)(g))$, where $Y(δ)(g)$ is an indecomposable element of degree $k-1$ or less, 
hence belongs to $Y_{k-1}$ by assumption.

Choose a single representative from each $\Sigma_k$-orbit in~$Z_k$ and denote the resulting set of representatives by~$S_k$.
We have a commutative diagram
$$\xymatrix{
\coprod_{g\in S_k}∂⟨k⟩\ar[r]\ar[d] & Y_{k-1}\ar@{^{(}->}[d]^{i}\\
\coprod_{g\in S_k}\Yo{⟨k⟩}\ar[r]^-{[g]_{g\in S_k}} & Y_k,}\eqlabel{push.sq}$$
which we now show to be cocartesian, which establishes that the right map is a weak equivalence if $k≥2$.
(The case $k=1$ was treated above: $Y_0=Y_1$ by assumption.)

\medskip 
\paragraph{\bf Step 3} (The pushout property.)
Fix $⟨\ell⟩\in \Gamma$ and evaluate at this object. 
We show that the bottom map in \cref{push.sq} induces a bijection of sets 
$$\coprod_{g\in S_k}\Yo{⟨k⟩}(⟨\ell⟩)∖∂⟨k⟩(⟨\ell⟩)→Y_k(⟨\ell⟩)∖Y_{k-1}(⟨\ell⟩).\eqlabel{cobase.change}$$
Elements in the left side are precisely pairs $(g,σ)$, where $g\in S_k$ and the map of sets $σ:⟨k⟩→⟨\ell⟩$ is active, i.e., satisfies $σ^{-1}\{\ast\}=\{\ast\}$. 
The displayed map sends $(g,σ)\mapsto Y(σ)(g)$.
If $Y(σ)(g)∈Y_{k-1}(⟨\ell⟩)$, then by the modified Eilenberg–Zilber property of~$Y_{k-1}$
we must have $g∈Y_{k-1}(⟨k⟩)$, contradicting the definition of~$g$.
Thus, $Y(σ)(g)∈Y_k(⟨\ell⟩)∖Y_{k-1}(⟨\ell⟩)$, which establishes the existence of the map \cref{cobase.change}.

It remains to show that any element $y\in Y_k(⟨\ell⟩)∖Y_{k-1}(⟨\ell⟩)$
comes from a unique pair $$(g,σ)\in \coprod_{g\in S_k}\Yo{⟨k⟩}(⟨\ell⟩)∖∂⟨k⟩(⟨\ell⟩).$$ 

As observed above, $y=Y(γ)(x)=Y(a)(Y(i)(x))$, for some indecomposable~$x∈Z_k$ and a morphism $γ=ia$ in~$Γ$,
presented as the composition of an active morphism~$a$ and an inert morphism~$i$.
If $i$ is not an isomorphism, then the indecomposable element $Y(i)(x)$ belongs to $Y_{k-1}$, hence $y∈Y_{k-1}$, contradicting the definition of~$y$.
Thus, we can assume $i=\id$, so that $y=Y(a)(x)$, where $x∈Z_k$.
Furthermore, $x$ is in the $Σ_k$-orbit of some unique $g∈S_k$, hence $x=Y(q)(g)$ for some isomorphism~$q$ in~$Γ$.
Now the pair $(g,σ)$ with $σ=qa$ establishes surjectivity.

Suppose $(g',σ')$ is another such pair.
By the uniqueness part of the modified Eilenberg–Zilber property,
there is a unique isomorphism $δ$ in~$Γ$ such that $g=Y(δ)(g')$.
Furthermore, $δσ=σ'$.
Thus, $g∈S_k$ and $g'∈S_k$ belong to the same $Σ_k$-orbit in~$Z_k$, hence $g=g'$.
Now $σ=σ'$ by \cref{indecomposable.free}.
\end{proof}

\section{Smooth bordism categories}
\label{bordcts}

In this section, we give a precise definition of the smooth bordism categories as smooth symmetric monoidal $(\infty,d)$-categories.
Specifically, for every geometric structure $\gs\in \Struct_d$ (\cref{geometric.structure}),
we define a corresponding smooth symmetric monoidal $(\infty,d)$-category $\Bord_{d}^\gs$ (\cref{bordstr}) that encodes bordisms equipped with geometric structures given by $\gs$.
The construction is manifestly functorial in $\gs$, which proves axiom \ax1.
In \cref{axioms.section}, we will show that $\Bord_d^\gs$ satisfies the additional axioms \ax2 and \ax3 in \cref{axioms}. 
Once this has been established, we will have proved \cref{axtheorem}.

In parallel, we also define a smooth symmetric monoidal $(∞,d)$-category $\BBord_{d}^\gs$ (\cref{enrichedbordstr})
that encodes bordisms with isotopies and uses enriched geometric structures in $\frakStruct_{d}$ (\cref{geometric.structure.isotopy}).
This proves axiom \frakax1 for $\BBord_d$.
In \cref{axioms.section}, we will prove that this version of the bordism category also satisfies axioms \frakax2 and \frakax3 in \cref{multiple.isotopy.axiom}.
This will imply \cref{axtheorem2}.

\begin{remark}
Our bordism categories $\Bord_d$ and $\BBord_{d}$ are not fibrant objects in $\smcat_{\infty,d}$,
although they are local with respect to some of the morphisms in \cref{globular.model.structure}.
We provide a summary of the fibrancy properties of our bordism categories. 
\begin{itemize}
\item For a geometric structure $\gs$ that satisfies homotopy descent in $\Struct_{d}$ (respectively $\frakStruct_{d}$), the categories $\Bord_d^\gs$ and $\Bord_{d,\uple}^\gs$ (respectively $\BBord_d^\gs$ and $\BBord_{d,\uple}^\gs$) are local with respect to the Čech morphisms  \cref{homotopy.descent} for finite open covers, Segal $\Gamma$-maps \cref{monoidal1}, Segal $\Delta$-maps \cref{segal1}. 
\item If the geometric structure $\gs$ does not satisfy homotopy descent, then all three of the above conditions may be violated (see \cref{one.dimensional.bordisms}). 
\item For any geometric structure $\gs$, $\Bord_d^\gs$ is also local with respect to the Segal completion maps \cref{segal2}.  If $d\geq 5$, and $\gs=\ast$, then $\BBord_d^\gs$ is not local with respect to the completion maps, due to the existence of nontrivial $h$-cobordisms.
\item Neither bordism category is local with respect to the map \cref{monoidal2}, since the ``trash bin'' need not be contractible. 
\item Both $\Bord_d$ and $\BBord_d$ are local with respect to the globular maps in \cref{globular.maps}, while $\Bord_{d,\uple}$ and $\BBord_{d,\uple}$ are not. 
\end{itemize} 
If desired, one could make simple adjustments to the definitions of bordism categories,
requiring the trash bin to be empty (with appropriate adjustment of the $Γ$-maps),
and (in the case of $\BBord_d$)
adding h-cobordisms to the simplicial set of objects (as in Lurie \cite[Definition~2.2.10]{Lurie.TFT} and Calaque–Scheimbauer \cite[Definition~5.24]{CalaqueScheimbauer}).
These simple adjustments produce an equivalent bordism category, which satisfies all locality conditions of \cref{globular.model.structure} if $\gs$ satisfies homotopy descent.
We remark that the remaining fibrancy condition, that of \emph{injective} fibrancy, is typically false unless $\gs=∅$.
\end{remark}

\begin{remark}
Working with nonfibrant models for the bordism category gives rise to significant simplifications in the proof of locality.
In fact, the most important geometric structure in the proof of locality is the case of geometric framings,
which do not satisfy descent and give rise to a nonfibrant bordism category, as illustrated in \cref{one.dimensional.bordisms}.
Nonfibrancy should not be seen as a disadvantage in this case, but rather as an advantage that allows us to work with a much smaller model for the bordism category in our proofs. 
\end{remark}

\begin{remark}
\label{small}
To make all sites small and to make presheaves valued in sets, rather than proper classes, we require that the underlying set of any manifold is a subset of $\RR$.
We do not require any compatibility with the smooth structure on $\RR$. 
\end{remark}

\subsection{Geometric structures}
\label{geostr}
We now describe how to encode geometric structures on bordisms.
As pointed out in the introduction, our treatment generalizes the traditional treatment of tangential structures. 
Roughly, the passage from the traditional approach to our approach is given by taking the sheaf of sections of the homotopy pullback 
$$\xymatrix{
Y\times_{\tdeloop\GL(d)}M\ar[r]\ar[d] & Y\ar[d]^-{\xi}\\
M\ar[r]^-{\tau} & \tdeloop\GL(d),
}$$
where $\tau:M\to \tdeloop\GL(d)$ is the classifying map for the tangent bundle of~$M$. 
There are two subtleties which need to be accounted for. 
First, we need to encode geometric structures.
Hence, we need to replace the classifying space $\tdeloop\GL(d)$ with the classifying \emph{stack} of vector bundles of rank $d$. 
Second, the sheaf must live on a larger site than the open subsets of $M$, since we will need to pullback the structure along all open embeddings.
Moreover, we will need all geometric structures to vary in families, parametrized over cartesian spaces.
We resolve these issues by working with sheaves on the site of submersions (with $d$-dimensional fibers), with fiberwise open embeddings between them.

Recall the cartesian site $\stcart$ from \cref{def.cartsp}, including the definition of submersion.

\begin{definition}
\label{fembdef}
Fix $d≥0$. Recall \cref{def.cartsp} which introduces submersions of structured manifolds with $d$-dimensional fibers, pullbacks of submersions with $d$-dimensional fibers, and open embeddings of structured manifolds. Let $\FEmb_d$ be the site whose objects are morphisms $p:M\to U$ in $\stman$
such that $p$ is a submersion with $d$-dimensional fibers and $U\in\stcart$.
Morphisms $(p:M\to U)\to (q:N\to V)$ are pairs of morphisms $(f:M\to N,g:U\to V)$ in $\stman$ such that the diagram 
$$\xymatrix{
M\ar[r]^-{f}\ar[d]_-{p}  & N\ar[d]^-{q}\\
U\ar[r]^-{g} & V
}$$
commutes and such that $f$ is a fiberwise open embedding covering $g$, i.e., the map $M\to g^*N$ is an open embedding. 
Covering families are given by a collection of morphisms
$$\left\{\vcenter{
\xymatrix{
M_{\alpha}\ar[r]^-{i_{\alpha}}\ar[d] & M\ar[d]\\
U_{\alpha}\ar[r]_{j_{\alpha}} & U
}}
\right\}$$ 
such that both horizontal maps $i_{\alpha}$ and $j_{\alpha}$ are open embeddings 
and the collection $\{i_{\alpha}\}$ is a covering family of~$M$.
We do not require that $\{j_{\alpha}\}$ is a covering family of~$U$. 
\end{definition}

\begin{definition}
\label{geometric.structure}
A {\it fiberwise $d$-dimensional geometric structure\/} is a simplicial presheaf on $\FEmb_d$.
We denote the left Bousfield localization of the injective model structure at Čech nerves of the covers in \cref{fembdef} by 
$$\Struct_d≔\sPSh(\FEmb_d)_{\inj,\Cech}.$$
The $\sset$-enriched Bousfield localization exists by \cref{left.Bousfield.localization.exists} and  \cref{enrichedbous}.
\end{definition}

\begin{remark}
Alternatively, we could consider sheaves on manifolds and fiberwise etale maps. 
The canonical inclusion $\FEmb_d\into {\sf FEtale}_d$ exhibits $\FEmb_d$ as a dense subsite in the 1-categorical sense, i.e., it induces an equivalence of 1-categories of sheaves.
This continues to hold at the level of simplicial presheaves, which is a stronger property (not every 1-dense subsite is $\infty$-dense).
Hence, we are free to use either site. 
\end{remark}

\begin{remark}
Sheaves on the site $\FEmb_d$ were considered (in an equivalent reformulation) by Nijenhuis \cite{Nijenhuis},
who used them to define natural bundles and natural mappings in differential geometry.
The latter are closely related to the geometric cobordism hypothesis, as shown in Grady–Pavlov \cite{GradyPavlov.GCH}.
The case of nonfiberwise geometric structures (using the Quillen equivalent model category
of simplicial sheaves on the site of manifolds and etale maps) was considered by Freed–Teleman \cite[Appendix~A]{FreedTeleman}.
The special case of geometric structures valued in groupoids instead of simplicial sets is considered by Ludewig–Stoffel \cite{LudewigStoffel},
where sheaves on the site of submersions with fiberwise etale maps are considered.
By the above remarks, our geometric structures are a natural generalization of those considered in Ludewig–Stoffel \cite{LudewigStoffel}.
Ayala–Francis–Tanaka \cite{AyalaFrancisTanaka} investigate the enriched (\cref{geometric.structure.isotopy}) and unenriched variants of presheaves on $\Emb_d$.
\end{remark}

\begin{definition}
\label{fembdelta}
We define the enriched site $\frakFEmb_d$, with enrichment in $\smset$, as follows.
The objects are the same objects as in $\FEmb_d$.
Given two objects $W→U$ and $W'→U'$,
the corresponding hom-object between them is a smooth set
whose $L$-points (with $L\in \cart$) are given by
\begin{itemize}
\item morphisms in $\FEmb_d$
of the form $(f:W⨯L→W',u:U⨯L→U')$ with source $(W⨯L→U⨯L)=(W→U)⨯\id_L$ and target $W'→U'$,
where we implicitly embed $\cart→\stcart$ using \cref{def.cartsp}.
We require that the map $u$ is the composition of the projection $U⨯L\to U$ and a map $U\to U'$. 
\end{itemize}
The Grothendieck topology of $\frakFEmb_d$ is generated by the same coverage as for $\FEmb_d$.
\end{definition}

\begin{definition}
\label{geometric.structure.isotopy}
Given $d≥0$, a {\it fiberwise $d$-dimensional geometric structure with isotopies\/} is an enriched presheaf in the $\smsset$-enriched model category 
$$\frakStruct_{d}≔\PSh(\frakFEmb_d,\smsset)_{\inj,\Cech},$$
where we promote the $\smset$-enrichment of $\frakFEmb_d$ to a $\smsset$-enrichment, by taking the constant smooth simplicial set and the model structure on $\smsset$ is defined in \cref{smsets}.
The enriched left Bousfield localization exists by \cref{left.Bousfield.localization.exists} and \cref{enrichedbous}. 
\end{definition}

Our definition of geometric structures is extremely versatile and captures all significant geometric structures we can think of, 
including metrics, topological structures, tangential structures, etc.

\begin{example}
Let $\stcart=\cart$. Let $X$ be a smooth manifold of any dimension.
We can regard $X$ as a geometric structure via the sheaf that sends 
$$(M\to U)\mapsto \sm(M,X),$$
for a submersion $M\to U$.
This is clearly a sheaf on $\FEmb_d$ since the total space functor $\FEmb_d\to \Man$, which maps $M\to U$ to $M$, sends covering families to covering families.  
This sheaf is not representable.
\end{example}

\begin{example}
The previous example can be generalized easily to all simplicial presheaves on the site $\stman$.
Indeed, the total space functor $T:\FEmb_d\to \stman$ induces a restriction functor
$$T^*:\sPSh(\stman)_{\Cech}\to \Struct_d,$$
which manifestly preserves the homotopy descent property.
\end{example}

\begin{example}
\label{riemmet}
Let $\stcart=\cart$. An example of a geometric structure that does not come from a simplicial presheaf on smooth manifolds and smooth maps
is given by the presheaf of fiberwise Riemannian metrics. 
Let
$$\FRiem_d:\FEmb_d^{\op}\to \set$$
be the presheaf that sends a submersion $p:M\to U$ with $d$-dimensional fibers to the set of metrics on the fiberwise tangent bundle $\ker Tp$ over~$M$.
A morphism $(f,g):(p:M\to U)\to (q:N\to V)$ is sent to the function 
$$(f,g)^*:\FRiem_d(q)\to \FRiem_d(p)$$
that sends a metric $g$ on the fiberwise tangent bundle $\ker Tq$ to the pullback metric $$(f^*g)_{x}(v,w)=g_{f(x)}(T_x f(v),T_x f(w)),$$ for all $x\in M$.
This is a well-defined metric on $\ker Tq$ since $f$ is a fiberwise etale map and for all $v\in \ker(T_x p)$, we have
$$(T_{f(x)}q\circ T_x f)(v)=(T_{p(x)}g\circ T_x p)(v)=0\qquad\Longrightarrow\qquad T_x f(v)\in \ker(T_{f(x)}q).$$

We can also consider Riemannian metrics with restrictions on sectional curvature (e.g., positive, negative, nonpositive, nonnegative),
since these properties are preserved by pullbacks along etale maps.
For example, we define the subobject of positive sectional curvature metrics $\FRiem_{d,{\rm sc}>0}\subset \FRiem_d$
as the functor that sends $p:M\to U$ to the subset of metrics on the fiberwise tangent bundle such that for all $u\in U$,
the metric $g\vert_{p^{-1}(u)}$ on $T_u M\cong \ker(Tp\vert_{p^{-1}(u)})\to M_u$ is a metric of positive scalar curvature.
Such metrics can again be pulled back by fiberwise embeddings and the property of having positive sectional curvature is preserved under such pullbacks.
\end{example}

\begin{example}
\label{friemex}
Given a simplicial presheaf $\gs\in \Struct_d$,
we can obtain an object in $\frakStruct_d$ by enriched left Kan extension along the enriched functor $\FEmb_d\into \frakFEmb_d$
(where $\FEmb_d$ is enriched by taking the constant smooth set).
In this way, all the previous examples yield corresponding geometric structures in $\frakStruct_d$. 

For example, applying the enriched left Kan extension to the sheaf $\FRiem_d$,
we get an enriched presheaf $\frakFRiem_d\in \frakStruct_d$
such that $\frakFRiem_d(M\to U)$ is the smooth set whose $L$-points are given by
equivalence classes of pairs $(g,φ)$, where $g$ is a fixed fiberwise Riemannian metric on some $(N\to U)∈\FEmb_d$,
$φ$ is an $L$-point of $\frakFEmb_d(M→U,N→U)$
given by a family of fiberwise embeddings $\varphi:M\times L\to N$ covering the projection $U\times L\to U$,
and the equivalence relation is generated by $(ψ^*g,φ)\sim(g,ψφ)$, where $ψ:(P→U)→(N→U)$ is a morphism in $\FEmb_d$
and $φ∈\frakFEmb_d(M→U,P→U)(L)$.
This smooth set is different from the smooth set of $U⨯L$-parametrized families of Riemannian metrics on $M⨯L$.
In particular, the $U⨯L$-families of metrics in $\frakFRiem_d$ are ``$d$-thin''
in that the parametrizing map to the moduli stack must factor through a fixed $d$-dimensional manifold.

For $d=1$, it turns out that we can identify the smooth set $\frakFRiem_1(\RR\times U\to U)$ (up to $\RR$-local weak equivalence)
with $\tdeloop \sm(U, \RR)$ (the constant smooth simplicial set),
where $\tdeloop \sm(U, \RR)$ is the delooping of smooth real-valued functions (with the additive group structure).
For details on this calculation, see \gcref{shapemetrics}.
\end{example}

\subsection{Cuts, cut tuples, and cut grids}

In this section, we define the notions of a cut, cut tuple, and cut grid on an object $M\to U\in \FEmb_d$,
which are responsible for implementing the structure of a $d$-category for bordisms.
The notions of a cut (\cref{cut}) and cut tuple (\cref{cut.m.tuple})
are inspired by Stolz–Teichner \cite[Definition~2.21]{StolzTeichner.SUSY}.
The notion of a cut grid (\cref{cutgrid}) roughly resembles the constructions of Lurie \cite[Definition~2.2.9]{Lurie.TFT}
and Calaque–Scheimbauer \cite[Definition~5.1]{CalaqueScheimbauer}.
The notion of a globular cut grid is inspired by Henriques \cite[\S2.2]{Henriques}.

\begin{remark}
Below we work with families of structured manifolds (such as supermanifolds) indexed by $\stcart$ (\cref{fembdef}). 
We will need to use concepts such as open subsets, smooth functions, submersions, and regular values in this context. By these terms, we will always mean the corresponding notion on the reduced manifold obtained by applying the functor $\stu:\stman\to \man$. 
For the reader who is not concerned with the general setting of structured manifolds, simply take $\stu=\id$. In this case, all manifolds are ordinary smooth manifolds and terms such as open set, submersion, etc., have the usual meaning. 
\end{remark}

\begin{definition}
\label{cut}
A \emph{cut} of an object $p:M\to U$ in $\FEmb_d$ is a triple $(C_{<},C_=,C_>)$ of subsets of $\stu(M)$ such that there is a smooth map $h:\stu(M)\to \RR$ (called the \emph{height function}) whose fiberwise-regular values (i.e., points that are regular values of the restriction to each fiber) form an open neighborhood of $0$. 
Moreover, $h^{-1}(-\infty,0)=C_{<}$, $h^{-1}(0)=C_=$ and $h^{-1}(0,\infty)=C_>$. 
We set 
$$C_\le=C_<\cup C_=, \qquad C_{\ge}=C_>\cup C_=.$$
We equip the set of cuts with a natural ordering $\leq$, with $C\leq C'$ if and only if $C_{\leq}\subset C'_{\leq}$. 
There is a functor $\Cut:\FEmb_d^\op\to \poset$ that associates to an object $p:M\to U$ its poset of cuts, and to a morphism the induced map of posets that takes preimages of cuts.
\end{definition}

\begin{definition}
\label{cut.m.tuple}
Fix $d\geq 0$, a simplex $[m]\in \Delta$, and an object $p:M\to U$ in $\FEmb_d$. 
A \emph{cut} $[m]$-\emph{tuple} $C$ for $p:M\to U$ is a collection of cuts $C_j=(C_{<j},C_{=j},C_{>j})$ of $p:M\to U$ indexed by vertices $j\in [m]$ such that 
$$C_0\leq C_1\leq \cdots \leq C_m.$$
Given $j≤j'$, we set 
$$C_{(j,j')}=C_{>j}\cap C_{<j'}, \qquad C_{[j,j']}=C_{\ge j}\cap C_{\le j'}.$$
We also denote by
$$C_{⟨j,j'⟩}$$
the cut $[j'-j]$-tuple obtained from~$C$ by removing the cuts $C_k$ with $k<j$ or $k>j'$.

The functor $\Cut$ naturally extends to a functor 
$$\Cut:\Delta^\op\times \FEmb_d^\op\to \set$$
as follows. To an object $([m],p:M\to U)$, we associate the set of cut $[m]$-tuples of $p$. 
\begin{itemize}
\item For fixed $[m]\in \Delta$ and a morphism $(f,g):(M\to U)\to (N\to V)$, the corresponding structure map
$$\Cut([m],N\to V)\to \Cut([m],M\to U)$$
is given by sending a cut $[m]$-tuple on $N\to V$ to the cut $[m]$-tuple $(f,g)^*C$ on $M\to U$, where $$(f,g)^*C_j=(\stu(f)^{-1}(C_{<j}),\stu(f)^{-1}(C_{=j}), \stu(f)^{-1}(C_{>j}))$$ for all $j\in [m]$.
That this is a well defined cut tuple follows from the fact that $f$ is a fiberwise embedding.
In particular, if $h_j$ is a height function for the cut $C_j$, then $h_j\circ \stu(f)$ is a height function for the cut $(f,g)^*C_j$.

\item For fixed $M\to U\in \FEmb_d$ and a coface map $d^j:[m-1]\to [m]\in \Delta$, the corresponding face map 
$$d_j:\Cut([m],M\to U)\to \Cut([m-1],M\to U)$$
removes the $j$th cut in the cut $[m]$-tuple $C$.
A codegeneracy $s^j:[m+1]\to [m]$ is send to the map that duplicates the $j$th a cut in the tuple. 
\end{itemize}
\end{definition}

The above definition implies that
$$C_{>0}\supset C_{>1}\supset \cdots \supset C_{>j}\supset \cdots\supset C_{>m},$$
as well as the analogous chains for $C_{\le}$ and $C_{\ge}$.

\begin{figure}[ht]
\begin{center}
\begin{tikzpicture}[scale=.65]
\draw (0,0) to [out=20,in=120] (5,2) to [out=-60, in=110] (9,1);
\draw (0,2) to [out=10,in=115] (5,2) to [out=300, in=145] (9,3);
\draw (0,2.5) to [out=0, in=125] (6,3) to [out=305, in=110] (9,3.5);
\node at (0-.7,0) {$C_{=0}$};
\node at (0-.7,2) {$C_{=1}$};
\node at (0-.7,2.5) {$C_{=2}$};
\end{tikzpicture}
\end{center}
\caption{A cut $[2]$-tuple $C_{j}=\{(C_{<j},C_{=j},C_{>j})\}_{j\in [2]}$ on $\RR^2$.
The cuts $C_{=0}$ and $C_{=1}$ intersect and the region $C_{\leq 0}\subset C_{\leq 1}$.
Cuts are not allowed to intersect transversally, as this would violate the ordering in \cref{cut.m.tuple}.}
\end{figure}

Next, we will extend the functor $\Cut$ further to a functor 
$$\tCut:(\Delta^{\times d})^\op\times \FEmb_d^\op\to \set,$$
where the subscript notation will become apparent in a moment.
Morally, $\tCut$ sends a multisimplex~$[m]$ to the collection of cut $[m_i]$-tuples, where $i\in \{1,\ldots, d\}$.
However, we do not want arbitrary collections of cut tuples.
Instead, we only take those cut tuples that intersect transversally in directions indexed by different elements of $\{1,\ldots,d\}$.
This motivates the following definition.

\begin{definition}
\label{cutgrid}
Fix $d\geq 0$.
We define the \emph{cut grid functor} 
$$\tCut:(\Delta^{\times d})^\op\times \FEmb_d^\op\to \set,$$
as follows. Let $p:M\to U$ be an object in $\FEmb_d$ and let ${\bf m}=([m_1],\ldots,[m_d])\in \Delta^{\times d}$ be a multisimplex.
The set $\tCut({\bf m},p)$ has elements:
\begin{itemize}
\item For each $1\leq i\leq d$, a cut $[m_i]$-tuple $C^i$ on $p:M\to U$,
\end{itemize}
which satisfy the transversality property: 
\begin{enumerate}
\item[$\pitchfork$.] For every subset $S\subset \{1,\ldots,d\}$ and for any $j:S→\ZZ$ such that $0\leq j_i\leq m_i$ for all $i∈S$, there is a smooth map $h_j:\stu(M)\to \RR^S$ such that for any $i\in S$, the map $$\pi_i\circ h_j:\stu(M)\to  \RR,$$
where $\pi_i:\RR^S\to \RR$ is the $i$th projection, yields the $j_i$-th cut $C_{j_i}^i$ in the cut tuple $C^i$, as in \cref{cut}. 
We require that the fiberwise-regular values of $h_j$ form an open neighborhood of $0$ in $\RR^S$. 
\end{enumerate}
The structure maps for the $i$th factor of $\Delta$ in $\Delta^{\times d}$ are given by applying a simplicial map to the corresponding cut $[m_i]$-tuple. We observe that the transversality property is still satisfied after applying a simplicial map, since simplicial maps simply remove cuts or duplicate cuts,
i.e., the function $h_j$ in the above definition can be left unchanged. The structure map for the factor $\FEmb_d$ is given by applying the corresponding structure map to every cut tuple.
We call the elements of $\tCut({\bf m},p)$ \emph{cut ${\bf m}$-grids}.
\end{definition}

\begin{notation}
\label{coreandbead}
Recall the notation in \cref{cut.m.tuple}. For a cut ${\bf m}$-grid $C$, a subset $S\subset \{1,\hdots,d\}$ and $j,j':S→\ZZ$ satisfying $0\leq j_i\leq j'_i\leq m_i$ for all $i∈S$,
we use the following notation.
\begin{itemize}
\item We define $$C_{[j,{j'}]}≔\bigcap_{i\in S} C_{[j_i,j'_i]}\subset \stu(M)\qquad C_{(j,{j'})}≔\bigcap_{i\in S} C_{(j_i,j'_i)}\subset \stu(M).$$
These subsets are the regions between corresponding cuts in various directions.
\item We define $C_{⟨j,{j'}⟩}$ to be the cut ${\bf m}'$-grid,
whose $i$th cut tuple is $C^i$ if $i∉S$
and $C^i_{⟨j_i,j'_i⟩}$ if $i∈S$,
with ${\bf m}'$ defined accordingly.
\end{itemize}
\end{notation}

\begin{definition}
\label{compactgrid}
Let $C$ be a cut ${\bf m}$-grid for $p:M\to U$.
We say that $C$ is \emph{compact} if for $S=\{1,2,\hdots,d\}$, the restriction of $p$ to $C_{[j,j']}$
is proper, for all $j,j':S\to \ZZ$.  
\end{definition}

\cref{cutgrid} does indeed implement transversality of cuts in different directions, as illustrated in \cref{transcut}. 

\begin{figure}[ht]
\begin{tikzpicture}[scale=.60]
\draw (0,0) to [out=20,in=120] (5,2) to [out=-60, in=110] (9,1);
\draw (0,2) to [out=10,in=115] (5,2) to [out=300, in=145] (9,3);
\draw (0,2.5) to [out=0, in=125] (6,3) to [out=305, in=110] (9,3.5);
\node at (0-.7,0) {$C^1_{=0}$};
\node at (0-.7,2) {$C^1_{=1}$};
\node at (0-.7,2.8) {$C^1_{=2}$};
\draw (1,6) to [out=270, in=110] (2,3) to [out=290, in=110] (1,-1);
\draw (1.5,6) to [out=-45, in=110] (2,3) to [out=290, in=110] (1,-1);
\draw (3,6) to [out=-45, in=110] (5,3) to [out=290, in=110] (1,-1); 
\draw (4,6) to [out=-90, in=110] (6,3) to [out=290, in=90] (4,-1);  
\draw (5,6) to [out=200, in =180] (9,4);
\node at (.6,6.5+.3) {$C_{=0}^2$};
\node at (1.6,6.5+.3) {$C_{=1}^2$};
\node at (3,6.5+.3) {$C_{=2}^2$};
\node at (4,6.5+.3) {$C_{=3}^2$};
\node at (5.2,6.5+.3) {$C_{=4}^2$};
\end{tikzpicture}
\caption{A cut $([2],[4])$-grid on $\RR^2$. }
\end{figure}
The notion of a cut ${\bf m}$-grid will allow us to define a higher categorical structure on bordisms which is the $d$-fold generalization of a double category. To get the correct generalization of a bicategory, we will need a \emph{globular} version of the above cut ${\bf m}$-grids. Our cut ${\bf m}$-grids are general enough that we can simply extract those cut grids that satisfy the globular condition. 

\begin{definition}
\label{globgrid}
Fix $d\geq 0$. Let $p:M\to U$ be an object in $\FEmb_d$ and let ${\bf m}=([m_1],\ldots,[m_d])\in \Delta^{\times d}$ be a multisimplex. Let $m:\{1,\hdots,d\}\to \ZZ$ be a function. Let $0:\{1,\hdots,d\}\to \ZZ$ be the function $0(i)=0$. Let $v^j_i:[0]\to [m_i]$, $j\in [m_i]$, be the map that picks out the $j$th vertex and let $v_j^i$ denote the induced structure map. Recall the notation of \cref{cutgrid}. Let 
$$\tCutglob({\bf m},p:M\to U)\subset \tCut({\bf m},p:M\to U)$$
be the subset whose cut ${\bf m}$-grids $C$ satisfy the following property:
\begin{itemize}
\item Let ${\bf m}'$ be the multisimplex whose simplices are the same as those of ${\bf m}$ in all directions except the $i$th direction, where we set $[m'_i]=[0]$. For all $j\in [m_i]$, the cut grid
$$v^i_jC\in \tCut({\bf m}',p:M\to U),$$ 
admits an open neighborhood $N$ of $C_{[0,{\bf m}']}⊂\stu(M)$ such that the restriction of $v^i_jC$ to $p\vert_{N}:N\to U$  
is a simplicial degeneration of a cut grid in 
$$\tCut(([m_1],\ldots,[m_{i-1}],[0],\ldots,[0]), p\vert_{N}:N\to U),$$
in the simplicial directions $i+1,i+2,\ldots,d$.  
\end{itemize}
The above property is preserved under pullback of cut grids, so that $\tCutglob$ defines a functor on $(\Delta^{\times d})^{\op}\times \FEmb_d$.  
We call an element of $\tCutglob({\bf m},p:M\to U)$ a \emph{globular cut ${\bf m}$-grid}. If a globular cut ${\bf m}$-grid is compact, we call it a \emph{compact globular ${\bf m}$-grid}. 
\end{definition}

The next example illustrates cut ${\bf m}$-grids and globular cut ${\bf m}$-grids in the case $d=2$. 

\begin{example}
Set $d=2$.
The following images depict cut $([1],[0])$- and $([1],[1])$-grids on a $2$-manifold given by the gray sheet (as an object in $\FEmb_2$, we take the base of the submersion to be trivial). The image on the left is a cut $([1],[0])$-grid and the image in the center is a cut $([1],[1])$-grid. The image on the right depicts a globular cut $([1],[1])$-grid. 
\begin{center}
\begin{tikzpicture}[scale=.40]
\fill[fill=gray!25, rounded corners] (0+.3,-1+.3) -- (0+.3,6+.3) -- (9+.3,6+.3) -- (9+.3,-1+.3) -- (0+.3,-1+.3);
\draw[ blue] (1+.3,6+.3) to [out=270, in=110] (2+.3,3+.3);
\draw[dashed, rounded corners]  (0+.3,-1+.3) -- (0+.3,6+.3) -- (9+.3,6+.3) ;
\fill[fill=gray!25, rounded corners] (0,-1) -- (0,6) -- (9,6) -- (9,-1) -- (0,-1);
\draw[dashed, rounded corners] (9.7,-.5) to[out=270,in=0, looseness=.8] (9,-1) -- (0,-1) -- (0,6) -- (9,6);
\fill[fill=gray!25, dashed] (9,6) -- (9,-1) to[out=0,in=270, looseness=1.2] (9.7,.5) -- (9.7,6.2);
\filldraw[fill=gray!25, dashed] (9+.3,6+.3) to [out=0, in=0, looseness=4] (9,6);
\draw [purple] (0,0) to [out=20,in=120] (5,2) to [out=-60, in=130] (9,1) to [out=310, in=-90] (9.7,1.5);
\draw [purple] (0,2.5) to [out=0, in=125] (6,3) to [out=305, in=110] (9,3.5) to [out=300, in=-130] (9.7,3.7);
\node at (0-.7,0) {$\scriptstyle C^2_{=0}$};
\node at (0-.7,2.8) {$\scriptstyle C^2_{=1}$};
\draw[ blue] (1,6) to [out=270, in=110] (2,3) to [out=290, in=110] (1,-1);
\node at (.6,6.5+.3) {$\scriptstyle C_{=0}^1$};
\draw[fill=white] (2.65-.2,4.6) to [out=57, in=118] (3.9-.2,4.6);
\draw[thick] (2.5-.2,4.5) to [out=60, in=115] (4-.2,4.5);
\draw[thick, fill=white] (2.65-.2,4.6) to [out=-45, in=-135] (3.9-.2,4.6);
\draw[fill=white]  (6.5+.2,.5+.2) to [out=46, in=130] (8.5-.2,.5+.2);
\draw[thick] (6.5,.5) to [out=60, in=115] (8.5,.5);
\draw[thick, fill=white] (6.5+.2,.5+.2)  to [out=-45, in=-135] (8.5-.2,.5+.2);
\end{tikzpicture}
\begin{tikzpicture}[scale=.40]
\fill[fill=gray!25, rounded corners] (0+.3,-1+.3) -- (0+.3,6+.3) -- (9+.3,6+.3) -- (9+.3,-1+.3) -- (0+.3,-1+.3);
\draw[ blue] (1+.3,6+.3) to [out=270, in=110] (2+.3,3+.3);
\draw[ blue] (3+.3,6+.3) to [out=-45, in=110] (5+.3,3+.3); 
\draw[dashed, rounded corners]  (0+.3,-1+.3) -- (0+.3,6+.3) -- (9+.3,6+.3) ;
\fill[fill=gray!25, rounded corners] (0,-1) -- (0,6) -- (9,6) -- (9,-1) -- (0,-1);
\draw[dashed, rounded corners] (9.7,-.5) to[out=270,in=0, looseness=.8] (9,-1) -- (0,-1) -- (0,6) -- (9,6);
\fill[fill=gray!25, dashed] (9,6) -- (9,-1) to[out=0,in=270, looseness=1.2] (9.7,.5) -- (9.7,6.2);
\filldraw[fill=gray!25, dashed] (9+.3,6+.3) to [out=0, in=0, looseness=4] (9,6);
\draw [purple] (0,0) to [out=20,in=120] (5,2) to [out=-60, in=130] (9,1) to [out=310, in=-90] (9.7,1.5);
\draw [purple] (0,2.5) to [out=0, in=125] (6,3) to [out=305, in=110] (9,3.5) to [out=300, in=-130] (9.7,3.7);
\node at (0-.7,0) {$\scriptstyle C^2_{=0}$};
\node at (0-.7,2.8) {$ \scriptstyle C^2_{=1}$};
\draw[ blue] (1,6) to [out=270, in=110] (2,3) to [out=290, in=110] (1,-1);
\draw[ blue] (3,6) to [out=-45, in=110] (5,3) to [out=290, in=110] (1,-1); 
\node at (.6,6.5+.3) {$\scriptstyle C_{=0}^1$};
\node at (3,6.5+.3) {$ \scriptstyle C_{=1}^1$};
\draw[fill=white] (2.65-.2,4.6) to [out=57, in=118] (3.9-.2,4.6);
\draw[thick] (2.5-.2,4.5) to [out=60, in=115] (4-.2,4.5);
\draw[thick, fill=white] (2.65-.2,4.6) to [out=-45, in=-135] (3.9-.2,4.6);
\draw[fill=white]  (6.5+.2,.5+.2) to [out=46, in=130] (8.5-.2,.5+.2);
\draw[thick] (6.5,.5) to [out=60, in=115] (8.5,.5);
\draw[thick, fill=white] (6.5+.2,.5+.2) to [out=-45, in=-135] (8.5-.2,.5+.2);
\end{tikzpicture}
\begin{tikzpicture}[scale=.40]
\fill[fill=gray!25, rounded corners] (0+.3,-1+.3) -- (0+.3,6+.3) -- (9+.3,6+.3) -- (9+.3,-1+.3) -- (0+.3,-1+.3);
\draw[ blue] (1+.3,6+.3) to [out=270, in=110] (2+.3,3+.3);
\draw[ blue] (1.5+.3,6+.3) to [out=-45, in=110] (2+.3,3+.3);
\draw[dashed, rounded corners]  (0+.3,-1+.3) -- (0+.3,6+.3) -- (9+.3,6+.3) ;
\fill[fill=gray!25, rounded corners] (0,-1) -- (0,6) -- (9,6) -- (9,-1) -- (0,-1);
\draw[dashed, rounded corners] (9.7,-.5) to[out=270,in=0, looseness=.8] (9,-1) -- (0,-1) -- (0,6) -- (9,6);
\fill[fill=gray!25, dashed] (9,6) -- (9,-1) to[out=0,in=270, looseness=1.2] (9.7,.5) -- (9.7,6.2);
\filldraw[fill=gray!25, dashed] (9+.3,6+.3) to [out=0, in=0, looseness=4] (9,6);
\draw [purple] (0,0) to [out=20,in=120] (5,2) to [out=-60, in=130] (9,1) to [out=310, in=-90] (9.7,1.5);
\draw [purple] (0,2.5) to [out=0, in=125] (6,3) to [out=305, in=110] (9,3.5) to [out=300, in=-130] (9.7,3.7);
\node at (0-.7,0) {$\scriptstyle C^2_{=0}$};
\node at (0-.7,2.8) {$\scriptstyle C^2_{=1}$};
\draw[ blue] (1,6) to [out=270, in=110] (2,3) to [out=290, in=110] (1,-1);
\draw[ blue] (1.5,6) to [out=-45, in=110] (2.05,2.9) to [out=-45, in =110] (3,2.8) to [out=-90, in=80] (1.3,.5) to [out=-90,in=90] (3,-1);
\node at (.6,6.5+.3) {$\scriptstyle C_{=0}^1$};
\node at (2,6.5+.3) {$\scriptstyle C_{=1}^1$};
\draw[fill=white] (2.65-.2,4.6) to [out=57, in=118] (3.9-.2,4.6);
\draw[thick] (2.5-.2,4.5) to [out=60, in=115] (4-.2,4.5);
\draw[thick, fill=white] (2.65-.2,4.6) to [out=-45, in=-135] (3.9-.2,4.6);
\draw[fill=white]  (6.5+.2,.5+.2) to [out=46, in=130] (8.5-.2,.5+.2);
\draw[thick] (6.5,.5) to [out=60, in=115] (8.5,.5);
\draw[thick, fill=white] (6.5+.2,.5+.2) to [out=-45, in=-135] (8.5-.2,.5+.2);
\end{tikzpicture}
\end{center}
\end{example}

\subsection{Categories of bordisms}

In this section, we define both the uple and globular $d$-categories of bordisms,
in the terminology of Calaque–Scheimbauer \cite{CalaqueScheimbauer}.
The $d$-uple bordism category is \emph{not} local with respect to the globular morphisms \eqref{glob}
and hence there are distinguished composition directions, which need to be accounted for,
just as a double category has distinguished composition directions that are put on equal footing. 
The formalization of smoothness for bordism categories is due to Stolz–Teichner \cite{StolzTeichner.Elliptic, StolzTeichner.SUSY}.

\begin{definition}[$d$-uple/globular bordisms]
\label{bord}
Given $d≥0$, we specify an object $\Bord_{d,\uple}$, respectively $\Bord_{d}$, in the category
$\smcatuple_{\infty,d}$ (\cref{multiple.model.structure}), respectively $\smcat_{\infty,d}$ (\cref{globular.model.structure}), as follows.
For an object $(U,\langle \ell\rangle,{\bf m})\in \stcart\times \Gamma\times \Delta^{\times d}$,
the simplicial set $\Bord_{d,\uple}(U,\langle \ell\rangle,{\bf m})$ (respectively $\Bord_{d}(U,\langle \ell\rangle,{\bf m})$) 
is the nerve of the following category, which is small by \cref{small}. 

\medskip
\paragraph{\it Objects}
An object of the category is a {\it bordism\/} given by the following data.
\begin{enumerate}
\item[(1)] A $d$-dimensional smooth manifold $M$ (possibly open).
\item[(2)] A compact (and globular, in the case of $\Bord_d$) cut ${\bf m}$-grid $C$ (\cref{cutgrid}) for the projection map $p: M\times U\to U$.
\item[(3)] A choice of map $P:M\times U\to \langle \ell\rangle$.
This gives a partitions $M\times U$, and therefore $M$, into $\ell$ disjoint subsets
and another subset corresponding to the basepoint (the \emph{trash bin}).
\end{enumerate}
We also set
$$\core(M,C,P,j\leq j')= C_{[j,j']} \setminus \stu(P)^{-1}\{\ast\}.$$
In the case where $S=\{1,\ldots,d\}$, $j=0$ and $j'(i)=m_i$, we set 
$$\core(M,C,P)= C_{[j,j']} \setminus \stu(P)^{-1}\{\ast\} \eqlabel{core}$$
and call it the \emph{core} of the bordism $(M,C,P)$.

\medskip 
\paragraph{\it Morphisms}
We define a morphism $\varphi:(M,C,P)\to (M',C',P')$ (called a \emph{cut-respecting embedding}) as a morphism $\varphi:(M\times U\to U)\to (M'\times U\to U)$ in $\FEmb_d$ covering $\id_{U}$ such that
\begin{enumerate}
\item[(m)] $\varphi^*C'=C$, $\varphi^*P'=P$, and the image of $\stu(\varphi)$ contains $\core(M',C',P')$,
\end{enumerate} 
where the pullback of cut grids was defined in \cref{cutgrid}. 

\medskip
\paragraph{\it Presheaf structure maps}
The structure maps corresponding to morphisms in $\stcart$, $Γ$, and $Δ$
are given by nerves of functors of categories specified as follows.
For a given map $f:U'\to U$ in $\stcart$, we pull back cut tuples and $P$ via the corresponding map $\id\times f:M\times U'\to M\times U$.
Likewise for~$\varphi$. 
For $Γ$, a map $\langle\ell\rangle\to\langle\ell'\rangle$ is simply composed with the given map $P:M\times U\to\langle\ell\rangle$.
For $Δ^{\times d}$, the structure maps are induced by the maps in \cref{cutgrid}.
A face map removes the corresponding cut and a degeneracy map duplicates a cut.
\end{definition}

This definition of the extended bordism category is different from the definition of Lurie \cite{Lurie.TFT} and Calaque–\hskip0pt Scheimbauer \cite{CalaqueScheimbauer}.
(See \cref{comparison.conjecture}, though.)
One difference is that we do not keep track of the height function that implements cuts.
Instead, we simply keep track of the cuts.
We also allow cuts to overlap, but we do not allow cuts to intersect transversally.
These features allow us to implement globularity in our bordism categories which, in the case of $\Bord_d$, would otherwise collapse to just degenerate bordisms.

\begin{figure}
\begin{center}
\begin{tikzpicture}[yscale=1.3]
\node at (0,1.3) (t) {};
\node at (0,-1.3) (s) {};
\node at (0,-1.4) (t') {$C_{=0}^1$};
\draw[purple] (t) -- (s);
\node at (-1,-.2) (q) {};
\node at (1,-.2) (p) {};
\node at (1.5,-.4) (q') {$C_{=0}^2$};
\draw[blue] (q.center) -- (p.center);
\node at (-1,.4) (q) {};
\node at (1,.4) (p) {};
\node at (1.5,.6) (q') {$C_{=1}^2$};
\node at (0,-.2) (int') { $\bullet$};
\draw[blue] (q.center) to [out=0,in=180] (int');
\draw[blue] (int') to [out=0,in=180] (p.center);
\node at (-2,.4) (1) {};
\node at (-1,.4) (2) {};
\node at (-2,-.2) (3) {};
\node at (-1,-.2) (4) {};
\node at (1,.4) (5) {};
\node at (2,.4) (6) {};
\node at (1,-.2) (7) {};
\node at (2,-.2) (8) {};
\draw[blue] (1.center) -- (2.center);
\draw[blue] (3.center) -- (4.center);
\draw[blue] (5.center) -- (6.center);
\draw[blue] (7.center) -- (8.center);
\draw[dashed] (0.2,1.2) -- (0.2,-1.2);
\draw[dashed] (-0.2,1.2) -- (-0.2,-1.2);
\node at (2.5,0) {$\cong$};
\end{tikzpicture}
\begin{tikzpicture}[yscale=1.3]
\node at (0,1.3) (t) {};
\node at (0,-1.3) (s) {};
\node at (0,-1.4) (t') {$C_{=0}^1$};
\draw[purple] (t) -- (s);
\node at (-1,-.2) (q) {};
\node at (1,-.2) (p) {};
\node at (2,-.4) (q') {$C_{=0}^2=C^2_{=1}$};
\draw[blue] (q.center) -- (p.center);
\node at (-1,.4) (q) {};
\node at (1,.4) (p) {};
\node at (0,-.2) (int') { $\bullet$};
\node at (-2,-.2) (3) {};
\node at (-1,-.2) (4) {};
\node at (1,-.2) (7) {};
\node at (2,-.2) (8) {};
\draw[blue] (3.center) -- (4.center);
\draw[blue] (7.center) -- (8.center);
\end{tikzpicture}

\begin{tikzpicture}[yscale=1.3]
\node at (-.2,1.3) (t) {};
\node at (-.2,-1.3) (s) {};
\node at (-.3,-1.4) (t') {$C_{=0}^1$};
\draw[purple] (t) -- (s);
\node at (.4,1.3) (tt) {};
\node at (.4,-1.3) (ss) {};
\node at (.5,-1.4) (tt') {$C_{=1}^1$};
\draw[purple] (tt) -- (ss);
\node at (-2,-.2) (q) {};
\node at (2,-.2) (p) {};
\node at (1.5,-.4) (q') {$C_{=0}^2$};
\draw[blue] (q) --(p);
\node at (-1,.4) (qq) {};
\node at (1,.4) (pp) {};
\node at (1.5,.6)  {$C_{=1}^2$};
\node at (-.2,-.2) (int) {$\bullet$};
\node at (.4,-.2) (int') {$\bullet$};
\node at (.1,.2) (i) {};
\draw[blue]  (qq.center) to [out=0,in=180] (int.center);
\draw[blue]  (int.center) to [out=0,in=180] (i.center);
\draw[blue]  (i.center) to [out=0,in=180] (int'.center);
\draw[blue]  (int'.center) to [out=0,in=180]  (pp.center);
\node at (-2,.4) (3) {};
\node at (-1,.4) (4) {};
\node at (1,.4) (7) {};
\node at (2,.4) (8) {};
\draw[blue]  (3.center) -- (4.center);
\draw[blue]  (7.center) -- (8.center);
\node at (2.5,0) {$\not\cong$};
\draw[dashed] (-.25,-1.2) -- (-.25,1.2);
\draw[dashed] (-.15,-1.2) -- (-.15,1.2);
\draw[dashed] (.45,-1.2) -- (.45,1.2);
\draw[dashed] (.35,-1.2) -- (.35,1.2);
\end{tikzpicture}
\begin{tikzpicture}[yscale=1.3]
\node at (-.2,1.3) (t) {};
\node at (-.2,-1.3) (s) {};
\node at (-.3,-1.4) (t') {$C_{=0}^1$};
\draw[purple] (t) -- (s);
\node at (.4,1.3) (tt) {};
\node at (.4,-1.3) (ss) {};
\node at (.5,-1.4) (tt') {$C_{=1}^1$};
\draw[purple] (tt) -- (ss);
\node at (-2,-.2) (q) {};
\node at (2,-.2) (p) {};
\node at (2,-.4) (q') {$C_{=0}^2=C^2_{=1}$};
\draw[blue]  (q) -- (p);
\node at (-1,.4) (qq) {};
\node at (1,.4) (pp) {};
\node at (-.2,-.2) (int) {$\bullet$};
\node at (.4,-.2) (int') {$\bullet$};
\node at (.1,.2) (i) {};
\end{tikzpicture}
\caption{Take $d=2$. As an ambient manifold for a bordism in $\Bord_2$, we take $\RR^2$.
The picture illustrates isomorphic and non-isomorphic objects in bidegree $([0],[1])$ and $([1],[1])$ (respectively).
The dashed lines represent open neighborhoods of the cores of the simplicial faces. The top isomorphism is provided by the embedding which includes the dashed strip containing the vertical line into $\RR^2$. This defines 1-simplex in $\Bord_d([0],[1])$ connecting the two bordisms. 
The bottom image demonstrates that our definition of globularity does not force the bordism category to collapse into only degenerate multisimplices. The bordism on the left is globular, but it is not isomorphic to a degenerate multisimplex.}
\end{center}
\end{figure}

A second difference between our bordism category and that of Lurie \cite{Lurie.TFT} and Calaque--Scheimbauer \cite{CalaqueScheimbauer} is in the $(d+1)$-morphisms in our bordism category (encoded by cut-respecting embeddings in \cref{bord}).
The core of our bordisms may be contained in a large manifold with nontrivial topology outside of the core.
The portions of the manifold outside the core are irrelevant and we would really like to just keep track of a small neighborhood of the core.
This is secretly built in to our definition.
Indeed, Kan fibrant replacement of the nerve of the category in \cref{bord} formally inverts cut-respecting embeddings.
In particular, every bordism is in the same connected component as the bordism obtained by taking an arbitrarily small neighborhood of the core.
One might be concerned that we have added bizarre higher homotopy groups to our bordism category, but this is not the case, as we now show.  

\begin{definition}(germy bordism category)
\label{germybord}
Fix $d≥0$.
We define the \emph{germy bordism category} as the object $\gBord_d\in \smcat_{\infty,d}$ whose value on $(U,\langle \ell\rangle, {\bf m})$
is the nerve of the following groupoid~$\cD$, constructed using the category~$\cC$ of bordisms and open embeddings from \cref{bord}.
\begin{itemize}
\item Objects are the same as objects in~$\cC$.
\item Morphisms $A→B$ are given by equivalence classes of spans $A←G→B$.
Two spans $A←G_1→B$ and $A←G_2→B$ are identified if there is a third span $A←H→B$ that admits a morphism to both spans.
\item Spans are composed using pullbacks, which always exist in~$\cC$.
The composition manifestly respects the equivalence relation.
\end{itemize}
\end{definition}

Two remarks are in order.
First, a span $A←G→B$ should be interpreted as restricting to a smaller open neighborhood~$G$ of the core of~$A$, then embedding this neighborhood into~$B$.
Second, the equivalence relation on spans says that two spans coincide if their restrictions to an even smaller open neighborhood~$H$ coincide.
That is to say, morphisms $A→B$ in~$\cD$ are precisely germs of open embeddings of an open neighborhood of the core of~$A$ into~$B$.
Requiring the maps (i.e., open embeddings) to be morphisms forces the remainder of the structure of a bordism to be preserved by the embeddings.
The category $\cD$ is indeed a groupoid: the inverse of a morphism $A←G→B$ is the opposite span $B←G→A$.
Composing the spans $A←G→B$ and $B←G→A$ yields the span $A←(G⨯_B G)→A$, which is equivalent to the identity span $A←A→A$
via the map $G⨯_B G→A$. The other composition is equivalent to the identity by the same argument. 

We now define a canonical morphism in $\smcat_{\infty,d}$
\begin{equation}\label{fmap}f:\Bord_d\to \gBord_d.\end{equation}
Fix $U\in \stcart$, $\langle \ell\rangle\in \Gamma$, and ${\bf m}\in \Delta^{\times d}$.
After evaluating at $(U,\langle \ell\rangle,{\bf m})$, the map \cref{fmap} is given by taking the nerve of the following functor,
which we denote by $f_{(U,\langle \ell\rangle,{\bf m})}:\cC→\cD$, with $\cC$ and $\cD$ as in \cref{germybord}.
Send an object of~$\cC$ to itself.
Send a morphism $A→B$ in~$\cC$ to the span $A←A→B$, where the first leg is identity.
Composition is preserved by construction.

\begin{proposition}
\label{embtodiff}
The map $f$ in \cref{fmap} yields an objectwise weak equivalence of simplicial sets. 
\end{proposition}

\begin{proof}
For each $U\in \stcart$, $\langle \ell\rangle\in \Gamma$, and ${\bf m}\in \Delta^{\times d}$,
we apply Quillen's Theorem A to the corresponding functor $f_{(U,\langle \ell\rangle,{\bf m})}:\cC→\cD$.
Fix an arbitrary object $B∈\cD$.
We claim that the comma category $\cC/B$ is cofiltered, hence contractible.
Objects in $\cC/B$ are spans of the form $A←G→B$.
The identity span $B←B→B$ shows that $\cC/B$ is nonempty.
Given two objects of $\cC/B$ presented by spans $A_1←G_1→B$ and $A_2←G_2→B$,
the span $G_1⨯_B G_2←G_1 ⨯_B G_2→B$ is another object of~$\cC/B$ that maps to both of the above objects via the obvious projection maps.

Finally, suppose we have two morphisms $f,g:A_1→A_2$ in $\cC$ that make two triangles (whose edges are spans) with vertices $A_1$, $A_2$, $B$ commute in~$\cD$.
The commutativity is witnessed by a subobject $h:A_0→A_1$, which can be concretely described as the intersection of two subobjects of $G_1∩f^*G_2∩g^*G_2$
supplied by the definition of equality of spans that make the two triangles commute.
Then $h$ defines a morphism from the span $A_0←A_0→B$ to the span $A_1←G_1→B$ such that $fh=gh$.
This show that the comma category $\cC/B$ is cofiltered for any $B∈\cD$, completing the proof.
\end{proof}

\begin{example}
Let $\stcart=\cart$. We provide a simple example of a bordism in the uple version of our bordism category. This example illustrates that manifolds with corners are examples of morphisms in our bordism category. More complicated examples are depicted in \cref{figure2}. 

Let $d=2$ and let $U=\RR^0$, so that the family direction is trivial.
We will construct a bordism in bisimplicial degree $([2],[2])$.
The ambient manifold is given by the torus $T=S^1\times S^1$.
Let $e:T\into \RR^3$ be the embedding of the torus given by $e(\theta,\phi)=((\cos\theta+3)\cos\phi,(\cos\theta+3)\sin\phi,\sin\theta)$.
Take the cut tuple on $T$ in the first simplicial direction to be the preimage (under~$e$)
of the codimension~1 manifolds obtained by intersecting the hyperplanes $x=-3$, $x=0$, and $x=3$ with $e(T)$.
Similarly, take the cut tuple in the second simplicial direction to be the preimage of the intersection of $y=-3$, $y=0$ and $y=3$ with $e(T)$.
A sketch of the resulting bordism is given below.
\begin{center}
\begin{tikzpicture}
\draw (0,0) to[out=90, in=180] (2,1); 
\draw (2,1) to[out=0, in=90] (4,0); 
\draw (0,0) to[out=-90, in=180] (2,-1); 
\draw (2,-1) to[out=0, in=-90] (4,0); 
\draw (3.5,-.75) to[out=90+10, in=-90-10] (3.5,.75); 
\draw[dashed] (3.5,-.75) to[out=90-10, in=-90+10] (3.5,.75); 
\draw (.5,-.75) to[out=90+10, in=-90-10] (.5,.75); 
\draw[dashed] (.5,-.75) to[out=90-10, in=-90+10] (.5,.75); 
\draw (1.5,.25) to[out=40, in=180] (2,.4); 
\draw (2,.4) to[out=0, in=140] (2.5,.25); 
\draw (1.44,.3) to[out=-60, in=180] (2,0); 
\draw (2,0) to[out=0, in=-90-30] (2.6,.3); 
\draw (2,-1) to[out=100,in=-90-10] (2,0);
\draw[dashed] (2,-1) to[out=90-10,in=-90+10] (2,0);
\draw (2,.4) to[out=100,in=-90-10] (2,1);
\draw[dashed] (2,.4) to[out=90-10,in=-90+10] (2,1);
\draw (.05,-.3) to[out=10,in=180] (2,-.2);
\draw (2,-.2) to[out=0,in=170] (4-.05,-.3);
\draw[dashed] (.05,-.3) to[out=-10,in=180] (2,-.4);
\draw[dashed] (2,-.4) to[out=0,in=190] (4-.05,-.3);
\draw (.45,1-.3) to[out=10,in=180] (2,1-.2);
\draw (2,1-.2) to[out=0,in=170] (4-.45,1-.3);
\draw[dashed] (.45,1-.3) to[out=-10,in=180] (2,1-.4);
\draw[dashed] (2,1-.4) to[out=0,in=190] (4-.45,1-.3);
\draw (.05,.2) to[out=10,in=170] (1.5,.2);
\draw[dashed] (.05,.2) to[out=-10,in=-170] (1.5,.2);
\draw[dashed] (.05,-.3) to[out=-10,in=180] (2,-.4);
\draw[dashed] (2,-.4) to[out=0,in=190] (4-.05,-.3);
\draw (2.55,.2) to[out=10,in=170] (4,.2);
\draw[dashed] (2.55,.2) to[out=-10,in=-170] (4,.2);
\end{tikzpicture}
\end{center}
In this case, the torus is the ambient manifold. The region lying between the planes $x=-3$, $x=3$, $y=-3$ and $y=3$ is a manifold with corners.
This manifold with corners is the core of the bordism.
\end{example}

\begin{remark}
We have chosen to work exclusively with trivial bundles $M\times U\to U$ in \cref{bord} in order to ensure functoriality of pullbacks of bundles.
If we are working over the cartesian site, i.e., $\stcart=\cart$, all bundles are trivializable and our presheaf of $(\infty,d)$-categories satisfies descent on cartesian spaces.
Alternatively, we could include all bundles, but then we would need to work with Grothendieck fibrations instead of presheaves.
\end{remark}

\begin{figure}[ht]
\begin{center}
\begin{tikzpicture}[scale=.60]
\filldraw[fill=gray!25, rounded corners, dashed] (0+.3,-1+.3) -- (0+.3,6+.3) -- (9+.3,6+.3) -- (9+.3,-1+.3) -- (0+.3,-1+.3);
\draw[purple] (5+.3,2+.3) to [out=-60, in=110] (9+.3,1+.3);
\draw[purple]  (5+.3,2+.3) to [out=300, in=145] (9+.3,3+.3);
\draw[purple] (6+.3,3+.3) to [out=305, in=110] (9+.3,3.5+.3);
\draw[blue] (1+.3,6+.3) to [out=270, in=110] (2+.3,3+.3);
\draw[blue] (1.5+.3,6+.3) to [out=-45, in=110] (2+.3,3+.3);
\draw[blue] (3+.3,6+.3) to [out=-45, in=110] (5+.3,3+.3); 
\draw[blue] (4+.3,6+.3) to [out=-90, in=110] (6,3);  
\draw[orange] (5+.3,6+.3) to [out=200, in =180] (9+.3,4+.3);
\filldraw[fill=gray!25, rounded corners, dashed] (0,-1) -- (0,6) -- (9,6) -- (9,-1) -- (0,-1);
\draw[purple] (0,0) to [out=20,in=120] (5,2) to [out=-60, in=110] (9,1);
\draw[purple] (0,2) to [out=10,in=115] (5,2) to [out=300, in=145] (9,3);
\draw[purple] (0,2.5) to [out=0, in=125] (6,3) to [out=305, in=110] (9,3.5);
\node at (0-.7,0) {$C^1_{=0}$};
\node at (0-.7,2) {$C^1_{=1}$};
\node at (0-.7,2.8) {$C^1_{=2}$};
\draw[blue] (1,6) to [out=270, in=110] (2,3) to [out=290, in=110] (1,-1);
\draw[blue] (1.5,6) to [out=-45, in=110] (2,3) to [out=290, in=110] (1,-1);
\draw[blue] (3,6) to [out=-45, in=110] (5,3) to [out=290, in=110] (1,-1); 
\draw[blue] (4,6) to [out=-90, in=110] (6,3) to [out=290, in=90] (4,-1);  
\draw[orange] (5,6) to [out=200, in =180] (9,4);
\node at (.6,6.5+.3) {$C_{=0}^2$};
\node at (1.6,6.5+.3) {$C_{=1}^2$};
\node at (3,6.5+.3) {$C_{=2}^2$};
\node at (4,6.5+.3) {$C_{=3}^2$};
\node at (5.2,6.5+.3) {$D_{=}$};
\draw[fill=white] (2.65-.2,4.6) to [out=57, in=118] (3.9-.2,4.6);
\draw[thick] (2.5-.2,4.5) to [out=60, in=115] (4-.2,4.5);
\draw[thick, fill=white] (2.65-.2,4.6) to [out=-45, in=-135] (3.9-.2,4.6);
\draw[fill=white]  (6.5+.2,.5+.2) to [out=46, in=130] (8.5-.2,.5+.2);
\draw[thick] (6.5,.5) to [out=60, in=115] (8.5,.5);
\draw[thick, fill=white] (6.5+.2,.5+.2) to [out=-45, in=-135] (8.5-.2,.5+.2);
\end{tikzpicture}
\begin{tikzpicture}[scale=.60]
\fill[fill=gray!25, rounded corners] (0+.3,-1+.3) -- (0+.3,6+.3) -- (9+.3,6+.3) -- (9+.3,-1+.3) -- (0+.3,-1+.3);
\draw[blue] (1+.3,6+.3) to [out=270, in=110] (2+.3,3+.3);
\draw[blue] (1.5+.3,6+.3) to [out=-45, in=110] (2+.3,3+.3);
\draw[blue] (3+.3,6+.3) to [out=-45, in=110] (5+.3,3+.3); 
\draw[blue] (4+.3,6+.3) to [out=-90, in=110] (6,3);  
\draw[blue] (5+.3,6+.3) to [out=200, in =180] (9+.3,4+.3);
\draw[dashed, rounded corners]  (0+.3,-1+.3) -- (0+.3,6+.3) -- (9+.3,6+.3) ;
\fill[fill=gray!25, rounded corners] (0,-1) -- (0,6) -- (9,6) -- (9,-1) -- (0,-1);
\draw[dashed, rounded corners] (9.7,-.5) to[out=270,in=0, looseness=.8] (9,-1) -- (0,-1) -- (0,6) -- (9,6);
\fill[fill=gray!25, dashed] (9,6) -- (9,-1) to[out=0,in=270, looseness=1.2] (9.7,.5) -- (9.7,6.2);
\filldraw[fill=gray!25, dashed] (9+.3,6+.3) to [out=0, in=0, looseness=4] (9,6);
\draw[purple] (0,0) to [out=20,in=120] (5,2) to [out=-60, in=130] (9,1) to [out=310, in=-90] (9.7,1.5);
\draw[purple] (0,2) to [out=10,in=115] (5,2) to [out=300, in=145] (9,3) to [out=325, in=-110] (9.7,3.3);
\draw[purple] (0,2.5) to [out=0, in=125] (6,3) to [out=305, in=110] (9,3.5) to [out=300, in=-130] (9.7,3.7);
\node at (0-.7,0) {$C^1_{=0}$};
\node at (0-.7,2) {$C^1_{=1}$};
\node at (0-.7,2.8) {$C^1_{=2}$};
\draw[blue] (1,6) to [out=270, in=110] (2,3) to [out=290, in=110] (1,-1);
\draw[blue] (1.5,6) to [out=-45, in=110] (2,3) to [out=290, in=110] (1,-1);
\draw[blue] (3,6) to [out=-45, in=110] (5,3) to [out=290, in=110] (1,-1); 
\draw[blue] (4,6) to [out=-90, in=110] (6,3) to [out=290, in=90] (4,-1);  
\draw[blue] (5,6) to [out=200, in =180] (9,4) to [out=0, in=-130] (9.7,4.2);
\node at (.6,6.5+.3) {$C_{=0}^2$};
\node at (1.6,6.5+.3) {$C_{=1}^2$};
\node at (3,6.5+.3) {$C_{=2}^2$};
\node at (4,6.5+.3) {$C_{=3}^2$};
\node at (5.2,6.5+.3) {$D_{=}$};
\draw[fill=white] (2.65-.2,4.6) to [out=57, in=118] (3.9-.2,4.6);
\draw[thick] (2.5-.2,4.5) to [out=60, in=115] (4-.2,4.5);
\draw[thick, fill=white] (2.65-.2,4.6) to [out=-45, in=-135] (3.9-.2,4.6);
\draw[fill=white]  (6.5+.2,.5+.2) to [out=46, in=130] (8.5-.2,.5+.2);
\draw[thick] (6.5,.5) to [out=60, in=115] (8.5,.5);
\draw[thick, fill=white] (6.5+.2,.5+.2) to [out=-45, in=-135] (8.5-.2,.5+.2);
\end{tikzpicture}
\end{center}
\caption{Let $\stcart=\cart$. Two bordisms with cut tuples $C^1$ and $C^2$ for $d=2$.
Everything is parametrized by elements $x\in U$ for some fixed cartesian space $U$ and the figures depict the fiber over some point $x\in U$.
The gray region is the ambient 2-dimensional smooth manifold: two parallel sheets connected by two “tunnels” on the left and one folded sheet connected by two tunnels on the right.
The purple cuts are given by the cut $[2]$-tuple $C^1$.
On the left the blue cuts are given by the cut $[3]$-tuple $C^2$.
On the right, the blue cuts are given by the cut $[4]$-tuple $C^2\cup D$.
The orange cut $D$ is not permitted in the left picture, as this would cause the core to be noncompact.}
\label{figure2}
\end{figure}

\begin{example}
\label{transcut}
Again, let $\stcart=\cart$.
We provide an example of a manifold with cuts that is \emph{not} an object in our bordism category for $d=2$.
We let $U=\RR^0$, so that the family direction is trivial. The ambient manifold $M$ is $\RR^2$.
Consider the following two cut $[0]$-tuples (i.e., cuts) $C^1$ and $C^2$.
The set $C^1_{=}$ is given by the $y$-axis in $\RR^2$, $C_{<}^1=\{(x,y)\in \RR^2 \mid x<0\}$, $C_{>}^1=\{(x,y)\in \RR^2 \mid x>0\}$.
The set $C^2_{=}$ is given by the parabola $x=y^2$, $C^2_{<}=\{(x,y)\in \RR^2 \mid x<y^2\}$, $C^2_{>}=\{(x,y)\in \RR^2 \mid x>y^2\}$.
In this case, the two cuts $C^1_{=}$ and $C^2_{=}$ are tangent at the origin (see \cref{figure3}).

This is not an admissible bordism for the following reason.
The condition $\pitchfork$ in \cref{cutgrid} implies that for $S=\{1,2\}$ there is a smooth function $h_j:\RR^2\to \RR^2$,
corresponding to the function $j:\{1,2\}\to \ZZ$ defined by $j(1)=0$, $j(2)=0$, such that $h_j$ has zero as a regular value
and $(\pi_1\circ h_j)^{-1}(0)=C^1_{=}$, $(\pi_2\circ h_j)^{-1}(0)=C^2_{=}$.

The standard basis vector $e_2=(0,1)$ is tangent to both curves $C^1_{=}$ and $C^2_{=}$. 
Therefore, since $\pi_1\circ h_j$ vanishes on~$C^1_{=}$ and $\pi_2\circ h_j$ vanishes on $C^2_{=}$, 
we must have $d\pi_1dh_j(0,0)e_2=d\pi_2dh_j(0,0)e_2=0$. Thus, $dh_j(0,0)e_2=0$, which implies that $(0,0)$ is not a regular value of $h_j$. 

\begin{figure}[ht]
\begin{center}
\begin{tikzpicture}
\fill[blue, opacity=.3] (2,1) to[out=-70-90,in=180-90, looseness=.8] (1,0) to[out=-90, in=250-90, looseness=.8] (2,-1) -- (2,1); 
\draw (0,0) -- (2,0);
\draw[purple] (1,-1) -- (1,1);
\fill[red, opacity=.3] (1,1) -- (1,-1) -- (2,-1) -- (2,1) -- (1,1);
\draw[blue] (2,1) to[out=-70-90,in=180-90, looseness=.8] (1,0) to[out=-90, in=250-90, looseness=.8] (2,-1);
\draw[->] (2.5,0) -- (3.2,0);
\node at (2.88,.3) {$h_j$};
\begin{scope}[xshift=1.5in]
\fill[red, opacity=.3] (1,1) -- (1,-1) -- (2,-1) -- (2,1) -- (1,1);
\fill[blue, opacity=.3] (0,1) -- (0,0) -- (2,0) -- (2,1) -- (0,1);
\draw[blue] (0,0) -- (2,0);
\draw[purple] (1,-1) -- (1,1);
\end{scope}
\node at (1.9,1.3) {$C^2_{= }$};
\node at (1, 1.3) {$C^1_{= }$};
\end{tikzpicture}
\end{center}
\caption{The map $h_j$ sends the blue region $C_{>}^2$ into the blue region $y>0$ and red region $C_{>}^1$ into red region $x>0$.
It sends the blue curve $C_{=}^2$ into the blue line $x=0$ and the red line $C^1_{=}$ into the red line $y=0$.
The two curves $C^2_{=}$ and $C^1_{=}$ are tangent at the origin.}
\label{figure3}
\end{figure}
\end{example}

\begin{remark}
\cref{transcut} shows that cuts in \emph{different} simplicial directions must intersect transversally. This is forced by the requirement that $h_j$ have the origin as a regular value and by the compatibility of $h_j$ with the cuts. We remind the reader that cuts in the \emph{same} simplicial direction are allowed to intersect and even overlap, but are not allowed to intersect transversally (see Figure 2).
\end{remark}

\subsection{Categories of bordisms with geometric structure}

In \cref{bord}, the bordisms are equipped with no additional structure. 
We would like to equip bordisms with geometric structure.
Recall from \cref{geometric.structure} that these are given simply by simplicial presheaves on the site $\FEmb_d$. 

\begin{definition}[$d$-uple/globular bordisms with structure]
\label{bordstr}
Fix $d≥0$.
We define functors $$\Bord_{d,\uple}:\Struct_d→\smcatuple_{∞,d},\qquad \Bord_d:\Struct_d→\smcat_{∞,d}$$
by specifying their value at $\gs\in \Struct_d$ (\cref{geometric.structure})
as an object $\Bord_{d,\uple}^\gs$, respectively $\Bord_d^\gs$, in the category 
$\smcatuple_{\infty,d}$ (\cref{multiple.model.structure}), respectively $\smcat_{\infty,d}$ (\cref{globular.model.structure}),
as follows.
Fix $(U,\langle \ell\rangle,{\bf m})\in \stcart\times \Gamma\times \Delta^{\times d}$.
We define $\Bord_{d,\uple}^\gs(U,\langle \ell\rangle,{\bf m})$, respectively $\Bord_d^\gs(U,\langle \ell\rangle,{\bf m})$, by taking the diagonal of the nerve of the following category internal to simplicial sets.

\begin{enumerate}
\item
The simplicial set of objects is given by 
$$\Ob_\gs≔\coprod_{(M,C,P)}\gs(M\times U\to U),$$
where the coproduct ranges over the objects of \cref{bord}.
\item
The simplicial set of morphisms is given by
$$\Mor_\gs≔\coprod_{\varphi:(M,C,P)\to ( M', C', P')}\gs( M'\times U\to U),$$
where the coproduct is taken over the morphisms in \cref{bord}.
\item The target map of the category structure sends the component indexed by a cut-respecting embedding $\varphi:(M,C,P)\to ( M', C', P')$
to the component indexed by the bordism $( M', C', P')$, via the identity map on $\gs( M'\times U\to U)$.
The source map pulls back the structure using the morphism $\varphi$, which is a morphism in $\FEmb_d$, then maps to the component indexed by the bordism $(M,C,P)$.
\item The composition map $\Mor_\gs\times_{\Ob_\gs}\Mor_\gs\to \Mor_\gs$ is given as follows. An $l$-simplex in the pullback is a pair of morphisms $\varphi:(M,C,P)\to ( M', C',  P')$ and $\psi:( M', C', P')\to (M'',C'',P'')$, along with a pair of $l$-simplices $\tau\in \gs(M'\times U\to U)$ and $\sigma\in \gs(M''\times U\to U)$, such that $\psi^*(\sigma)=\tau$. The composition map sends this data to $\sigma$ in the summand $\gs(M''\times U\to U)$ indexed by the composition $\psi\circ \varphi$.
\end{enumerate}
The presheaf structure maps are specified as follows.
For $\Gamma\times \Delta^{\times d}$, we simply use the presheaf structure maps on objects $(M,C,P)$ as specified in \cref{bord}.
For $\stcart$, we use the structure map from \cref{bord}, along with the structure map for the presheaf $\gs\in \sPSh(\FEmb_d)$.
The resulting objects $\Bord_{d,\uple}^\gs$ and $\Bord_d^\gs$ are manifestly functorial in $\gs$.
\end{definition}

\begin{remark}
Analogous to \cref{germybord}, one can define an alternative presentation for $\Bord_d^\gs(U,\langle \ell\rangle,{\bf m})$ which includes germs. The resulting bordism category can again be shown to be equivalent to \cref{bordstr} by a simple modification of the argument in \cref{embtodiff}. Since we will not need this alternative version of the bordism category, we will not provide the details of the construction or proof. 
\end{remark}

\begin{example}
Taking $\gs$ to be the terminal object in $\Struct_d$, we get
$$\Bord_d^\gs=\Bord_d,$$
where $\Bord_d$ is as in \cref{bord}.
\end{example}

\begin{example}
Let $X$ be a smooth manifold, viewed as a sheaf on $\FEmb_d$ via
$$X(N\to U)=C^{\infty}(N,X).$$
For $N=M\times U\to U$, the corresponding geometric structure on a bordism $(M,C,P)$ is a smooth function $f:M\times U\to X$. 
A morphism $\varphi$ in \cref{bordstr} acts simply by pulling back $f$ by the cut-respecting embedding $\varphi: M'\times U\to M\times U$. 
\end{example}

\begin{example}
\label{metricsregb}
(See \cref{metricsregb.figure}.)
Consider the sheaf $\FRiem_1\in \Struct_1$ of fiberwise Riemannian metrics (\cref{riemmet}).
A vertex in  $\Bord_1^{\FRiem_1}(\RR^0,\langle 1\rangle,[1])$ is a quadruple $(M,C,P,g)$,
where $M$ is a 1-dimensional manifold, $C$ is a cut $[1]$-tuple on $M$, and $g$ is a Riemannian metric on $M$. Assuming the trash bin is empty, $P$ sends everything to the element $1\in \langle 1\rangle$. 
In the case that $M$ is connected, the metric gives $M$ a Riemannian length $t$. 

A 1-simplex $\varphi:(M,C,P,g)\to (M',C',P',g')$ in $\Bord_1^{\FRiem_1}(\RR^0,\langle 1\rangle,[1])$ is given by the following data. The map $\varphi:M\to M'$ is a cut-respecting embedding, $C$ and $ C'$ are cut $[1]$-tuples on $M$ respectively $M'$ such that $\varphi^*C'=C$, $g'$ is a Riemannian metric on $M'$ such that $\varphi^*(g')=g$, and $P$ and $P'$ are maps to $\langle 1\rangle$ such that $\varphi^*P'=P$.
The two endpoints of the 1-simplex are given by $d_0(\varphi)=(M', C', P',g)$ and $d_1(\varphi)=(M, C, P, g)$. 
\end{example}

\begin{figure}[ht]
\begin{center}
\begin{tikzpicture}[scale=.65]
\draw (0,0) to [out=20,in=120] (5,2) to [out=-60, in=110] (9,1);
\filldraw (2,1.28) circle (2pt);
\filldraw (7,1.52) circle (2pt);
\node at (2,2) {$C_{=0}$};
\node at (7,2) {$C_{=1}$};
\node at (4.5,3) {$t$};
\node at (0,1) {$M$};
\end{tikzpicture}
\end{center}
\caption[]{A vertex in $\Bord_1^{\FRiem_1}(\RR^0,\langle 1\rangle,[1])$.
The 1-manifold $M$ is equipped with a Riemannian metric $g$.
The core is the 1-manifold with boundary lying between the cuts $C_{=0}$ and $C_{=1}$.
The Riemannian length of the core is $t$.}
\label{metricsregb.figure}
\end{figure}

\begin{remark}
For an example where the higher simplices of the geometric structure are present, take the $\Delta^1$-family of cuts in \cref{gbunexample} to be constant.
This gives an example of a 1-simplex in $$\Bord_2^{\deloop_{\nabla}G}(\RR^0,\langle 1\rangle, ([1],[1])).$$
\end{remark}

\subsection{Categories of bordisms with isotopies and geometric structure}
\label{isotopybordsec}

We have the following variant of the bordism category,
which incorporates isotopies between cuts.

\begin{definition}[$d$-uple/globular bordisms with isotopies]
\label{frakbord}
Given $d≥0$, we specify objects $\BBord_{d,\uple}$ and $\BBord_d$ respectively in the categories $\fraksmcatuple_{\infty,d}$ and $\fraksmcat$ by pulling back $\Bord_{d,\uple}$ and $\Bord_{d}$ (\cref{bord}) along the functor 
$$ \stcart\times \Gamma\times \Delta^{\times d}\times \cart\to \stcart\times \Gamma\times \Delta^{\times d}$$
that sends $(U,\langle \ell\rangle,{\bf m},L)\mapsto (U⨯L,\langle \ell\rangle, {\bf m})$, where $\times$ denotes the tensoring in \cref{def.cartsp}. Explicitly, 
$$\BBord_{d,\uple}(U,\langle \ell\rangle,{\bf m})(L)=\Bord_{d,\uple}(U⨯L,\langle \ell\rangle,{\bf m}), \qquad\BBord_{d}(U,\langle \ell\rangle,{\bf m})(L)=\Bord_{d}(U⨯L,\langle \ell\rangle,{\bf m}).$$
\end{definition}

We are now ready to add geometric structures to the isotopy bordism category. Examples will be given subsequently.

\begin{definition}[$d$-uple/globular bordisms with geometric structures and isotopies]
\label{enrichedbordstr}
Fix $d≥0$.
We define $\smsset$-enriched functors $$\BBord_{d,\uple}:\frakStruct_d→\fraksmcatuple_{∞,d},\qquad \BBord_d:\frakStruct_d→\fraksmcat_{∞,d}$$
by specifying their value at $\gs\in \frakStruct_d$ (\cref{geometric.structure.isotopy})
as an object $\BBord_{d,\uple}^\gs$, respectively $\BBord_d^\gs$, in the category 
$\fraksmcatuple_{\infty,d}$ (\cref{multiple.model.structure.smooth}), respectively $\fraksmcat_{\infty,d}$ (\cref{globular.model.structure.smooth}),
as follows.
Fix $(U,\langle \ell\rangle,{\bf m})\in \stcart\times \Gamma\times \Delta^{\times d}$.
We define the smooth simplicial sets $\BBord_{d,\uple}^\gs(U,\langle \ell\rangle,{\bf m})$ and $\BBord_d^\gs(U,\langle \ell\rangle,{\bf m})$ by taking the diagonal of the nerve of the following category internal to smooth simplicial sets (see \cref{bordstr} for the unenriched case).

\begin{enumerate}
\item The $L$-points of the smooth simplicial set of objects are given by 
$$\Ob_\gs(L)≔\coprod_{(M,C,P)}\gs(M\times U\to U)(L),$$
where the coproduct ranges over the objects in \cref{bord}, with $U⨯L$ replacing $U$. The structure maps for the smooth simplicial set are provided by those of $\gs(M\times U\to U)$ and the $\cart$ structure maps for the objects in \cref{bord}.

\item
The $L$-points of the smooth simplicial set of morphisms are given by
$$\Mor_\gs(L)≔\coprod_{\varphi:(M,C,P)\to ( M', C', P')}\gs( M'\times U\to U)(L),$$
where the coproduct is taken over the morphisms in \cref{bord}, with $U⨯L$ replacing $U$. The structure maps for the smooth simplicial set are provided by those of $\gs(M'\times U\to U)$ and the $\cart$ structure maps for the morphisms in \cref{bord}.

\item The target map of the category structure sends the component indexed by an $U⨯L$-family of cut-respecting embeddings $\varphi:(M,C,P)\to ( M', C', P')$ 
to the component indexed by the $U⨯L$-family of bordisms $( M', C', P')$, via the identity map on $\gs( M'\times U\to U)(L)$.
The source map pulls back the structure using the morphism $\varphi$, which is an $L$-point of $\frakFEmb_d(M\times U\to U,M'\times U\to U)$, and then maps to the component indexed by the $U⨯L$-family of bordisms $(M,C,P)$.
\item The $L$-points of the composition map $\Mor_\gs\times_{\Ob_\gs}\Mor_\gs\to \Mor_\gs$ are given as follows. An $l$-simplex in the pullback is a pair $\varphi,\psi\in \Mor(U⨯L)$, where $\Mor(U⨯L)$ denotes the set of morphisms in \cref{bord}, with $U⨯L$ replacing $U$, along with a pair of $l$-simplices $\tau\in \gs(M'\times U\to U)(L)$ and $\sigma\in \gs(M''\times U\to U)(L)$, such that $\psi^*(\sigma)=\tau$. The composition map sends this data to $\sigma$ in the summand $\gs(M''\times U\to U)(L)$ indexed by the composition $\varphi\circ \psi$.
\end{enumerate}
The presheaf structure maps are specified as follows.
For $\Gamma\times \Delta^{\times d}$, we simply use the presheaf structure maps on objects $(M,C,P)$ as specified in \cref{frakbord}.
For $\stcart$, we use the structure map from \cref{frakbord}, along with the structure map for the enriched presheaf $\gs\in \frakStruct_d$.
Finally, to turn $\BBord_{d,\uple}^\gs$ and $\BBord_d^\gs$ into $\smsset$-enriched functors in $\gs$,
observe that the enriched structure maps
$$\frakStruct_d(\gs_1,\gs_2)⨯\BBord_{d,\uple}^{\gs_1}→\BBord_{d,\uple}^{\gs_2}$$
can be defined by moving the first factor past the diagonal and nerve functors (both of which preserve the $\smsset$-tensoring),
and then acting on $\Ob_{\gs_1}(L)$ and $\Mor_{\gs_1}(L)$ by the elements of $\frakStruct_d(\gs_1,\gs_2)(L)$.
Likewise for $\BBord_d$.
\end{definition}

\begin{remark}
Just as in the case of the bordism category without isotopies, the bordism category with isotopies is also equivalent to its germy variant (see \cref{embtodiff}).
Since we will not use the germy variant in the present work, we omit the proof.  
\end{remark}

\begin{conjecture}
\label{comparison.conjecture}
Evaluating $\BBord_d^\gs$ at the terminal object of $\stcart$ yields a symmetric monoidal $(∞,d)$-category
that is weakly equivalent to the symmetric monoidal $(∞,d)$-categories of bordisms
of Calaque–Scheimbauer \cite[Definition~9.10]{CalaqueScheimbauer} and Schommer-Pries \cite[Definition~5.8]{SchommerPries},
taken for the topological structure given by evaluating~$\gs$ at the terminal object of $\stcart$.
(For comparison to other models of topological structures, such as spaces with an action of the orthogonal group, see Grady–Pavlov \cite{GradyPavlov.Str}.)
\end{conjecture}

\begin{remark}
\label{frakbordl0}
From \cref{frakbord}, it follows at once that evaluation of $\BBord_d$ at $L=\RR^0$ yields precisely $\Bord_d$. In the case of geometric structure (\cref{enrichedbordstr}), taking $L=\RR^0$ yields $\Bord_d^{\ev_{\RR^0}\gs}$, where $\ev_{\RR^0}\gs\in \Struct_d$ is the presheaf of simplicial sets obtained by evaluating the components of the enriched presheaf $\gs\in \frakStruct_d$ at $\RR^0\in \cart$. 
\end{remark}

The previous definition is fairly complicated, so some examples are in order. 

\begin{example}
\label{1driem}
Consider the sheaf $\FRiem_1\in \Struct_1$ of fiberwise Riemannian metrics (see \cref{riemmet}).
In \cref{metricsregb}, we described the data of a vertex and 1-simplex in $ \Bord_1^{\FRiem_1}(\RR^0,\langle 1\rangle,[1])$.
We promote $\FRiem_1$ to an object in $\frakStruct_{1}$ as in \cref{friemex}.
An $\RR^0$-point of $\BBord_1^{\frakFRiem_1}(\RR^0,\langle 1\rangle,[1])$ is the same as a vertex in 
$\Bord_1^{\FRiem_1}(\RR^0,\langle 1\rangle,[1])$,
i.e., a quadruple $(M, C, P, g)$ where $M$ is a 1-dimensional manifold, $C$ is a cut $[1]$-tuple on $M$ and $g$ is a Riemannian metric on $M$.
In general, an $L$-point is an $L$-family of bordisms $(M,C,P)$ (\cref{frakbord}) together with an equivalence class of pairs $(g,\varphi)$, where $g$ is a fiberwise Riemannian metric on $M$ and $\varphi:M\times L\to M$ is an $L$-family of embeddings. 

A legitimate $\RR$-point is given by taking $\varphi(x,t)=x$, for $(x,t)\in M\times \RR$ 
and $C_t$ to be the $\RR$-family of cut $[1]$-tuples given by keeping the cut $C_{=1}$ fixed and moving the cut $C_{=0}$ in \cref{metricsregb} in the direction of $C_{=1}$,
with the two cuts being equal at $t=1$.
This deformation defines an isotopy that appears to collapse the data of the Riemannian length. 
However, this is not quite the case, since the same isotopy also moves the source 0-bordism along the same interval.
Thus, the isotopy of 1-bordisms constructed above translates the data of a 1-bordism given by the interval to the data of an isotopy of 0-bordisms given by the sources of 1-bordisms.

\begin{figure}[ht]
\begin{center}
\begin{tikzpicture}[scale=.6]
\draw (0,0) to [out=20,in=120] (5,2) to [out=-60, in=110] (9,1);
\filldraw (2,1.28) circle (2pt);
\filldraw (7,1.52) circle (2pt);
\node at (3,.8) {$(C_{=0})_{t=0}$};
\node at (7,2) {$(C_{=1})_{t=0}$};
\node at (4.5,3) {$t$};
\node at (0,1) {$M$};
\end{tikzpicture}
\qquad
\begin{tikzpicture}[scale=.6]
\draw (0,0) to [out=20,in=120] (5,2) to [out=-60, in=110] (9,1);
\filldraw (2,1.28) circle (2pt);
\filldraw (7,1.52) circle (2pt);
\filldraw (4,2.55) circle (2pt);
\node at (7.1,2) {$(C_{=1})_{t=.5}$};
\node at (4,3) {$(C_{=0})_{t=.5}$};
\node at (6,1) {$t$};
\node at (2.5,2.2) {$t'$};
\node at (0,1) {$M$};
\end{tikzpicture}
\end{center}
\caption{The family of cut $[1]$-tuples $C_t$ evaluated at $t=0$ and $t=.5$.
At $t=1$ $C_{=0}$ and $C_{=1}$ coincide.
The Riemannian length of the core of each individual cut tuple changes from $t$ to $t'$ as $t$ changes from $t=0$ to $t=.5$.
However, the length of the \emph{entire $\RR$-family} is recorded throughout the deformation.
}
\end{figure}
\end{example}

In the next example, we consider a geometric structure with higher simplicial data. 

\begin{example}\label{gbunexample}
Let $G$ be a Lie group and let $\gs=\deloop_{\nabla}G∈\Struct_d$ be the moduli stack of principal $G$-bundles with fiberwise connection,
which sends an object $(M→U)∈\FEmb_d$ to the nerve of the groupoid whose objects are
trivial principal $G$-bundles $G⨯M→M$ with a fiberwise connection on $M→U$ and morphisms are connection-preserving isomorphisms,
with the structure maps for $\FEmb_d$ defined by pulling back fiberwise connection 1-forms.
We used trivial bundles in the construction since they are easy to pull back.
To encode nontrivial principal $G$-bundles, we can rectify the Grothendieck fibration of principal $G$-bundles over manifolds to a strict presheaf of groupoids,
then proceed as in the above construction.
Alternatively, we can replace $\deloop_\nabla G$ with its associated $\infty$-sheaf by applying the fibrant replacement functor for $\Struct_d$.
We again promote $\deloop_{\nabla}G$ to an object in $\frakStruct_d$ by taking the enriched left Kan extension as described in \cref{friemex}.

Set $d=2$.
A vertex in $\BBord_2^{\deloop_{\nabla}G}(\RR^0,\langle 1\rangle, ([1],[1]))(\RR^0)$ is a 2-dimensional manifold $M$, equipped with two cut $[1]$-tuples $C^1$ and~$C^2$, and a map $P:M\to \langle 1\rangle$.
Moreover, the manifold is equipped with a trivial principal $G$-bundle $G\times M\to M$ with connection $\nabla$.
A 1-simplex in $\BBord_2^{\deloop_{\nabla}G}(\RR^0,\langle 1\rangle, ([1],[1]))(\RR^0)$ is a cut respecting embedding $\varphi:M \to M'$, along with a pair of cut $[1]$-tuples $C^1$ and $C^2$ on $M'$ and a map $P:M'\to \langle 1\rangle$.
Moreover, the manifold~$M'$ is equipped with a morphism of principal $G$-bundles with connection $q:(P,\nabla)\to (P',\nabla')$. 

As an example of how the simplicial presheaf structure maps operate on the above 1-simplex, consider the two coface maps $d^0:[0]\to [1]$ and $d^1:[0]\to [1]$. The corresponding face maps $d_0$ and $d_1$ operate on the above data as follows. The face map $d_0$ acts by
$$d_0(\varphi,  (C^1,C^2), q)=(M', (C^1,C^2),(P',\nabla)).$$
The face map $d_1$ acts by 
$$d_1(\varphi,  (C^1,C^2), q)=(M, (\varphi^*C^1,\varphi^*C^2), (P,\nabla)).$$
\end{example}

We also include the example given by taking the geometric structure to be a representable object, which is not a sheaf on $\frakFEmb_d$. We call this geometric structure on bordisms a \emph{geometric framing}. 

\begin{example}
\label{one.dimensional.bordisms}
Consider the representable enriched presheaf $(\RR\times U\to U)\in \frakFEmb_d$. The bordism category $\BBord_1^{\RR\times U\to U}$ has the following explicit description.

Fix $V\in \stcart$, $\langle 1\rangle\in \Gamma$, $[0]\in\Delta$, and $L\in \cart$.
A vertex in the simplicial set $\BBord_1^{\RR\times U\to U}(V,\langle 1\rangle,[0])(L)$ is given by a smooth $d$-manifold $M$, a smooth $L$-family of cuts $C=(C_{<},C_{=},C_{>})$ on $M\times V\to V$,
along with an $L$-family of fiberwise embeddings
$$i_t:M\times V\to \RR\times V, \quad t\in L,$$
and a smooth map $f:V\to U$.
For all $t\in L$, the restriction of $(C_=)_{t}\subset M\times V$ to each fiber $v\in V$ is embedded in $\RR\times \{f(v)\}$, via the smooth map $f$ and the embedding $i_t:M\into \RR$.  

\begin{figure}[ht]
\begin{center}
\tikzset{->-/.style={decoration={markings,mark=at position #1 with {\arrow[scale=1.5]{>}}},postaction={decorate}}}
\tikzset{-<-/.style={decoration={markings,mark=at position #1 with {\arrow[scale=-1.5]{>}}},postaction={decorate}}}
\begin{tikzpicture}
\draw [decorate,decoration={brace,amplitude=5pt, raise=4pt},yshift=0pt] (-.6,.2) -- (.6,.2); 
\node at (0,.7) {germ};
\draw (-2,0) -- (2,0);
\node at (0,0) {$\bullet$};
\node at (0,-.4) {$C_{=}$};
\draw (0,-.2) -- (0,.2);
\draw (-.4,.2) arc (90:270:.2);
\draw (.4,.2) arc (90:-90:.2);
\node at (-1.5,-.4) {$C_{<}$};
\node at (1.5,-.4) {$C_{>}$};
\draw[->] (0,-.7)--(0,-1.7);
\draw[->-=.7] (-2,-2) -- (2,-2);
\node at (.2,-1) {$h$};
\node at (0,-2.4) {$0$};
\draw (0,-1.8) -- (0,-2.2);
\node at (1.5,-1.7) {$\RR$};
\node at (1.5,.3) {$M\subset \RR$};
\end{tikzpicture}
\end{center}
\caption{A point $+_s$ ($s∈\RR$) in the bordism category with $V=U=\ast$ and $L=\RR^0$. 
The cut tuple is provided by a smooth map $h:\RR\to \RR$ ($t↦t-s$) and is defined by $C_<=h^{-1}(-\infty,0)=(-∞,s)$, $C_==h^{-1}(0)=\{s\}$ and $C_>=h^{-1}(0,\infty)=(s,∞)$.
The core of this bordism is precisely $C_==\{s\}$.
Likewise, the height function $t↦s-t$ produces another point, denoted by~$-_s$.}
\label{embedded.points}
\end{figure}
1-simplices in $\BBord_1^{\RR\times U\to U}(V,\langle 1\rangle,[0])(L)$ are given by an $L$-family of commutative triangles 
\begin{equation}\label{commtriangframing}
\xymatrix{
M\times V\ar[rr]^-{g_t}\ar[dr]_-{i_t} &&  M'\times V\ar[dl]^-{i_t}
\\
&\RR\times V,&
}
\end{equation}
where $t\in L$ and $g_t$ is the value of an element $g\in \frakFEmb_d(M\times V\to V,M\times V\to V)(L)$ at $t:\RR^0\to L$.
The map~$g_t$ is required to send the cut tuple $C_t$ on $M\times V\to V$ to the corresponding cut tuple $C'_t$ on $M'\times V\to V$.
Later, we will show that this category is actually discrete: we can replace $M$ by its embedded image in $\RR$, fiberwise over $f:V\to U$. 

For all $t\in L$ and $v\in V$, the embedding $i_{t,v}:M\to \RR$ induces a framing on $M$.
Since both $V$ and $L$ are connected, the family of framings is in the same connected component of the space of framings on $M$.
We call a $L$-point \emph{negative} (respectively \emph{positive})
if the flow induced by the framing flows in the direction of $(C_{t,v})_{<}$ (respectively $(C_{t,v})_{>}$), for some (hence all) $t\in L$ and $v\in V$.
Signs of $L$-points are not allowed to change in families. 
Moreover, $L$-points related by a 1-simplex must have the same sign, which is forced by the commutativity of the triangle \eqref{commtriangframing}.

\cref{embedded.points} provides examples of 0-bordisms given by positive and negative points.
More generally, given $s:U→\RR$, we have a 0-bordism~$+_s$ in $\BBord_1^{\RR\times U\to U}(V,\langle 1\rangle,[0])$
given by $M=\RR^1$, the fiberwise embedding $M⨯U→\RR⨯U$ given by the identity map,
the cut $[0]$-tuple given by the $U$-family of cuts induced by the height function $(m,u)↦m-s(u)$,
and the map $P:M⨯U→⟨1⟩$ that sends everything to $1∈⟨1⟩$.
We also have a 0-bordism~$-_s$ given by the same data except for the height function, which is now $(m,u)↦s(u)-m$.

We construct two $L=\RR$-parametrized isotopies in the smooth simplicial set $\BBord_1^{\RR\times U\to U}(V,\langle 1\rangle,[0])$.
The first isotopy $ρ_+$ has $M=\RR^1$,
the fiberwise embedding $M⨯U⨯L→\RR⨯U$ sends $(m,u,l)↦(m,u)$,
the cut $[0]$-tuple is given by the $U⨯L$-family of cuts induced by the height function $(m,u,l)↦m-l$,
and the map $P:M⨯U→⟨1⟩$ sends everything to $1∈⟨1⟩$.
Pulling back along the maps $\{-2\}→L$ and $\{2\}→L$ yields the 0-bordisms $+_{-2}$ and~$+_2$ respectively,
where the subscripts $-2$ and $2$ refer to constant functions $s:U→\RR$ with indicated values.
The second isotopy $ρ_-$ is given by pulling back the 0-bordism $-_0$ along the map $L→\RR^0$.

We can join the 0-bordisms $+_{-2}$ and $-_{0}$ in $\BBord_1^{\RR\times U\to U}$ by an elbow~$\epsilon$ as follows.
Take $M=\RR$, with the fiberwise open embedding $M⨯V→\RR⨯V$ given by the identity map.
Set $C^0_<=((-∞,-2)∪(0,∞))⨯V$, $C^0_==\{-2,0\}⨯V$, $C^0_>=(-2,0)⨯V$.
The cut $C^0$ is exhibited by the height function $(t,v)↦-t(t+2)$.
Set $C^1_<=M⨯V$ and $C^1_==C^1_>=∅$.
The cut $C^1$ is exhibited by the height function $(t,v)↦-1$.
Set $P:M⨯V→⟨1⟩$ to the map that sends $M⨯V$ to $1∈⟨1⟩=\{*,1\}$.
Altogether we have a 1-bordism from the disjoint union of 0-bordisms $+_{-2}$ and $-_0$,
defined by the height functions $(t,v)↦t+2$ and $(t,v)↦-t$ respectively, to the empty 0-bordism.

We can also encode the other elbow~$\eta$ using the same data with a different cut tuple~$C$.
Set $C^0_>=M⨯V$ and $C^0_==C^0_<=∅$.
The cut $C^0$ is exhibited by the height function $(t,v)↦1$.
Set $C^1_>=((-∞,0)∪(2,∞))⨯V$, $C^1_==\{0,2\}⨯V$, $C^1_<=(0,2)⨯V$.
The cut $C^1$ is exhibited by the height function $(t,v)↦t(t-2)$.
Altogether we have a 1-bordism from the empty 0-bordism to the 0-bordism given by the disjoint union
of 0-bordisms $-_0$ and $+_2$ defined by the height functions $(t,v)↦-t$ and $(t,v)↦t-2$ respectively.

\begin{figure}[ht]
\begin{center}
\tikzset{->-/.style={decoration={markings,mark=at position #1 with {\arrow[scale=1.2]{>}}},postaction={decorate}}}
\tikzset{-<-/.style={decoration={markings,mark=at position #1 with {\arrow[scale=-1.2]{>}}},postaction={decorate}}}
\begin{tikzpicture}
\begin{scope}[rotate=90]
\draw[->-=.7,thick] (6,1.4) -- (6,3.4);
\node at (6,2.4) {$\bullet$};
\node at (6.3,3) {$\RR$};
\node at (6.4,2.4) {$-_s$};
\node at (5.7, 2.4) {$s$};

\draw[-<-=.3,thick] (8,1.4) -- (8,3.4);
\node at (8,2.4) {$\bullet$};
\node at (8.3,3) {$\RR$};
\node at (8.4,2.4) {$+_{s}$};
\node at (7.7,2.4) {$s$};

\node at (8.3,-1) {$\RR$};

\draw[-<-=.7, ->-=.35, thick] (6,-2.4) -- (6,.4);
\node at (6,0) {$\bullet$};
\node at (6,-2) {$\bullet$};
\node at (6.4,0) {$+_{-2}$};
\node at (6.4,-2) {$-_{0}$};
\node at (5.7, -2) {$0$};
\node at (5.7, 0) {$-2$};
\node at (6.3,-1) {$\RR$};
\node at (5.7,-1) {$\epsilon$};

\draw[->-=.75, -<-=.25, thick] (8,-2.4) -- (8,.4);
\node at (8,0) {$\bullet$};
\node at (8,-2) {$\bullet$};
\node at (8.4,0) {$-_0$};
\node at (8.4,-2) {$+_2$};
\node at (7.7, -2) {$2$};
\node at (7.7, 0) {$0$};
\node at (7.7, -1) {$\eta$};
\end{scope}
\end{tikzpicture}
\end{center}
\caption{The 0-bordisms $+_s$ and $-_s$, for $s\in \RR$, and the 1-bordisms $\epsilon$ and $\eta$ in $\BBord_1^{\RR\times U\to U}$.}
\label{1bordisms.framed}
\end{figure}
\end{example}

\begin{example}
\label{d1bords}
The bordisms depicted in \cref{1bordisms.framed} exhibit $+_0$ as the dual of $-_0$. Observe that all $+$ points are isomorphic and all $-$ points are isomorphic in $\BBord_1^{\RR\times U\to U}$, via the isotopy $\rho_+$, respectively $\rho_-$ (\cref{one.dimensional.bordisms}). Hence, $+_0$ can be identified with both $+_{-2}$ and $+_{2}$. One of the triangle identities for the duality is given by composing the two 1-bordisms $\epsilon \sqcup +_{2}$ and $+_{-2}\sqcup \eta$:
\begin{center}
\begin{tikzpicture}
\tikzset{->-/.style={decoration={markings,mark=at position #1 with {\arrow[scale=1.2]{>}}},postaction={decorate}}}
\tikzset{-<-/.style={decoration={markings,mark=at position #1 with {\arrow[scale=-1.2]{>}}},postaction={decorate}}}

\draw[->-=.4, -<-=.8, thick] (-4.5,-3)--(-2,-3);
\draw[-<-=.15, ->-=.7, thick] (-2,-3)--(.5,-3);

\node at (-4,-3) {$\bullet$};
\node at (-2,-3) {$\bullet$};
\node at (0,-3) {$\bullet$};
\node at (-3,-3.3) {$\epsilon$};
\node at (-1,-3.3) {$\eta$};

\node at (-4,-2.6) {$+_{-2}$};
\node at (-2, -2.6) {$-_0$};
\node at (0,-2.6) {$+_2$};
\node at (-4,-3.3) {$-2$};
\node at (-2,-3.3) {$0$};
\node at (0,-3.3) {$2$};
\end{tikzpicture}
\end{center}
The composition is the interval
\begin{center}
\begin{tikzpicture}
\tikzset{->-/.style={decoration={markings,mark=at position #1 with {\arrow[scale=1.2]{>}}},postaction={decorate}}}
\tikzset{-<-/.style={decoration={markings,mark=at position #1 with {\arrow[scale=-1.2]{>}}},postaction={decorate}}}

\draw[->-=.4, thick] (-4.5,-3)--(-2,-3);
\draw[ ->-=.7, thick] (-2,-3)--(.5,-3);

\node at (-4,-3) {$\bullet$};
\node at (0,-3) {$\bullet$};

\node at (-4,-2.6) {$+_{-2}$};
\node at (0,-2.6) {$+_2$};
\node at (-4,-3.3) {$-2$};
\node at (0,-3.3) {$2$};
\end{tikzpicture}
\end{center} 
which is isomorphic to the identity via the isotopy of cut tuples that fixes $+_{-2}$ and moves $+_{2}$ to $+_{-2}$ by the isotopy~$\rho_+$.
For the other triangle identity, replace $\epsilon$ by the isotopic bordism obtained by shifting $\epsilon$ by 4 units to the right and then form the composition. 
\end{example}

\begin{example}
\label{frcircle}
\def\fBBord{R\BBord}
Continuing \cref{d1bords},
we can also compute the trace of $+_0$ by composing the 1-bordism~$\eta$,
the appropriate braiding 1-bordism from the target of~$\eta$ to the source of~$\epsilon$, and the 1-bordism $\epsilon$.
To construct the braiding 1-bordism, consider the isotopies $\rho_+$ (from $+_{-2}$ to $+_2$) and $\rho_-$ (from $-_0$ to itself).
Their disjoint union would give us the desired isotopy from $-_0⊔+_2$ to $+_{-2}⊔-_0$.
However, the disjoint union does not map to $\RR⨯U$ via a fiberwise embedding, but merely via a fiberwise etale map,
since $\rho_+$ and $\rho_-$ have the same embedding as a geometric structure. 
Thus, the disjoint union of $\rho_+$ and $\rho_-$ does not give an isotopy in $\BBord_1^{\RR\times U\to U}$.
Once we map $\BBord_1$ to its fibrant replacement $\fBBord_1$, this problem disappears
and we can use the Čech descent property for the cover $\RR\sqcup \RR\to \RR$ in $\cart$
to glue the isotopies $\rho_+$ and $\rho_-$ into a single element $\rho=ρ_+⊔ρ_-\in \fBBord_1^{\RR⨯U→U}(\RR^0,⟨1⟩,[0])(\RR^1)_0$.

Use \cref{companion} to convert $\rho$ into a 1-bordism $\hat \rho$ using the point $2\in \RR$ as a source and $-2\in \RR$ as a target.
The target of $\eta$ is not quite the same as the source of $\hat \rho$.
The former is two $U$-families of points cut out from $M_0⨯U=\RR⨯U$, while the latter looks like two points cut out from $M_2⨯U=(\RR\sqcup \RR)⨯U$.
Setting $M_1=(-1,1)\sqcup (1,3)$, we have an embedding $M_1\into M_0$ whose image is the union $(-1,1)\cup (1,3)\subset \RR$.
We also have an embedding $M_1\into M_2$ whose image is the union of the interval $(-1,1)$ in the second copy of $\RR$ in the disjoint union,
and the interval $(1,3)$ in the first copy of $\RR$. 
We induce a cut tuple $C$ on $M_1⨯U→U$ by pulling back the cut tuple on $M_0⨯U→U$ along the embedding.
The map $P$ sends everything to $1\in \langle 1\rangle$.
The resulting triple $(M_1,C,P)$ with the inclusion $M_1⨯U→\RR⨯U$ as a geometric structure is a 0-bordism.
The embeddings $M_1→M_0$ and $M_1→M_2$ induce 1-simplices, i.e., elements of $\BBord_1^{\RR⨯U→U}(\RR^0,⟨1⟩,[0])(\RR^0)_1$, between the corresponding bordisms.
Thus, the target of~$\eta$ is connected to the source of $\hat\rho$ via a zigzag of such 1-simplices.
Use \cref{companion} to convert the resulting zigzag of 1-simplices to two 1-bordisms.
Perform an analogous construction for the target of~$\hat\rho$ and source of~$\epsilon$, which yields two more 1-bordisms.

Having converted all entities to seven elements of $\fBBord_1^{\RR⨯U→U}(\RR^0,⟨1⟩,[1],\RR^0)_0$,
we use the Segal property in the direction of the first~$Δ$ to compose these 1-bordisms into a single element of $\fBBord_1^{\RR⨯U→U}(\RR^0,⟨1⟩,[1],\RR^0)_0$.
The resulting element is a presentation of the framed circle.

Consider the geometric structure of topological framings,                                     
which sends $M→U$ to the smooth set of trivializations of the fiberwise tangent bundle of $M→U$.
Consider also the map from the enriched geometric structure represented by $\RR⨯U→U$ to the weakly equivalent geometric structure of topological framings $F$,
which sends a fiberwise embedding of $N→V$ into $\RR⨯U→U$ to the induced framing on $N→V$.
Applying the functor $\BBord_1^{\RR⨯U→U}→\BBord_1^F$ to the maps $η$ and~$ε$,
makes composing them much easier, since after discarding the data of the embedding itself,
the codomain of~$η$ is simplicially homotopic to the domain of~$ε$, where the homotopy is implemented as a zigzag of 1-simplices similar to $M_0←M_1→M_2$.
Applying \cref{companion} as before, we now have to compose only four entities instead of seven.
Here $\BBord_1^F$ can be seen as a partial fibrant replacement of $\BBord_1^{\RR⨯U→U}$,
since $\BBord_1^F$ is local with respect to all maps in \cref{bousloc} except for completion maps.
\end{example}

\begin{remark}
The bordism category defined above does encode ordinary closed manifolds equipped with a geometric structure.
This is apparent in all the previous examples, except possibly the example of geometric framings (\cref{one.dimensional.bordisms}, \cref{d1bords}, \cref{frcircle}).
In the case of geometric framings, the geometric structure does not satisfy descent on the enriched site $\frakFEmb_1$.
When the geometric structure does not satisfy descent,
the Segal gluing condition will be violated and the corresponding bordism category will not be local with respect to the Segal maps.
In such cases, we must first sheafify the geometric structure in order to see an object in the fibrant replacement of the bordism category.
More precisely, any element in the fibrant replacement of the isotopy bordism category with geometric structure $\gs$ is homotopic to an honest bordism (no fibrant replacement) with geometric structure in the sheafification of $\gs$. 
\end{remark}

\section{Verification of the axioms}
\label{axioms.section}

In this section, we will prove that the bordism category without isotopies $\Bord_{d}$ satisfies the axioms \ax1--\ax3 (\cref{axioms})
and that the bordism category with isotopies $\BBord_{d}$ satisfies \frakax1--\frakax3 (\cref{multiple.isotopy.axiom}).

\subsection{Proof of the first two axioms}
We begin by proving that our bordism category satisfies the axiom \ax2 of \cref{axioms} and the axiom \frakax2 of \cref{multiple.isotopy.axiom}.
Axiom \ax1 was established in \cref{bordcts}. 

\begin{proposition}
\label{a1a2}
Fix $d\geq 0$.
Recall the definition of $\Struct_d$ (\cref{geometric.structure}) and $\frakStruct_{d}$ (\cref{geometric.structure.isotopy}).
Denote by $\Struct_{d,\inj}$ and $\frakStruct_{d,\inj}$ the injective model structure on $\sPSh(\FEmb_d)$ and $\sm \sPSh(\frakFEmb_d)$, respectively. 
The functors
$$\Bord_{d,\uple}:\Struct_{d,\inj}\to \smcatuple_{\infty,d}, \quad \gs\mapsto \Bord_{d,\uple}^\gs$$
and
$$\Bord_{d}:\Struct_{d,\inj}\to \smcat_{\infty,d}, \quad \gs\mapsto \Bord_{d}^\gs$$
are $\sset$-enriched left Quillen functors that preserve all weak equivalences. 
In particular, they are homotopy cocontinuous.
Similarly, the functors
$$\BBord_{d,\uple}:\frakStruct_{d,\inj}\to \fraksmcatuple_{\infty,d}, \quad \gs\mapsto \BBord_{d,\uple}^\gs$$
and
$$\BBord_{d}:\frakStruct_{d,\inj}\to \fraksmcat_{\infty,d}, \quad \gs\mapsto \BBord_{d}^\gs$$
are $\smsset$-enriched left Quillen functors that preserve all weak equivalences.
In particular, they are homotopy cocontinuous.
\end{proposition}

\begin{proof}
We first claim that all four functors preserve monomorphisms and weak equivalences, and are enriched cocontinuous.
The enrichment is given by simplicial sets in the case of $\Bord_d$ and smooth simplicial sets in the case of $\BBord_d$.

Analyzing \cref{bordstr}, we first observe that taking the nerve yields a simplicial diagram whose simplicial set of $n$-simplices is given by 
$$\coprod_{(M_0,C_0,P_0)\to \cdots \to (M_n,C_n,P_n)}\gs(M_n\times U\to U),$$ 
where the disjoint union runs over composable chains of morphisms in \cref{bord}. For $\BBord_d$, the simplicial diagram is similar after evaluating at $L\in \cart$, except that $\gs(M_n\times U\to U)$ is replaced by $\gs(M_n\times U\to U)(L)$ and $U$ is replaced by $U⨯L$ in \cref{bord}.

Given the above formula, the simplicial diagram associated to a geometric structure $\gs$ manifestly preserves objectwise weak equivalences, monomorphisms, and tensoring, for both $\Bord_d$ and $\BBord_d$. Since colimits of simplicial presheaves can be computed objectwise and colimits commute with coproducts, colimits in $\gs$ can be passed outside the disjoint union at each simplicial level $n$. The case of $\BBord_d$ is analogous, with simplicial sets replaced by smooth simplicial sets. 

Next, we observe that taking the diagonal of a bisimplicial set preserves colimits, objectwise weak equivalences, tensoring, and monomorphisms.
From this, it follows that both $\Bord_d$ and $\BBord_d$ are enriched left Quillen functors.
The right adjoint to $\Bord_d$ is given by the functor that send $X\in \smcat_{\infty,d}$ to field theories with values in~$X$: 
$$X\mapsto \left((W\to U)\mapsto \map(\Bord_d^{W\to U},X)\right).$$
Here $\map(-,-)$ denotes the simplicial enrichment of $\smcat_{\infty,d}$.
For $\BBord_d$ the right adjoint is defined in a similar way, with $\BBord_d$ replacing $\Bord_d$ and the the $\smsset$-enrichment replacing the $\sset$-enrichment.

The case of $\Bord_{d,\uple}$ and $\BBord_{d,\uple}$ is proved completely analogously, by dropping the globularity condition (\cref{globgrid}) on cut ${\bf m}$-grids. 
\end{proof}

It remains to prove that $\Bord_{d}$ satisfies axiom \ax3 and that $\BBord_{d}$ satisfies axiom \frakax3.
This will occupy the the next two subsections.

\begin{remark}
\label{uple.globular}
In the rest of \cref{axioms.section} and in \cref{codescent}, we will work in the globular case only.
The multiple versions are defined by removing the word globular from globular cut ${\bf m}$-grids.
We remark that the proof of the main theorem also works in the uple case.
In more detail, there is an obvious axiomatization of the uple category of bordisms with the same \ax1 and \ax2 axioms and with the \ax3 axiom obtained by removing the word globular everywhere below.
The proof of axiom \ax3 in this section carries through verbatim in the uple case.
All the constructions in \cref{codescent} also work in the uple case, by just forgetting that the cut tuples are globular.
\end{remark}

\subsection{Embedded bordism categories}

We will define functors
$$\EBord_d:\FEmb_d\to \smcat_{\infty,d}, \qquad \frakEBord_d:\frakFEmb_d\to \fraksmcat_{\infty,d},$$
which serve as a model for the bordism category (respectively bordism category with isotopies) in the case of representable geometric structures.
Axiom \ax3 (respectively \frakax3 asserts that $\Bord_{d}^{p:W\to U}$ and $\BBord_{d}^{p:W\to U}$ are weakly equivalent to $\EBord_d^{p:W\to U}$ and $\frakEBord_d^{p:W\to U}$, respectively.
Moreover, the equivalences are natural in $p:W\to U$. 

To begin, let us recall that if $p:W\to U$ is a submersion with $d$-dimensional fibers and $f:V\to U$ is a smooth function, we can form the pullback 
$$
\xymatrix{
W\times_{U}V\ar[r]^-{j}\ar[d]^-{\pi} & W\ar[d]^-{p}
\\
V\ar[r]^-{f} & U.
}
$$
The pullback square defines a morphism in $\FEmb_d$. 

\begin{definition}
\label{embeddedbordcat.unenriched}
Fix $d≥0$.
We define a functor 
$$\EBord_d:\FEmb_d\to \smcat_{\infty,d}$$ 
as follows. 
Let $p:W\to U\in \FEmb_d$.
We define the value of $\EBord_d$ at $p$ as follows.
Fix $V\in \stcart$, $\langle \ell\rangle\in \Gamma$, and ${\bf m}\in \Delta^{\times d}$. The set $\EBord_d^{p}(V,\langle \ell\rangle,{\bf m})$, whose elements $[f,N,C,P]$ we will call \emph{embedded bordisms}, is defined as the quotient of the set of quadruples:
\begin{enumerate}
\item[(1)] a morphism $f:V\to U$ in $\stcart$;
\item[(2)] an open subobject $N\subset W\times_{U}V$ such that the restriction $\pi:N\to V$ is a (locally trivial) fiber bundle;
\item[(3)] a compact (globular) cut ${\bf m}$-grid $C$ for the projection $\pi:N\to V$;
\item[(4)] a morphism $P:N\to \langle \ell\rangle$ in $\stman$;
\end{enumerate}
by the equivalence relation
\begin{itemize}
\item  Two quadruples $(f,N,C,P)$ and $(f',N',C',P')$ are equivalent if and only if
\begin{enumerate}
\item[(5)] $f=f'$;
\item[(6)] There is an open subobject $N''\subset N\cap N':=N\times_{W\times_{U}V}N'$ and a compact (globular) cut ${\bf m}$-grid $C''$ on $\pi:N''\to V$ such that $\stu(N'')$ contains the core of both $C$ and $C'$ and $C''=C=C'$ and $P=P'$, after restricting to $N''$.
\end{enumerate}
\end{itemize} 
The structure maps for the presheaf are defined the same way as in \cref{bord} in the $\Delta^{\times d}$ and $\Gamma$ directions.
For the $\stcart$ direction, given a smooth map $g:V\to V'$ we pullback the submersion along the composition $f\circ g$.
Cut grids can be pulled back along the induced fiberwise diffeomorphism $h:W\times_{U}V\to W\times_{U}V'$.
The equivalence relation is also preserved by such pullbacks (by pulling back the data in (6) along the fiberwise diffeomorphism $h$). 

Finally, to turn $\EBord_d^p$ into a functor in $p$,
observe that the structure maps
$$\FEmb_d(p_1,p_2)⨯\EBord_d^{p_1}(V,\langle \ell\rangle,{\bf m})→\EBord_d^{p_2}(V,\langle \ell\rangle,{\bf m})$$
can be defined by sending a morphism $(φ:W_1→W_2,g:U_1→U_2)$ in $\FEmb_d$
together with an equivalence class $[f,N,C,P]$
to the equivalence class of $[gf,(φ⨯_g V)N,C,P]$,
where the $N$ is interpreted as an open embedding that is composed with the open embedding $$\varphi⨯_g V:W_1⨯_{U_1}V→W_2⨯_{U_2}V.$$
\end{definition}

\begin{remark}
In Condition~(2) of \cref{embeddedbordcat.unenriched}, the requirement that $\pi$ is a (locally trivial) fiber bundle is not needed when $\stcart=\cart$. 
A modification of Ehresmann's theorem implies that for a sufficiently small open neighborhood of the core (recall the properness condition in \cref{compactgrid}), the map $\pi$ is locally trivial. 
However, in the general case (for example, supermanifolds) this condition is needed because we chose to work with trivial bundles in  \cref{bord}. 
\end{remark}

\begin{remark}
\label{smallercore}
Any open subobject of $N$ whose reduction contains the core contains a smaller open subobject satisfying Condition~(2) whose reduction also contains the core.
Thus, $N''$ in Condition~(6) of the equivalence relation in \cref{embeddedbordcat.unenriched} can be chosen so that $\pi:N''\to V$ is a (locally trivial) fiber bundle. 
\end{remark}

Morally, the elements in $\EBord_d^p(V,\langle \ell\rangle,{\bf m})$ are $V$-families of germs of cores embedded in $W$; one imagines a family of diced cubes sitting in $W$ with each diced cube having a germy neighborhood in which the cuts extend. Of course, we allow for more than just diced cubes, but this picture is helpful in that it emphasizes the essential aspects: the combinatorics of cut tuples, the core of the bordism, and the germ.

\begin{definition}
\label{bracketcore}
Every embedded bordism $[f,N,C,P]\in \EBord_d^p(V,\langle\ell\rangle,{\bf m})$ has an associated \emph{core} 
$$\core[f,N,C,P]=C_{[0,{\bf m}]}\setminus \stu(P)^{-1}(\ast)\subset \stu(W\times_{U}V)$$
(see \cref{cutgrid} for $C_{[0,{\bf m}]}$), which does not depend on the representative.
In other words, we have a well defined function 
$$\core:\EBord_d^p(V,\langle \ell\rangle,{\bf m})\to \mathcal{P}(\stu(W\times_{U}V)),$$
where $\mathcal{P}$ denotes the power set. 
\end{definition}

Next we turn to the model for the bordism category with isotopies. The definition is almost verbatim the same as the the model without isotopies, except that we now have $L$-families of cut tuples.

\begin{definition}
\label{embeddedbordcat.enriched}
Fix $d≥0$. We define a $\smsset$-enriched functor 
$$\frakEBord_d:\frakFEmb_d\to \fraksmcat_{\infty,d}$$ 
as follows.
Let $p:W\to U\in \frakFEmb_d$.
We define the value of $\frakEBord_d$ at $p$ as the object in $\fraksmcat_{\infty,d}$ obtained by pulling back $\EBord_d^p$ along the functor 
$$\stcart\times \Gamma\times \Delta^{\times d}\times \cart\to \stcart\times \Gamma\times \Delta^{\times d}$$
that sends $(U,\langle \ell\rangle,{\bf m},L)\mapsto (U⨯L, \langle \ell\rangle, {\bf m})$, and then taking the subobject whose value at $(U,\langle \ell\rangle,{\bf m},L)$ is given by the set in \cref{embeddedbordcat.unenriched}, with $V⨯L$ replacing $V$, but with Condition~(1) replaced by
\begin{enumerate}
\item[(1)] A morphism $f:V⨯L\to U$ in $\stcart$ that factors as a composition of the projection $V⨯L\to V$ and a morphism $V\to U$. 
\end{enumerate}
The structure maps for the presheaf are defined the same way as in \cref{bord} in the $\Gamma$ and $\Delta^{\times d}$ directions.
For the $\stcart$ direction, given a smooth map $g:V\to V'$ we pullback the submersion along the composition $f\circ g$.
$L$-families of cut grids can be pulled back along the induced fiberwise diffeomorphism $h:W\times_{U}V\to W\times_{U}V'$.
The equivalence relation is also preserved by such pullbacks (by pulling back the data in Condition~(6) along the fiberwise diffeomorphism~$h$).

Finally, to turn $\frakEBord_d^p$ into $\smsset$-enriched functor in $p$,
observe that the enriched structure maps
$$\frakFEmb_d(p_1,p_2)⨯\frakEBord_d^{p_1}(V,\langle \ell\rangle,{\bf m})→\frakEBord_d^{p_2}(V,\langle \ell\rangle,{\bf m})$$
can be defined by applying the structure map of \cref{embeddedbordcat.unenriched} in $L$-families, for every $L∈\cart$.
The resulting operation preserves Condition~(1).
\end{definition}

\begin{remark}
\label{reduction.to.enriched}
It is immediate from the definition of $\frakEBord_d^p$ that evaluation at $L=\RR^0$ yields precisely $\EBord_d^p$, since an $\RR^0$-family of cut grids is precisely an ordinary cut grid. 
\end{remark}

\subsection{Stalks of embedded bordism categories}
We now analyze the stalk of the isotopy embedded bordism category $\frakEBord_d^p$ at a point $\mathfrak{p}$ of $\stcart$ (\cref{stalks}) and a point $\mathfrak{l}$ of $\cart$.
The case of $\EBord_d^p$ then follows formally by taking $\mathfrak{l}=(\RR^0,\id)$ (\cref{reduction.to.enriched}).
This analysis leads to simplifications of the proof, by eliminating parametrizing families and isotopies of cut grids.

\begin{proposition}
\label{bordstalkisot}
Fix $d\geq 0$, $p:W\to U\in \frakFEmb_d$, a point $\mathfrak{p}=(T,\rho)$ of $\stcart$ (\cref{stalks}), a point $\mathfrak{l}=(\RR^n,\id)$ of $\cart$ (\cref{stalks}), $\langle \ell\rangle\in \Gamma$, and ${\bf m}\in \Delta^{\times d}$.
Then the corresponding double stalk $\mathfrak{l}^*\mathfrak{p}^*\frakEBord_d^p(\langle\ell\rangle,{\bf m})$ is the quotient of the set of quadruples in \cref{embeddedbordcat.enriched},
taken for all possible $V=B^{\mathfrak{p}}_{\delta}\subset T\in \stcart$ and $L=B^{\mathfrak{l}}_{\delta}\subset\RR^{n}$ (\cref{stalks}), by the following equivalence relation.
\begin{itemize}
\item  Two quadruples $(f,N\to  B^{\mathfrak{p}}_{\delta}\times B^{\mathfrak{l}}_{\delta} ,C,P)$ and $(f',N'\to B^{\mathfrak{p}}_{\delta'}\times B^{\mathfrak{l}}_{\delta'},C',P')$ are equivalent if there is $\epsilon<\min\{\delta,\delta'\}$ such that the restriction of the two quadruples to $B^{\mathfrak{p}}_{\epsilon}\times B^{\mathfrak{l}}_{\epsilon}$ are equivalent in the sense of \cref{embeddedbordcat.unenriched}. 
\end{itemize}
Moreover, in  \cref{embeddedbordcat.unenriched}, we may optionally require the following extra conditions.
\begin{enumerate}
\item  $N=M\times V⨯L$ and $\pi:M\times V⨯L\to V⨯L$ is the projection map, for some $d$-manifold $M$,
\item $N''$ and $\pi$ in the equivalence relation of \cref{embeddedbordcat.unenriched} is a product and projection, respectively.
\item $N$ can be chosen so that the inclusion of the fiber of the core
$$\core[f,N,C,P]_0≔\core[f,N,C,P]\cap \stu(W)_{\stu(f)(0,\rho^{-1}(0))}\into \stu(N)$$
(\cref{bracketcore})
is bijective on $\pi_0$.
\end{enumerate}
Condition (3) above gives a well defined $\pi_0$ for an element $[f,N,C,P]\in \mathfrak{l}^*\mathfrak{p}^*\frakEBord_d^p(\langle\ell\rangle,{\bf m})$, given by 
$$\pi_0[f,N,C,P]:=\pi_0(\core[f,N,C,P]_0)\overset{\cong}{\to} \pi_0(N)$$ 
with $N$ is an arbitrary open subset satisfying (3). This partitions the class $[f,N,C,P]$ into a family of classes $[f,N_j,C_j,P_j]_{j\in \pi_0[f,N,C,P]}$, where $N_j$ is the connected component of $N$ corresponding to $j\in \pi_0(N)$, and $C_j$ and $P_j$ are pulled back along the inclusion $N_j\into N$.
\end{proposition}

\begin{proof}
This follows immediately from the calculation of the double stalk at $\mathfrak{p}$ and $\mathfrak{l}$ in \cref{doublestalk} and \cref{embeddedbordcat.enriched}. The optional requirement (1) follows immediately from local triviality condition on $\pi$. The optional requirement (2) is \cref{smallercore}. That $N$ can be chosen as in (3) follows by passing to an open subobject of $N$ given by the lifting (\cref{def.cartsp}) of a union of disjoint connected open subsets, each of which contains a single connected component of the fiber of the core and restricting the data of $C$ and $P$ to this open subset. Such subsets exist by compactness of the fiber of the core. The local triviality condition on $\pi:N\to V⨯L$ can be ensured by passing to a smaller open subobject of $N$ using \cref{smallercore}.
\end{proof}

\begin{remark}
\label{0stalkbord}
Continuing \cref{reduction.to.enriched}, observe that for $\mathfrak{l}=(\RR^0,\id)$, we have 
$$\mathfrak{l}^*\mathfrak{p}^*\frakEBord_d^{p}=\mathfrak{p}^*\EBord_d^p,$$
since $B_{\delta}^{\mathfrak{l}}=\RR^0$ and by \cref{reduction.to.enriched}, evaluating $\frakEBord_d^p$ at $L=\RR^0$ yields $\EBord_d^p$.
This will allow us to prove statements about $\mathfrak{p}^*\EBord_d^p$ in \cref{codescent},
by proving a claim for $\mathfrak{l}^*\mathfrak{p}^*\frakEBord_d^p$ and then taking $\mathfrak{l}=(\RR^0,\id)$. 
\end{remark}

\subsection{Comparison of bordisms and embedded bordisms}
The goal of this subsection is to construct a comparison map $\frake_d$ from $\BBord_d$ to $\frakEBord_d$, which will be shown to be a weak equivalence in \cref{a3axiom}, proving axiom \frakax3.
The corresponding axiom \ax3 for $\Bord_d$ will follow formally, by taking the stalk of the map $\frake_d$ at $\mathfrak{l}=(\RR^0,\id)$, using \cref{reduction.to.enriched} and \cref{frakbordl0}.

\begin{construction}
\label{emapconstr}
We will define a natural transformation of $\smsset$-enriched functors
\begin{equation}\label{embmans}\frake_d:\BBord_d\to \frakEBord_d,\end{equation}
where we regard $\BBord_d$ (\cref{enrichedbordstr}) as a functor $\frakFEmb_d\to \fraksmcat_{\infty,d}$, by restricting along the $\smsset$-enriched Yoneda embedding $\frakFEmb_d\into \frakStruct_d$. We proceed in stages.

\medskip
\paragraph{\bf Step 1}
For a fixed $(p:W\to U)\in \frakFEmb_d$ and $(V,\langle \ell\rangle,{\bf m},L)\in \stcart\times \Gamma\times \Delta^{\times d}\times \cart$, the $(V,\langle \ell\rangle,{\bf m},L)$-component of the corresponding map $\frake^p_d$ is a map of simplicial sets
$$\frake_d^p(V,\langle \ell\rangle,{\bf m},L):\BBord_d^{p:W\to U}(V,\langle \ell\rangle,{\bf m},L)\to \frakEBord_d^p(V,\langle \ell\rangle,{\bf m},L),$$ 
whose domain is the nerve of a category and whose codomain is the nerve of a discrete category.
Hence, to define the map $\frake_d^p$ it suffices to define functor 
between the corresponding categories. 

\medskip
\paragraph{\bf Step 2} We analyze the simplicial set $\BBord_d^{p:W\to U}(V,\langle \ell\rangle,{\bf m},L)$ (\cref{enrichedbordstr}), which is given by the nerve of the following category.

An object is given by a quadruple $(M,C,P,(i,f))$, where $(M,C,P)$ is a bordism in the isotopy bordism category (\cref{frakbord}) and $(i,f):(M\times V⨯L\to V⨯L)\to (p:W\to U)$ is an $L$-point in $\frakFEmb_d(M\times V\to V, W\to U)$.

A morphism is a commutative triangle
\begin{equation}\label{commtriagemb}
\xymatrix{
M\times V⨯L\ar[rr]^{g}\ar[dr]_{i} && M'\times V⨯L\ar[dl]^-{i'}
\\
& W, &
}
\end{equation}
where $g:M\times V⨯L\to M'\times V⨯L$ is a fiberwise embedding covering $\id_{(V⨯L)}$ satisfying Condition (m) in \cref{bord}, i.e., $g$ is a cut-respecting fiberwise embedding.
Since $i'$ is fiberwise injective over $f'$, this implies that there is at most one $g$ making the diagram commute. It follows that this category is a preorder.

\medskip
\paragraph{\bf Step 3} By definition of an $L$-point in $\frakFEmb_d(M\times V\to V,W\to U)$, the pair $(i,f)$ gives a diagram 
$$\xymatrix{M\times V⨯L\ar[r]^-i\ar[d] & W\ar[d]^-{p}\\
V⨯L\ar[r]^f & U,}$$
where $i$ is a fiberwise embedding and $f$ is the composition of the projection $V⨯L\to V$ and a map $V\to U$.
By the universal property of the pullback, we have a canonical map 
$$u:M\times V⨯L\to W\times_{U}(V⨯L),$$
which is an open embedding.
Thus, $u$ determines an open subobject of $W\times_U (V⨯L)$.

We define the functor whose nerve is the map $\frake_d^p(V,\langle \ell\rangle,{\bf m},L)$ as follows. An object given by a quadruple $(M,C,P,(i,f))$ is sent to the object in $\frakEBord_d^p(V,\langle \ell\rangle,{\bf m},L)$ (\cref{embeddedbordcat.enriched}) given by quadruple $[f,N,C,P]$, where $N$ is the subobject determined by the open embedding $u$ in Step 2. A morphism $g:M\times V⨯L\to M'\times V⨯L$ is sent to the identity on $[f,N,C,P]$. This functor is well defined by the equivalence relation in \cref{embeddedbordcat.enriched}, the commutativity of the diagram \cref{commtriagemb}, and the fact that $g$ is a cut-respecting open embedding. 
\end{construction}

\subsection{Proof of the third axiom}
\label{a3axiom}
In this subsection, we prove that the comparison map in \cref{emapconstr} is a weak equivalence. The basic idea is to work stalkwise to remove parametrizing families of bordisms. By Step~1 in \cref{emapconstr}, the map $\frake_d$ is induced by taking the nerve of a functor between categories. We apply Quillen's Theorem A to show the comma categories are contractible. Roughly speaking, the comma categories are contractible because the embedded image of a bordism gives a terminal object in the comma categories. 

\begin{proposition}
\label{emaplociso}
Fix $d\geq 0$.
The natural transformation $$\frake_d:\BBord_d\to \frakEBord_d,$$ (\cref{emapconstr}) is an objectwise weak equivalence.
\end{proposition}

\begin{proof}
Fix $(p:W\to U)\in \frakFEmb_d$, $\langle \ell\rangle\in \Gamma$, ${\bf m}\in \Delta^{\times d}$, a point $\mathfrak{l}$ of $\cart$, and a point $\mathfrak{p}$ of $\stcart$. By \cref{stalks}, it suffices to prove that  the simplicial map
$$\mathfrak{l}^*\mathfrak{p}^*\frake^p_{d}(\langle \ell\rangle,{\bf m}):\mathfrak{l}^*\mathfrak{p}^*\BBord_d^{p:W\to U}(\langle \ell\rangle,{\bf m})\to \mathfrak{l}^*\mathfrak{p}^*\frakEBord_d^{p}(\langle \ell\rangle,{\bf m})$$
is a weak equivalence. 

Taking the $\mathfrak{l}$ and $\mathfrak{p}$ stalks of the functor from Steps 2 and 3 of \cref{emapconstr} gives a functor $e$ whose nerve is the map $\mathfrak{l}^*\mathfrak{p}^*\frake^p_{d}(\langle \ell\rangle,{\bf m})$.
We apply Quillen's Theorem A to the functor $e$ 
and verify that the comma categories are contractible. Fix an object $[f,N,C,P]$ in the codomain of $e$. By \cref{bordstalkisot}, we can assume $V=B^{\mathfrak{p}}_{\epsilon}$, $L=B^{\mathfrak{l}}_{\epsilon}$, $N=M\times V⨯L$, and $\pi:N\to V⨯L$ is the projection. By Step 1 in \cref{emapconstr}, the target category is discrete. Therefore, the comma category $[f,N,C,P]\downarrow e$ is precisely the full subcategory of the domain that maps to $[f,N,C,P]$. We denote this full subcategory by ${\sf D}$. 

Consider the quadruple $(M,C,P,(j\circ \iota, f))$, where $M$ is the factor of $N$ indicated above and
$$j\circ \iota: N\overset{\iota}{\into} W\times_{U}V⨯L\overset{j}{\to} W,$$
where $\iota$ is the open embedding in \cref{embeddedbordcat.unenriched} and $j$ is the projection. 
By definition of the functor $e$,
$$e(M,C,P,(j\circ \iota, f))=[f,N,C,P],$$
which shows that ${\sf D}$ is nonempty.

By Step 2 in \cref{emapconstr}, ${\sf D}$ is a  preorder. To prove that ${\sf D}$ is cofiltered, it suffices to show that for any two objects $d$ and $d'$, there is an object $c$ and morphisms $c\to d$ and $c\to d'$. Assume $d=(M,C,P,(\varphi,f))$, $d'=(M',C',P',(\varphi',f'))$, and $e(M,C,P,(\varphi,f))=[f,N,C,P]=e(M',C',P',(\varphi',f'))$. Then the quadruples $(f,N=M\times V⨯L, C,P)$ and $(f',N'=M'\times V⨯L, C',P')$ define equivalent elements of the set in \cref{embeddedbordcat.unenriched}. Therefore, 
\begin{itemize}
\item $f=f'$; 
\item there is $N''$, which we can take to be of the form $M''\times V⨯L$ by \cref{bordstalkisot};
\item there are open embeddings $i:N''\into N$ and $i':N''\into N'$ that commute with the maps to $W$;
\item there is a compact (globular) cut ${\bf m}$-grid $C''$ on the projection $\pi:N''\to V⨯L$ such that $C$ and $C'$ restrict to $C''$;
\item there is a map $P'':N''\to \langle \ell\rangle$ such that $P$ and $P'$ restrict to $P''$.
\end{itemize}
Set $c=(M'',C'',P'',(\varphi\circ i,f))$. Then $c$ maps to both $d$ and $d'$ by construction and we are done.
\end{proof}

This completes the verification of axiom \frakax3.
We are now ready to complete the proof of \cref{existencebord2} and deduce \cref{existencebord} from it.

\begin{theorem}
\label{proofof106}
\label{bordisms.with.isotopies.representable}
There exists a smooth symmetric monoidal $(\infty,d)$-category of bordisms with isotopies.
That is, the bordism category $\BBord_d$ satisfies the three axioms \frakax1--\frakax3 of \cref{axioms,multiple.isotopy.axiom}.
Likewise for the case of $\BBord_{d,\uple}$.
\end{theorem}

\begin{proof}
By \cref{emaplociso}, the map $\frake_d$ in \cref{embmans} is an equivalence.
This proves axiom \frakax3.
By \cref{a1a2}, axiom \frakax1 and \frakax2 are also satisfied. 
For the uple case, apply \cref{uple.globular}.
\end{proof}

\begin{theorem}
\label{proofof104}
There exists a smooth symmetric monoidal $(\infty,d)$-category of bordisms.
That is, the bordism category $\Bord_d$ satisfies the three axioms \ax1--\ax3 of \cref{axioms}.
Likewise for $\Bord_{d,\uple}$.
\end{theorem}

\begin{proof}
By \cref{a1a2}, axioms \ax1 and \ax2 are satisfied. 
The axiom \ax3 follows from \cref{emaplociso} by evaluating at $L=\RR^0$
and applying \cref{0stalkbord} and \cref{reduction.to.enriched}.
In the uple case, use \cref{uple.globular}.
\end{proof}

We have completed the verification of the axioms.
The remainder of the paper will be concerned with proving \cref{axtheorem} and \cref{axtheorem2}.

\section{Codescent for bordism categories}
\label{codescent}

In this section, we will complete the proof of \cref{axtheorem} and \cref{axtheorem2}.
These results formally follow from \cref{codescent.embedded},
which establishes the codescent property for the embedded bordism category.

\begin{theorem}
\label{axioms.contractible.codescent}
Fix $d\geq 0$.
Consider the following relative category.
\begin{itemize}
\item Objects are $\BBord'_d$ satisfying the axioms of \cref{multiple.isotopy.axiom},
with a fixed choice of a weak equivalence $ι^*\BBord'_d→R\frakEBord_d$,
where $ι:\frakFEmb_d→\frakStruct_d$ is the $\smsset$-enriched Yoneda embedding
and $R$ is the fibrant replacement functor for the injective model structure on $\smsset$-enriched functors $\frakFEmb_d→\fraksmcat_{∞,d}$.

\item Morphisms are natural transformations $\BBord'_d→\BBord''_d$ such that after applying $ι^*$ the resulting triangle with $R\frakEBord_d$ commutes.
\end{itemize}
Then the map from the resulting relative category to the terminal relative category is a Dwyer–Kan equivalence.
That is to say, the axioms of \cref{multiple.isotopy.axiom}
determine the functor $\BBord'_d$ uniquely up to a weakly contractible choice.
Furthermore, the resulting smooth symmetric monoidal $(\infty,d)$-category of bordisms
$$\BBord'_d:\frakStruct_{d}→\fraksmcat_{∞,d},$$
where $\frakStruct_{d}$ has the Čech local model structure (\cref{geometric.structure.isotopy}), is homotopy cocontinuous.

The analogous statements hold for the axioms of \cref{axioms}, replacing $\BBord_d$ with $\Bord_d$, $\frakEBord_d$ with $\EBord_d$, $\frakFEmb_d$ with $\FEmb_d$,
$\frakStruct_d$ with $\Struct_d$, $\fraksmcat_{∞,d}$ with $\smcat_{∞,d}$, and $\smsset$-enrichments with $\sset$-enrichments.
\end{theorem}

\begin{proof}
Given that $\frakStruct_{d,\inj}$ is the homotopy cocompletion of $\frakFEmb_d$,
the restriction functor~$ι^*$ along the $\smsset$-enriched Yoneda embedding $$ι:\frakFEmb_d→\frakStruct_{d,\inj}$$ induces a Dwyer–Kan equivalence
from the relative category of $\smsset$-enriched homotopy cocontinuous functors $$\frakStruct_{d,\inj}→\fraksmcat_{∞,d}$$
to the relative category of $\smsset$-enriched functors $$\frakFEmb_d→\fraksmcat_{∞,d}.$$
(Some care must be exercised when dealing with size issues caused by the fact that the involved categories are large.)
This establishes the contractibility statement.

The full subcategory of $\smsset$-enriched functors $\frakStruct_{d,\inj}→\fraksmcat_{∞,d}$
that are also homotopy cocontinuous as $\smsset$-enriched functors $\frakStruct_d→\fraksmcat_{∞,d}$
is mapped by $ι^*$ to the full subcategory
of $\smsset$-enriched functors $\frakFEmb_d→\fraksmcat_{∞,d}$
that satisfy Čech homotopy codescent for covering families in $\frakFEmb_d$.
By axiom \frakax3, the category $ι^*\BBord'_d$ is weakly equivalent to the isotopy embedded bordism category $\frakEBord_d$ (\cref{embeddedbordcat.enriched}).
The category $\frakEBord_d$ satisfies Čech homotopy codescent by \cref{codescent.embedded}.
Therefore, $\BBord'_d$ satisfies Čech homotopy codescent and is homotopy cocontinuous on $\frakStruct_d$.
The statement for the axioms of \cref{axioms} is proved in the same way.
\end{proof}

\cref{axioms.contractible.codescent} proves \cref{axtheorem} and \cref{axtheorem2} for any bordism category that satisfies the axioms \ax1--\ax3 and \frakax1--\frakax3, respectively.
In the case the bordism categories defined in \cref{bordcts}, we have the following stronger statement. 

\begin{theorem}
\label{a1a2a3}
Fix $d≥0$.
Recall the definition of $\Bord_{d,\uple}$ and $\Bord_d$ (\cref{bordstr}),
as well as $\BBord_{d,\uple}$ and $\BBord_d$ (\cref{enrichedbordstr}). The smooth symmetric monoidal $(\infty,d)$-category of bordisms with isotopies $\BBord_d$
satisfies the axioms of \cref{multiple.isotopy.axiom},
and, therefore, satisfies Čech homotopy codescent.
Likewise for the smooth symmetric monoidal $(\infty,d)$-category of bordisms $\Bord_d$. In fact, we have the following stronger statement.
The functors
$$\Bord_{d,\uple}:\Struct_d\to \smcatuple_{\infty,d}, \quad \gs\mapsto \Bord_{d,\uple}^\gs$$
and
$$\Bord_{d}:\Struct_d\to \smcat_{\infty,d}, \quad \gs\mapsto \Bord_{d}^\gs$$
are $\sset$-enriched left Quillen functors that preserve all weak equivalences. 
In particular, they are homotopy cocontinuous.
Similarly, the functors
$$\BBord_{d,\uple}:\frakStruct_d\to \fraksmcatuple_{\infty,d}, \quad \gs\mapsto \BBord_{d,\uple}^\gs$$
and
$$\BBord_{d}:\frakStruct_d\to \fraksmcat_{\infty,d}, \quad \gs\mapsto \BBord_{d}^\gs$$
are $\smsset$-enriched left Quillen functors that preserve all weak equivalences.
In particular, they are homotopy cocontinuous.
\end{theorem}

\begin{proof}
The axioms hold by \cref{proofof106} and \cref{proofof104}, respectively.
Čech homotopy codescent now follows from \cref{axioms.contractible.codescent}.
For the stronger statement, observe that by the universal property of left Bousfield localizations (\cref{left.Bousfield.localization}),
it suffices to show that the functors are left Quillen functors prior to left Bousfield localizations (i.e., using the injective model structure),
which is done in \cref{a1a2},
and then show that the Čech codescent morphisms are weak equivalences,
which follows from \cref{emaplociso} combined with \cref{codescent.embedded}. 
\end{proof}

Throughout the remainder of the paper, we will use the following notation.

\begin{notation}
\label{intnotation}
Fix $d≥0$, an object $(p:W\to U)\in \frakFEmb_d$, and a covering family $\{p_a:W_a\to U_a\}_{a\in A}$ of~$p$.
For a multi-index $\alpha:[n]\to A$ (where $[n]∈Δ$), we let $p_{\alpha}$ denote the restriction of $p$ to $W_α→U_α$,
where $$W_{\alpha} ≔ W_{\alpha(0)}\cap W_{\alpha(1)}\cap W_{\alpha(2)}\cap \cdots\cap W_{\alpha(n)},
\qquad U_{\alpha} ≔ U_{\alpha(0)}\cap U_{\alpha(1)}\cap U_{\alpha(2)}\cap \cdots\cap U_{\alpha(n)}.$$
\end{notation}

We now turn to the proofs of the supporting propositions needed to establish \cref{axioms.contractible.codescent} and \cref{a1a2a3}. We begin with the main proposition, which we break up into smaller claims that are proved subsequently. 

\begin{proposition}
\label{codescent.embedded}
Fix $d≥0$, an object $(p:W\to U)\in \frakFEmb_d$, and a covering family $\{p_a:W_a\to U_a\}_{a\in A}$ of~$p$.
Recall $\EBord_d^p$ (\cref{embeddedbordcat.unenriched}) and $\frakEBord_d^p$ (\cref{embeddedbordcat.enriched}).   
Then the canonical homotopy codescent map
$$\hocolim_{[n]\in \Delta^\op} \coprod_{\alpha:[n]\to A}\frakEBord_d^{p_\alpha:W_{\alpha}\to U_{\alpha}}\to \frakEBord_d^{p:W\to U},$$
induced by the inclusions $p_α→p$, is a weak equivalence in $\fraksmcat_{\infty,d}$.

The homotopy codescent map is a weak equivalence also for $\EBord_d$, using $\FEmb_d$ instead of $\frakFEmb_d$ and $\smcat_{∞,d}$ instead of $\fraksmcat_{∞,d}$.
\end{proposition}

\begin{proof}
By \cref{hostr}, the canonical map from the homotopy colimit to the strict colimit
$$\hocolim_{[n]\in \Delta^\op} \coprod_{\alpha:[n]\to A}\frakEBord_d^{p_{\alpha}}\to \colim_{[n]\in \Delta^\op} \coprod_{\alpha:[n]\to A}\frakEBord_d^{p_{\alpha}}$$
is a weak equivalence in $\fraksmcat_{\infty,d}$.
Thus, it remains to show that the map
$$\colim_{[n]\in \Delta^\op} \coprod_{\alpha:[n]\to A}\frakEBord_d^{p_{\alpha}} → \frakEBord_d^p\eqlabel{rem.map.ebord}$$
is a stalkwise weak equivalence (\cref{stalks}), which implies that it is a weak equivalence in $\fraksmcat_{\infty,d}$.
We factor the stalk of \cref{rem.map.ebord} at every point $\mathfrak{l}$ of $\cart$ and $\mathfrak{p}$ of $\stcart$
as a composition of monomorphisms introduced in \cref{subcov}:
$$\colim_{[n]\in \Delta^\op}\coprod_{\alpha:[n]\to A}\mathfrak{l}^*\mathfrak{p}^*\frakEBord_d^{p_{\alpha}}=\frakB_{-1}\overset{{\bf 1}}{\into} \frakB_0\overset{{\bf 2}}{\into} \frakB_d\overset{{\bf 3}}{\into} \frakB=\mathfrak{l}^*\mathfrak{p}^*\frakEBord_d^p.$$
Here we moved the stalk functors, which are computed as filtered colimits, inside the colimit over $Δ^\op$ on the left side.
The map ${\bf 1}$ is a stalkwise weak equivalence in $\smcat_{\infty,d}$ by \cref{gammawkeq}.  
The map ${\bf 2}$ is a stalkwise weak equivalence in $\smcat_{\infty,d}$ by \cref{induction.i}. 
The map ${\bf 3}$ is an equality by \cref{subcov}.
Hence the composition is a stalkwise equivalence, which establishes the case of $\frakEBord_d$.

The case of $\EBord_d$ now follows by setting $\mathfrak{l}=(\RR^0,\id)$
and invoking \cref{0stalkbord} to get $\mathfrak{l}^*\mathfrak{p}^*\frakEBord_d≃\mathfrak{p}^*\EBord_d$,
which completes the proof.
\end{proof}

\subsection{Reduction of homotopy colimits to strict colimits}

\begin{proposition}
\label{hostr}
Fix $d≥0$, an object $(p:W\to U)\in \frakFEmb_d$, and a covering family $\{p_a:W_a\to U_a\}_{a\in A}$ of~$p$.
Recall $\frakEBord_d^p$ (\cref{embeddedbordcat.enriched}) and the notation $p_α$ (\cref{intnotation}) for a multi-index $\alpha:[n]\to A$, where $[n]∈Δ$.
Then the canonical map
$$\hocolim_{[n]\in \Delta^\op} \coprod_{\alpha:[n]\to A}\frakEBord_d^{p_{\alpha}}\lto3{} \colim_{[n]\in \Delta^\op} \coprod_{\alpha:[n]\to A}\frakEBord_d^{p_{\alpha}}$$
is a weak equivalence.
\end{proposition}

\begin{proof}
\def\cP{{\cal P}}
\def\cQ{{\cal Q}}
First, we combine the double (homotopy) colimit into a single (homotopy) colimit using the Grothendieck construction:
$$\hocolim_{(α:[n]→A)∈(\Delta/A)^\op}\frakEBord_d^{p_{\alpha}} \lto3{} \colim_{(α:[n]→A)∈(\Delta/A)^\op}\frakEBord_d^{p_{\alpha}},$$
where $Δ/A$ is the following category.
Objects are maps of sets $α:[m]→A$.
Morphisms $(α:[m]→A)→(α':[m']→A)$ are maps of simplices $σ:[m]→[m']$ such that $α=α'σ$.

Pick a total order on~$A$ and consider the full subcategory~$Δ'/A$ of $Δ/A$ on objects given by injective order-preserving maps $λ:[n]→A$.
Equivalently, objects of $Δ'/A$ are finite nonempty subsets of~$A$ and morphisms are inclusions.
Consider the following functor $\cP:Δ/A→Δ'/A$.
Given an object $α:[m]→A$ in $Δ/A$, we set $\cP(α)$ to the image of~$α$.
Given a morphism $σ:(α:[m]→A)→(α':[m']→A)$, we set $\cP(σ)$ to the corresponding inclusion of images.

The functor~$\cP$ is a homotopy initial functor by Quillen's Theorem~A.
To this end, pick an arbitrary object $Λ∈Δ'/A$ given by a finite nonempty subset $Λ⊂A$
and consider the comma category $\cP/Λ$, which can be described concretely as the full subcategory of $Δ/A$ on objects $α:[m]→A$ such that the image of~$α$ is a subset of~$Λ$.
The category $\cP/Λ$ is the category of simplices of the nerve of the codiscrete groupoid on~$Λ$, i.e.,
the set of objects is~$Λ$ and there is exactly one morphism between any pair of objects.
Therefore, the nerve of $\cP/Λ$ is the barycentric subdivision of the nerve of the codiscrete groupoid on~$Λ$, and the latter is contractible.

Since the Čech diagram $$C:(Δ/A)^\op→\smcat_{∞,d}, \qquad α↦\frakEBord_d^{p_α}$$
factors as the functor $\cP^\op:(Δ/A)^\op→(Δ'/A)^\op$ followed by the restriction~$C'$ of~$C$ to~$(Δ'/A)^\op$
and the functor $\cP^\op$ is homotopy final,
it remains to show that the map $\hocolim C'→\colim C'$ is a weak equivalence.

Next, we reduce to the case when $A$ is finite.
Consider the following category $Δ''/A$.
Objects are pairs of finite nonempty sets $(Λ,J)$, where $Λ⊂J⊂A$.
Morphisms $(Λ,J)→(Λ',J')$ are inclusions $Λ⊂Λ'$, $J⊃J'$.
There is a functor $\cQ:Δ''/A→Δ'/A$ that sends $(Λ,J)↦Λ$, which defines it uniquely also on morphisms.
Given $Σ∈Δ'/A$, the comma category $\cQ/Σ$ has objects $(Λ,J)$ such that $Λ⊂Σ$ and morphisms induced from $Δ''/A$.
The projection functor $$π:\cQ/Σ→Δ''/Σ,\qquad (Λ,J)↦Λ$$ is a Grothendieck fibration,
whose fiber over~$Λ$ is the poset of all~$J$ such that $J⊃Λ$.
This poset has a terminal object $J=Λ$, so the nerve of every fiber is contractible.
This implies that the nerve of~$π$ is a weak equivalence,
i.e., the nerve of $\cQ/Σ$ is weakly equivalent to the nerve of $Δ''/Σ$.
Since $Σ$ is finite, the latter poset has a terminal object~$Σ$.
Hence, the nerve of $\cQ/Σ$ is contractible and the functor $\cQ$ is homotopy initial.

Thus, the restricted Čech diagram $C':(Δ'/A)^\op→\smcat_{∞,d}$ can be precomposed with $\cQ^\op$ without changing its (homotopy) colimit.
The projection functor $$(Δ'/A)^\op→Δ'/A, \qquad (Λ,J)↦J$$ is a Grothendieck fibration,
whose fiber over some $J∈Δ'/A$ is the category $(Δ'/J)^\op$.
Therefore, the (homotopy) colimits of $C'$ over $(Δ'/A)^\op$ can be replaced by the double (homotopy) colimits:
$$\hocolim_{J∈Δ'/A} \hocolim_{Λ∈(Δ'/J)^\op}\frakEBord_d^{p_Λ}\lto3{} \colim_{J∈Δ'/A} \colim_{Λ⊂(Δ'/J)^\op}\frakEBord_d^{p_Λ}.$$

The outer (homotopy) colimit is indexed by the filtered category of finite subsets of~$A$.
In $\smcat_{∞,d}$, like in any localized category of simplicial presheaves,
the map from a filtered homotopy colimit to the corresponding colimit is always a weak equivalence.
Therefore, it remains to show that for any finite subset $J⊂A$ the inner map
$$\hocolim_{Λ∈(Δ'/J)^\op}\frakEBord_d^{p_Λ}\lto3{} \colim_{Λ∈(Δ'/J)^\op}\frakEBord_d^{p_Λ}$$
is a weak equivalence, i.e., we can assume $A$ to be finite.
To this end, it suffices to show that the diagram $$(Δ'/J)^\op→\smcat_{∞,d}, \qquad Λ↦\frakEBord_d^{p_Λ}$$ is projectively cofibrant.
Since $J$ is finite, the indexing category $(Δ'/J)^\op$
can be turned into a Reedy category using $Λ↦|J|-|Λ|$ as the degree function,
declaring all arrows to be positive and only identity arrows to be negative.
For such a category, the projective model structure coincides with the Reedy model structure.
Thus, we need to show that the latching map $$\colim_{K⊃Λ,K≠Λ}\frakEBord_d^{p_K}→\frakEBord_d^{p_Λ}$$
is a monomorphism for any $Λ∈(Δ'/J)^\op$, i.e., a finite nonempty subset $Λ⊂J$.

At this point we need to invoke a property that is specific to the diagram in question:
we have $$\frakEBord_d^{p_K}∩\frakEBord_d^{p_{K'}}=\frakEBord_d^{p_{K∪K'}}$$
for any finite nonempty subsets $K,K'⊂J$.
Indeed, a bordism embedded into $p_K$ and $p_{K'}$ is necessarily embedded into $p_K∩p_{K'}=p_{K∪K'}$.

Since in the category $(Δ'/J)^\op$ we have morphisms $K∪K'→K$ and $K∪K'→K'$,
any element of $\frakEBord_d^p$ that is present in several stages of the latching diagram is already present at a stage that precedes them all,
and therefore all of its copies will be identified in the latching object.
Thus, the colimit that defines the latching object can be computed simply as the union of subobjects.
This immediately implies that the latching map is a monomorphism, which completes the proof.
\end{proof}

\subsection{Filtration on the bordism category}
\label{filtrationsec}

Recall from \cref{bordstalkisot} the $(\mathfrak{l},\mathfrak{p})$-stalk of the embedded bordism category with isotopies:
$$\mathfrak{l}^*\mathfrak{p}^*\frakEBord_d^p\in \PSh(\Gamma\times \Delta^{\times d}).$$
This simplicial presheaf sends $(\langle \ell\rangle,{\bf m})$ to the set of germs of embedded bordisms $[f,N,C,P]$,
equipped with a germ of $L$-families of cut ${\bf m}$-grids and a map $P:N\to \langle \ell\rangle$ (see \cref{bordstalkisot} for a more detailed description). 

\begin{remark}
\label{constsimpd}
From this point onward, we will work only with presheaves of sets on $\Gamma\times \Delta^{\times d}$.
Occasionally we will need to talk about weak equivalences, cofibrations, acyclic cofibrations between such objects.
In such instances, we will always implicitly promote the presheaf of sets to a presheaf of simplicial sets, by taking the constant simplicial set.
\end{remark}

\begin{notation}
\label{filtration.notation}
We fix the following data:
\begin{itemize}
\item A natural number $d\in \NN$, indexing the dimension.
\item A point $\mathfrak{l}$ of $\cart$, indexing a stalk.
\item A point $\mathfrak{p}$ of $\stcart$, indexing a stalk.
\item An object $p:W\to U$ of $\frakFEmb_d$.
\item A covering $\mathcal{W}=\{p_a:W_a\to U_{a}\}_{a\in A}$ of $p:W\to U$.
\end{itemize}
Since all statements in this subsection will use the above data, we do not need to carry around the notation $\mathfrak{l}$, $\mathfrak{p}$, $p:W\to U$, $p_a:W_a\to U_a$.
In this section, we will use the much more compact notation
$$\frakB≔\mathfrak{l}^*\mathfrak{p}^*\frakEBord_d^p, \quad \frakB_{-1}≔\colim_{[n]\in \Delta}\coprod_{\alpha:[n]\to A}\mathfrak{l}^*\mathfrak{p}^*\frakEBord_d^{p_\alpha},$$
where $p_{\alpha}$ is defined in \cref{intnotation}.
\end{notation}

Next, we will define the filtration on $\frakB$ that we used in \cref{codescent.embedded}. The meaning of the filtration can be described roughly as follows.
\begin{itemize}
\item Passing to stalks on both $\stcart$ and $\cart$ give tiny perturbations of chopped embedded bordisms, where the chopping is provided by the cut ${\bf m}$-grid. Perturbations of both the grid (encoded by the site $\cart$) and the embedding (provided by the site $\stcart$) are allowed. 
\item Embedded bordisms in filtration degree $-1$ have an embedding into $p$ that factors through $p_a$, for some $a\in A$ indexing an element of the cover  (\cref{filtration.notation}). 
\item Embedded bordisms in filtration degree 0 are disjoint unions of bordisms in degree $-1$. 
\item Embedded bordisms in filtration degree $i>0$ are bordisms that can be decomposed as a composition of two bordisms in the $i$th direction such that the first is in filtration degree $0$ and the second is an $i$-dimensional bordism. 
\end{itemize}

\begin{definition}
\label{icore}
Assume \cref{filtration.notation}. Let $[f,N,C,P]\in \frakB$ and let $i\geq -1$. If $i\geq 0$, we say that $[f,N,C,P]$ is an \emph{embedded $i$-bordism} if it simplicially degenerate in all simplicial directions $k>i$. An embedded $(-1)$-bordism is the empty bordism. Recall the definition of the core in a single fiber in  \cref{bordstalkisot}. Working in the same fiber, we define the \emph{$i$-core} $\core_i[f,N,C,P]_0$ as follows. 
\begin{itemize}
\item  Consider the collection of open subsets $V\subset \stu(N)_{\stu(f)(0,\rho^{-1}(0))}$, where $V=\stu(U)_{\stu(f)(0,\rho^{-1}(0))}$ and $U\subset N$ is an open subobject such that the embedded bordism $[f,U,C\vert_{U},P\vert_{U}]$ is an $i$-bordism. This set is ordered by inclusion and has a maximal element $M$. We define
$$\core_i[f,N,C,P]_0=\stu(N)_{\stu(f)(0,\rho^{-1}(0))}\setminus M.$$
\end{itemize}
We observe that $\core_{-1}[f,N,C,P]_0=[f,N,C,P]_0$ and $\core_d[f,N,C,P]_0=\emptyset$.
\end{definition}

\begin{remark}
If $i<i'$ then $\core_{i}[f,N,C,P]_0\supset \core_{i'}[f,N,C,P]_0$, since every embedded $i$-bordism is an embedded  $i'$-bordism. 
\end{remark}

\begin{definition}
\label{subcov} 
Assume \cref{filtration.notation}.
We define a filtration 
$$\frakB_{-1} \into \frakB_0\into \frakB_1\into \cdots \into \frakB_d=\frakB$$
of presheaves of sets on the category $Γ⨯Δ^{⨯d}$
inductively as follows.
Objectwise on $\Gamma\times \Delta^{\times d}$, we define $\frakB_i(\langle \ell\rangle,{\bf m})$ for $i\geq 0$ as the subset of $\frakB(\langle \ell \rangle,{\bf m})$ whose elements $[f,N,C,P]$ (\cref{bordstalkisot}) satisfy the following. 
\begin{itemize}
\item Every connected component of $\core_i[f,N,C,P]_0$ (\cref{icore}) is contained in $W_a$ for some $W_a\in {\cal W}$. 


\end{itemize}
Observe that $\frakB_d=\frakB$ by definition. 
\end{definition}

\begin{figure}[ht]
\begin{tikzpicture}[scale=.4]
	\begin{pgfonlayer}{nodelayer}
		\node [style=none] (4) at (0.5, -2.25) {};
		\node [style=none] (5) at (2.5, -2.25) {};
		\node [style=none] (6) at (-1.5, -0.5) {};
		\node [style=none] (7) at (0.5, -0.5) {};
		\node [style=none] (8) at (-4, -2.75) {};
		\node [style=none] (9) at (-2, -2.75) {};
		\node [style=none] (10) at (-6, -0.75) {};
		\node [style=none] (11) at (-4, -0.75) {};
		\node [style=none] (12) at (3.5, -4.25) {};
		\node [style=none] (13) at (-7.75, -4.25) {};
		\node [style=none] (14) at (-5.25, 1.5) {};
		\node [style=none] (15) at (5.75, 1.5) {};
		\node [style=none] (16) at (-0.75, -4.25) {};
		\node [style=none] (17) at (2.25, 1.5) {};
		\node [style=none] (18) at (-14, -4.25) {};
		\node [style=none] (19) at (-14, -1.25) {};
		\node [style=none] (20) at (-12.25, -3) {};
		\node [style=none] (21) at (-10.75, -4.25) {};
		\node [style=none] (23) at (-12.25, -4.75) {1};
		\node [style=none] (24) at (-13, -2.75) {};
		\node [style=none] (25) at (-12.25, -4.75) {};
		\node [style=none] (26) at (-13, -2.75) {2};
		\node [style=none] (27) at (-14.75, -2.25) {3};
		\node [style=none] (29) at (-5.5, -4.25) {};
		\node [style=none] (30) at (-4.5, 1.5) {};
		\node [style=none] (31) at (3.75, -3) {};
		\node [style=none] (32) at (-7.5, -2.25) {};
		\node [style=none] (33) at (4.5, -1.5) {};
		\node [style=none] (34) at (0.75, -1) {};
		\node [style=none] (35) at (1.25, -4) {};
		\node [style=none] (36) at (3, -0.5) {};
		\node [style=none] (37) at (5.5, -4) {};
		\node [style=none] (38) at (5.5, -4) {$W_a$};
	\end{pgfonlayer}
	\begin{pgfonlayer}{edgelayer}
		\draw [style=new edge style 2, in=75, out=105, looseness=1.75] (5.center) to (4.center);
		\draw [style=new edge style 3, bend left=90] (5.center) to (4.center);
		\draw [style=new edge style 2, in=75, out=105, looseness=1.75] (7.center) to (6.center);
		\draw [style=new edge style 3, bend left=90] (7.center) to (6.center);
		\draw [style=new edge style 2, in=75, out=105, looseness=1.75] (9.center) to (8.center);
		\draw [style=new edge style 3, bend left=90] (9.center) to (8.center);
		\draw [style=new edge style 2, in=75, out=105, looseness=1.75] (11.center) to (10.center);
		\draw [style=new edge style 3, bend left=90] (11.center) to (10.center);
		\draw [style=new edge style 2, in=0, out=-180] (12.center) to (13.center);
		\draw [style=new edge style 2, in=-90, out=90] (12.center) to (15.center);
		\draw [style=new edge style 2] (15.center) to (14.center);
		\draw [style=new edge style 2, in=105, out=-105] (14.center) to (13.center);
		\draw [style=new edge style 0] (18.center) to (19.center);
		\draw [style=new edge style 0] (18.center) to (20.center);
		\draw [style=new edge style 0] (18.center) to (21.center);
		\draw [style=new edge style 2, in=-105, out=90] (16.center) to (17.center);
		\draw [style=new edge style 2, in=-45, out=75, looseness=2.00] (29.center) to (30.center);
		\draw [style=new edge style 2, in=15, out=-180, looseness=2.50] (31.center) to (32.center);
		\draw [style=new edge style 3, in=0, out=0, looseness=3.00] (35.center) to (36.center);
		\draw [style=new edge style 3, in=-180, out=-180, looseness=2.50] (35.center) to (36.center);
	\end{pgfonlayer}
\end{tikzpicture}
\caption{For $d=3$, an image of the core of a bordism in $\frakB_{2}$. The axes on the left indicate the composition directions. The dashed disk is an element $W_a$ in the cover ${\cal W}$. The $2$-core is the 3-dimensional region given by the disjoint union of the bumps (including the dashed curves on the boundary). Each bump is contained in some open set in ${\cal W}$. Removing the $2$-core yields a 2-dimensional bordism.}
\end{figure}

\subsection{$\Gamma$-direction}

We are now ready to prove that each inclusion in the above filtration is a weak equivalence.
We begin with the $\Gamma$-direction. 

\begin{proposition}[$\Gamma$-direction]
\label{gammawkeq}
The map
$$\frakB_{-1}→\frakB_0$$
(\cref{subcov})
is a weak equivalence in the model category $\cat_{\infty,d}$ (\cref{globular.multiple.model.structure}).
\end{proposition}

\begin{proof}
By \cref{multiple.single},
it suffices to show that after evaluating at some ${\bf m}∈Δ^{⨯d}$
we get a weak equivalence in the model category $\sPSh(\Gamma)_{\local}$ (\cref{gammasp}).

To this end, we apply \cref{EZ.criterion}.
We use the description of stalks in \cref{bordstalkisot},
including the special form of~$N$ in Condition~(3) that guarantees a convenient description of~$π_0$:
elements of $\frakB_0(⟨\ell⟩)$ are equivalence classes $[f,N,C,P]$ as described there,
with the map $P:N→⟨\ell⟩$ factoring as
$$N \lto7{Q} \{*\}⊔π_0[f,N,C,P] \lto7{γ} ⟨\ell⟩.$$

An equivalence class $[f,N,C,P]$ is indecomposable (\cref{indecomposable.element}) if and only if the map~$γ$ is an isomorphism.
Inert maps throw some components in the trash bin and do nothing else.
Therefore, inert maps preserve indecomposable elements.

The existence part of the modified Eilenberg–Zilber property (\cref{EZ.property}) for $\frakB_0$
is established by the relation $$[f,N,C,P]=\frakB_0(γ)[f,N,C,Q].$$
supplied by the factorization above.
This decomposition is unique up to a unique permutation since the decomposition into connected components is unique up to a unique permutation.
Finally, the symmetric group acts freely on indecomposable elements.
This establishes the modified Eilenberg–Zilber property for $\frakB_0$.
The modified Eilenberg–Zilber property for $\frakB_{-1}$ now follows because $[f,N,C,P]∈\frakB_{-1}$ implies $[f,N,C,Q]∈\frakB_{-1}$.

Finally, indecomposable elements of $\frakB_0$ in degree~$⟨1⟩$ belong to $\frakB_{-1}$.
Indeed, by \cref{icore}, we have $\core_0[f,N,C,P]_0=\core[f,N,C,P]_0∖M$,
where $M$ is the core of the embedded 0-bordism given by all embedded point bordisms, i.e., embedded bordisms degenerate in all simplicial directions, in $[f,N,C,P]$.
Since every embedded point bordism is automatically subordinate to the covering family $\mathcal{W}$,
embedded bordisms in $\frakB_0$ can be equivalently described as embedded bordisms all of whose connected components are subordinate to the covering family~$\mathcal{W}$.
Thus, indecomposable embedded bordisms in $\frakB_0$ are connected embedded bordisms in $\frakB$ that are subordinate to~$\mathcal{W}$,
i.e., precisely the embedded bordisms in $\frakB_{-1}$.
\end{proof}

\subsection{$\Delta$-direction}

We now turn to the most technical part of the proof.
Supporting propositions will be proved subsequently. 

\begin{proposition}
\label{induction.i}
For all~$i$ such that $0<i≤d$,
the monomorphism $\rho_i:\frakB_{i-1}\into \frakB_i$ in \cref{subcov} is a weak equivalence in $\cat_{\infty,d}$.
\end{proposition}

\begin{proof}
By \cref{multiple.single}, it suffices to show that the monomorphism of presheaves of sets
$$\rho_i(\langle \ell\rangle,{\bf m}):\frakB_{i-1}(\langle \ell\rangle,{\bf m})\into \frakB_{i}(\langle \ell\rangle,{\bf m}),$$ 
obtained by evaluating $\rho_i$ on an arbitrary multisimplex $\langle \ell\rangle\in \Gamma$ and ${\bf m}\in \Delta^{\times d-1}$, via restriction along the functor  
$$\Delta\to \Gamma\times \Delta^{\times d}, \quad [\omega]\mapsto (\langle \ell\rangle,(m_1,\ldots,m_{i-1},\omega, m_i,\ldots,m_{d-1})),$$
is an equivalence in the Rezk model structure (\cref{rezkdelta}) on simplicial objects in simplicial sets (recall \cref{constsimpd}, which promotes sets to discrete simplicial sets).
We evaluate on the objects $\langle \ell\rangle$ and ${\bf m}$ and omit them in the notation.
By \cref{necklaces.prop}, it suffices to show that $\Cnec(ρ_i)$ is a Dwyer–Kan equivalence of simplicial categories.
We have $\frakB_i([0])=\frakB_{i-1}([0])$ by definition of~$\frakB_i$, so $\Cnec(ρ_i)$ is essentially surjective.
It remains to observe that by \cref{weak.eq.necklaces}, for each pair of vertices $x$ and $y$ in $\frakB_i([0])$, the induced map
$$\Cnec(\rho_i):\Cnec(\frakB_{i-1})(x,y)\to \Cnec(\frakB_i)(x,y)$$
is a weak equivalence of simplicial sets.
\end{proof}

\begin{notation}
\label{simplex.notation}
In addition to \cref{filtration.notation}, we also fix the following data:
\begin{itemize}
\item an object $\langle \ell\rangle\in \Gamma$;
\item an integer $i$ such that $1\leq i \leq d$;
\item a multisimplex ${\bf m}\in \Delta^{\{1,\hdots, \hat \imath,\hdots d\}}=\Delta^{\times d-1}$;
\item elements $x,y\in \frakB_i(\langle \ell\rangle,{\bf m})([0])=\frakB_{i-1}(\langle \ell\rangle,{\bf m})([0])$.
\end{itemize}
Since all statements in this subsection will use the above data, we do not need to carry around the notation $\langle \ell\rangle$ and~${\bf m}$.
In this section, we will use the much more compact notation
$$\frakB^i_i≔\frakB_i(\langle \ell\rangle,{\bf m}) \qquad \frakB^i_{i-1}≔\frakB_{i-1}(\langle \ell\rangle,{\bf m})$$
with $\frakB_i(\langle \ell\rangle,{\bf m})$ and $\frakB_{i-1}(\langle \ell\rangle,{\bf m})$ as in \cref{subcov}. Hence, both $\frakB^i_i$ and $\frakB^i_{i-1}$ are presheaves of sets on $\Delta$, where $\Delta$ is interpreted as the $i$th factor in $\Delta^{\times d}$. 
We emphasize that $i$ is fixed throughout this section and we are studying compositions of bordisms in the $i$th direction for both $\frakB^i_i$ and $\frakB^i_{i-1}$, as indicated by the superscript. 
\end{notation}

Next, we will need an explicit description of the Dugger–Spivak necklace categories. We begin by defining a category which we will show is equivalent to the necklace categories.

\begin{definition}
\label{etalemapsnow}
Assume \cref{simplex.notation}.
We define a category $\frakB_{i-1}^{i,x,y}$, respectively $\frakB_i^{i,x,y}$, whose objects are called \emph{etale necklaces}.
The set of objects is the quotient of the following set of tuples,
taken for all possible $V=B^{\mathfrak{p}}_{\delta}\subset T\in \stcart$ and $L=B^{\mathfrak{l}}_{\delta}\subset\RR^{n}$ (\cref{stalks}):
\begin{enumerate}
\item a simplex $[m_i]\in \Delta$. We let ${\bf n}=([m_1],\hdots,[m_i],\hdots,[m_d])$, where ${\bf m}=([m_1],\hdots,\widehat{[m_i]},\hdots,[m_d])\in \Delta^{\times (d-1)}$ is the fixed multisimplex ${\bf m}$ in \cref{simplex.notation}; 
\item a morphism $f: V\times L\to U$ in $\stcart$;
\item an etale map $h:N\to  W\times_{U}V\times L$ such that the submersion $\pi\circ h:N\to V\times L$ is a (locally trivial) fiber bundle, where $\pi$ is the projection;
\item a compact (globular) cut ${\bf n}$-grid $C$ on $\pi\circ h:N\to V\times L$;
\item a morphism $P:N\to \langle \ell\rangle$ in $\stman$;
\item a subset $\Upsilon\subset [m_i]$ containing $0$ and $m_i$ that satisfies the following property. 
\end{enumerate}
For all pairs of vertices $(\upsilon,\upsilon')$ in $Υ$, we have the following.
\begin{enumerate}
\item[(6a)] If $\upsilon=\upsilon'$ or $\upsilon<\upsilon'$ are consecutive in $\Upsilon$, then $h$ restricts to an open embedding $h_{\langle \upsilon,\upsilon'\rangle}:N_{\langle \upsilon,\upsilon'\rangle}\to W\times_{U}V\times L$ for some open subobject $N_{\langle \upsilon,\upsilon'\rangle}\subset N$ whose reduction contains the core of $C_{\langle \upsilon,\upsilon'\rangle}$ (see \cref{coreandbead} with $S=\{i\}$).
This gives rise to an embedded bordism 
$$A_{\upsilon,\upsilon'}=[f,N_{\langle \upsilon,\upsilon'\rangle},h_{\langle \upsilon,\upsilon'\rangle}^*C_{\langle \upsilon,\upsilon'\rangle},h_{\langle \upsilon,\upsilon'\rangle}^*P].$$
The definition of the embedded bordism $A_{υ,υ'}$ continues to make sense for any $υ,υ'∈[m_i]$ such that $υ≤υ'$ and there is no $υ''∈Υ$ such that $υ<υ''<υ'$.
\item[(6b)] If $\upsilon<\upsilon'$ are consecutive, the embedded bordism $A_{\upsilon,\upsilon'}$ is in $\frakB^i_{i-1}$, respectively $\frakB^i_i$.
\item[(6c)] $A_{0,0}=x$ and $A_{m_i,m_i}=y$.
\end{enumerate}
The equivalence relation on the set of tuples is given by: two tuples $(m_i,f,h, C,P,\Upsilon)$ and $(m_i',f',h',C',P',\Upsilon')$ are equivalent if there is $\epsilon<\min\{\delta,\delta'\}$ such that after restricting to $B_{\epsilon}^{\mathfrak{p}}\times B_{\epsilon}^{\mathfrak{l}}$, we have
\begin{enumerate}
\item[(e1)] $m_i=m_i'$;
\item[(e2)] $f=f'$; 
\item[(e3)] There are open embeddings $N''\into N$ and $N''\into N'$ such that the square
$$\xymatrix{
N''\ar[r]\ar[d] & N\ar[d]^-{h}
\\
N'\ar[r]^-{h'} & W\times_{U}V\times L
}$$
commutes and there is a compact (globular) cut ${\bf n}$-grid $C''$ on $\pi\circ h:N''\to B_{\varepsilon}^{\mathfrak{p}}\times B^{\mathfrak{l}}_{\varepsilon}$ such that: 
\begin{itemize}
\item the image of $\stu(N'')$ under the two embeddings $\stu(N'')\into \stu(N)$ and $\stu(N'')\into \stu(N')$ contains the core of $C$ and $C'$, respectively;
\item $C''=C=C'$ and $P=P'$, after restricting to $\stu(N'')$;
\end{itemize}
\item[(e4)] $\Upsilon=\Upsilon'$.
\end{enumerate} 
A morphism $[m_i,f,h,C,P,\Upsilon]\to [m_i',f',h',C',P',\Upsilon']$ is given by
\begin{enumerate}
\item[(m)] a morphism of simplices $\sigma:[m_i]\to [m_i']$ preserving the initial and terminal vertices such that $\Upsilon'\subset \sigma(\Upsilon)$ and $[m_i,f,h,C,P,\Upsilon]=[m_i,f',h',\sigma^*C',P',\Upsilon]$, where $\sigma^*C'$ is defined using the simplicial structure map for cut grids in \cref{cutgrid}, applied in the $i$th direction.
\end{enumerate}
\end{definition}

\begin{remark}
\label{taumap}
The category $\frakB_{i-1}^{i,x,y}$ is a subcategory of $\frakB_i^{i,x,y}$, since Condition (6b) for $\frakB^i_{i-1}$ is a stronger condition than (6b) for $\frakB^i_i$. We let
$$\tau_{i}^{x,y}:\frakB_{i-1}^{i,x,y}\into \frakB_{i}^{i,x,y}$$
denote the inclusion. 
\end{remark}

\begin{proposition}
\label{neckequiv}
Assume \cref{simplex.notation}. There is a functor 
$$\mathfrak{u}^{x,y}_{i}:\frakB_{i}^{i,x,y}\to \Nec/(\frakB^i_{i},x,y),$$
which is an equivalence of categories.
Moreover, the functor $\mathfrak{u}^{x,y}_i$ restricts along the inclusion $\tau_i$ (\cref{taumap}) to a functor $\mathfrak{u}_{i-1}^{x,y}:\frakB_{i-1}^{i,x,y}\to \Nec/(\frakB^i_{i-1},x,y)$, which is also an equivalence of categories.
\end{proposition}

\begin{proof}
The functor $\mathfrak{u}^{x,y}_{i}$ sends an object in $\frakB^{i,x,y}_{i}$
to the necklace whose joint cuts are indexed by $\Upsilon$ and whose beads are indexed by consecutive pairs of elements in $\Upsilon$,
as constructed in \cref{alternative.nec}.
The embedded bordism corresponding to a bead $\upsilon<\upsilon'\in \Upsilon$ is given by $A_{\upsilon,\upsilon'}$ in \cref{etalemapsnow}. 	

By \cref{alternative.nec}, the data of $\sigma$ in \cref{etalemapsnow} uniquely corresponds to a morphism $\sigma'$ of necklaces.
The condition that $[m_i,f,h,C,P,\Upsilon]=[m_i,f',h',\sigma^*C',P',\Upsilon]$ implies that the resulting morphism of necklaces is a morphism in $\Nec/(\frakB^i_{i},x,y)$.
We define $\mathfrak{u}^{x,y}_{i}(\sigma)=\sigma'$. That this assignment respects composition again follows from \cref{alternative.nec}, so that $\mathfrak{u}^{x,y}_{i}$ defines a functor. 

Faithfulness of $\mathfrak{u}$ is clear.
To prove fullness, let $A=\mathfrak{u}[m_i,f,h,C,P,\Upsilon]$ and $B=\mathfrak{u}[m'_i,f',h',C',P',\Upsilon']$
and let $\sigma':A\to B$ be a morphism in $\Nec/(\frakB^i_{i},x,y)$.
Let $\sigma:([m_i],\Upsilon)\to ([m_i'],\Upsilon')$ be the morphism in $\Nec'$
associated to $\sigma'$ under the equivalence of categories \cref{alternative.nec}.
To see that $\sigma$ yields a morphism in $\frakB_{i-1}^{i,x,y}$ (\cref{etalemapsnow}),
we must verify the equality $[m_i,f,h,C,P,\Upsilon]=[m_i,f',h',\sigma^*C',P',\Upsilon]$.
Conditions (e1) ($m_i=m'_i$) and (e4) ($Υ=Υ'$) hold by construction.
It remains to construct $ε<\min\{δ,δ'\}$ and $N''$ such that $f=f'$ and Condition~(e3) holds, after restricting to $B^{\mathfrak{p}}_ε⨯B^{\mathfrak{l}}_ε$.
The morphism $\sigma':A\to B$ provides for every consecutive pair $υ<υ'$ in~$Υ$
the corresponding $ε_{υ,υ'}$ and $N''_{υ,υ'}$ that witness the equality of the bead $A_{υ,υ'}$ and the embedded bordism $A'_{σ(υ),σ(υ')}$.
Take $ε=\min_{υ<υ'}ε_{υ,υ'}$ and restrict to $B^{\mathfrak{p}}_ε⨯B^{\mathfrak{l}}_ε$.
Now $f=f'$.
Take $N''=⋃_{υ<υ'}N''_{υ,υ'}$, where the union is taken among subobjects of $W⨯_U B^{\mathfrak{p}}_ε⨯B^{\mathfrak{l}}_ε$. The resulting map $N''\to W⨯_U B^{\mathfrak{p}}_ε⨯B^{\mathfrak{l}}_ε$ factors uniquely through the subobjects $N$ and $N'$ because $N''_{\upsilon,\upsilon'}$ satisfies this property, for each consecutive $\upsilon<\upsilon'\in \Upsilon$. Replace $N''$ with a smaller open subobject
whose reduction contains the core of $[m_i,f,h,C,P,\Upsilon]$
and the map $N''→W⨯_U B^{\mathfrak{p}}_ε⨯B^{\mathfrak{l}}_ε$ is a (locally trivial) fiber bundle,
which is possible by \cref{smallercore}.
Furthermore, $C=σ^*C'$, $P=P'$ after restricting to $N''$,
since these conditions hold on a covering family of~$N''$ given by $N''_{υ,υ'}∩N''$ for all consecutive $υ<υ'$ in~$Υ$.

It is essentially surjective since any object in the comma category $\Nec/(\frakB^i_{i},x,y)$ can be glued to obtain an object in \cref{etalemapsnow} as follows.
Pick a representative for every bead in the necklaces.
Use Condition~(e2) in \cref{etalemapsnow} to shrink the radius of the open balls $V$ and $L$ so that all beads have the same~$f$.
Use Condition~(e3) in \cref{etalemapsnow} to shrink the neighborhoods~$N$ so that the $C$- and $P$-components of all beads are compatible with each other.
We are going to construct open neighborhoods $Z_{υ,υ'}\subset \stu(N_{\langle \upsilon,\upsilon'\rangle})$ of $\core(A_{υ,υ'})_0$,
where $υ=υ'$ or $υ<υ'$ is a consecutive pair of elements in~$Υ$
such that $Z_{υ,υ}⊂Z_{υ,υ'}$ and $Z_{υ',υ'}⊂Z_{υ,υ'}$.
To this end, pick $ε>0$ such that the $ε$-neighborhood of $\core(A_{υ,υ'})_0$ is a subset of $\stu(N_{υ,υ'})$, for all $υ≤υ'$.
Denote by $\zeta_{υ,υ'}$, $\zeta_{υ',υ''}$, and $\zeta_{υ',υ'}$
the height functions of the cut~$υ'$, defined on $\stu( N_{υ,υ'})$, $\stu( N_{υ',υ''})$, and $\stu( N_{υ',υ'})$ respectively.
Now take
$$E_{υ,υ'}=\zeta_{υ,υ'}^{-1}(-∞,-ε)∪\zeta_{υ',υ'}^{-1}(-2ε,2ε)⊂\stu(N_{υ,υ'}),$$
$$F_{υ',υ''}=\zeta_{υ',υ''}^{-1}(ε,∞)∪\zeta_{υ',υ'}^{-1}(-2ε,2ε)⊂\stu(N_{υ',υ''}),$$
$$Z_{υ',υ'}=(Z_{υ,υ'}⨯_{\stu(W\times_{U}V\times L)} Z_{υ',υ''})∩\zeta_{υ',υ'}^{-1}(-ε,ε),$$
for consecutive $υ<υ'<υ''$
and set $Z_{υ,υ'}=E_{υ,υ'}∩F_{υ,υ'}$.
This guarantees that for every consecutive triple $υ<υ'<υ''$ the map $$Z_{υ',υ'}→Z_{υ,υ'}⨯Z_{υ',υ''}$$ is a closed map,
hence the pushout
$$Z_{υ_0,υ_1}⊔_{Z_{υ_1,υ_1}}⊔⋯⊔_{Z_{υ_{|Υ|-1},υ_{|Υ|-1}}}Z_{υ_{|Υ|},υ_{|Υ|}}$$
is a (Hausdorff) manifold equipped with the induced etale map to~$\stu(W\times_{U}V\times L)$. Lift this map to a unique morphism of structured manifolds $h:N\to W\times_{U}V\times L$.

We can now restrict the beads to the $Z$-neighborhoods constructed above and glue the $C$- and $P$-components together,
which form a compatible family by the above construction.
The $f$-component is the same for the chosen representatives of all beads.
This yields the desired etale bordism, which show the essential surjectivity of~$\mathfrak{u}$.

The claim that the restriction of $\mathfrak{u}_i^{x,y}$ along $\tau_i$ factors through $\Nec/(\frakB^i_{i-1},x,y)$ follows immediately by Condition (6b) for $\frakB^i_{i-1}$.
That it is an equivalence follows by the same argument as for $\frakB^i_i$. 
\end{proof}

The above category may seem complicated, however it turns out that the map $\tau_i$ splits as a disjoint union whose summands are indexed by composed bordisms, where the composition happens in the $i$th direction. This allows us to reduce the analysis of the map $\tau_i$ in \cref{taumap} to each summand. 

\begin{proposition}
\label{decompbi}
Consider the category $G$ defined exactly as in \cref{etalemapsnow} for $\frakB_i^{i,x,y}$, but with Conditions (6a) and (6b) dropped.
The category $\frakB_{i}^{i,x,y}$ is a full subcategory of $G$.
Denote by $D$ the union of connected components of $G$ that intersect with $\frakB_i^{i,x,y}$.
Denote by $\frakA$ the full subcategory of~$D$ on objects that have $[m_i]=[1]$.
Then $\frakA$ is a discrete category.
Furthermore, the inclusion $\tau_{i}$ in \cref{taumap} splits as a disjoint union 
$$\tau_i^{x,y}:\frakB_{i-1}^{i,x,y}=\coprod_{\fraka\in \frakA}\frakB_{i-1}^{\fraka}
\lto9{\coprod_{\fraka\in \frakA} \tau_i^\fraka}
\coprod_{\fraka\in \frakA}\frakB_i^{\fraka}=\frakB_i^{i,x,y},$$
where each summand is obtained by intersecting with the connected component of $\fraka$ inside $G$. 
\end{proposition}

\begin{proof}
That $\frakA$ is a discrete category follows from the fact that $σ:[1]→[1]$ must be the identity map for any morphism in~$\frakA$
and in case $σ=\id_{[1]}$, the relation $[m_i,f,h,C,P,\Upsilon]=[m_i,f',h',\sigma^*C',P',\Upsilon]$ implies that the source object is equal to the target object.

Suppose two objects $\fraka,\fraka'∈\frakA$ are in the same connected component of $G$,
meaning that are connected by a chain of morphisms going in either direction.
Applying the simplicial structure map to all intermediate objects in the chain
to remove all intermediate cuts in the $i$th direction (i.e., other than the initial and terminal cuts corresponding to $x$ and~$y$)
produces another chain of morphisms connecting $\fraka$ to $\fraka'$ that lies entirely within~$\frakA$.
Since $\frakA$ is a discrete category, this implies $\fraka=\fraka'$.
This induces a decomposition of the inclusion $\tau_i$ into a disjoint union $\coprod_{\fraka\in \frakA} \tau_i^\fraka$.
\end{proof}

The summands $\tau^\fraka_i$ in \cref{decompbi} are relatively simple. A choice of summand $\fraka$ fixes most of the data in \cref{etalemapsnow}. The remaining data is given by necklaces and morphisms of necklaces on a fixed bordism. Such morphisms are given by compositions of two elementary operations: insertion of a non-joint cut and changing a joint cut into a non-joint cut.

\begin{definition}
\label{etalebordism}
Assume the notation of \cref{decompbi}. We call an element $\fraka\in \frakA$ an \emph{etale bordism}. This terminology is justified by the fact that such an element corresponds to an equivalence class in \cref{etalemapsnow} of the form $[[1],f,h,C,P,\{0,1\}]$, with Conditions (6a) and (6b) dropped. In particular, the map $h:N\to W\times_{U}V\times L$ is only required to be an etale map.
\end{definition}
\begin{figure}[ht]
\begin{tikzpicture}
\draw[line width=.5mm] (0,0) -- (0,-2);
\draw (.25,0) -- (.25,-2);
\draw (.5,0) -- (.5,-2);
\draw (.75,0) -- (.75,-2);
\draw[line width=.5mm] (1,0) -- (1,-2);
\draw (1.25,0) -- (1.25,-2);
\draw[line width=.5mm] (1.5,0) -- (1.5,-2);
\draw (1.75,0) -- (1.75,-2);
\draw (2,0) -- (2,-2);
\draw (2.25,0) -- (2.25,-2);
\draw (2.5,0) -- (2.5,-2);
\draw (2.75,0) -- (2.75,-2);
\draw[line width=.5mm]  (3,0) -- (3,-2);
\draw (3.25,0) -- (3.25,-2);
\draw (3.5,0) -- (3.5,-2);
\draw (3.75,0) -- (3.75,-2);
\draw[line width=.5mm] (4,0) -- (4,-2);
\draw (0,0) -- (4,0);
\draw (0,-2) -- (4,-2);
\draw[decoration={brace,mirror},decorate] (0,-2.2) -- (1,-2.2);
\node at (.5,-2.75) {\text{first bead} $C_1$};
\draw[->] (4.5,-1) -- (7,-1);
\begin{scope}[xshift=2cm,yshift=0cm]
\draw[line width=.5mm] (5.5,0) -- (5.5,-2);
\draw (5.75,0) -- (5.75,-2);
\draw (6,0) -- (6,-2);
\draw (6.125,0) -- (6.125,-2);
\draw (6.25,0) -- (6.25,-2);
\draw (6.5,0) -- (6.5,-2);
\draw (6.75,0) -- (6.75,-2);
\draw[line width=.5mm] (7,0) -- (7,-2);
\draw (7.25,0) -- (7.25,-2);
\draw (7.375,0) -- (7.375,-2);
\draw (7.5,0) -- (7.5,-2);
\draw (7.625,0) -- (7.625,-2);
\draw (7.75,0) -- (7.75,-2);
\draw (8,0) -- (8,-2);
\draw (8.25,0) -- (8.25,-2);
\draw[line width=.5mm]  (8.5,0) -- (8.5,-2);
\draw (8.75,0) -- (8.75,-2);
\draw (9,0) -- (9,-2);
\draw (9.25,0) -- (9.25,-2);
\draw[line width=.5mm] (9.5,0) -- (9.5,-2);
\draw (5.5,0) -- (9.5,0);
\draw (5.5,-2) -- (9.5,-2);
\draw[decoration={brace,mirror},decorate] (5.5,-2.2) -- (7,-2.2);
\node at (6,-2.75) {\text{first bead} $C_1$};
\end{scope}
\end{tikzpicture}
\caption{An composition of two elementary morphisms in the necklace category. The morphism turns a joint cut into a non-joint cut and inserts a non-joint cut}
\end{figure}

\begin{proposition}
\label{weak.eq.necklaces}
For each pair of vertices $x$ and $y$ in $\frakB^i_i([0])$, the induced map
$$\Cnec(\rho_i):\Cnec(\frakB^i_{i-1})(x,y)\to \Cnec(\frakB^i_i)(x,y)$$
is a weak equivalence of simplicial sets.
\end{proposition}

\begin{proof}
By \cref{neckequiv}, we have a commutative diagram 
$$\xymatrix{
\frakB_{i-1}^{i,x,y}\ar[r]^-{\mathfrak{u}_{i}^{x,y}}\ar[d]_-{\tau_i^{x,y}} & \Nec/(\frakB^i_{i-1},x,y)\ar[d]\\
\frakB_{i}^{i,x,y}\ar[r]^-{\mathfrak{u}_{i-1}^{x,y}} & \Nec/(\frakB^i_{i},x,y),
}$$
where the right vertical map is induced by the monomorphism $\rho_i:\frakB^i_{i-1}\into \frakB^i_i$.  By \cref{neckequiv}, the horizontal maps in the diagram are equivalences of categories. By definition of $\Cnec$ (\cref{Cnec.def}), the nerve of the right vertical map is precisely $\Cnec(\rho_i)$. Thus, it suffices to prove that the nerve of $\tau_i$ is a weak equivalence. 

By \cref{decompbi}, the map $\tau_i$ splits as a disjoint union $\coprod_{\fraka\in \frakA}\tau^\fraka_i$, with $\frakA$ defined in \cref{decompbi}. Therefore, it suffices to show that the nerve of each $\tau^\fraka_i:\frakB_{i-1}^{\fraka}\to \frakB_{i}^{\fraka}$ is a weak equivalence. By \cref{contractible.cuts}, the nerve of  the comma categories $c\downarrow \tau^\fraka_i$ are weakly contractible, for all $c\in \frakB_{i}^{\fraka}$. Hence the nerve of $\tau^\fraka_i$ is a weak equivalence, by Quillen's Theorem A, for all $\fraka\in \frakA$. This implies that the nerve of $\tau_i$ is also a weak equivalence.
\end{proof}

\subsection{Contractibility of categories of necklaces}

We are ready to prove the claim about contractibility of the summands of necklace categories.

\begin{notation}
\label{etale.bordism}
Recall \cref{filtration.notation} and \cref{simplex.notation}.
In addition to these, we also fix the following data,
taken for some fixed $V=B^{\mathfrak{p}}_{\delta}\subset T\in \stcart$ and $L=B^{\mathfrak{l}}_{\delta}\subset\RR^{n}$ (\cref{stalks}):
\begin{itemize}
\item a morphism $f:V\times L\to U$ in $\stcart$;
\item an etale map $h:N\to  W\times_{U}V\times L$ such that the submersion $\pi\circ h:N\to V\times L$ is a (locally trivial) fiber bundle, where $\pi$ is the projection;
\item a compact (globular) cut ${\bf n}$-grid $C$ on $\pi\circ h:N\to V\times L$,
where ${\bf n}=([m_1],\ldots,[1],\ldots,[m_d])$, where ${\bf m}=([m_1],\hdots,\widehat{[m_i]},\hdots,[m_d])\in \Delta^{\times (d-1)}$
is the fixed multisimplex ${\bf m}$ in \cref{simplex.notation}; 
\item a morphism $P:N\to \langle \ell\rangle$ in $\stman$;
\item an object $c\in \frakB_i^{\fraka}$. 
\end{itemize}
We require the class $\fraka=[[1],f,h,C,P,\{0,1\}]$ to belong to the subcategory~$\frakA$ defined in \cref{decompbi}. The latter condition forces the source and target in the $i$th direction to be $x$ and~$y$ respectively.
We also use the following notation.
\begin{itemize}
\item Let $\phi\leq \chi$ be cuts of $c$. We denote by $c_{[\phi, \chi]}$ the etale necklace obtained by removing all cuts of $c$ that are greater than $\chi$ or less than $\phi$ (possibly including the source and target cut). Thus, $\phi$ is the source cut of $c_{[\phi, \chi]}$ and $\chi$ is the target cut. 
\item We set
$$\frakD_{[\phi,\chi]}:= c_{[\phi, \chi]} \downarrow \tau_i^{\fraka_{[\phi, \chi]}},$$ 
where $\tau_i^{\fraka_{[\phi, \chi]}}$ is the summand defined in \cref{decompbi}, by taking $x=\phi$ and $y=\chi$ (including replacing the source and target cuts of $\fraka$). In particular, 
 $$\frakD_{[x, y]}:= c \downarrow \tau_i^{\fraka}$$
\item For an object $g:c\to d\in \frakD_{[x,y]}$, we also define an object $g_{[\phi, \chi]}\in \frakD_{[\phi, \chi]}$ as the morphism of necklaces obtained by removing all cuts less than $\phi$ or greater than $\chi$ in both $c$ and $d$. 
\end{itemize}
%
\end{notation}

\begin{remark}\label{necklacecutsd}
By definition of the morphisms in the etale necklace category (\cref{etalemapsnow}), an object in $\frakD_{[x,y]}$ is an etale necklace with the property that removing some non-joint cuts and changing some non-joint cuts into joint cuts yields the fixed necklace $c\in \frakB_i^{\fraka}$. A morphism in the comma category $d\to d'\in \frakD_{[x,y]}$ is a morphism of necklaces from the target of $d$ to the target of $d'$ so that the corresponding maps of simplices commute. For this reason, the forgetful functor $\frakD_{[x,y]}\to \frakB_{i-1}^{\fraka}$ is faithful. 

The comma categories $\frakD_{[\phi,\chi]}$ have the same description, except that the source and target cuts change from $x$ and $y$ to $\phi$ and $\chi$ (respectively). When working with objects $g\in \frakD_{[\phi,\chi]}$, we will always work with the necklace given by the target of $g$. In particular, cuts are only inserted or removed in the target and we only perform these operations when they yield an object in $\frakD_{[\phi,\chi]}$. 
\end{remark}

The following notation will be used in \cref{contractible.cuts}.

\begin{proposition}
\label{move.cut} 
Assume \cref{etale.bordism} as well as the conventions of \cref{necklacecutsd}. Let $\phi\leq \chi$ be cuts of $c$. Recall that for a cut tuple~$C$ its \emph{inner cuts} are defined as all cuts of~$C$ except for the initial and terminal cut. Let $g:Q→\frakD_{[\phi,\chi]}$ be a functor and suppose $A$ and $B$ are cuts in the $i$th direction such that $\phi≤A≤B≤\chi$. For every $q∈Q$ define the following etale necklaces (\cref{etalemapsnow}).
\begin{itemize}
\item
Denote by $E$ the set of indices of the cut tuple in the $i$th direction of $g(q)$ that are equal to~$B$ and are not initial.
Then $g^2(q)$ is the etale necklace obtained by inserting exactly $|E|$ copies of the cut~$A$ (in the same status as joint or non-joint cut)
immediately before~$E$ in the $i$th direction. 
\item $g^1(q)$ is obtained by turning the newly added cuts in $g^2(q)$ into non-joint cuts.
\item $\hat g^1(q)$ is obtained from $g(q)$ by converting all inner cuts equal to $B$ into non-joint cuts.
\item $\hat g^2(q)$ is obtained from $g(q)$ by removing all inner cuts equal to~$B$. Optionally, we may keep the last copy of $B$.
\end{itemize}
Then the following hold. 
\begin{enumerate}
\item If for all $q\in Q$, the cut tuple in the $i$th direction~$C^i$ of $g(q)$ and the etale necklace $g^2(q)$ satisfy
\begin{itemize}
\item if $C^i_k=B$ for some~$k$,
then for every index $l$ we have $C^i_l≤A$ or $C^i_l≥B$;
\item $g^2(q)$ belongs to $\frakD_{[\phi,\chi]}$;
\end{itemize}
then there is a zigzag of natural transformations of functors $Q→\frakD_{[\phi,\chi]}$ 
from the functor~$g$ to a functor~$g'$ such that $g'(q)=g^2(q)$.

\item If for every $q∈Q$, the necklace $\hat g^1(q)$ belongs to $\frakD_{[\phi,\chi]}$ 
then there is a zigzag of natural transformations of functors $Q→\frakD_{[\phi,\chi]}$ to a functor $g'$ such that $g'(q)=\hat g^2(q)$. 
\end{enumerate}
\end{proposition}

\begin{proof}
Parts (1) and (2) use the same construction, except that in Part~(2) we reverse the direction of the zigzag.
Thus, the next two paragraphs talk about Part~(1) only, constructing a zigzag of the form $g=g^0→g^1←g^2=g'$. 

Fix an arbitrary element $q\in Q$.
Denote by $([m],Υ)$ the necklace corresponding to the etale necklace $g^0(q)$.
Denote by $F⊂[m]$ the subset of (consecutive) inner vertices of~$[m]$ for which the corresponding inner cut in the etale necklace of $g^0(q)$ equals~$B$.
Denote by $[m']$ the simplex given by the disjoint union $[m]_{<F}⊔F⊔F⊔[m]_{>F}$, ordered as indicated.
Denote by $ι_1,ι_2:[m]→[m']$ the two embeddings that map $[m]_{<F}$ and $[m]_{>F}$ via the identity map
and send $F⊂[m]$ to the first respectively second copy of~$F$ in~$[m']$.
%

Next, fix an arbitrary morphism $σ:q→\bar q$ in~$Q$,
which is given by a morphism of necklaces $σ:([m],Υ)→([\bar m],\bar\Upsilon)$.
Denote by $F⊂[m]$ and $\bar F⊂[\bar m]$ the subsets defined above.
Denote by $σ':[m']→[\bar m']$ the map of simplices (with $[m']$, $[\bar m']$ defined above)
defined as the map
$$σ':[m']=[m]_{<F}⊔F⊔F⊔[m]_{>F}→[\bar m]_{<\bar F}⊔\bar F⊔\bar F⊔[\bar m]_{>\bar F}$$
that is induced by $σ$ on all four summands, mapping the $i$th summand on the left to the $i$th summand on the right.
Set the morphisms $g^1(σ)$ and $g^2(σ)$ to the map~$σ'$.
This turns $g^1$ and $g^2$ into functors.

For Part~(1), for all~$q$ the object $g^2(q)$ belongs to $\frakD_{[\phi,\chi]}$, by assumption.
Turning some of the newly added cuts into non-joint cuts 
produces $g^1(q)$, which therefore also belongs to $\frakD_{[\phi,\chi]}$,
since this operation does not change the $(i-1)$-core. The maps of necklaces $$([m],Υ)\lto3{ι_2}([m'],ι_2(Υ))\lgets3{\id_{[m']}}([m'],ι_1(Υ)∪ι_2(Υ))$$
yield the components of natural transformations $g^0→g^1←g^2$. 

For Part~(2), the object $\hat g^1(q)$ belongs to $\frakD_{[\phi,\chi]}$ by assumption.
The object $\hat g^2(q)$ belongs to $\frakD_{[\phi,\chi]}$ since the operation of removing non-joint cuts does not change the $(i-1)$-core. 
\end{proof}

\begin{figure}
\def\jcut{\color{blue}\bf}
\label{processing.cuts2}
\begin{tikzpicture}
\draw[->] (0,1.3) -- (0,1.7);
\draw[->] (0,2.3) -- (0,2.7);
\draw[->] (0,3.3) -- (0,3.7);
\draw[->] (0,4.3) -- (0,4.7);

\node at (0,1) {$\jcut C_0$};
\node at (0,2) {$C_1$};
\node at (0,3) {$B$};
\node at (0,4) {$\jcut B$};
\node at (0,5) {$\jcut C_4$};
\draw[->] (1+1,1.3) -- (1+1,1.7);
\draw[->] (1+1,2.3) -- (1+1,2.7);
\draw[->] (1+1,3.3) -- (1+1,3.7);
\draw[->] (1+1,4.3) -- (1+1,4.7);
\draw[->] (2,5.3) -- (2,5.7);
\draw[->] (2,6.3) -- (2,6.7);

\node at (1+1,1) {$\jcut C_0$};
\node at (1+1,2) {$C_1$};
\node at (1+1,3) {$A$};
\node at (1+1,4) {$A$};
\node at (1+1,5) {$B$};
\node at (2,6) {$\jcut B$};
\node at (2,7) {$\jcut C_4$};

\draw[->] (4,1.3) -- (4,1.7);
\draw[->] (4,2.3) -- (4,2.7);
\draw[->] (4,3.3) -- (4,3.7);
\draw[->] (4,4.3) -- (4,4.7);
\draw[->] (4,5.3) -- (4,5.7);
\draw[->] (4,6.3) -- (4,6.7);

\node at (4,1) {$\jcut C_0$};
\node at (4,2) {$C_1$};
\node at (4,3) {$A$};
\node at (4,4) {$\jcut A$};
\node at (4,5) {$ B$};
\node at (4,6) {$\jcut B$};
\node at (4,7) {$\jcut C_4$};

\draw[->] (6,1.3) -- (6,1.7);
\draw[->] (6,2.3) -- (6,2.7);
\draw[->] (6,3.3) -- (6,3.7);
\draw[->] (6,4.3) -- (6,4.7);
\draw[->] (6,5.3) -- (6,5.7);
\draw[->] (6,6.3) -- (6,6.7);

\node at (6,1) {$\jcut C_0$};
\node at (6,2) {$ C_1$};
\node at (6,3) {$A$};
\node at (6,4) {$\jcut A $};
\node at (6,5) {$ B$};
\node at (6,6) {$B$};
\node at (6,7) {$\jcut C_4 $};

\node at (8,1) {$\jcut C_0$};
\node at (8,2) {$ C_1$};
\node at (8,3) {$A$};
\node at (8,4) {$\jcut A$};
\node at (8,5) {$\jcut C_4 $};

\draw[->] (8,1.3) -- (8,1.7);
\draw[->] (8,2.3) -- (8,2.7);
\draw[->] (8,3.3) -- (8,3.7);
\draw[->] (8,4.3) -- (8,4.7);

\draw[->] (.2,5) -- (1.7,7);
\draw[->] (.2,4) -- (1.7,6);
\draw[->] (.2,3) -- (1.7,5);
\draw[->] (.2,2) -- (1.7,2);
\draw[->] (.2,1) -- (1.7,1);

\draw[->] (3.7,7) -- (2.3,7);
\draw[->] (3.7,6) -- (2.3,6);
\draw[->] (3.7,5) -- (2.3,5);
\draw[->] (3.7,4) -- (2.3,4);
\draw[->] (3.7,3) -- (2.3,3);
\draw[->] (3.7,2) -- (2.3,2);
\draw[->] (3.7,1) -- (2.3,1);

\draw[->] (4.3,7) -- (5.7,7);
\draw[->] (4.3,6) -- (5.7,6);
\draw[->] (4.3,5) -- (5.7,5);
\draw[->] (4.3,4) -- (5.7,4);
\draw[->] (4.3,3) -- (5.7,3);
\draw[->] (4.3,2) -- (5.7,2);
\draw[->] (4.3,1) -- (5.7,1);

\draw[->] (7.7,5) -- (6.3,7);
\draw[->] (7.7,4) -- (6.3,4);
\draw[->] (7.7,3) -- (6.3,3);
\draw[->] (7.7,2) -- (6.3,2);
\draw[->] (7.7,1) -- (6.3,1);
\end{tikzpicture}
\caption{An example of the zigzags of natural transformations (1) and (2) of \cref{move.cut}. The etale necklace $g(q)$ is the rightmost necklace. The blue cuts are the joint cuts. The first zigzag in the composition is the $q$ component of the zigzag of natural transformations (1), which inserts duplicate $A$ cuts immediately behind duplicate $B$ cuts. The second zigzag is the $q$ component of the zigzag of natural transformations in (2), which removes $B$ cuts.}
\end{figure}

In the next proposition, we will apply \cref{new.whitehead} to prove contractibility of $\frakD_{[x,y]}$. Starting with a full subcategory with finitely many objects $Q$, we need to construct a natural transformation of zigzags connecting the inclusion $Q\into \frakD_{[x,y]}$ to a constant functor.  The idea is to use the cuts $\Psi$ constructed in \cref{morse.decomposition} as an anchor and slowly move all the cuts in the $i$th direction in all beads in all necklaces of $Q$ to the left of the $\Psi$ cuts until all that remains is the $\Psi$ cuts. The zigzag in \cref{move.cut} is what allows us to move cuts. The size of the region between two successive $\Psi$ cuts is controlled by a parameter $\varepsilon$. By choosing a small enough $\varepsilon$, we can ensure that each time we move a cut past a $\Psi$ cut, the connected components of the $(i-1)$-cores possibly shrink and expand at most by $\varepsilon$. This ensures that the zigzag remains in $\frakD_{[x,y]}$. 

\begin{figure}[ht]
\begin{center}
\begin{tikzpicture}
\draw[fill=gray] (-2,2.5) parabola bend (0,1) (2,2.5);
\draw[fill=gray] (-2,-2.5) parabola bend (0,-1) (2,-2.5);
\draw[thick] (-2.5,2) parabola bend (0,.7) (2.5,2);
\draw[thick, rotate=180] (-2.5,2) parabola bend (0,.7) (2.5,2);
\draw[rotate=90, fill=gray!25] (-2,2.5) parabola bend (0,1) (2,2.5);
\draw[rotate=90, fill=gray!25] (-2,-2.5) parabola bend (0,-1) (2,-2.5);
\draw[thick] (-2.5,2) .. controls (-2.37,2) and (-1.8,1.5) .. (-1.4,1.1); 
\begin{scope}[xscale=1,yscale=-1]
\draw[thick] (-2.5,2) .. controls (-2.37,2) and (-1.8,1.5) .. (-1.4,1.1); 
\end{scope}
\begin{scope}[xscale=-1,yscale=1]
\draw[thick] (-2.5,2) .. controls (-2.37,2) and (-1.8,1.5) .. (-1.4,1.1); 
\end{scope}
\begin{scope}[xscale=-1,yscale=-1]
\draw[thick] (-2.5,2) .. controls (-2.37,2) and (-1.8,1.5) .. (-1.4,1.1); 
\end{scope}
\draw (-3,0) -- (3,0);
\draw (0,-3) -- (0,3);
\draw[thick, dashed] (-3,1) .. controls (4,0) and (-4,0) .. (3,-1);
\draw[thick] (-3,1) .. controls (-1.3,.7) and (-.4,.5) .. (-1.5,-1.2);
\begin{scope}[xscale=-1,yscale=-1]
\draw[thick] (-3,1) .. controls (-1.3,.7) and (-.4,.5) .. (-1.5,-1.2);
\end{scope}
\node at (0,0) {$\bullet$};
\node at (-3,2) {$\Psi_{j+1}$};
\node at (-2.2,.5) {$D_k'$};
\node at (-2,1.1) {$D_k$};
\node at (-.5,-.5) {$H$};
\node at (-.2,.2) {$c$};
\node at (4.5,0) {$(u_1,\ldots,u_i)$-axis};
\node at (.5,3.3) {$(u_{i+1},\ldots,u_d)$-axis};
\end{tikzpicture}
\end{center}
\caption{
The cut $\Psi$ is the two parabolas that bound the light gray region.
This picture illustrates the procedure of \cref{contractible.cuts}, when the cut $\Psi$ is given by Milnor's construction \cite{Milnor} at a critical point of a Morse function.
The cut $\Psi_{j+1}$ is the bold cut that follows the two parabolas (bounding the light gray region) and then passes through the ``bridge''~$H$.
The white region labeled $H$ (the bridge) is the closure of $\Psi_{(j,j+1)}$, which is contained in an open ball centered at $c$.
The solid cut $D'_k$ is the result of the surgery, and $D_k'\leq \Psi_j$, while the dashed cut $D_k$ passes through $\Psi_{(j,j+1)}$.}
\end{figure}

\begin{figure}[ht]
\label{processing.cuts}
\begin{tikzpicture}
\draw[line width=.5mm] (0,0) -- (0,-2);
\draw (.25,0) -- (.25,-2);
\draw (.5,0) -- (.5,-2);
\draw (.75,0) -- (.75,-2);
\draw[line width=.5mm] (1,0) -- (1,-2);
\draw (1.25,0) -- (1.25,-2);
\draw[line width=.5mm] (1.5,0) -- (1.5,-2);
\draw (1.75,0) -- (1.75,-2);
\draw (2,0) -- (2,-2);
\draw (2.25,0) -- (2.25,-2);
\draw (2.5,0) -- (2.5,-2);
\draw (2.75,0) -- (2.75,-2);
\draw[line width=.5mm]  (3,0) -- (3,-2);
\draw (3.25,0) -- (3.25,-2);
\draw (3.5,0) -- (3.5,-2);
\draw (3.75,0) -- (3.75,-2);
\draw[line width=.5mm] (4,0) -- (4,-2);
\draw (0,0) -- (4,0);
\draw (0,-2) -- (4,-2);
\draw[very thick] (2,0) .. controls (2.5,0) and (2,-2) .. (3.5,-2);
\node at (2.65,-1) {$\Psi$};
\node at (2.4,-2.5) {$(C_q)_3$};
\draw[->] (4.3,-1) -- (7,-1);
\node at (5.5,-.7) {\cref{contractible.cuts}};
\end{tikzpicture}
\begin{tikzpicture}
\draw[line width=.5mm] (0,0) -- (0,-2);
\draw (.25,0) -- (.25,-2);
\draw (.5,0) -- (.5,-2);
\draw (.75,0) -- (.75,-2);
\draw[line width=.5mm] (1,0) -- (1,-2);
\draw (1.25,0) -- (1.25,-2);
\draw[line width=.5mm] (1.5,0) -- (1.5,-2);
\draw (1.75,0) -- (1.75,-2);
\draw (2,0) -- (2,-2);

\draw (0,0) -- (4,0);
\draw (4,0) -- (4,-2);
\draw (4,-2) -- (0,-2);
\draw (2.25,-2) .. controls (2.25,.3) and (2,0) .. (2,0);
\draw (2.5,-2) .. controls (2.5,-1.3) and (2.42,-1.2) .. (2.42,-1.2);
\draw (2.75,-2) .. controls (2.75,-1.8) and (2.69,-1.7) .. (2.69,-1.7);
\draw[line width=.5mm] (3,-2) .. controls (3,-1.9) and (2.89,-1.845) .. (2.89,-1.845);
\draw[very thick] (2,0) .. controls (2.5,0) and (2,-2) .. (3.5,-2);
\node at (2.65,-1) {$\Psi$};
\node at (2.4,-2.5) {$(C'_q)_3$};
\end{tikzpicture}
\caption{An illustration of zigzag in \cref{contractible.cuts}.
The outer box containing the cuts represents the core of the bordism.
The cut $\Psi$ is the curved bold cut.
Joint cuts are represented by bold vertical cuts. The cut tuple labeled $(C_q)_3$ is the third bead in the necklace $C_q$. 
Processing the cuts with respect to a fixed total ordering, we cut and glue cuts one by one to the cut $\Psi$.
}
\end{figure}

\begin{proposition}
\label{contractible.cuts}
Assume \cref{etale.bordism}. 
The nerve of the comma category $\frakD_{[x,y]}$ is weakly contractible.
\end{proposition}

\begin{proof}
We use \cref{new.whitehead}, which requires us to construct
for every finite full subcategory~$Q$ of~$\frakD_{[x,y]}$,
a zigzag of natural transformations from the inclusion functor $g:Q→\frakD_{[x,y]}$
to a composition $\Psi:Q→1→\frakD_{[x,y]}$, where $1$ denotes the terminal category.  

Recall that objects in $\frakD_{[x,y]}$ are morphisms of necklaces $c→d$.
Such a morphism can add new cuts in~$c$ and it can also convert the existing joint cuts of~$c$ into non-joint cuts.

\medskip
\paragraph{\bf Making all cuts in the image of $c$ joint.}

The very first natural transformation $g←g'$ in the zigzag of natural transformations that connects $g$ to~$\Psi$
ensures that for every $q∈Q$ the object $g'(q)$ is given by a morphism of necklaces $c→d$
in which all cuts of~$c$ (joint or not) map to joint cuts of~$d$.
Specifcally, given $q∈Q$, we define the component $g(q)←g'(q)$ as follows.
The object $g'(q)$ is given by the same data $c→d$ as $g(q)$,
except that all cuts of~$d$ in the image of~$c$ are turned into joint cuts.
The map $g(q)←g'(q)$ is given by the identity map of corresponding simplices and removes some elements of the subset $\Upsilon$, corresponding to the non-joint cuts in the image of $c$.
The resulting morphisms $g(q)←g'(q)$ are morphisms of necklaces that yield a natural transformation $g←g'$.

From now on, we assume without loss of generality that for all $q∈Q$ the object $g(q)$
is given by a morphism of necklaces $c→d$ that maps all cuts of~$c$ to joint cuts of~$d$.
This property will be preserved by all constructions below.


 The zigzag will be constructed using three nested inductive constructions.
 
 \medskip
 \paragraph{\bf Outer inductive layer: induction on the cuts of $c$.}
 
The outer induction will be on the cuts~$\chi$ of~$c$, in decreasing order.
After the corresponding inductive step for~$\chi$ is complete,
we will have constructed an object $\Psi_{\chi}\in \frakD_{[\chi, y]}$, whose cuts are all joint cuts, and a zigzag of homotopies connecting the map $g:Q→\frakD$
to a map $g_\chi:Q→\frakD$ such that for every $q∈Q$, the object $g_\chi(q)$
has following properties:
\begin{itemize}
\item  $g_{\chi}(q)_{[\chi,y]}=\Psi_{\chi}$; 
\item  $g_{\chi}(q)_{[x,\chi]}=g(q)_{[x, \chi]}$.
\end{itemize}
The base of the induction occurs when $\chi=y$, in which case both conditions are tautologically satisfied, taking $g_{y}=g$.
Below, we fix~$χ$ and concentrate our attention on a single embedded bordism between two consecutive cuts $\phi$ and $\chi$ of $c$, the terminal cut being~$\chi$.
That is to say, we can assume that $c$ has precisely two cuts $\phi$ and $\chi$, both of which are joint cuts.
Once the embedded bordism $\Psi$ and the zigzag of natural transformations $g$ is constructed for such~$c$,
we set $\Psi_{\phi}=\Psi\vee \Psi_{\chi}$ and replace the object $g_j(q)$ for $q∈Q$ and all~$j$
with the object $g(q)_{[x,\phi]}∨g_j(q)∨Ψ_{\chi}$,
where $∨$ is the monoidal structure on morphisms of necklaces (\cref{alternative.nec}) given by joining necklaces at their endpoints. Thus, from now on we assume $c$ has exactly two cuts $\phi≤\chi$, hence $c$ also has exactly one bead. From the definition of $\frakB_i^{i,x,y}$, this implies that $c$ is an embedded bordism in $\frakB^i_{i}$.  We truncate the data of $g_{\chi}:Q\to \frakD$ by removing all cuts in all necklaces that are less than $\phi$ or greater than $\chi$ to obtain a functor $(g_{\chi})_{[\phi,\chi]}:Q\to \frakD_{[\phi,\chi]}$, which we denote by $g$ below. We also adjust $\fraka$ accordingly (by removing $x$ and $y$ and inserting $\phi$ and $\chi$). 
\medskip

\paragraph{\bf Middle inductive layer: induction on $\Psi$ from \cref{morse.decomposition}.} Next, we will construct an embedded bordism~$Ψ$ using \cref{morse.decomposition}, taking $x=\phi$ and $y=\chi$, which will be used in the second layer of the induction.
To use the proposition, we must specify $g$, which is already provided (the functor $g=(g_{\chi})_{[\phi,\chi]}$ in the preceding paragraph), and a cover $\Omega$ of $\stu(W)_{\stu(f)(0,\rho^{-1}(0))}$.
We construct $\Omega$ as follows.
Fix a Riemannian metric on $\stu(W)_{\stu(f)(0,\rho^{-1}(0))}$. 
Choose $ε>0$ so that the $3ε$-neighborhood of the $(i-1)$-core of every bead in every necklace of $g(Q)$ is subordinate to ${\cal W}$, which is possible by finiteness of $Q$, compactness of the $(i-1)$-core, and the definition of $\frakB_{i-1}^{\fraka}$ (\cref{decompbi}).
Further reduce $\varepsilon$ so that every open subset of $\stu(W)_{\stu(f)(0,\rho^{-1}(0))}$ of diameter less that $2\varepsilon$ that intersects the core of $\fraka$ is subordinate to ${\cal W}$, which is possible because the core of $\fraka$ is compact. Now choose $Ω$ to consist of all open subsets of $\stu(W)_{\stu(f)(0,\rho^{-1}(0))}$ with diameter less than~$ε$.
Applying \cref{morse.decomposition} to $g$ and $\Omega$ yields an etale necklace $\Psi$. Let $\vert \Psi\vert$ denote the number of cuts of $\Psi$ in the $i$th direction. \cref{morse.decomposition} guarantees that $\Psi$ is in $\frakD_{[\phi,\chi]}$. Set $\varepsilon'=\varepsilon/\vert \Psi\vert$.

We now construct the zigzag of natural transformations connecting $g$ and the constant functor $\Psi:Q\to 1\to \frakD_{[\phi,\chi]}$. By induction on $j\in [\vert \Psi\vert]$, in decreasing order, we construct functors
$g_j:Q\to \frakD_{[\phi,\chi]}$
and a zigzag chain of natural transformations connecting $g_{j+1}$ to $g_j$
such that the constructed functor $g_j$, for $j<\vert \Psi\vert$,
has the following property.
\begin{itemize}
\item For all $q\in Q$, the last $\vert \Psi\vert-j$ beads of $g_j(q)$ coincide with the last $\vert \Psi\vert-j$ beads of~$\Psi$.
\item The number of the remaining beads of $g_j(q)$ is equal to the number of beads in $g_{\vert \Psi\vert}(q)$.
Furthermore, the $i$-core of every remaining bead is empty and the $(i-1)$-core of every remaining bead of $g_j(q)$ is contained in the $ε'(\vert \Psi\vert-j)$-neighborhood of the $(i-1)$-core of the corresponding bead in $g_{\vert \Psi\vert-1}(q)$. Only the last bead of $g_j(q)$ has nonempty $i$-core. 
\end{itemize}
For the base of the induction, we set $g_{\vert \Psi\vert}=g$.
At the end of the induction, the last $\vert \Psi\vert$ beads of $g_{0}(q)$ coincide with the beads of~$Ψ$,
meaning the remaining beads of~$g_{0}(q)$ must be all degenerate, with their cuts equal to the initial cut~$\phi$.
A simple argument at the end of the proof shows how to connect $g_{0}$ to the constant functor on~$Ψ$.

For the inductive step, suppose $g_{j+1}$ has been constructed.
We will define the zigzag chain from $g_{j+1}$ to $g_j$ as follows.
In the proof below, we will talk about inserting and removing cuts into the necklaces given by the target of objects in $\frakD_{[\phi,\chi]}$, using \cref{move.cut}.
Cuts will always be assumed to be in the $i$th direction.

\medskip
\paragraph{\bf Inner inductive layer: induction on cuts of $g(Q)$} Consider the set of all cuts~$G$ in all etale necklaces in $g_{j+1}(Q)$ such that $G≤Ψ_{j+1}$, i.e., excluding the cuts $Ψ_{j'}$ for $j'>j+1$.
Extend the natural partial order $<$ on this set (\cref{cut}) to a total order~$<_{T}$.
We identify the resulting finite nonempty totally ordered set with some $[K]\in \Delta$ and denote the cut corresponding to an element $k\in [K]$ by~$D_k$.
We construct functors $g_{j,k}$ for all $k\in [K]$ together with a double zigzag of natural transformations connecting $g_{j,k-1}$ and $g_{j,k}$,
by induction on $k\in [K]$ in increasing order.
For the base of the induction, set $g_{j,0}=g_{j+1}$.
Suppose $g_{j,k'}$ and the double zigzags connecting $g_{j,k'-1}$ to $g_{j,k'}$ have been constructed for all $k'<_{T}k$.
The double zigzag connecting $g_{j,k-1}$ to $g_{j,k}$ is constructed as follows.

We apply \cref{gluing.lemma} to the cuts $D_k$ and $\Psi_j$ to obtain a new cut $D'_k$ such that $D'_{k}\leq D_k$ and $D'_{k}\leq  \Psi_j$. 
We take the neighborhood $M$ in \cref{gluing.lemma} to be
$$M=\bigcap_{\{k'\mid D_{k'}<D_k\}}(D'_{k'})_{>0}\cap Z_{\varepsilon'}\subset \stu(W)_{\stu(f)(0,\rho^{-1}(0))},$$ 
where  $Z_{\varepsilon'}$ is the $\varepsilon'$-neighborhood of $(D_k)_{=0}$.
This ensures that if $k'<_Tk$ and $D_{k'}<D_k$, then $D_{k'}'<D_k'$.

Now if $D'_k≠D_k$ (hence $D'_k<D_k$),
we invoke \cref{move.cut} to construct a double zigzag that replaces $D_k$ by~$D'_k$ in all etale necklaces in~$g_{j,k-1}(Q)$.
First, we add $D'_k$ using Part~(1), taking $A=D'_k$, $B=D_k$.
Second, we remove~$D_k$ using Part~(2),
unless $D_k=Ψ_{j+1}$, in which case we keep the last copy of $D_k$ in every necklace (using the optional $\hat g^2(q)$ in \cref{move.cut}).

To apply Part~(1) of \cref{move.cut}, we must show that the necklace $g^2(q)$ obtained by inserting $D'_k$ into $g_{j,k-1}(q)$ is a legitimate object in $\frakD_{[\phi,\chi]}$. The result of the insertion is an embedded  bordism $B$ by \cref{gluing.lemma}. To show that the resulting bordism $B$ is globular, suppose we extract a vertex in the $l$th direction from $B$. The resulting bordism is degenerate in directions $i+1,\hdots,d$ because the original bordism is an $i$-bordism (since its $i$-core is empty, by the second bullet point above). The resulting bordism is degenerate in directions $l+1,\hdots,i$, since the original bordism was globular. 

Finally, the $(i-1)$-cores are not changed by the insertion of non-joint cuts, which implies that the insertion yields a necklace in $\frakD_{[\phi,\chi]}$. The $i$-core is empty after the first stage of the middle induction (that is, after inserting $\Psi_{\vert\Psi\vert-1}$), since the $i$-core is concentrated in the last bead of $\Psi$.

To apply Part~(2) of \cref{move.cut}, we must show that the necklace $\hat g^1(q)$
obtained by changing the joint cuts equal to~$D_k$ in the target of $g^2(q)$ to non-joint cuts yields an object in $\frakD_{[\phi,\chi]}$.
The $(i-1)$-core~$c_1$ of the new bead in $\hat g^1(q)$
is contained in the union of the $(i-1)$-cores of the two beads that were joined together to form the new bead,
namely, the $(i-1)$-core~$c_2$ of the bead in $g^2(q)$ with initial cut~$D_k$
and the $(i-1)$-core~$c_3$ of the bead in $g^2(q)$ with initial cut $D'_k$ and terminal cut $D_k$.
By definition of the open neighborhood $M$ above and the construction of \cref{gluing.lemma}, the $(i-1)$-core~$c_3$ is in the $ε'$-neighborhood of the $(i-1)$-core~$c_4$ of the bead in~$Ψ$ with initial cut~$Ψ_j$ and terminal cut~$Ψ_{j+1}$. 
By definition of~$Ω$, the diameter of $c_4$ is less than~$ε$.
Therefore, the diameter of~$c_3$ is less than $ε+ε'≤2ε$.
Hence, $c_4$ is subordinate to ${\cal W}$.

There are now two cases, which we handle separately.
If $c_2$ and $c_3$ are disjoint, then $c_1$ is the disjoint union of the connected components of $c_2$ and $c_3$.
By assumption, every connected component of $c_2$ is subordinate to ${\cal W}$.
As shown above, $c_4$ is also subordinate to ${\cal W}$.
Therefore, every connected component of $c_1$ is subbordinate to ${\cal W}$ and hence $\hat g^1(q)\in \frakD_{[\phi,\chi]}$. 

In the second case, the intersection of $c_2$ and $c_3$ is nonempty. 
Therefore, $c_1$ is contained in the $ε'+ε$-neighborhood of~$c_2$ (the $(i-1)$-core of the bead in $g^2(q)$ with initial cut~$D_k$).
By inductive assumption in the second bullet point, the $(i-1)$-core~$c_2$ is contained in the $ε'(|Ψ|-j)$-neighborhood
of the $(i-1)$-core~$c_5$ of the corresponding bead in $g_{|Ψ|-1}(q)$.
Combined together, we see that the $(i-1)$-core~$c_1$ of the new bead in $\hat g^1(q)$ is contained in the $ε'+ε+ε'(|Ψ|-j)<3\varepsilon$-neighborhood
of the $(i-1)$-core~$c_5$ of the corresponding bead in $g_{|Ψ|-1}(q)$.
By construction of~$ε$, the latter neighborhood is subordinate to~$\mathcal{W}$.
This shows that each bead in the target of $\hat g^1(q)$ is in $\frakB_{i-1}^i$, hence $\hat g^1(q)\in \frakD_{[\phi,\chi]}$.

Next, we have to show that after the inner inductive layer is complete, the inductive assumption in the second bullet point continues to hold.
That is to say, the $(i-1)$-core of every bead of $g_j(q)$ except for the last bead is contained in the $ε'(\vert \Psi\vert-j)$-neighborhood of the $(i-1)$-core of the corresponding bead in $g_{\vert \Psi\vert-1}(q)$.
Consider a bead of $g_{j+1}(q)$ with initial and terminal cuts $D_k$ and $D_l$ that is not the last bead of $g_{j+1}(q)$.
The $(i-1)$-core $d_1$ of this bead is contained in the $ε'(\vert \Psi\vert-(j+1))$-neighborhood of the $(i-1)$-core~$d_0$ of the corresponding 
bead in $g_{\vert \Psi\vert-1}(q)$.
After the inner inductive layer is complete, the newly modified bead has initial and terminal cuts $D'_k$ and $D'_l$ as constructed in the inductive step.
By \cref{gluing.lemma}, the $(i-1)$-core $d_2$ of the new bead is contained in the $ε'$-neighborhood of~$d_1$
and therefore in the $ε'+ε'(\vert \Psi\vert-(j+1))=ε'(\vert \Psi\vert-j)$-neighborhood of~$d_0$,
which proves the inductive assumption for~$j$. The last natural transformation in the zigzag chain constructed so far is $g_0$.

\medskip
\paragraph{\bf Removing leftover degenerate cuts.}
The underlying necklaces of all objects in the image of $g_0$ are identical to~$Ψ$ except that the 0th cut (i.e., $\phi$) can be repeated several times, marked as a joint or non-joint cut.
Apply Part~(2) of \cref{move.cut} (taking $B=\phi$) to $g_0(q)$ for all $q\in Q$ to remove the inner cuts equal to~$\phi$. The resulting etale necklace is equal to $\Psi$ by construction. 
This completes the construction of the zigzag chain.

\end{proof}

\subsection{Butchering bordisms}

In this section, we prove the two remaining propositions used in \cref{contractible.cuts}. \cref{morse.decomposition} provides a way to cut bordisms into very small pieces.
The cutting is done in a completely globular way and only uses elementary Morse theory.
In particular, we do not need Morse theory for manifolds with corners. \cref{gluing.lemma} provides a way to cut and glue two transversally compatible cuts.

\begin{figure}[ht]
\begin{center}
\begin{tikzpicture}
\draw[thick, blue] (0,2) -- (0,-2);
\draw[purple] (-2,-1.8) -- (2,-1.8);
\draw[purple] (-2,1.8) -- (2,1.8);
\draw[thick, blue] (0,1.8) to [out=-90, in=90] (2,0) to[out=-90,in=90] (0,-1.8);
\foreach \x in {1,...,7}
\draw[orange,dashed] (2*\x*.1428,2) -- (2*\x*.1428,-2);
\draw (2*3*.1428,2) -- (2*3*.1428,1) to [out=-90,in=90] (2*3*.1428+.1428,.5) to [out=-90,in=90] (2*3*.1428,0) -- (2*3*.1428,-2);
\draw (2*3*.1428,2) -- (2*3*.1428,-2);
\draw (2*3*.1428,2) -- (2*3*.1428,1.3) to [out=-90,in=90] (2*3*.1428+.2,1) to [out=-90,in=90] (2*3*.1428+.1482,.5) to [out=-90,in=90] (2*3*.1428+.2,0) to [out=-90,in=90] (2*3*.1428,-.5);
\draw[->] (3,0) -- (5,0);
\node at (4,.5) {\cref{gluing.lemma}};
\foreach \y in {1,...,7}
\node at (2*\y*.1428,-2.3) {$\y$};
\node at (0,-2.3) {$0$};
\end{tikzpicture}
\begin{tikzpicture}
\draw[thick, blue] (0,2) -- (0,-2);
\draw[purple] (-2,-1.8) -- (2,-1.8);
\draw[purple] (-2,1.8) -- (2,1.8);
\draw[thick, blue] (0,1.8) to [out=-90, in=90] (2,0) to [out=-90,in=90] (0,-1.8);
\foreach \x in {1,...,7}
\draw[orange,dashed] (2*\x*.1428,2) -- (2*\x*.1428,-2);
\draw[blue] (2*3*.1428,.95) to [out=-45,in=90] (2*3*.1428+.1428,.5) to [out=-90,in=90] (2*3*.1428,0) -- (2*3*.1428,-.65) to [out=-80,in=30,looseness=1.2] (2*3*.1428,-.95);
\draw[blue] (2*3*.1428+.2,.88) to[out=-45,in=90] (2*3*.1428+.1428,.5) to [out=-90,in=90] (2*3*.1428+.2,0) to [out=-90,in=90] (2*3*.1428,-.5);
\foreach \y in {1,...,7}
\node at (2*\y*.1428,-2.3) {$\y$};
\node at (0,-2.3) {$0$};
\end{tikzpicture}
\end{center}
\caption{An illustration of \cref{morse.decomposition} to modify the cuts $\Phi_{2j+1,k,b}$ for $k=3$ and $b=1,2,3$ to obtain new cuts $\Psi_{2j+1,3,b}$.
The orange dashed vertical cuts are the preimages $g^{-1}(s_k)$ for $k\in \{1,\ldots,7\}$.
The black cuts around line 3 are the cuts $\Phi_{2j+1,3,b}$ for $b=1,2,3$.
Each blob between the black cuts (excluding the ``tails'') is contained in some element of the cover $\Omega$.
On the right, a surgery is performed to cut and glue the black cuts to the target cut $y$, depicted as the bold solid blue curve. }
\label{morse.cuts}
\end{figure}

\begin{remark}\label{critpntsbord}
In the next proposition, we construct the $\Psi$ cuts using a locally defined Morse flow in an open neigborhood $N'$ of the core with compact closure. For the sake of clarity, we observe the following.
\begin{itemize}
\item If $i<d$, then $\core(\frakb)_0$ is an $i$-bordism.
By perturbing $g$, we can push all critical points out of the core and into $N'$.
Then we can take $N''$ small enough so that $g$ has no critical points in $N''$. 
\item When $i=d$, then $\core(\frakb)_0$ is an $d$-bordism.
In this case, we cannot push critical points out of the core, since $g$ has to be compatible with the source cut and the core has the same dimension as $N'$.
However, we can still shrink $N'$ so that all critical points of~$g$ are in the interior of the core and there is only finitely many of them. 
\item By perturbing $g$, we can also assume that all critical values of $g$ on $N'$ are distinct.
\end{itemize}
We emphasize that all that is being used in \cref{morse.decomposition} is a locally defined Morse flow and the Morse lemma. In particular, there is no globally defined Morse flow since $N'$ is an open manifold.
\end{remark}

\begin{proposition}
\label{morse.decomposition}
Assume \cref{etale.bordism}. Let $\phi<\chi$ be successive cuts of $c$. Suppose $g:Q\to \frakD_{[\phi,\chi]}$ and $Ω$ is an open cover of~$\stu(W)_{\stu(f)(0,\rho^{-1}(0))}$ that is subordinate to ${\cal W}$.
Then there is an object $Ψ\in \frakD_{[\phi,\chi]}$ such that the following hold.
\begin{enumerate}
\item All cuts of $Ψ$ in the $i$th direction are joint cuts.
\item Every bead of $\Psi$ (except for the last bead in globular case) has an empty $i$-core and its $(i-1)$-core is subordinate to $\Omega$. 
\item In the globular case, every connected component of the $(i-1)$-core of the last bead is subordinate to ${\cal W}$.
\item In the globular case, we denote by $\Xi$ the joint cut of $\Psi$ preceding $y$ in the last bead. Then every cut in every necklace of $g(Q)$ that is not equal to $y$ precedes $\Xi$.
\item The cuts of $\Psi$ in the $i$th direction are transversal to every cut in $g(Q)$ in every direction $l\neq i$ and are transversally compatible (\cref{openntrans}) with every cut in $g(Q)$ in the $i$th direction.
\end{enumerate}
\end{proposition}

\begin{proof}
In the globular case, we construct the cut $\Xi$ that cuts off the $i$-core of $\fraka_{[\phi,\chi]}$ (\cref{subcov}). In the uple case, we skip the construction of $\Xi$ and set $\frakb=\fraka_{[\phi,\chi]}$ below.

Use the Whitney extension theorem to construct a smooth map $\zeta:\stu(N)→[0,∞)$ (\cref{icore}), whose zero set coincides with~$(\core_i\fraka_{[\phi,\chi]})_0$. For a generic $\eta>0$ the preimage $\zeta^{-1}(\eta)$ will be a cut on~$\fraka_{[\phi,\chi]}$
that is transversal to all cuts of $\fraka_{[\phi,\chi]}$ and all cuts in all etale necklaces in $g(Q)$.
For a sufficiently small $\eta>0$,
every connected component of $\zeta^{-1}[0,\eta]$ is compact and subordinate to~${\cal W}$ because every connected component of $(\core_i\fraka_{[\phi,\chi]})_0$
is compact and subordinate to ${\cal W}$. By globularity, there is a sufficiently small $\eta>0$ such that
the preimage $\zeta^{-1}(\eta)$ is disjoint from all cuts of~$\fraka_{[\phi,\chi]}$ and all cuts of all necklaces in~$g(Q)$, in the first $i$ directions. Set 
$$\Xi=(\zeta^{-1}(\eta,∞)\cap \chi_{<0},\zeta^{-1}(\eta)\cup \chi_{=0},\zeta^{-1}[0,\eta)\cup \chi_{>0}).$$
By construction of $\eta$, every connected component of the $(i-1)$-core of the resulting bead between $\Xi$ and~$\chi$ is subordinate to~${\cal W}$.

Denote by $\frakb$ the etale bordism $\fraka_{[\phi,\Xi]}$.
In particular, $\frakb$ is an embedded bordism whose $i$-core is empty.
By construction, in the globular case, an arbitrary cut $C\neq \chi$ in the $i$th direction of an arbitrary etale necklace~$g(q)$ ($q\in Q$) satisfies $C\leq \Xi\leq \chi$.
Thus, we can consider $g(q)$ as an object of $\frakD_{[\phi,\Xi]}$. Below we work with $\frakb$ and construct an object $\Psi$ in $\frakD_{[\phi,\Xi]}$. We promote $\Psi$ to an object in $\frakD_{[\phi,\chi]}$ by inserting $\chi$.


Recall that by the optional Condition (1) of \cref{bordstalkisot}, $N\to V\times L$ (\cref{etale.bordism}) can be chosen to be a trivial bundle.
Hence, we can extend the construction of cuts in a single fiber to the germ of the fiber (by taking them to be constant in the family direction).
If the cuts satisfy Property (5) (transversality property) in a single fiber, they also satisfy Property~(5) for a small family, since the core is compact and transversality is a generic condition.
Below we work in a single fiber $\stu(N)_{\stu(f)(0,\rho^{-1}(0))}$, which we denote simply by $N$, and we implicitly extend cuts to the germ as indicated above.

Fix some Morse function $g$ defined on $N$ that vanishes on the source cut $\phi\in \frakb$ and whose critical points are disjoint from all cuts in $\frakb$.
By compactness of the core, there is a neighborhood $N'\subset N$ of $\core(\frakb)_0$ such that only finitely many critical points of $g$ are in $N'$, $g$ is bounded on~$N'$, all critical points of $g$ are contained in the interior of the core, and all critical values of $g$ are distinct (see \cref{critpntsbord}). 
Since we take the germ of the core, we can replace $N$ by $N'$ and restrict all data to $N'$ accordingly.

When $i=d$, there may be critical points contained in the interior of the core. In this case, order the critical points of $g$ according to corresponding critical values, in increasing order.
For all critical points~$c$ of $g$, the construction in Milnor \cite[Theorem 3.2]{Milnor} modifies $g$ to obtain a Morse function $G$ that coincides with $g$ everywhere, except for an arbitrarily small neighborhood of~$c$.
The neighborhoods can be chosen to be disjoint since there are only finitely many critical points.
Choose an $\epsilon>0$ such that for all $j$ indexing the critical points of~$g$,
the cut tuple with endcuts $\Psi_{2j}=g^{-1}(g(c_j)+\epsilon)$ and $\Psi_{2j-1}=G^{-1}(g(c_j)+\epsilon)$ satisfies the following property.
\begin{itemize}
\item The closure of $\Psi_{(2j-1,2j)}$ is contained in $\omega$, for some $\omega\in \Omega$, and is also contained in the interior of the core.
\item The cuts $\Psi_{2j}$ and $\Psi_{2j-1}$ intersect all cuts in $\frakb$ and $g(Q)$ in the directions other than $i$ transversally.
\end{itemize}
Since $\Omega$ is subordinate to ${\cal W}$, this implies that the bead with endcuts $\Psi_{2j-1}$ and $\Psi_{2j}$ is in $\frakB_{i-1}$
and its $(i-1)$-core, which is a subset of $\overline{\Psi}_{(2j-1,2j)}$, is subordinate to~$Ω$.

We set $\Psi_0=\phi$ and $\Psi_{2J+1}=\chi$, where $J$ indexes the last critical point. We now concentrate our attention on the bead with joint cuts $\Psi_{2j}$ and $\Psi_{2j+1}$, where $0\leq j\leq J$. Let $s=1+\sup_{x\in N}g(x)$ (see \cref{morse.cuts} for an illustration).   
Partition $[g(c_j)+\epsilon,s]$ into intervals $I_{k}=[s_{k-1},s_k]$, $k=1,\ldots, N$, where $g(c_j)+\epsilon=s_0\leq s_1\leq \cdots \leq s_N=s$ (with $\epsilon>0$ defined above),
so that for any $x\in g^{-1}[s_{k-1},s_k]\cap (\core \frakb)_0$
the Morse flow line $\varphi_t(x)$ is defined for all $t\in [s_{k-1}-g(x),s_k-g(x)]$ and is subordinate to the cover $\Omega$. 
We cut $g^{-1}[g(c_j)+\epsilon,s]$ into thin strips by the cuts $g^{-1}(s_k)$.
We can assume that for each $1\leq k\leq N$, the cut $g^{-1}(s_k)$ intersects all cuts in $\frakb$ and $g(Q)$ in the directions other than $i$ transversally.
If not, we perturb the $s_k$ for $1\leq k\leq N$ until the cut $g^{-1}(s_k)$ is transversal. 

We now concentrate our attention on a single strip $g^{-1}([s_{k-1},s_k])$, $1\leq k\leq N$.
For any $\omega\in \Omega$, we define $$\omega'=\{x\in g^{-1}(s_{k-1})\mid \varphi_t(x)\in \omega,\quad t\in [0,s_{k}-s_{k-1}]\},$$
where $\varphi_t$ is the Morse flow.
Then $\{\omega'\}$ covers $(\core\frakb)_0∩g^{-1}(s_{k-1})$, which is compact.
Choose a finite subcover $\{X_b\}_{1\leq b\leq B_k}$ of $\{\omega'\}$. 
Let $X=\bigcup_{b=1}^{B}X_b$. Choose a partition of unity $\{\psi_b:X\to \RR\}$ subordinate to the cover $\{X_b\}$.
Rescale $\psi_b$ so that their sum is $s_k-s_{k-1}$.
Define $\Phi_{2j,k,b}$ to be the the diffeomorphic image of $X$ under the map 
$$x\mapsto \varphi_{(\psi_1(x)+\cdots +\psi_{b}(x))}(x).$$
Here $\Phi_{2j,k,b}$ is almost a cut, except that it need not extend across the entire ambient manifold~$N$,
but merely across some smaller open neighborhood~$N'$ of the core.
After shrinking $N$ further to~$N'$, we can assume $\Phi_{2j,k,b}$ to be a cut.
After perturbing $ψ_b$ generically,
we can assume that for each $b\in \{1,\ldots,B_k\}$,
the cut~$\Phi_{2j,k,b}$ intersects all cuts in~$g(Q)$ transversally, in all directions.

Now for each $k∈\{1,…,N\}$ and $b\in \{1,\ldots,B_k\}$,
we apply \cref{gluing.lemma} to the cuts $\Phi_{2j,k,b}$ and $Ψ_{2j+1}$, producing a cut $Ψ_{2j,k,b}$,
which is transversal to all cuts in $\frakb$ and all cuts of $\frakb$ and all cuts of all etale necklaces in~$g(Q)$ in directions $l\neq i$ and transversally compatible to all cuts in direction~$i$.
When invoking \cref{gluing.lemma}, we choose the open neighborhood $M$ in \cref{gluing.lemma}
small enough so that $Ψ_{2j,k,b'}<Ψ_{2j,k,b}$, for all $b'<b$, and the $(i-1)$-core of the resulting bead is subordinate to $\Omega$.
The cuts 
$$Ψ_{2j},Ψ_{2j,1,1},…,Ψ_{2j,1,B_1},…,Ψ_{2j,N,1},…,Ψ_{2j,N,B_N}=Ψ_{2j+1}$$ 
form a cut tuple by definition.
Concatenating all such cut tuples for different~$j$ yields the desired cut tuple~$Ψ$.

Declare all cuts in $\Psi$ to be joint cuts. To see that this defines a necklace in $\frakB_{i-1}^{\frakb}$, observe that the cuts of $\Psi$  are transversal to all the cuts of $\frakb$ and $g(Q)$ in all directions $l\neq i$. By construction, each bead is in $\frakB_{i-1}$, since the $(i-1)$-core of every bead of $\Psi$ is subordinate to $\Omega$, and therefore to ${\cal W}$.
\end{proof}

\begin{definition}
\label{openntrans}
Recall \cref{etale.bordism}, including the etale map $h:N\to  W\times_{U}V\times L$,
and the projection map $N→V⨯L$, which is an object in $\frakFEmb_d$.
Suppose $C$ and $C'$ are representatives of germs of cuts on $N→V⨯L$.
We say that $C$ and $C'$ are \emph{transversally compatible} if, after passing to some smaller representatives of the germs,
there is an open subobject $N'\subset N$ such that 
\begin{itemize}
\item $C=C'$ on the complement of $\stu(N')$ in $\stu(N)$. 
\item The cuts $C$ and $C'$ intersect (fiberwise) transversally in $\stu(N')$. 
\end{itemize}
\end{definition}

\begin{figure}[ht]
\begin{tikzpicture}[scale=-1]
\draw (-2,0) .. controls (0,0) and (1,-1) .. (2.5,-1);
\draw (1,-2)-- (-1,2);
\draw[very thick] (-1,2) .. controls (0,-.1) and (.6,-.9) .. (2.5,-1);
\draw[dashed] (.2,-.4)  circle (40pt);
\node at (.7,.4) {$M$};
\node at (1.5,1.5) {$C_{\leq 0}\cap C'_{\leq 0}$};
\end{tikzpicture}
\caption{Gluing together two transversal cuts. The two transversally intersecting thin lines are the cuts $C_{=0}$ and $C'_{=0}$. The construction $(C\wedge C')_{=0}$ joins the two cuts together. The region lying above the thin horizontal cut $C_{=0}$ is $C_{\leq 0}$. The region lying to the right of the straight vertical thin cut $C'_{=0}$ is $C'_{\leq 0}$. The region lying above and to the right of the thick black curved cut $(C\wedge C')_{=0}$ is $(C\wedge C)_{\leq 0}$. The dashed circle is the open neighborhood $M$. }
\end{figure}

\begin{proposition}
\label{gluing.lemma}
Assume the notation of \cref{openntrans}.
Suppose $C$ and $C'$ are transversally compatible.
Let $N'$ be an open neighborhood satisfying the conditions of \cref{openntrans}.
Then for any open neighborhood $M\subset \stu(N')$ of the intersection $C_{=0}\cap C'_{=0}\cap \stu(N')$,
there is a (nonunique) germ of cuts $C'\wedge C$ on $N→V⨯L$ having the following properties:
\begin{enumerate}
\item Outside of $M$, $$(C'\wedge C)_{=0}=(C'_{=0}\cap C_{\leq 0})\cup (C_{=0}\cap C'_{\leq 0}).$$  
Moreover, $$(C'\wedge C)_{\leq 0}\subset C_{\leq 0}\cap C'_{\leq 0}.$$  
\item $(C'\wedge C)$ is compatible with both $C$ and $C'$, in the sense of \cref{cut.m.tuple}. 
\end{enumerate}
\end{proposition}

\begin{proof}
Fix maps $h,k:\stu(N)\to \RR$ defining the $L$-families of cuts $C$ and $C'$. Restrict both maps to $\stu(N')$ and consider the map
$$(h,k):\stu(N') \to \RR^2.$$
Then 
$$(h,k)^{-1}(0)=C_{=0}\cap C'_{=0},$$
and the preimages of the two axes are $C_{=0}$ and $C'_{=0}$. Since the cuts are transversal, the map $(h,k):\stu(N')\to \RR^2$ is a submersion. Choose an open ball of radius $\delta$ around $0\in \RR^2$. 
Fix a fiberwise embedding $\iota:\RR\times V\times L\to \RR^2$ (fiberwise over $V\times L$)  that satisfies the following properties: 
\begin{enumerate} 
\item $\iota(s,x,t)=(0,-s)$ for all $(s,x,t)\in \RR\times V\times L$, with $s<-1$;
\item $\iota(s,x,t)=(s,0)$, for all $(s,x,t)\in \RR\times V\times L$, with $s>1$;
\item $\iota(s,x,t)\in B_{\delta}$, where $B_{\delta}$ is the open ball of radius $\delta$ in $\RR^2$, for all $(s,x,t)\in \RR\times V\times L$, with $-1<s<1$. 
\end{enumerate}
Then $(C'\wedge C)_{=0}≔(h,k)^{-1}(\iota(\RR\times U))$ gives rise to the cut with the desired properties. 
Since a tubular neighborhood of $h^{-1}(0)\cap k^{-1}(0)$ is given by $(h,k)^{-1}(B_{\delta})$, we can get an arbitrarily small radius by letting $\delta\to 0$. Transversality of the cuts ensures that $(C'\wedge C)_{\leq 0}\subset C_{\leq 0}\cap C'_{\geq 0}$ holds in $\stu(N')$. Finally, since $C=C'$ on the complement of $\stu(N')$ in $\stu(N)$, the cuts $C'\wedge C$ extends to the whole of $\stu(N)$, by defining $C'\wedge C=C'=C$ on $\stu(N)\setminus \stu(N')$. 
\end{proof}

\section{Applications}
\label{applications}

\subsection{Shapes and concordances}
\label{shapebordism}

In this section, we recall the notion of the \emph{shape} of a simplicial presheaf on cartesian spaces (or smooth manifolds). 
We first recall the notion of a \emph{cohesive} $\infty$-topos of Schreiber \cite{Schreiber}. 
These are certain toposes which participate in a quadruple adjunction $(\csp  \dashv \Delta \dashv {\it \Gamma} \dashv \nabla)$, with $\Delta$ and $\nabla$ fully faithful and $\csp$ preserving products. 
In our context the adjunction takes the form of a quadruple Quillen adjunction
$$\xymatrix{\csp, {\it \Gamma} \colon \sPSh(\cart)_{\proj,\Cech}\ar@<.225cm>[r]\ar@<-.075cm>[r] & \sset \colon \Delta,\nabla,\ar@<-.075cm>[l]\ar@<.225cm>[l]}$$
with left adjoints depicted above their corresponding right adjoints. 
This adjoint quadruple is induced by an adjunction between $\cart$ and the terminal category $1$. 
The right adjoint functor $1\to \cart$ picks $\RR^0$.
The further left and right adjoints are given by Kan extensions. 
Thus, the functor ${\it \Gamma}$ evaluates at $\RR^0$ and the functor $\Delta$ takes the locally constant stack associated to a simplicial set.
Using the fact that the inclusion $i:\cart\into \Man$ is the inclusion of a full $\infty$-dense subsite,
one shows that the restriction functor $i^*:\sPSh(\Man)\to \sPSh(\cart)$ is part of a Quillen equivalence
(in fact, two Quillen equivalences, corresponding to the injective and projective model structures).

\begin{remark}
The constant presheaf functor $\sset→\sPSh(\Man)$ does not extend to a quadruple Quillen adjunction when $\sPSh(\Man)$ is equipped with local weak equivalences.
However, one can compute explicitly the shape functor via derived functors
by first applying the restriction functor $i^*:\sPSh(\Man)\to \sPSh(\cart)$
and then the left derived functor of $\csp$ in the projective model structure.
However, a more preferable formula will be presented below. 
\end{remark}

These quadruple adjoint functors give rise to idempotent monads 
\begin{equation}\label{monads}\smallint≔\Delta\csp\dashv \flat≔\Delta{\it \Gamma}\dashv \sharp≔\nabla{\it \Gamma},\end{equation}
each of which reflects a different nature of a smooth stack (see Schreiber \cite{Schreiber} for detailed discussion). 
We will focus on one of these idempotent functors, namely the leftmost adjoint $\smallint ≔\Delta\csp$, which is also known as the \emph{shape} functor.
Abusing terminology, we also call the left adjoint $\csp$ the \emph{shape} functor.
We call the value of $\csp$ on a sheaf $X\in \sPSh(\Man)_\Cech$, the \emph{shape of} $X$. 

\begin{remark}
Berwick-Evans, Boavida de Brito, and the second author \cite{BEBdBP}
(see also \cite{Pavlov.Diffeo} for the structured case)
established an explicit formula for the shape functor~$\smallint$ in \eqref{monads}.
More precisely, \cite[Theorem~1.1]{BEBdBP} proves that if $X\in \sPSh(\Man)_\Cech$,  
then we have an equivalence 
\begin{equation}\label{dimaform}\smallint X\simeq \hocolim_{[n]\in \Delta^\op}\hom(\gsim^n,X),\end{equation}
where the homotopy colimit is taken in the category of presheaves
and $\gsim^n$ is the smooth $n$-simplex, viewed as a smooth stack via its diffeological space of plots.
Although it is fairly easy to see that the colimit on the right is invariant under concordance, it is highly nontrivial to show that it satisfies descent.
See \cite{Pavlov.Diffeo} for a version that treats sheaves valued in algebraic ∞-categories, like $Γ$-spaces. 
\end{remark}

\begin{remark}
The formula for $\smallint$ also gives rise to a formula for $\csp$.
Denote by $\map(-,-)$ the simplicial enriched hom for simplicial presheaves on~$\cart$
and by $\hom(-,-)$ the internal hom.
We have 
$$\hocolim_{[n]\in \Delta^\op}\map(\gsim^n,X)
\simeq {\it \Gamma} \hocolim_{[n]\in \Delta^\op}\hom(\gsim^n,X)
\simeq {\it \Gamma} \smallint (X)
\simeq {\it \Gamma} \Delta\csp (X)
\simeq \csp(X),$$
where the last equivalence follows from the fact that $\Delta$ is homotopically fully-faithful
and the (derived) counit $\epsilon:{\it \Gamma} \Delta\lto3{\simeq} \id$ is an equivalence.
In particular, if $X$ is a smooth manifold, viewed as a representable smooth stack, then 
$$\csp(X)\simeq \sing(X),$$
where on the right we have the singular simplicial set of the underlying topological space.
Taking $\pi_0$ of both sides of the above identifies concordance classes with the topological connected components of the underlying space.
\end{remark}

Formula \eqref{dimaform} has some striking consequences. 
For example, if $\simcon$ denotes the equivalence relation given by concordance, then for any smooth manifold $M$, we have the formula
$$X(M)/{\simcon}\cong [M,\hocolim_{[n]\in \Delta^\op}X(\gsim^n)],$$
where on the right we take homotopy classes of maps and on the left we mod out by the relation of concordance.
This immediately gives a classifying space construction for~$X$.
We illustrate with the following example.

\begin{example}
Let $X=\Vect$ be the sheaf of groupoids\footnote{The functor which sends a smooth manifold $M$ to the groupoid of vector bundles $\Vect(M)$ is only a pseudofunctor.
Hence, here we mean its strictification to an honest functor.} given by vector bundles with isomorphisms between them.
Because every vector bundle is locally trivial, we have an equivalence $\Vect(\gsim^n)\simeq\deloop\GL(\gsim^n)$
on every smooth $n$-simplex, where $\deloop\GL(S)=\coprod_{d≥0}\deloop\GL_d(S)$ is the disjoint union of deloopings of smooth general linear groups,
with $\GL_d(S)=\sm(S,\GL(\RR^d))$.
The classifying space of~$\Vect$ can be computed (for a fixed summand in dimension $d≥0$) as 
\begin{align*}
\csp \Vect_d  &\simeq \hocolim_{[n]\in \Delta^\op}\Vect_d(\gsim^n)
\simeq\hocolim_{[n]\in \Delta^\op}\deloop\GL_d(\gsim^n)
\simeq\hocolim_{[n]\in \Delta^\op}\hocolim_{[m]\in \Delta^\op} \sm(\gsim^n,\GL_d^m)\\
&\simeq\hocolim_{[m]\in \Delta^\op}\hocolim_{[n]\in \Delta^\op} \sm(\gsim^n,\GL_d)^m
\simeq\hocolim_{[m]\in \Delta^\op} \sing(\GL_d)^m
\simeq\tdeloop\sing(\GL_d),\\
\end{align*}
where we have used the fact that $\Delta^{\rm op}$ is sifted to distribute over the product. Here $\tdeloop\sing(\GL_d)$ denotes the classifying space of the simplicial group $\sing(\GL_d)$,
or, equivalently, the singular simplicial set of the topological classifying space of the topological group~$\GL_d$.
Thus, we recover the classifying space construction
$$\Vect(M)/{\simcon}\cong [M,\coprod_{d≥0}\tdeloop\GL_d].$$
\end{example}

\subsection{Smooth field theories and their classifying spaces}
\label{classifying.space}

In this section, we prove \cref{classcnc}.
This provides an affirmative answer to a long-standing conjecture of Stolz and Teichner.
We begin by defining field theories.
Recall the conventions and notations of \cref{presheaf.notation} and \cref{globular.hom}.

\begin{remark}
\label{two.variants}
The results in the section apply equally well to both bordism categories, $\Bord_d$ and $\BBord_d$,
since they only rely on the codescent property.
We formulate them for $\Bord_d$ only, leaving the other case implicit.
\end{remark}

\begin{definition}
Fix $d\geq 0$, $\gs\in \sPSh(\FEmb_d)$, and $T\in \smcat_{∞,d}$.
The smooth symmetric monoidal $(\infty,d)$-category of $d$-dimensional field theories, with geometric structure $\gs$ and target $T$, is the globular functor object (\cref{globular.hom})
$$\FFT_{d,T}(\gs)≔\Funmonglob(\Bord_d^\gs,T).$$
The assignment $\gs\mapsto \FFT_{d,T}(\cal S)$ defines a functor $\sPSh(\FEmb_d)\to \smcat_{\infty,d}$.
We will mostly work with the invertible part (i.e., the core) of $\FFT_{d,T}(\gs)$,
which we denote by 
$$\FFT^{\times}_{d,T}(\gs)≔(\Funmonglob(\Bord_d^\gs,T))^⨯,$$
where the superscript~$⨯$ denotes the core of a smooth symmetric monoidal $(∞,d)$-category (defined by evaluating its fibrant replacement at ${\bf0}∈Δ^{⨯d}$, which amounts to discarding noninvertible $k$-morphisms for $1≤k≤d$).
Thus, we have a functor $$\FFT^⨯_{d,T}:\sPSh(\FEmb_d)\to \sPSh(\stcart⨯Γ).$$
A $d$-dimensional \emph{functorial field theory} with geometric structure $\gs$ and target~$T$ is a vertex in the simplicial set
$$\FFT^⨯_{d,T}(\gs)(\RR^0,⟨1⟩).$$
\end{definition}

By \cref{local} (combined with the fact that the functor $(-)^⨯$ preserves homotopy limits), the functor $\FFT^⨯_{d,T}$
sends local weak equivalences in $\sPSh(\FEmb_d)$ to (local) weak equivalences of sheaves of $Γ$-spaces.
If we regard a manifold as a sheaf on $\FEmb_d$ via its functor of smooth points, then Čech nerves of covers of manifolds are local weak equivalences in $\sPSh(\FEmb_d)$. 
Hence if we restrict the above functor along the embedding $\Man\into \sPSh(\FEmb_d)$ we get a well defined sheaf of smooth $Γ$-spaces on the site of smooth manifolds.
For a fixed geometric structure $\gs$, we can also consider the functor on manifolds that sends 
$$X\mapsto \Funmonglob(\Bord_d^{X\times \gs},T)^⨯=\FFT^{\times}_{d,T}(X\times \gs).$$ 
Since this functor (as a functor from $\Man^\op$) is homotopy continuous, \cref{local} implies that this presheaf of smooth $Γ$-spaces satisfies descent.
More generally, to allow for twisted situations, we consider an arbitrary $\infty$-cosheaf
$$F:\Man\to \sPSh(\FEmb_d).$$

\begin{example}
Fix $d\geq 0$ and $\gs\in \sPSh(\FEmb_d)$.
The assignment $X\mapsto X\times \gs$ is a cosheaf~$F$ on $\Man$, with values in $\sPSh(\FEmb_d)$.
In the resulting bordism category $\Bord_d^{X\times \gs}$, a bordism $M$ will be equipped with a smooth map to $X$ and a geometric structure, i.e., a vertex in $\gs(M)$.
\end{example}

\begin{example}
We consider the case where $\stcart$ is the category of supercartesian spaces and $\stu$ is the reduction functor $\stu:\SMan\to \Man$ (\cref{def.supercart}).
Thus, families of bordisms in the corresponding bordism category are parametrized by supercartesian spaces.

Consider the cosheaf $F:\Man\to \sPSh(\FEmb_d)$ that sends a smooth manifold~$X$ to the presheaf of groupoids on $\FEmb_d$ that sends $(M\to U)\in\FEmb_d$
to the groupoid $F(M→U)$ whose objects are pairs $(f:G\to X,p:G\to U)$ of morphisms of supermanifolds
with $\stu(p:G\to U)=(M\to U)$ and $p$ a map of relative superdimension~$d|1$.
The morphisms in $F(M→U)$ are isomorphisms $q:G\to G'$ that make the triangles with $f$ and $p$ commute and $\Re(q)=\id_M$.

We also have a variant where $G$ is equipped with a super-Euclidean structure in the sense of Stolz–Teichner \cite{StolzTeichner.SUSY,StolzTeichner.Elliptic}.
This recovers the $0|1$, $1|1$, and $2|1$-dimensional bordism categories of Stolz and Teichner.
Hence, \cref{concordance.twisted} produces a classifying space of $d|1$-Euclidean field theories.
\end{example}

\begin{definition}
Fix $d\geq 0$, $T\in \smcat_{\infty,d}$ and an $\infty$-cosheaf $F:\Man\to \sPSh(\FEmb_d)$.
The functor 
$$\FFT^⨯_{d,T,F}:\Man^\op\to \sPSh(\stcart⨯\Gamma)_{\local}$$
sends 
$$X\mapsto \Funmonglob(\Bord_d^{F(X)},T)^⨯=\FFT^{\times}_{d,T}(F(X)).$$ 
\end{definition}

We apply the construction of Berwick-Evans–Boavida de Brito–Pavlov \cite{BEBdBP,Pavlov.Diffeo}
to the presheaf of smooth $\Gamma$-spaces $\FFT^{\times}_{d,T,F}$,
which we convert to a presheaf of $Γ$-spaces by evaluating at $\RR^0∈\stcart$.

\begin{definition}
\label{csp.cosheaf}
Fix $d\geq 0$, $T\in \smcat_{\infty,d}$ and an $\infty$-cosheaf $F:\Man\to \sPSh(\FEmb_d)$.
We define the \emph{classifying $Γ$-space of field theories} as
$$\csp\FFT^{\times}_{d,T,F}≔\hocolim_{[n]\in \Delta^\op}\FFT^⨯_{d,T,F}(\gsim^n)(\RR^0).$$
\end{definition}

\begin{theorem}
\label{concordance.twisted}
Fix $d\geq 0$, $T\in \smcat_{\infty,d}$ and an $\infty$-cosheaf $F:\Man\to \sPSh(\FEmb_d)$.
We have a natural bijection 
$$\FFT^{\times}_{d,T,F}[X]\cong [X,\csp \FFT^{\times}_{d,T,F}],$$
where $[X]$ denotes concordance classes with respect to~$X$ (after evaluating on $\RR^0∈\stcart$)
and $[X,-]$ denotes homotopy classes of maps from a smooth manifold~$X$ to a $Γ$-space.
More generally, we have a weak equivalence 
$$\hocolim_{[n]\in \Delta^\op}\FFT^{\times}_{d,T,F}(X\times \gsim^n)(\RR^0)\simeq \map(X,\csp \FFT^{\times}_{d,T,F}),$$
where $\map$ on the right side denotes the powering of a $Γ$-space by the singular complex of~$X$.
\end{theorem}

\begin{proof}
By \cref{a1a2a3}, the functor $X\mapsto\Bord_d^{F(X)}$ is an ∞-cosheaf on $\Man$.
Thus, the functor $\FFT^⨯_{d,T,F}$ is an ∞-sheaf on $\Man$ valued in $Γ$-spaces.
Hence, the main theorems of \cite{BEBdBP} and \cite{Pavlov.Diffeo} immediately imply the claim.
Taking $F(X)=X\times \gs$ yields \cref{classcnc}.
\end{proof}

In the case when all objects of $T$ are invertible with respect to the monoidal product,
the resulting $Γ$-spaces are group-like,
so can be identified with connective spectra.
In particular, $\csp\FFT^⨯_{d,T,F}$ is a connective spectrum
whose associated cohomology theory allows for a concrete model for cohomology classes over~$X$:
these are precisely concordance classes of field theories over~$X$.

Other constructions on cohomology classes, such as cup products and cohomology operations
also allow for such geometric models in terms of field theories.
In the next section, we explore one such model for power operations.

\subsection{Power operations in the Stolz–Teichner program}
\label{power.operations}

As explained at the end of the previous section, we can expect geometric models
for various operations on cohomology classes in terms of field theories.
In this section, we explore one such model for power operations,
building upon the work of Barthel–Berwick-Evans–Stapleton \cite{BarthelBerwickEvansStapleton}.

\begin{remark}
\cref{two.variants} continues to apply in this section: all results below are formulated for $\Bord_d$, but apply equally well to $\BBord_d$.
\end{remark}

\begin{definition}
\label{power.cooperation}
Fix $d\geq 0$. The $n$th power cooperation ($n≥0$) is a morphism
$$\Bord_d^{X^{⨯n}\hq Σ_n}→\Bord_d^{X},$$
where $X$ is a simplicial presheaf on $\FEmb_d$.
The notation ${-}\hq Σ_n$ denotes the homotopy quotient of simplicial presheaves presented as a bar construction, where $Σ_n$ acts on $X^{⨯n}$ by permutations.
Fixing an object in $\stcart⨯Γ⨯Δ^{⨯d}$, the cooperation on individual bordisms is defined using a pull-push construction along
$$X^{⨯n}\hq Σ_n←(X^{⨯n}⨯\bar{n})\hq Σ_n→X,$$
where $\bar{n}$ is the set on which $\Sigma_n$ acts, the left map is induced by the projection $\bar n→1$, and the right map is given by evaluation 
$$(X^{⨯n}⨯\bar{n})\hq Σ_n=(X^{\bar n}⨯\bar{n})\hq Σ_n\lto3\ev X.$$ 
The left map is an $n$-fold covering, so pulling back a map $M→X^{⨯n}\hq Σ_n$
produces an $n$-fold covering $\tilde M→M$ (in our strict model $\tilde M=M\times \bar n\to M$ is the projection map, but nontrivial $n$-fold coverings emerge after fibrantly replacing the bordism category) together with a map $\tilde M\to (X^{⨯n}⨯\bar n)\hq Σ_n$ defining the geometric structure.
Postcomposing the latter map with the evaluation map $(X^{⨯n}⨯\bar n)\hq Σ_n→X$ yields the desired bordism.
\end{definition}

\begin{definition}
\label{fft.power.operations}
Fix $d\geq 0$, $T\in \smcat_{\infty,d}$ and a simplicial presheaf $X\in \sPSh(\FEmb_d)$. The $n$th geometric power operation ($n≥0$)
$$\gpop_n:\FFT_{d,T}(X)→\FFT_{d,T}(X^{⨯n}\hq Σ_n)$$
is the morphism of smooth symmetric monoidal $(∞,d)$-categories defined via precomposition with the power cooperation
$\Bord_d(X^{⨯n}\hq Σ_n)→\Bord_d(X)$ from \cref{power.cooperation}.
\end{definition}

We now explain how to recover power operations in cohomology via \cref{fft.power.operations}.

Fix $d\geq 0$, $T\in \smcat_{\infty,d}$ and an $\infty$-cosheaf $F:\Man\to \sPSh(\FEmb_d)$ (\cref{csp.cosheaf}).
Then the functor $\FFT^{\times}_{d,T,F}$ is an $\infty$-sheaf on the site of smooth manifolds with values in $\sPSh(\stcart⨯Γ)_\local$, i.e., smooth $\Gamma$-spaces.
More generally, for an arbitrary $\infty$-sheaf of $\Gamma$-spaces $X\in \sPSh(\Man\times \Gamma)$, we set 
$$\FFT^{\times}_{d,T,F}(X)≔\Funmonglob(X,\FFT^{\times}_{d,T,F})^⨯\in \sPSh(\stcart⨯\Gamma).$$
Denote by $\hom(-,-)$ the (noncartesian) internal hom functors on $\sPSh(Γ)_\local$ and $\sPSh(\stcart⨯Γ)_\local$.
For two $\infty$-sheaves of $\Gamma$-spaces $X$ and $Y$ (i.e., fibrant objects in $\sPSh(\stcart⨯Γ)_\local$), we have a comparison map 
$$\iota:\hom(X,Y)(\RR^0)\to \hom(\csp(X),\csp(Y))\in \sPSh(\Gamma),$$
defined as the adjoint of the map 
$$\hom(X,Y)(\RR^0)\times \csp X\to \csp Y.$$
In expanded form, the latter map reads
$$\hom(X,Y)(\RR^0)\times\hocolim_{n∈Δ^\op}X(\gsim^n)\to\hocolim_{n∈Δ^\op}Y(\gsim^n),$$
which can be rewritten as the map
$$\hocolim_{n∈Δ^\op}(\hom(X,Y)(\RR^0)\times X(\gsim^n))\to\hocolim_{n∈Δ^\op}Y(\gsim^n).$$
The latter is defined objectwise for each $n$-separately, as the composition map.
The $\iota$ assemble together to form a natural transformation (in both $X$ and $Y$).
The $n$-fold monoidal product $(\FFT^⨯_{d,T,F})^{⨯n}→\FFT^⨯_{d,T,F}$ gives rise to a map
$$m_n:\hom((\csp X)^{\times n}, (\csp\FFT^⨯_{d,T,F})^{⨯n})\simeq \hom((\csp X)^{\times n}, \csp((\FFT^⨯_{d,T,F})^{⨯n}))\to \hom((\csp X)^{\times n}, \csp\FFT^⨯_{d,T,F}),$$
where we have used the fact that $\csp$ preserves finite products. We use the maps $\iota$ and $m_n$ in the following proposition, which gives refinements of power operations in cohomology. 

\begin{theorem}
The diagram
$$\xymatrix@C=2cm@R=1.2cm{
\FFT^{\times}_{d,T,F}(X)(\RR^0) \ar@/_{6pc}/[dd]_-{{\bf P}_n}\ar[r]^-{\iota_X}\ar[d]^-{\ev^*}&\hom(\csp(X),\csp(\FFT^\times_{d,T,F}))\ar[d]^-{\ev^*}\ar@/^8pc/[dd]^-{P_n}\\
\FFT^{\times}_{d,T,F}((X^{⨯n}⨯\bar n)\hq Σ_n)(\RR^0) \ar[d]\ar[r]^-{ι_{(X^{⨯n}⨯\bar n)\hq Σ_n}} & \hom((\csp(X)^{⨯n}⨯\bar n)\hq Σ_n,\csp(\FFT^{\times}_{d,T,F}))\ar[d]^-{m_n}\\
\FFT^{\times}_{d,T,F}(X^{⨯n}\hq Σ_n)(\RR^0) \ar[r]^-{ι_{X^{⨯n}\hq Σ_n}} & \hom(\csp(X)^{⨯n}\hq Σ_n,\csp(\FFT^{\times}_{d,T,F}))
}$$
commutes.
The map $\ev^*$ is defined by precomposition with the evaluation map $X^n\times \bar{n}\to X$.
The bottom square is the homotopy fixed points of the $Σ_n$-action on the commutative square adjoint to the naturality square of the $n$-fold monoidal product $(\FFT^⨯)^{⨯n}→\FFT^⨯$.
The left composition is the geometric power operation $\gpop_n$.
The right composition is the (homotopical) power operation $\tpop_n$.
\end{theorem}

\begin{proof}
The proof is simply an unwinding of the definitions. The map ${\bf P}_n$ is defined by first precomposing with the evaluation map and then precomposing with the pullback map
$$p:\Bord_d^{X^{\times n}\hq \Sigma_n}\to \Bord_d^{X^{\times n}\times \bar{n}\hq \Sigma_n},$$ which is well defined since $X^{\times n}\times \bar{n}\hq \Sigma_n\to X^{\times n}\hq \Sigma_n$ is an $n$-fold covering. The right vertical map $m_n$ is defined above, via the monoidal maps for $\FFT_{d,T,F}^{\times}$. The commutativity of the bottom square (without $\Sigma_n$) follows by the identification
$$
\Bord_d^{X^{\times n}\times \bar n}\simeq \Bord_d^{\coprod_{\bar n}X^{⨯n}}\simeq \coprod_{\bar n}\Bord_d^{X^{⨯n}},
$$
which gives rise to an identification $\FFT^{\times}_{d,T,F}(X^{⨯n}\times \bar n)\cong (\FFT^{\times}_{d,T,F}(X^{⨯n}))^{\times n}$.
The left vertical  map in the bottom square is defined by the monoidal map $(\FFT_{d,T,F}^{\times})^{\times n}\to (\FFT_{d,T,F}^{\times})$, evaluated on $X^{⨯n}$.
Then commutativity of the bottom square follows by passing to fixed points. The top square commutes for free, as it is the naturality square for $\iota$.  
\end{proof}

In particular, if the target symmetric monoidal $(\infty,d)$-category~$T$ has all objects and $k$-morphisms for $k≥1$ invertible,
then $\csp(\FFT^⨯_{d,T,F})$ is a sheaf of $\Ei$-ring spectra
and the bottom map~$\tpop_n$ is the $n$th power operation for the space $\csp(X)$
with a geometric refinement $\gpop_n$ on the level of field theories.

\begin{example}
If $X$ is a smooth manifold, then $\csp(X)$ is equivalent to the (singular simplicial set of) underlying space of $X$. Hence, the power operation $\tpop_n$ takes the form 
$$\tpop_n:\hom(X,\csp(\FFT^{\times}_{d,T,F}))\to \hom(X^{⨯n}\hq Σ_n,\csp(\FFT^{\times}_{d,T,F})),$$
which exhibits the $n$-th power operation for the connective spectrum $\csp(\FFT^{\times}_{d,T,F})$.
\end{example}


\def\doi#1{\href{https://doi.org/#1}{doi:#1}}
\def\arXiv#1{\href{https://arxiv.org/abs/#1}{arXiv:#1}}
\def\jstor#1{\href{https://www.jstor.org/stable/#1}{JSTOR:#1}}



\begin{thebibliography}{[9999.m]}

\rightskip0pt plus 1fil
\itemsep0pt plus 2pt

\bibitem[1958]{Nijenhuis}
Albert Nijenhuis.
Geometric aspects of formal differential operations on tensor fields.
Proceedings of the International Congress of Mathematicians 1958, 463–469.
\url{https://www.mathunion.org/fileadmin/ICM/Proceedings/ICM1958/ICM1958.ocr.pdf}.

\bibitem[1969]{Milnor}
John Milnor.
Morse Theory.
Annals of Mathematics Studies 51 (1969).
\jstor{j.ctv3f8rb6}.

\bibitem[1970]{Day}
Brian~J.~Day.
Construction of biclosed categories.
PhD thesis.
School of Mathematics, University of New South Wales (1970).
\doi{10.26190/unsworks/8048}.

\bibitem[1972]{BousfieldKan}
Aldridge~K.~Bousfield, Daniel~M.~Kan.
Homotopy limits, completions and localizations.
Lecture Notes in Mathematics 304 (1972).
\doi{10.1007/978-3-540-38117-4}.

\bibitem[1974]{Segal.Gamma}
Graeme Segal.
Categories and cohomology theories.
Topology 13:3 (1974), 293–312.
\doi{10.1016/0040-9383(74)90022-6}.

\bibitem[1978]{BousfieldFriedlander}
Aldridge~K.~Bousfield, Eric~M.~Friedlander.
Homotopy theory of $Γ$-spaces, spectra, and bisimplicial sets.
Lecture Notes in Mathematics 658 (1978), 80–130.
\doi{10.1007/bfb0068711}.

\bibitem[1983]{DwyerKan}
William~G.~Dwyer, Daniel~M.~Kan.
Function complexes for diagrams of simplicial sets.
Indagationes Mathematicae (Proceedings) 86:2 (1983), 139–147.
\doi{10.1016/1385-7258(83)90051-3}.

\bibitem[1988]{Atiyah}
Michael Atiyah.
Topological quantum field theories.
Publications Mathématiques de l’IHÉS, 68 (1988), 175–186.
\doi{10.1007/bf02698547}

\bibitem[1991]{MoerdijkReyes}
Ieke Moerdijk, Gonzalo~E.~Reyes.
Models for Smooth Infinitesimal Analysis.
Springer, 1991.
\doi{10.1007/978-1-4757-4143-8}.

\bibitem[1992]{Freed.Extended}
Daniel~S.~Freed.
Higher algebraic structures and quantization.
Communications in Mathematical Physics 159:2 (1994), 343–398.
\arXiv{hep-th/9212115v2}, \doi{10.1007/bf02102643}.

\bibitem[1993]{Lawrence}
Ruth~J.~Lawrence.
Triangulations, categories and extended topological field theories.
Quantum topology, 191–208.
World Scientific, Series on Knots and Everything 3 (1993).
\doi{10.1142/9789812796387\_0011}.

\bibitem[1995]{BaezDolan}
John~C.~Baez, James Dolan.
Higher‐dimensional algebra and topological quantum field theory.
Journal of Mathematical Physics 36:11 (1995), 6073–6105.
\arXiv{q-alg/9503002v2}, \doi{10.1063/1.531236}.

\bibitem[1998]{Rezk.CSS}
Charles Rezk.
A model for the homotopy theory of homotopy theory.
Transactions of the American Mathematical Society 353:3 (2000), 973–1007.
\arXiv{math/9811037v3}, \doi{10.1090/s0002-9947-00-02653-2}.

\bibitem[1999.a]{GoerssJardine}
Paul~G.~Goerss, John~F.~Jardine.
Simplicial homotopy theory.
Progress in Mathematics 174 (1999).
\doi{10.1007/978-3-0346-0189-4}.

\bibitem[1999.b]{Hovey}
Mark Hovey.
Model categories.
Mathematical Surveys and Monographs 63 (1999).
\doi{10.1090/surv/063}.

\bibitem[1999.c]{Schwede}
Stefan Schwede.
Stable homotopical algebra and $Γ$-spaces.
Mathematical Proceedings of the Cambridge Philosophical Society 126:2 (1999), 329–356.
\doi{10.1017/s0305004198003272}.

\bibitem[2000]{Toen}
Bertrand Toën.
Dualité de Tannaka supérieure I: Structures monoidales.
June 10, 2000.
\url{https://archive.mpim-bonn.mpg.de/id/eprint/912.}


\bibitem[2002.a]{Hirschhorn}
Philip~S.~Hirschhorn.
Model categories and their localizations.
American Mathematical Society.
Mathematical Surveys and Monographs 99 (2002).
\doi{10.1090/surv/099}.

\bibitem[2002.b]{DuggerIsaksen}
Daniel Dugger, Daniel~C.~Isaksen.
Weak equivalences of simplicial presheaves.
Contemporary Mathematics 346 (2004), 97–113.
\arXiv{math/0205025v1}, \doi{10.1090/conm/346/06292}.

\bibitem[2002.c]{DHI}
Daniel Dugger, Sharon Hollander, Daniel~C.~Isaksen.
Hypercovers and simplicial presheaves.
Mathematical Proceedings of the Cambridge Philosophical Society 136:1 (2004), 9–51.
\arXiv{math/0205027v2}, \doi{10.1017/S0305004103007175}.

\bibitem[2004.a]{Segal.CFT}
Graeme Segal.
The definition of conformal field theory.
London Mathematical Society Lecture Notes Series 308 (2004), 421–577.
\doi{10.1017/cbo9780511526398.019}.

\bibitem[2004.b]{StolzTeichner.Elliptic}
Stephan Stolz, Peter Teichner.
What is an elliptic object?
London Mathematical Society Lecture Note Series 308 (2004), 247–343.
\doi{10.1017/cbo9780511526398.013}.

\bibitem[2005]{Barwick.CSS}
Clark Barwick.
$(\infty,n)$-Cat as a closed model category.
Dissertation, University of Pennsylvania, 2005.
\url{https://repository.upenn.edu/dissertations/AAI3165639/}.


\bibitem[2006.a]{JoyalTierney}
André Joyal, Myles Tierney.
Quasi-categories vs Segal spaces.
In: Categories in Algebra, Geometry and Mathematical Physics.
Contemporary Mathematics 431 (2007), 277–326.
\arXiv{math/0607820v2}, \doi{10.1090/conm/431/08278}.

\bibitem[2006.b]{Shulman}
Michael Shulman.
Homotopy limits and colimits and enriched homotopy theory.
\arXiv{math/0610194v3}.

\bibitem[2007]{Barwick.Model}
Clark Barwick.
On left and right model categories and left and right Bousfield localizations.
Homology, Homotopy and Applications 12:2 (2010), 245–320.
\arXiv{0708.2067v2}, \doi{10.4310/hha.2010.v12.n2.a9}.


\bibitem[2009.a]{Rezk.Theta}
Charles Rezk.
A cartesian presentation of weak $n$-categories.
Geometry \& Topology 14:1 (2010), 521–571.
\arXiv{0901.3602v3}, \doi{10.2140/gt.2010.14.521}.

\bibitem[2009.b]{Lurie.TFT}
Jacob Lurie.
On the classification of topological field theories.
Current Developments in Mathematics 2008, 129–280.
\arXiv{0905.0465v1}, \doi{10.4310/cdm.2008.v2008.n1.a3}.


\bibitem[2009.c]{GorchinskiyGuletskii}
Sergey Gorchinskiy, Vladimir Guletskiĭ.
Symmetric powers in abstract homotopy categories.
Advances in Mathematics 292 (2016), 707–754.
\arXiv{0907.0730v4}, \doi{10.1016/j.aim.2016.01.011}.

\bibitem[2009.d]{Dugger.Spivak}
Daniel Dugger, David~I.~Spivak.
Rigidification of quasi-categories.
Algebraic \& Geometric Topology 11:1 (2011), 225–261.
\arXiv{0910.0814v1}, \doi{10.2140/agt.2011.11.225}.



\bibitem[2010]{FiorenzaSchreiberStasheff}
Domenico Fiorenza, Urs Schreiber, Jim Stasheff.
Čech cocycles for differential characteristic classes: an ∞-Lie theoretic construction.
Advances in Theoretical and Mathematical Physics 16:1 (2012), 149–250.
\arXiv{1011.4735v2}, \doi{10.4310/atmp.2012.v16.n1.a5}.

\bibitem[2011.a]{StolzTeichner.SUSY}
Stephan Stolz, Peter Teichner.
Supersymmetric field theories and generalized cohomology.
Proceedings of Symposia in Pure Mathematics 83 (2011), 279–340.
\arXiv{1108.0189v1}, \doi{10.1090/pspum/083/2742432}.

\bibitem[2011.b]{BarwickSchommerPries}
Clark Barwick, Christopher Schommer-Pries.
On the unicity of the homotopy theory of higher categories.
Journal of the American Mathematical Society 34:4 (2021), 1011–1058.
\arXiv{1112.0040v6}, \doi{10.1090/jams/972}.



\bibitem[2012]{FreedTeleman}
Daniel S. Freed, Constantin Teleman.
Relative quantum field theory.
Communications in Mathematical Physics 326:2 (2014), 459–476.
\arXiv{1212.1692v3}, \doi{10.1007/s00220-013-1880-1}.



\bibitem[2013.a]{Henriques}
André Henriques.
Three-tier CFTs from Frobenius algebras.
Topology and Field Theories.
Contemporary Mathematics 613 (2014), 1–40.
\arXiv{1304.7328v2}, \doi{10.1090/conm/613/12233}.


\bibitem[2013.b]{BunkeNikolausVolkl}
Ulrich Bunke, Thomas Nikolaus, Michael Völkl.
Differential cohomology theories as sheaves of spectra.
Journal of Homotopy and Related Structures 11:1 (2016), 1–66.
\arXiv{1311.3188v1}, \doi{10.1007/s40062-014-0092-5}.

\bibitem[2014.a]{Riehl}
Emily Riehl.
Categorical homotopy theory.
New Mathematical Monographs 24 (2014), Cambridge University Press.
\doi{10.1017/cbo9781107261457}.


\bibitem[2014.b]{BergnerRezk}
Julia~E.~Bergner, Charles Rezk.
Comparison of models for $(∞,n)$-categories, II.
Journal of Topology 13:4 (2020), 1554-1581.
\arXiv{1406.4182v3}, \doi{10.1112/topo.12167}.

\bibitem[2014.c]{AyalaFrancisTanaka}
David Ayala, John Francis, Hiro Lee Tanaka.
Local structures on stratified spaces.
Advances in Mathematics 307 (2017), 903–1028.
\arXiv{1409.0501v6}, \doi{10.1016/j.aim.2016.11.032}.

\bibitem[2015.a]{Jardine}
John~F.~Jardine.
Local Homotopy Theory.
Springer, 2015.
\doi{10.1007/978-1-4939-2300-7}.

\bibitem[2015.b]{CalaqueScheimbauer}
Damien Calaque, Claudia Scheimbauer.
A note on the $(∞,n)$-category of cobordisms.
Algebraic \& Geometric Topology 19:2 (2019), 533–655.
\arXiv{1509.08906v2}, \doi{10.2140/agt.2019.19.533}.


\bibitem[2017.a]{Lurie.HTT}
Jacob Lurie.
Higher Topos Theory.
April 9, 2017.
\url{https://www.math.ias.edu/~lurie/papers/HTT.pdf}.

\bibitem[2017.b]{Schreiber}
Urs Schreiber.
Differential cohomology in a cohesive ∞-topos.
Version 2 (August 11, 2017).
\url{https://ncatlab.org/schreiber/files/dcct170811.pdf}.

\bibitem[2017.c]{Lurie.HA}
Jacob Lurie.
Higher Algebra.
September 18, 2017.
\url{https://www.math.ias.edu/~lurie/papers/HTT.pdf}.

\bibitem[2017.d]{SchommerPries}
Christopher Schommer-Pries.
Invertible topological field theories.
\arXiv{1712.08029v1}.

\bibitem[2018]{Bergner}
Julia~E.~Bergner.
The homotopy theory of $(∞,1)$-categories.
London Mathematical Society Student Texts 90 (2018), Cambridge University Press.
\doi{10.1017/9781316181874}.


\bibitem[2019]{BEBdBP}
Daniel Berwick-Evans, Pedro Boavida de Brito, Dmitri Pavlov.
Classifying spaces of infinity-sheaves.
\arXiv{1912.10544v2}.

\bibitem[2020.a]{LudewigStoffel}
Matthias Ludewig, Augusto Stoffel.
A framework for geometric field theories and their classification in dimension one.
Symmetry, Integrability and Geometry: Methods and Applications (SIGMA) 17:072 (2021), 1–58.
\arXiv{2001.05721v2}, \doi{10.3842/SIGMA.2021.072}.

\bibitem[2020.b]{BarthelBerwickEvansStapleton}
Tobias Barthel, Daniel Berwick-Evans, Nathaniel Stapleton.
Power operations in the Stolz–Teichner program.
Geometry \& Topology 26:4 (2022), 1773–1848.
\arXiv{2006.09943v2}, \doi{10.2140/gt.2022.26.1773}.

\bibitem[2020.c]{SatiSchreiber.POC}
Hisham Sati, Urs Schreiber.
Proper Orbifold Cohomology.
\arXiv{2008.01101v2}.

\bibitem[2021.a]{Balchin}
Scott Balchin.
A Handbook of Model Categories.
Springer, 2021.
\doi{10.1007/978-3-030-75035-0}.

\bibitem[2021.b]{AmabelDebrayHaine}
Araminta Amabel, Arun Debray, Peter~J.~Haine.
Differential Cohomology: Categories, Characteristic Classes, and Connections.
\arXiv{2109.12250v2}.

\bibitem[2021.c]{GradyPavlov.GCH}
Daniel Grady, Dmitri Pavlov.
The geometric cobordism hypothesis.
\arXiv{2111.01095v3}.

\bibitem[2022.a]{Pavlov.Numerable}
Dmitri Pavlov.
Numerable open covers and representability of topological stacks.
Topology and its Applications 318:108203 (2022), 1–28.
\arXiv{2203.03120v2}, \doi{10.1016/j.topol.2022.108203}.

\bibitem[2022.b]{Pavlov.Diffeo}
Dmitri Pavlov.
Projective model structures on diffeological spaces and smooth sets and the smooth Oka principle.
\arXiv{2210.12845v1}.

\bibitem[2023.a]{nLab.Day}
nLab.
Day convolution.
Revision 72 (January 11, 2023).
\url{https://ncatlab.org/nlab/revision/Day+convolution/72}.

\bibitem[2023.b]{GradyPavlov.Str}
Daniel Grady, Dmitri Pavlov.
Geometric structures on bordisms.
In preparation.

\end{thebibliography}
\end{document}